\documentclass[a4paper, 12pt, oneside]{book}
\usepackage{listings}
\usepackage[ps2pdf, bookmarks=true,
            bookmarksnumbered=true, breaklinks=true,
            pdfstartview=FitH, hyperfigures=false,
            plainpages=false, naturalnames=true,
            colorlinks=true, linkcolor=Blue,
            citecolor=BrickRed, urlcolor=Blue,
            filecolor=Blue, pagecolor=Blue,
            pdfpagelabels, pagebackref=false]{hyperref}
\usepackage[all]{xy}
\usepackage[usenames]{color}
\usepackage{setspace}
\usepackage[compact]{titlesec}
\usepackage{makeidx}
\input{xypic}
\usepackage{geometry, enumerate, amsmath, amssymb,
            latexsym, theorem, graphicx, mathdots}
\usepackage{algorithmic}
\usepackage[ruled]{algorithm}
\hypersetup{
pdfauthor   = {Gareth Alun Evans},
pdftitle    = {Noncommutative Involutive Bases},
pdfsubject  = {Mathematics},
pdfkeywords = {AMS MSC 2000: 13P10 (Polynomial ideals,
Gr\"obner bases); 16S15 (Finite generation, finite
presentability, normal forms (diamond lemma, term rewriting));
16Z05 (Computational aspects of associative rings).},
pdfcreator  = {LaTeX with hyperref package},
pdfproducer = {dvips + ps2pdf}
}

\makeindex

{\theorembodyfont{\rmfamily}\newtheorem{example}{Example}[section]}
{\theorembodyfont{\rmfamily}\newtheorem{defn}[example]{Definition}}
{\theorembodyfont{\rmfamily}}
\newtheorem{prop}[example]{Proposition}
\newtheorem{thm}[example]{Theorem}
\newtheorem{conj}[example]{Conjecture}
{\theorembodyfont{\rmfamily}\newtheorem{remark}[example]{Remark}}
{\theorembodyfont{\rmfamily}\newtheorem{openquestion}{Open Question}}
\newtheorem{cor}[example]{Corollary}
\newtheorem{lem}[example]{Lemma}

\newenvironment{pf}{\noindent \textbf{Proof:} }{\hfill $\Box$}

\newcommand{\init}{\mathrm{in}}
\DeclareMathAlphabet{\mathpzc}{OT1}{pzc}{m}{it}

\newcommand{\LC} {\mathrm{LC}}
\newcommand{\LM} {\mathrm{LM}}
\newcommand{\LT} {\mathrm{LT}}
\newcommand{\lcm} {\mathrm{lcm}}
\newcommand{\val} {\mathrm{val}}
\newcommand{\PRE} {\mathrm{Prefix}}
\newcommand{\SUFF} {\mathrm{Suffix}}
\newcommand{\SUB} {\mathrm{Subword}}
\newcommand{\Sug} {\mathrm{Sug}}
\newcommand{\Rem} {\mathrm{Rem}}

\begin{document}

% SET PAGE NUMBERING
\pagestyle{plain}
\pagenumbering{Roman}

% TITLE PAGE
\begin{titlepage}
\vspace*{0.5in}
\begin{center}
{\Huge{Noncommutative Involutive Bases}} \\
\begin{large}
\vspace{0.75in}
Thesis submitted to the University of Wales in support of \\
the application for the degree of {Philosophi\ae} Doctor  \\[0.25in]
by \\[0.25in]
{\Large{Gareth Alun Evans}} \\
\end{large}
\end{center}
\vfill
\hfill
\parbox{2.8in}{School of Informatics \\
The University of Wales, Bangor \\
September 2005 \\[5mm]
\includegraphics[scale=0.2]{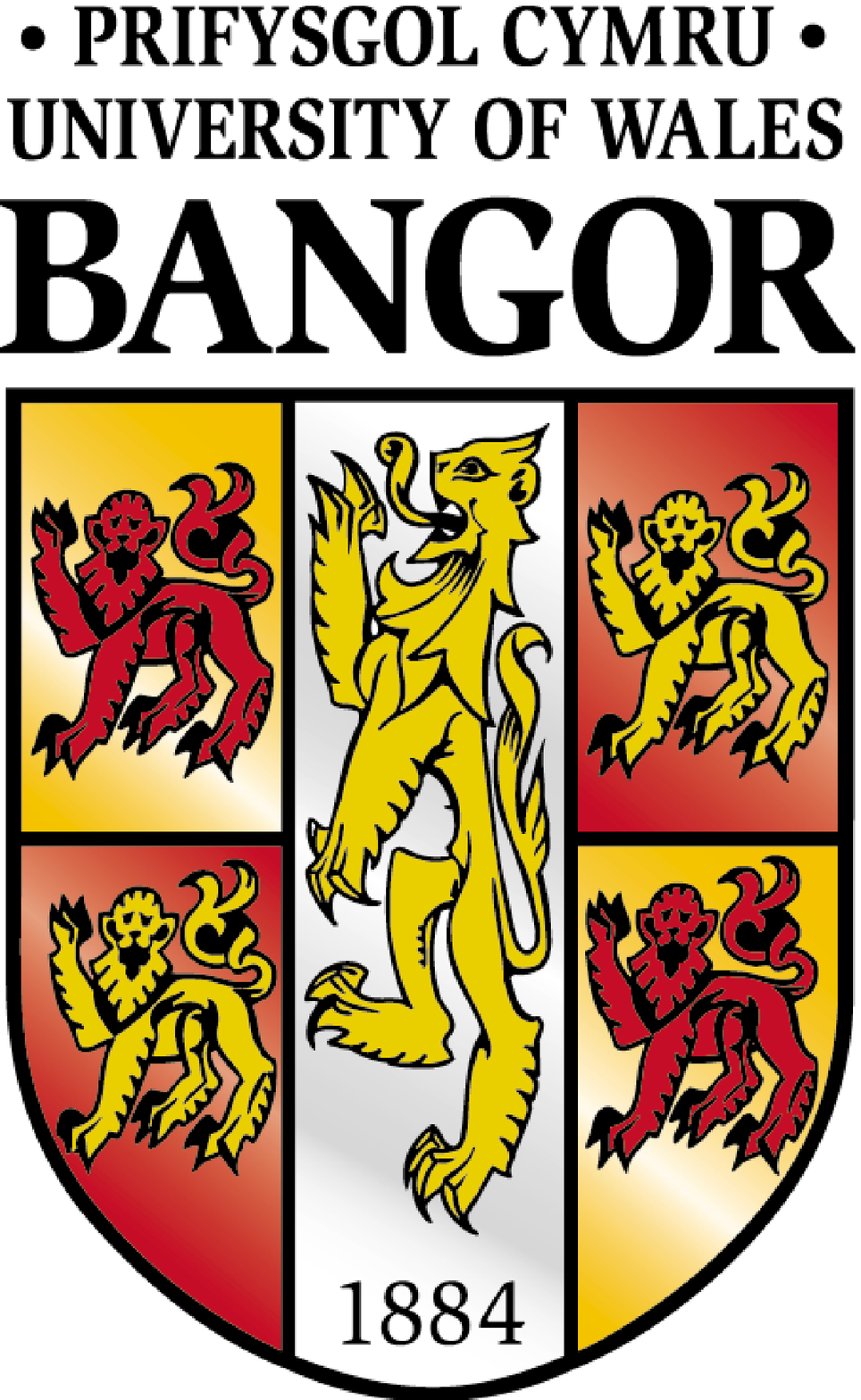}} \\[0.5cm]
\end{titlepage}

% DECLARATIONS
\clearpage
\pagenumbering{roman}
%
% Declarations
% Author: Gareth Evans
% Last Modified: 29th September 2005
%

{\large DECLARATION}

This work has not previously been accepted in substance for any degree and is
not being concurrently submitted in candidature for any degree.

Signed \dotfill\hspace*{1mm}(candidate)\hspace*{2.5cm} \\
Date \dotfill\hspace*{1mm}\phantom{(candidate)}\hspace*{2.5cm} \\[5mm]

{\large STATEMENT 1}

This thesis is the result of my own investigations, except where otherwise
stated. \\
Other sources are acknowledged by explicit references.
A bibliography is appended.

Signed \dotfill\hspace*{1mm}(candidate)\hspace*{2.5cm} \\
Date \dotfill\hspace*{1mm}\phantom{(candidate)}\hspace*{2.5cm} \\[5mm]

{\large STATEMENT 2}

I hereby give consent for my thesis, if accepted, to be available for
photocopying and for inter-library loan, and for the title and summary
to be made available to outside organisations.

Signed \dotfill\hspace*{1mm}(candidate)\hspace*{2.5cm} \\
Date \dotfill\hspace*{1mm}\phantom{(candidate)}\hspace*{2.5cm} \\[5mm]
\thispagestyle{empty}

% SUMMARY
\begin{spacing}{1}
%
% Summary
% Author: Gareth Evans
% Last Modified: 24th September 2005
%

\clearpage
\phantomsection
\addcontentsline{toc}{chapter}{Summary}
\chapter*{Summary}

%\vspace*{-3mm}
The theory of Gr\"obner Bases originated in the work
of Buchberger \cite{buch65} and is now considered
to be one of the most important and useful
areas of symbolic computation. A great deal of effort
has been put into improving Buchberger's algorithm
for computing a Gr\"obner Basis, and indeed in finding
alternative methods of computing Gr\"obner Bases. Two
of these methods include the Gr\"obner Walk method \cite{AGK97}
and the computation of Involutive Bases \cite{ZharBlink93}.

By the mid 1980's, Buchberger's work had been generalised 
for noncommutative polynomial rings by Bergman \cite{Bergman78}
and Mora \cite{mora86}. This thesis provides
the corresponding generalisation for Involutive Bases and
(to a lesser extent) the Gr\"obner Walk, with the main results 
being as follows.
% of the thesis being the following.
\vspace*{-2mm}
\begin{enumerate}[(1)]
\item
Algorithms for several new noncommutative
involutive divisions are given, including strong; weak;
global and local divisions.
\item
\vspace*{-2mm}
An algorithm for computing a
noncommutative Involutive Basis is given. When used with
one of the aforementioned involutive divisions,
it is shown that this algorithm returns a noncommutative
Gr\"obner Basis on termination.
\item
\vspace*{-2mm}
An algorithm for a noncommutative Gr\"obner
Walk is given, in the case of conversion between
two harmonious monomial orderings. It is shown that this algorithm
generalises to give an algorithm for performing
a noncommutative Involutive Walk, again in the case of
conversion between two harmonious monomial
orderings.
\item
\vspace*{-2mm}
Two new properties of commutative involutive divisions
are introduced (stability and extendibility), respectively
ensuring the termination of the Involutive Basis
algorithm and the applicability (under certain conditions)
of homogeneous methods of computing Involutive Bases.
\end{enumerate}
\vspace*{-2mm}
Source code for an initial implementation of an algorithm to compute
noncommutative Involutive Bases is provided in Appendix \ref{appB}. 
This source code, written using ANSI C and a series of libraries 
({\sf AlgLib}) provided by MSSRC \cite{MSSRC},
forms part of a larger collection of programs providing examples for the
thesis, including implementations of the commutative and noncommutative
Gr\"obner Basis algorithms \cite{buch65, mora86}; 
the commutative Involutive Basis algorithm
for the Pommaret and Janet involutive divisions \cite{ZharBlink93}; and
the Knuth-Bendix critical pairs completion algorithm for monoid
rewrite systems \cite{KB}.
% thesis, including implementations of:
% \vspace*{-1mm}
% \begin{enumerate}[(1)]
% \item
% The commutative Gr\"obner Basis algorithm \cite{buch65},
% incorporating theoretical
% improvements such as the sugar strategy \cite{GMNRT91};
% Buchberger's criteria \cite{buch79}; and the
% Gr\"obner Walk \cite{AGK97};
% \item
% The commutative Involutive Basis algorithm for the Pommaret and
% Janet divisions \cite{ZharBlink93};
% \item
% The noncommutative Gr\"obner Basis algorithm \cite{mora86}; 
% % , including options
% % allowing the selection of the wreath product monomial ordering and the
% % computation of Logged Gr\"obner Bases;
% \item
% The Knuth-Bendix critical pairs completion algorithm 
% for monoid rewrite systems \cite{KB}.
% \end{enumerate}

% ACKNOWLEDGEMENTS
%
% Acknowledgements
% Author: Gareth Evans
% Last Modified: 14th September 2005
%

\clearpage
\phantomsection
\addcontentsline{toc}{chapter}{Acknowledgements}
\chapter*{Acknowledgements}

Many people have inspired me to complete this thesis,
and I would like to take this opportunity to thank some
of them now.

I would like to start by thanking my family for their
constant support, especially my parents who have
encouraged me every step of the way. {\it Mae fy nyled
yn fawr iawn i chi.}

I would like to thank Prof. Larry Lambe from MSSRC,
whose software allowed me to test my theories in a way
that would not have been possible elsewhere.

Thanks to all the Mathematics Staff and Students
I have had the pleasure of working with over the past
seven years. Particular thanks go to Dr. Bryn Davies,
who encouraged me to think independently; to Dr. Jan Abas,
who inspired me to reach goals I never thought I could
reach; and to Prof. Ronnie Brown, who introduced me to
Involutive Bases.

I would like to finish by thanking my Supervisor
Dr. Chris Wensley. Our regular meetings kept the cogs in
motion and his insightful comments enabled me to
avoid wrong turnings and to get the little details right.
{\it Diolch yn fawr!}

This work has been gratefully supported by the EPSRC
and by the School of Informatics at the University 
of Wales, Bangor.
\vfill
\begin{flushright}
Typeset using \LaTeX,
{\sf XFig} and \XY-pic. \\
% AMS MSC 2000: 13P10, 16S15, 16Z05.
\end{flushright}

\end{spacing}

% QUOTE
\newpage
\begin{center}
{\it ``No one has ever done anything like this.'' \\
     ``That's why it's going to work.'' }
\begin{flushright}
\vspace{-4mm}
The Matrix \cite{Matrix}
\end{flushright}
\end{center}

% TABLE OF CONTENTS
\newpage
\begin{spacing}{0.9}
\pdfbookmark[0]{\contentsname}{toc}
\tableofcontents
\end{spacing}

% LIST OF ALGORITHMS
\clearpage
\phantomsection
\pdfbookmark[0]{List of Algorithms}{loa}
\begin{spacing}{1.5}
\listofalgorithms
\end{spacing}

% SET PAGE NUMBERING
\pagestyle{headings}
\clearpage
\pagenumbering{arabic}

% CHAPTERS
%
% Chapter 0
% Author: Gareth Evans
% Last Modified: 20th September 2005
%

\clearpage
\phantomsection
\addcontentsline{toc}{chapter}{Introduction}
\chapter*{Introduction}

\phantomsection
\addcontentsline{toc}{section}{Background}
\section*{Background}

\subsection*{Gr\"obner Bases}

During the second half of the twentieth century, one of the most
successful applications of symbolic computation was in the development
and application of    % commutative
{\it Gr\"obner Basis} theory for finding special bases of
ideals in commutative polynomials rings.
Pioneered by Bruno Buchberger in 1965 \cite{buch65},
the theory allowed an answer to the question
``What is the unique remainder when a polynomial is divided by a set
of polynomials?''. Buchberger's algorithm for computing a
Gr\"obner Basis was improved and refined over several decades
\cite{AGK97, buch79, FGLM, GMNRT91},
aided   % in no small part
by the development of powerful symbolic computation systems
over the same period.
% Research on improving the efficiency of the algorithm continues.
Today there is an implementation of Buchberger's algorithm
in virtually all general purpose symbolic computation systems,
including Maple \cite{Maple05} and Mathematica \cite{Mathematica04},
and many more specialised systems.
% One of the best implementations is currently
% generally acknowledged to be that in Singular \cite{Sing05}.

\subsubsection*{What is a Gr\"obner Basis?}

Consider the problem of finding the remainder when a
number is divided by a set of numbers.
If the dividing set contains just one number, then the problem
only has one solution. For example, ``$5$'' is the only possible
answer to the question ``What is $20 \div 4$?''.
If the dividing set contains
more than one number however, there may be several solutions,
as the division can potentially be performed in more than
one way.

{\bf Example.}
Consider a tank containing 21L of water. Given two 
empty jugs, one with a capacity of 2L and the other 5L,
is it possible to empty the tank using just the jugs,
assuming only full jugs of water may be
removed from the tank?
\begin{center}
\includegraphics{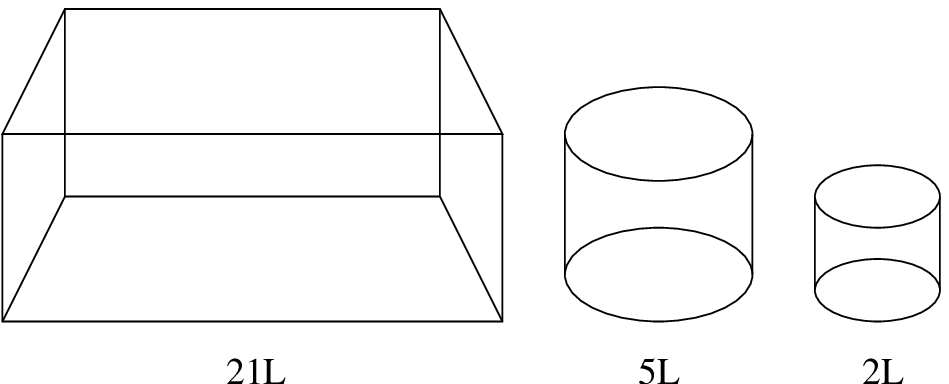}
\end{center}
Trying to empty the tank using the 2L jug only, we are
able to remove $10 \times 2 = 20$L of water
from the tank, and we are left
with 1L of water in the tank. Repeating with
the 5L jug, we are again left with 1L of water in the tank.
If we alternate between the jugs however
(removing 2L of water followed by 5L followed
by 2L and so on), the tank this time does become empty,
because $21 = 2+5+2+5+2+5$.

The observation that we are left with a different volume of water
in the tank dependent upon how we try to empty it
corresponds to the idea that the remainder obtained
when dividing the number 21 by the numbers 2 and 5 is
dependent upon how the division is performed.

This idea also applies when dividing polynomials by sets
of polynomials --- remainders here will also be dependent
upon how the division is performed. However, if we divide a
polynomial with respect to a set of polynomials that is
a Gr\"obner Basis, then we will always
obtain the same remainder no matter how the
division is performed. This fact, along with the fact that any set of
polynomials can be transformed into an equivalent set
of polynomials that is a Gr\"obner Basis, provides the
main ingredients of Gr\"obner Basis theory.

{\bf Remark.} The `Gr\"obner Basis' for our water tank
example would be just a 1L jug, allowing us to empty any
tank containing $n$L of water (where $n \in \mathbb{N}$).

\subsubsection*{Applications}

There are numerous applications of Gr\"obner Bases in all branches of
mathematics, computer science, physics and engineering \cite{Buch98}.
% both in academia and in the wider world.
Topics vary from geometric theorem proving
to solving systems of polynomial equations,
and from algebraic coding theory to the design of experiments in statistics.

{\bf Example.}
% Let $F := \{x+y+z=3, \, x^2+y^2+z^2=9, \, x^3+y^3+z^3=27\}$
% be a set of polynomial equations. One way of solving this
% set for $x$, $y$ and $z$ is to calculate a {\it lexicographic}
% Gr\"obner Basis for $F$. This yields the set
% $G := \{x+y+z=3, \, y^2+yz-3y+z^2-3z=0, \, z^3=3z^2\}$, which
% then allows us to deduce that $z = 0$ or $3$, that $y = 0$
% when $z = 3$, and so on.
Let $F := \{x+y+z=6, \, x^2+y^2+z^2=14, \, x^3+y^3+z^3=36\}$
be a set of polynomial equations. One way of solving this
set for $x$, $y$ and $z$ is to compute a {\it lexicographic}
Gr\"obner Basis for $F$. This yields the set
$G := \{x+y+z=6, \, y^2+yz+z^2-6y-6z=-11, \, z^3-6z^2+11z=6\}$, 
the final member of which is a univariate polynomial in $z$,
a polynomial we can solve to deduce that $z = 1$, $2$ or $3$.
Substituting back into the second member of $G$, when
$z = 1$, we obtain the polynomial $y^2-5y+6=0$, which
enables us to deduce that $y = 2$ or $3$; when
$z = 2$, we obtain the polynomial $y^2-4y+3=0$, which
enables us to deduce that $y = 1$ or $3$; and when
$z = 3$, we obtain the polynomial $y^2-3y+2=0$, which
enables us to deduce that $y = 1$ or $2$. Further
substitution into $x+y+z=6$ then enables us to deduce
the value of $x$ in each of the above cases, enabling us
to give the following table of solutions for $F$.
\begin{center}
\begin{tabular}{|c|cccccc|}  \hline
x & 3 & 2 & 3 & 1 & 2 & 1 \\ \hline
y & 2 & 3 & 1 & 3 & 1 & 2 \\ \hline
z & 1 & 1 & 2 & 2 & 3 & 3 \\ \hline
\end{tabular}
\end{center}

\subsection*{Involutive Bases}

As Gr\"obner Bases became popular, researchers noticed
a connection between Buchberger's ideas and ideas originating from
the Janet-Riquier theory of Partial Differential Equations % \cite{Janet29}
developed in the early 20th century (see for example \cite{Mansfield91}).
This link was completed for commutative polynomial rings
by Zharkov and Blinkov in the early 1990's \cite{ZharBlink93} when they
gave an algorithm to compute an {\it Involutive Basis} that
provides an alternative way of computing a Gr\"obner Basis.
Early implementations of this algorithm
(an elementary introduction to which can be found in \cite{CHS00})
% which differs considerably from Buchberger's algorithm,
compared favourably with the most advanced implementations of
Buchberger's algorithm, with results in \cite{Gerdt98a} showing the
potential of the Involutive method in terms of efficiency.
% The behaviour of examples is unpredictable: sometimes the
% involutive method completes far more quickly, and sometimes
% the Gr\"obner method is superior.

\subsubsection*{What is an Involutive Basis?}

Given a Gr\"obner Basis $G$, we know that the remainder
obtained from dividing a polynomial with respect to $G$
will always be the same no matter how the division is
performed. With an Involutive
Basis, the difference is that there is only one
way for the division to be performed, so that unique
remainders are also obtained uniquely.

This effect is achieved through assigning a set of
{\it multiplicative variables} to each polynomial in an
Involutive Basis $H$, imposing a restriction on how
polynomials may be divided by $H$ by only allowing
any polynomial $h \in H$ to be multiplied by its
corresponding multiplicative variables.
Popular schemes of assigning multiplicative variables
include those based on the work of Janet \cite{Janet29},
Thomas \cite{Thomas37} and Pommaret \cite{Pommaret78}.

{\bf Example.}
Consider the Janet Involutive Basis
$H := \{xy-z, \; yz+2x+z, \; 2x^2+xz+z^2, \;
2x^2z+xz^2+z^3\}$ with multiplicative
variables as shown in the table below.
\begin{center}
\begin{tabular}{c|c}
Polynomial & Janet Multiplicative Variables \\ \hline
$xy-z$ & $\{x, y\}$ \\
$yz+2x+z$ & $\{x, y, z\}$ \\
$2x^2+xz+z^2$ & $\{x\}$ \\
$2x^2z+xz^2+z^3$ & $\{x, z\}$ \\ \hline
\end{tabular}
\end{center}
To illustrate that any polynomial may only be
{\it involutively divisible} by at most one member of any Involutive
Basis, we include the following two diagrams,
showing which monomials are involutively divisible
by $H$, and which are divisible by the
corresponding Gr\"obner Basis
$G := \{xy-z, \; yz+2x+z, \; 2x^2+xz+z^2\}$.
\begin{center}
\begin{picture}(0,0)%
\includegraphics{ch0d2.pstex}%
\end{picture}%
\setlength{\unitlength}{2605sp}%
\begingroup\makeatletter\ifx\SetFigFont\undefined%
\gdef\SetFigFont#1#2#3#4#5{%
  \reset@font\fontsize{#1}{#2pt}%
  \fontfamily{#3}\fontseries{#4}\fontshape{#5}%
  \selectfont}%
\fi\endgroup%
\begin{picture}(10302,5352)(1141,-6586)
\put(9676,-2461){\makebox(0,0)[lb]{\smash{\SetFigFont{12}{14.4}{\familydefault}{\mddefault}{\updefault}{\color[rgb]{0,0,0}$y$}%
}}}
\put(3976,-2461){\makebox(0,0)[lb]{\smash{\SetFigFont{12}{14.4}{\familydefault}{\mddefault}{\updefault}{\color[rgb]{0,0,0}$y$}%
}}}
\put(1201,-6586){\makebox(0,0)[lb]{\smash{\SetFigFont{12}{14.4}{\familydefault}{\mddefault}{\updefault}{\color[rgb]{0,0,0}Gr\"obner Basis}%
}}}
\put(6901,-6586){\makebox(0,0)[lb]{\smash{\SetFigFont{12}{14.4}{\familydefault}{\mddefault}{\updefault}{\color[rgb]{0,0,0}Involutive Basis}%
}}}
\put(2641,-6136){\makebox(0,0)[lb]{\smash{\SetFigFont{12}{14.4}{\rmdefault}{\mddefault}{\updefault}{\color[rgb]{0,0,0}$x^2$}%
}}}
\put(8371,-6136){\makebox(0,0)[lb]{\smash{\SetFigFont{12}{14.4}{\rmdefault}{\mddefault}{\updefault}{\color[rgb]{0,0,0}$x^2$}%
}}}
\put(8386,-5536){\makebox(0,0)[lb]{\smash{\SetFigFont{12}{14.4}{\rmdefault}{\mddefault}{\updefault}{\color[rgb]{0,0,0}$x^2z$}%
}}}
\put(1396,-4276){\makebox(0,0)[lb]{\smash{\SetFigFont{12}{14.4}{\rmdefault}{\mddefault}{\updefault}{\color[rgb]{0,0,0}$yz$}%
}}}
\put(7126,-4276){\makebox(0,0)[lb]{\smash{\SetFigFont{12}{14.4}{\rmdefault}{\mddefault}{\updefault}{\color[rgb]{0,0,0}$yz$}%
}}}
\put(4861,-6511){\makebox(0,0)[lb]{\smash{\SetFigFont{12}{14.4}{\familydefault}{\mddefault}{\updefault}{\color[rgb]{0,0,0}$x$}%
}}}
\put(10561,-6511){\makebox(0,0)[lb]{\smash{\SetFigFont{12}{14.4}{\familydefault}{\mddefault}{\updefault}{\color[rgb]{0,0,0}$x$}%
}}}
\put(2281,-5206){\makebox(0,0)[lb]{\smash{\SetFigFont{12}{14.4}{\rmdefault}{\mddefault}{\updefault}{\color[rgb]{0,0,0}$xy$}%
}}}
\put(7981,-5191){\makebox(0,0)[lb]{\smash{\SetFigFont{12}{14.4}{\rmdefault}{\mddefault}{\updefault}{\color[rgb]{0,0,0}$xy$}%
}}}
\put(1141,-1441){\makebox(0,0)[lb]{\smash{\SetFigFont{12}{14.4}{\familydefault}{\mddefault}{\updefault}{\color[rgb]{0,0,0}$z$}%
}}}
\put(6841,-1426){\makebox(0,0)[lb]{\smash{\SetFigFont{12}{14.4}{\familydefault}{\mddefault}{\updefault}{\color[rgb]{0,0,0}$z$}%
}}}
\end{picture}

\end{center}
Note that the irreducible monomials of both bases
all appear in the set
$\{1, \; x, \; y^i, \; z^i, \; xz^i\}$,
where $i \geqslant 1$; and that the cube, the 2 planes
and the line shown in the right hand diagram do not
overlap.

\subsection*{Noncommutative Bases}

There are certain types of noncommutative algebra to which methods for
commutative Gr\"obner Bases may be applied.
Typically, these are algebras with generators $\{x_1,\ldots,x_n\}$
for which products $x_jx_i$ with $j>i$
may be rewritten as $(x_ix_j \: +$ other terms$)$.
For example, version 3-0-0 of Singular \cite{Sing05}
(released in June 2005)
allows the computation of Gr\"obner Bases for $G$-algebras.
% (which have a Poincar\'e-Birkhoff-Witt (PBW) basis).
% The related skew polynomial rings are 
% discussed by Reinert in \cite{Reinert95},
% a thesis which considers the application of Gr\"obner Bases
% to monoid and group rings.

To compute Gr\"obner Bases for ideals 
in free associative algebras however,
one must turn to the theory of
{\it noncommutative Gr\"obner Bases}.
Based on the work of Bergman \cite{Bergman78} 
and Mora \cite{mora86}, the theory answers the question
% as well as answering the question
``What is the remainder when a noncommutative polynomial is divided
by a set of noncommutative polynomials?'', and
allows us to find Gr\"obner Bases for such algebras
as path algebras \cite{Keller97}.

The final piece of the jigsaw is to mirror the application
of Zharkov and Blinkov's Involutive methods to the noncommutative case.
This thesis provides the first
extended attempt at accomplishing this task,
improving the author's first basic algorithms for computing
{\it noncommutative Involutive Bases} \cite{Evans04}
and providing a full theoretical foundation for these algorithms.

\phantomsection
\addcontentsline{toc}{section}{Structure and Principal Results}
\section*{Structure and Principal Results}

This thesis can be broadly divided into two parts:
Chapters \ref{ChPr} through \ref{ChCIB} survey the 
% theory providing the
building blocks required for the theory of noncommutative
Involutive Bases; the remainder of the thesis then
describes this theory together with different ways of
computing noncommutative Involutive Bases.

\subsection*{Part 1}

Chapter \ref{ChPr} contains accounts of some
necessary preliminaries for our studies -- a review of
both commutative and noncommutative polynomial
rings; ideals; monomial orderings; and polynomial
division.

We survey the theory of {\it commutative Gr\"obner Bases} in
Chapter \ref{ChCGB}, basing our account on many sources,
but mainly on the books \cite{becker93} and \cite{Froberg98}.
We present the theory from the viewpoint of S-polynomials
(for example defining a Gr\"obner Basis in terms
of S-polynomials), mainly because Buchberger's algorithm
for computing a Gr\"obner Basis deals predominantly
with S-polynomials. Towards the end of the Chapter, we
describe some of the theoretical improvements of
Buchberger's algorithm, including the usage of
selection strategies, optimal variable orderings and
Logged Gr\"obner Bases.

The viewpoint of defining Gr\"obner Bases in terms
of S-polynomials continues in Chapter \ref{ChNCGB},
where we encounter the theory of {\it noncommutative
Gr\"obner Bases}. We discover that the theory is
quite similar to that found in the previous chapter,
apart from the definition of an S-polynomial and the
fact that not all input bases will have finite
Gr\"obner Bases.

In Chapter \ref{ChCIB}, we acquaint ourselves with
the theory of {\it commutative Involutive Bases}. This is
based on the work of Zharkov and Blinkov \cite{ZharBlink93};
Gerdt and Blinkov \cite{Gerdt98a, Gerdt97};
Gerdt \cite{Gerdt02, Gerdt05}; Seiler \cite{Seiler02b, Seiler02c};
and Apel \cite{Apel95, Apel98a}, with the notation
and conventions taken from a combination of these papers.
For example, notation for involutive cones and
multiplicative variables is taken from \cite{Gerdt98a},
and the definition of an involutive division and the
algorithm for computing an Involutive Basis is taken
from \cite{Seiler02b}.

As for the content of Chapter \ref{ChCIB}, we
introduce the Janet, Pommaret and Thomas divisions
in Section \ref{4point1}; describe what is meant
by a prolongation and autoreduction in
Section \ref{4point2}; introduce the properties
of continuity and constructivity in Section
\ref{4point3}; give the Involutive Basis algorithm
in Section \ref{4point4}; and describe some
improvements to this algorithm in Section \ref{4point5}.
In between all of this, we introduce two new
properties of involutive divisions,
stability and extendibility, that ensure (respectively)
the termination of the Involutive Basis algorithm and
the applicability (under certain conditions) of homogeneous
methods of computing Involutive Bases.

\subsection*{Part 2}

The main results of the thesis are contained in
Chapter \ref{ChNCIB}, where we introduce the theory
of {\it noncommutative Involutive Bases}. In Section
\ref{5point1}, we define two methods of performing
noncommutative involutive reduction, the first of which
(using thin divisors) allows the mirroring of
theory from Chapter \ref{ChCIB}, and the second
of which (using thick divisors) allows
efficient computation of involutive remainders.
We also define what is meant by a
noncommutative involutive division, 
and give an algorithm for performing
noncommutative involutive reduction.

In Section \ref{5point2}, we generalise the notions of
prolongation and autoreduction to the noncommutative
case, introducing two different types of prolongation
(left and right) to reflect the fact that left and
right multiplication are different operations in
noncommutative polynomial rings. These notions are then
utilised in the algorithm for computing a noncommutative
Involutive Basis, which we present in Section \ref{5point4}.

In Section \ref{CoCo}, we introduce two properties of
noncommutative involutive divisions. Continuity
helps ensure that any Locally Involutive Basis is an
Involutive Basis; conclusivity
ensures that for any given input basis, a
finite Involutive Basis will exist if and only if
a finite Gr\"obner Basis exists. A third property
is also introduced for weak involutive divisions
to ensure that any Locally Involutive Basis is a Gr\"obner
Basis (Involutive Bases with respect to strong involutive
divisions are automatically Gr\"obner Bases).

Section \ref{5point5} provides several involutive divisions 
for use with the noncommutative Involutive Basis
algorithm, including two global
divisions and ten local divisions. The properties
of these divisions are analysed, with full proofs
given that certain divisions satisfy
certain properties. We also show that some
divisions are naturally suited for efficient involutive reduction,
and speculate on the existence of further
involutive divisions.

In Section \ref{5point6}, we briefly discuss the topic of
the termination of the noncommutative Involutive
Basis algorithm. In Section \ref{Ch6Ex},
we provide several examples showing how noncommutative
Involutive Bases are computed, including examples
demonstrating the computation of involutive complete
rewrite systems for groups. Finally, in Section
\ref{5point8}, we discuss improvements to the
noncommutative Involutive Basis algorithm,
including how to introduce efficient
involutive reduction and Logged Involutive Bases.

Chapter \ref{ChWalk} introduces and generalises the theory
of the {\it Gr\"obner Walk}, where a Gr\"obner Basis
with respect to one monomial ordering may be computed from a
Gr\"obner Basis with respect to another monomial
ordering. In Section \ref{6point1}, we summarise the theory
of the commutative Gr\"obner Walk (based on the
papers \cite{AGK97} and \cite{CKMWalk}), and we describe a
generalisation of the theory to the Involutive case due to
Golubitsky \cite{Golub}. In Section \ref{6point2}, we then go on
to partially generalise the theory to the noncommutative
case, giving algorithms to perform both Gr\"obner and
Involutive Walks between two harmonious monomial orderings.
% Examples demonstrating the use of both of these new
% techniques are given, and we speculate on the existence
% of noncommutative walks between two arbitrary monomial
% orderings.

After some concluding remarks in Chapter \ref{ChConc},
we provide full proofs for two Propositions from Section 
\ref{5point5} in Appendix \ref{appA}.
Appendix \ref{appB} then provides ANSI C
source code for an initial implementation of
the noncommutative Involutive Basis algorithm,
together with a brief description of the {\sf AlgLib}
libraries used in conjunction with the code.
Finally, in Appendix \ref{appC}, we provide 
sample sessions showing the program given in
Appendix \ref{appB} in action. % together with some
% timings for different examples.

%
% Chapter 1
% Author: Gareth Evans
% Last Modified: 2nd February 2006
%

\chapter{Preliminaries} \label{ChPr}

In this chapter, we will set out some algebraic concepts
that will be used extensively in the following chapters. In
particular, we will introduce polynomial rings and ideals,
the main objects of study in this thesis.

\section{Rings and Ideals}

\subsection{Groups and Rings}

\begin{defn}
A
\index{binary operation} {\it binary operation}
on a set $S$ is a function 
$\ast: \: S \times S \rightarrow S$
such that associated with each ordered pair $(a, b)$ of elements of $S$ is a
uniquely defined element $(a\ast b) \in S$.
\end{defn}

\begin{defn}
A
\index{group} {\it group}
is a set $G$, with a binary operation $\ast$, such that the following
conditions hold.
\begin{enumerate}[(a)]
\item
$g_1\ast g_2 \in G$ for all $g_1, g_2 \in G$ (closure).
\item
$g_1\ast (g_2\ast g_3) = (g_1\ast g_2)\ast g_3$ for all $g_1, g_2, g_3 \in G$
(associativity).
\item
There exists an element $e \in G$ such that for all $g \in G$,
$e\ast g = g = g\ast e$ (identity).
\item
For each element $g \in G$, there exists an element $g^{-1} \in G$
such that $g^{-1}\ast g = e = g\ast g^{-1}$ (inverses).
\end{enumerate}
\end{defn}

\begin{defn}
A group $G$ is
\index{abelian} {\it abelian}
if the binary operation of the
group is commutative, that is $g_1\ast g_2 = g_2\ast g_1$
for all $g_1, g_2 \in G$. The operation in an abelian group
is often written additively, as $g_1 + g_2$, with the
inverse of $g$ written $-g$.
\end{defn}

\begin{defn}
A
\index{rng} {\it rng}
is a set $R$ with two binary operations $+$ and $\times$, known as
addition and multiplication, such that addition has an
identity element $0$, called {\it zero}, and the following axioms hold.
\begin{enumerate}[(a)]
\item
$R$ is an abelian group with respect to addition.
\item
$(r_1 \times r_2) \times r_3 = r_1 \times (r_2 \times r_3)$ for all
$r_1, r_2, r_3 \in R$ (multiplication is associative).
\item
$r_1 \times (r_2 + r_3) = r_1 \times r_2 + r_1 \times r_3$ and
$(r_1 + r_2) \times r_3 = r_1 \times r_3 + r_2 \times r_3$ for all
$r_1, r_2, r_3 \in R$ (the distributive laws hold).
\end{enumerate}
\end{defn}

\begin{defn}
A rng $R$ is a
\index{ring} {\it ring}
if it contains a unique
element $1$, called the {\it unit} element, such that
$1 \neq 0$ and $1 \times r = r = r \times 1$
for all $r \in R$.
\end{defn}

\begin{defn}
A ring $R$ is
\index{ring!commutative} {\it commutative}
if \index{commutative ring}
multiplication (as well as
addition) is commutative, that is
$r_1 \times r_2 = r_2 \times r_1$ for all $r_1, r_2 \in R$.
\end{defn}

\begin{defn}
A ring $R$ is
\index{ring!noncommutative} {\it noncommutative}
if \index{noncommutative ring}
$r_1 \times r_2 \neq r_2 \times r_1$
for some $r_1, r_2 \in R$.
\end{defn}

\begin{defn}
If $S$ is a subset of a ring $R$ that is itself a ring under
the same binary operations of addition and multiplication,
then $S$ is a
\index{subring} {\it subring}
of \index{ring!sub-} $R$.
\end{defn}

\begin{defn}
A ring $R$ is a
\index{ring!division} {\it division ring}
\index{division ring}
if every nonzero element $r \in R$ has a
multiplicative inverse $r^{-1}$. A
\index{field} {\it field}
is a commutative division ring.
\end{defn}

\subsection{Polynomial Rings}

\subsubsection{Commutative Polynomial Rings}

A nontrivial
\index{polynomial} {\it polynomial}
$p$ in $n$ (commuting) variables $x_1, \hdots, x_n$
is usually written as a sum
\begin{equation}
p = \sum_{i = 1}^k a_ix_1^{e^1_i}x_2^{e^2_i}\hdots x_n^{e^n_i},
\end{equation}
where $k$ is a positive integer and
each summand is a
\index{term} {\it term}
made up of a nonzero
\index{coefficient} {\it coefficient}
$a_i$ from some ring $R$ and a
\index{monomial} {\it monomial}
$x_1^{e^1_i}x_2^{e^2_i}\hdots x_n^{e^n_i}$ in which the
exponents $e^1_i, \hdots, e^n_i$ are nonnegative integers.
It is clear that each monomial may be represented in terms of its
exponents only, as a
\index{multidegree} {\it multidegree}
$e_i = (e^1_i, e^2_i, \hdots, e^n_i)$, so that a monomial
may be written as a multiset $\mathbf{x}^{e_i}$
over the set $\{x_1, \hdots, x_n\}$.
This leads to a more elegant representation of a
nontrivial polynomial,
\begin{equation} \label{polysum}
p = \sum_{\alpha \in \mathbb{N}^n} a_{\alpha} \mathbf{x}^{\alpha},
\end{equation}
and we may think of such a
polynomial as a function $f$ from the set of all
multidegrees $\mathbb{N}^n$ to the ring $R$ with finite
support (only a finite number of nonzero images).

\begin{example}
Let $p = 4x^2y + 2x + \frac{19}{80}$ be a polynomial in two variables
$x$ and $y$ with coefficients in $\mathbb{Q}$. This polynomial
can be represented by the function $f: \: \mathbb{N}^2 \rightarrow
\mathbb{Q}$ given by
$$f(\alpha) = \left\{
\begin{array}{lll}
4,             &  \hspace{4mm} & \alpha = (2, 1) \\
2,             && \alpha = (1, 0) \\
\frac{19}{80}, && \alpha = (0, 0) \\
0              && \mathrm{otherwise.}
\end{array} \right.$$
\end{example}

\begin{remark}
The zero polynomial $p = 0$ is represented by the function 
$f(\alpha) = 0_R$ for all possible $\alpha$. The constant
polynomial $p = 1$ is represented by the function
$f(\alpha) = 1_R$ for $\alpha = (0, 0, \hdots, 0)$, and
$f(\alpha) = 0_R$ otherwise.
\end{remark}

\begin{remark}
The product $m_1\times m_2$ of two monomials
$m_1, m_2$ with corresponding multidegrees $e_1, e_2
\in \mathbb{N}^n$ is the monomial corresponding to the
multidegree $e_1+e_2$.
For example, if $m_1 = x_1^2x_2x_3^3$ and
$m_2 = x_1x_2x_3^2$ (so that $e_1 = (2, 1, 3)$ and
$e_2 = (1, 1, 2)$), then $m_1\times m_2 = x_1^3x_2^2x_3^5$
as $e_1+e_2 = (3, 2, 5)$.
\end{remark}

\begin{defn}
Let $R[x_1, x_2, \hdots, x_n]$ denote the set of all functions
$f: \: \mathbb{N}^n \rightarrow R$ such that each function $f$
represents a polynomial in $n$ variables $x_1, \hdots, x_n$ with 
coefficients over a ring $R$. Given two functions $f, g \in
R[x_1, x_2, \hdots, x_n]$, let us define the functions $f+g$ and
$f\times g$ as follows.
\begin{alignat*}{3}
(f+g)(\alpha) & = f(\alpha) + g(\alpha) 
& \quad \mbox{for all}~\alpha \in \mathbb{N}^n; \\
(f\times g)(\alpha) & =  \sum_{\beta+\gamma = \alpha} f(\beta)\times g(\gamma)
& \quad \mbox{for all}~\alpha \in \mathbb{N}^n.
\end{alignat*}
Then the set
\index{$R$@$R[x_1, x_2, \hdots, x_n]$} $R[x_1, x_2, \hdots, x_n]$
becomes a ring, known as the 
\index{polynomial ring!commutative}
{\it polynomial ring in $n$ variables over $R$},
\index{commutative polynomial ring}
with the functions corresponding to the zero and 
constant polynomials being the respective zero and unit elements
of the ring.
\end{defn}

\begin{remark}
In $R[x_1, x_2, \hdots, x_n]$, $R$ is known
as the 
\index{ring!coefficient} {\it coefficient ring}.
\index{coefficient ring}
\end{remark}

%\begin{remark}
%It is easily verified that the set $R[x_1, x_2, \hdots, x_n]$
%is a commutative ring, with the polynomials $0$ and $1$ as zero
%and unit, since the functional sum and product of
%$f, g \in R[x_1, x_2, \hdots, x_n]$
%are the functions representing the polynomials given
%by adding and multiplying the two polynomials represented by
%$f$ and $g$ using the distributive laws.
%as long as we use the usual operations
%$$\sum_{i = 0}^n a_ix^i
% + \sum_{i = 0}^n b_ix^i = \sum_{i = 0}^n (a_i+b_i)x^i$$
%and
%$$\sum_{i = 0}^n a_ix^i \times \sum_{i = 0}^n b_ix^i =
%\sum_{i = 0}^{2n}\left( \sum_{j = 0}^n a_jb_{n-j}\right) x^i$$
%for adding and multiplying polynomials.
%\end{remark}

\subsubsection{Noncommutative Polynomial Rings}

A nontrivial {\it polynomial} \index{polynomial}
$p$ in $n$ noncommuting variables
$x_1, \hdots, x_n$ is usually written as a sum
\begin{equation} \label{polysumNC}
p = \sum_{i = 1}^k a_iw_i,
\end{equation}
where $k$ is a positive integer and each summand is a
{\it term} \index{term}
made up of a nonzero {\it coefficient} \index{coefficient}
$a_i$ from some ring $R$ and a {\it monomial} \index{monomial}
$w_i$ that is a
word over the alphabet $X = \{x_1, x_2, \hdots, x_n\}$.
We may think of a noncommutative polynomial
as a function $f$ from the set of all words $X^{\ast}$
to the ring $R$.

\begin{remark}
The zero polynomial $p = 0$ is the polynomial $0_R\varepsilon$,
where $\varepsilon$ is the empty word in $X^{\ast}$. Similarly
$1_R\varepsilon$ is the constant polynomial $p = 1$.
\end{remark}

\begin{remark}
The product $w_1\times w_2$ of two monomials
$w_1, w_2 \in X^{\ast}$ is given by concatenation.
For example, if $X = \{x_1, x_2, x_3\}$, $w_1 = x_3^2x_2$ and
$w_2 = x_1^3x_3$, then $w_1\times w_2 = x_3^2x_2x_1^3x_3$.
\end{remark}

\begin{defn}
Let $R\langle x_1, x_2, \hdots, x_n\rangle$ denote the set of all functions
$f: \: X^{\ast} \rightarrow R$ such that each function $f$
represents a polynomial in $n$ noncommuting variables with coefficients
over a ring $R$. Given two functions $f, g \in
R\langle x_1, x_2, \hdots, x_n\rangle$, let us define the
functions $f+g$ and $f\times g$ as follows.
\begin{alignat*}{3}
(f+g)(w) & = f(w) + g(w) 
& \quad \mbox{for all}~w \in X^{\ast}; \\
(f\times g)(w) & =  \sum_{u\times v = w} f(u)\times g(v)
& \quad \mbox{for all}~w \in X^{\ast}.
\end{alignat*}
Then the set
\index{$R$@$R\langle x_1, x_2, \hdots, x_n\rangle$} 
$R\langle x_1, x_2, \hdots, x_n\rangle$
becomes a ring, known as the
\index{polynomial ring!noncommutative}
{\it noncommutative polynomial ring in $n$ variables over $R$},
\index{noncommutative polynomial ring}
with the functions corresponding to the zero and
constant polynomials being the respective zero and unit elements
of the ring.
\end{defn}

\subsection{Ideals}

\begin{defn}
Let $\mathcal{R}$ be an arbitrary commutative ring. An
\index{ideal} {\it ideal}
$J$ in $\mathcal{R}$ is a subring of $\mathcal{R}$ satisfying the following
additional condition: $jr \in J$ for all $j \in J$, $r \in \mathcal{R}$.
\end{defn}

\begin{remark}
In the above definition, if $\mathcal{R}$ is a polynomial ring in $n$ 
variables over a ring $R$ ($\mathcal{R} = R[x_1, \hdots, x_n]$),
the ideal $J$ is a {\it polynomial ideal}. We will only consider
polynomial ideals in this thesis.
\end{remark}

\begin{defn}
Let $\mathcal{R}$ be an arbitrary noncommutative ring.
\begin{itemize}
\item
A
\index{ideal!left} {\it left (right) ideal} 
\index{ideal!right} $J$ in $\mathcal{R}$
is \index{left ideal} a \index{right ideal} subring
of $\mathcal{R}$ satisfying the following additional condition:
$rj \in J$ ($jr \in J$) for all $j \in J$, $r \in \mathcal{R}$.
\item
A
\index{ideal!two-sided} {\it two-sided ideal} $J$ in $\mathcal{R}$ is a
\index{two-sided ideal}
subring of $\mathcal{R}$ satisfying the following additional condition:
$r_1jr_2 \in J$ for all $j \in J$, $r_1, r_2 \in \mathcal{R}$.
\end{itemize}
\end{defn}

\begin{remark}
Unless otherwise stated, all noncommutative ideals
considered in this thesis will be two-sided ideals.
\end{remark}

\begin{defn}
A set of polynomials $P = \{p_1, p_2, \hdots, p_m\}$
is a 
\index{basis} {\it basis} for an ideal $J$ of a
noncommutative polynomial ring $\mathcal{R}$ if
every polynomial $q \in J$ can be written as
% $J$ consists of all polynomials $q$ of the form
\begin{equation}\label{combination}
q = \sum_{i = 1}^k \ell_i p_i r_i \; \;
(\ell_i, r_i \in \mathcal{R}, \; p_i \in P).
\end{equation}
We say that $P$ generates $J$, written
$J = \langle P \rangle$.
\end{defn}

\begin{remark}
The above definition has an obvious generalisation
for left and right ideals of noncommutative polynomial rings and
for ideals of commutative polynomial rings.
\end{remark}

\begin{example}
Let $\mathcal{R}$ be the noncommutative polynomial ring
$\mathbb{Q}\langle x, y \rangle$, and let
$J = \langle P \rangle$ be an ideal in $\mathcal{R}$, where $P :=
\{x^2y + yx - 2, \, yxy - x + 4y\}$. Consider the
polynomial $q := 2x^3y + yx^2y + 2xyx - 4x^2y + x^3 - 2xy - 4x$,
and let us ask if $q$ is a member of the ideal. To answer this question, 
we have to find out if there is an expression for $q$ of the type shown 
in Equation (\ref{combination}).
In this case, it turns out that $q$ is indeed a member of the ideal
(because $q = 2x(x^2y+yx-2) + (x^2y+yx-2)xy - x^2(yxy-x+4y)$), 
but how would we answer the question in general?  
This problem is known as the
Ideal Membership Problem and is stated as follows.
\end{example}

\begin{defn}[The Ideal Membership Problem]
\index{ideal membership problem}
Given an ideal $J$ and a polynomial $q$, does $q \in J$?
\end{defn}

As we shall see shortly, the Ideal Membership Problem can
be solved by dividing a polynomial with respect to a
Gr\"obner Basis for the ideal $J$. But before
we can discuss this, we must first introduce the
notion of polynomial division, for which we require a
fixed ordering on the monomials in any given polynomial.

\section{Monomial Orderings}

A \index{monomial ordering} {\it monomial ordering}
\index{ordering!monomial}
is a bivariate function $\mathrm{O}$ which tells us
which monomial is the larger of any two
given monomials $m_1$ and $m_2$. We
will use the convention that
\index{$<$}
$\mathrm{O}(m_1, m_2) = 1$ if and only if $m_1 < m_2$, and
$\mathrm{O}(m_1, m_2) = 0$ if and only if $m_1 \geqslant m_2$.
We can use a monomial ordering to order an arbitrary
polynomial $p$ by inducing an order on the terms of $p$
from the order on the monomials associated with the terms.

\begin{defn}
A monomial ordering $\mathrm{O}$ is
\index{monomial ordering!admissible} {\it admissible} if the following
\index{admissible monomial ordering}
conditions are satisfied.
\begin{enumerate}[(a)]
\item
$1 < m$ for all monomials $m \neq 1$.
\item
$m_1 < m_2 \Rightarrow m_{\ell}m_1m_r
< m_{\ell}m_2m_r$ for all monomials\footnote{For a
commutative monomial ordering, we can ignore the
monomial $m_r$.}
$m_1, m_2, m_{\ell}, m_r$.
\end{enumerate}
\end{defn}

By convention, a polynomial is always written in
descending order (with respect to a given monomial ordering), so that the
\index{leading term} {\it leading term}
\index{term!leading}
of the polynomial (with associated
\index{leading coefficient} {\it leading coefficient}
\index{coefficient!leading} and
\index{leading monomial} {\it leading monomial})
\index{monomial!leading}
always comes first.

\begin{remark}
For an arbitrary polynomial $p$, we will use
\index{$LT$@$\mathrm{LT}(p)$}
$\mathrm{LT}(p)$, 
\index{$LM$@$\mathrm{LM}(p)$}
$\mathrm{LM}(p)$ and 
\index{$LC$@$\mathrm{LC}(p)$}
$\mathrm{LC}(p)$
to denote the leading term, leading monomial and
leading coefficient of $p$ respectively.
\end{remark}

\subsection{Commutative Monomial Orderings} \label{CMO}
\index{monomial ordering}

A monomial ordering 
\index{ordering!monomial} usually requires an ordering on the
variables in our chosen polynomial ring. Given such a
ring $R[x_1, x_2, \hdots, x_n]$, we will assume this order to be
$x_1 > x_2 > \cdots > x_n$.

We shall now consider the most frequently used monomial orderings, where
throughout $m_1$ and $m_2$ will denote arbitrary monomials
(with associated multidegrees $e_1 = (e^1_1, e^2_1, \hdots, e^n_1)$
and $e_2 = (e^1_2, e^2_2, \hdots, e^n_2)$),
\index{degree}
and $\deg(m_i)$ will denote the total degree of the monomial $m_i$
(for example $\deg(x^2yz) = 4$).
All orderings considered will be admissible.

\subsubsection{The Lexicographical Ordering ($\mathrm{Lex}$)}
\index{monomial ordering!lexicographic}

Define $m_1 < m_2$ if $e^i_1 < e^i_2$ for some
$1 \leqslant i \leqslant n$ and $e^j_1 = e^j_2$ for all
$1 \leqslant j < i$. In words, $m_1 < m_2$ if the first
variable with different exponents in $m_1$ and $m_2$
has lower exponent in $m_1$.
\index{lexicographic ordering}

\subsubsection{The Inverse Lexicographical Ordering ($\mathrm{InvLex}$)}
\index{monomial ordering!inverse lexicographic}

Define $m_1 < m_2$ if $e^i_1 < e^i_2$ for some
$1 \leqslant i \leqslant n$ and $e^j_1 = e^j_2$ for all
$i < j \leqslant n$. In words, $m_1 < m_2$ if the last
variable with different exponents in $m_1$ and $m_2$
has lower exponent in $m_1$.
\index{inverse lexicographic ordering}

% \subsubsection{The Reverse Lexicographical Ordering ($\mathrm{RevLex}$)}
% \index{monomial ordering!reverse lexicographic}
% 
% Define $m_1 < m_2$ if $e^i_1 > e^i_2$ for some
% $1 \leqslant i \leqslant n$ and $e^j_1 = e^j_2$ for all
% $i < j \leqslant n$. In words, $m_1 < m_2$ if the last
% variable with different exponents in $m_1$ and $m_2$
% has higher exponent in $m_1$.

\subsubsection{The Degree Lexicographical Ordering ($\mathrm{DegLex}$)}
\index{monomial ordering!degree lexicographic}

Define $m_1 < m_2$ if $\deg(m_1) < \deg(m_2)$ or if $\deg(m_1) = \deg(m_2)$
and $m_1 < m_2$ in the Lexicographic Ordering.
\index{degree lexicographic ordering}

\begin{remark}
The $\mathrm{DegLex}$ ordering is also known as the $\mathrm{TLex}$ ordering
(T for total degree).
\end{remark}

\subsubsection{The Degree Inverse Lexicographical Ordering
               ($\mathrm{DegInvLex}$)}
\index{monomial ordering!degree inverse lexicographic}

Define $m_1 < m_2$ if $\deg(m_1) < \deg(m_2)$ or if $\deg(m_1) = \deg(m_2)$
and $m_1 < m_2$ in the Inverse Lexicographical Ordering.
\index{degree inverse lexicographic ordering}

\subsubsection{The Degree Reverse Lexicographical Ordering
               ($\mathrm{DegRevLex}$)}
\index{monomial ordering!degree reverse lexicographic}

Define $m_1 < m_2$ if $\deg(m_1) < \deg(m_2)$ or if $\deg(m_1) = \deg(m_2)$
and $m_1 < m_2$ in the Reverse Lexicographical Ordering, where
$m_1 < m_2$ if the last
variable with different exponents in $m_1$ and $m_2$
has higher exponent in $m_1$ ($e^i_1 > e^i_2$ for some
$1 \leqslant i \leqslant n$ and $e^j_1 = e^j_2$ for all
$i < j \leqslant n$).
\index{degree reverse lex. ordering}

\begin{remark}
On its own, the Reverse Lexicographical Ordering (RevLex)
is not admissible, as $1 > m$ for any monomial $m \neq 1$.
\end{remark}

\begin{example}
With $x > y > z$, consider the monomials 
$m_1 := x^2yz$; $m_2 := x^2$ and $m_3 := xyz^2$,
with corresponding multidegrees $e_1 = (2, 1, 1)$; $e_2 = (2, 0, 0)$
and $e_3 = (1, 1, 2)$. The following table shows the order placed on
the monomials by the various monomial orderings defined above.
The final column shows the order induced on the polynomial
$p := m_1 + m_2 + m_3$ by the chosen monomial ordering.
\begin{center}
\begin{tabular}{l|c|c|c|c}
Monomial Ordering $\mathrm{O}$ & $\mathrm{O}(m_1, m_2)$
& $\mathrm{O}(m_1, m_3)$ & $\mathrm{O}(m_2, m_3)$ & $p$ \\
\hline
$\mathrm{Lex}$       & 0 & 0 & 0 & $x^2yz + x^2 + xyz^2$ \\
$\mathrm{InvLex}$    & 0 & 1 & 1 & $xyz^2 + x^2yz + x^2$ \\
% $\mathrm{RevLex}$    & 1 & 0 & 0 & $x^2 + x^2yz + xyz^2$ \\
$\mathrm{DegLex}$    & 0 & 0 & 1 & $x^2yz + xyz^2 + x^2$ \\
$\mathrm{DegInvLex}$ & 0 & 1 & 1 & $xyz^2 + x^2yz + x^2$ \\
$\mathrm{DegRevLex}$ & 0 & 0 & 1 & $x^2yz + xyz^2 + x^2$ \\ \hline
\end{tabular}
\end{center}
\end{example}

\subsection{Noncommutative Monomial Orderings} \label{NCMO}

In the noncommutative case, because we use words and not
multidegrees to represent monomials, our definitions
for the lexicographically based orderings will have
to be adapted slightly. All other definitions and
conventions will stay the same.

\subsubsection{The Lexicographic Ordering ($\mathrm{Lex}$)}
\index{monomial ordering!lexicographic}

Define $m_1 < m_2$ if, working
left-to-right, the first (say $i$-th) letter on which $m_1$ and $m_2$ differ
is such that the $i$-th letter of $m_1$ is lexicographically
{\it less} than the $i$-th letter of $m_2$ in the variable ordering.
Note: this ordering is {\it not} admissible
(counterexample: if $x > y$ is the variable ordering,
then $x < xy$ but $x^2 > xyx$).
\index{lexicographic ordering}

\begin{remark}
When comparing two monomials $m_1$ and $m_2$
such that $m_1$ is a proper prefix of
$m_2$ (for example $m_1 := x$ and $m_2 := xy$ as in
the above counterexample), a
problem arises with the above definition in that we
eventually run out of letters in the shorter word to compare with
(in the example, having seen that the first letter of both
monomials match, what do we compare the second letter
of $m_2$ with?). One answer is to introduce a
\index{padding symbol}
padding symbol $\$ $ to pad $m_1$ on the right to make sure it is the
same length as $m_2$, with the convention that any
letter is greater than the padding symbol (so that
$m_1 < m_2$). % In the example,
% we find that $\$ < y$ and hence $x < xy$). 
The padding symbol will not explicitly appear 
anywhere in the remainder of this thesis,
but we will bear in mind that it can be introduced to
deal with situations where prefixes and suffixes of
monomials are involved.
\end{remark}

\begin{remark}
The lexicographic ordering is also known as the dictionary
ordering since the words in a dictionary (such as the Oxford
English Dictionary) are ordered using the lexicographic
ordering with variable (or alphabetical) ordering
$a < b < c < \cdots$.
Note however that while a dictionary orders words in 
increasing order, we will write polynomials in decreasing order.
\end{remark}

\subsubsection{The Inverse Lexicographical Ordering ($\mathrm{InvLex}$)}
\index{monomial ordering!inverse lexicographic}

Define $m_1 < m_2$ if, working
left-to-right, the first (say $i$-th) letter on which $m_1$ and $m_2$ differ
is such that the $i$-th letter of $m_1$ is lexicographically
{\it greater} than the $i$-th letter of $m_2$.
Note: this ordering (like $\mathrm{Lex}$) is {\it not} admissible
(counterexample: if $x > y$ is the variable ordering,
then $xy < x$ but $xyx > x^2$).
\index{inverse lexicographic ordering}

\subsubsection{The Degree Reverse Lexicographical Ordering
($\mathrm{DegRevLex}$)}
\index{monomial ordering!degree reverse lexicographic}

Define $m_1 < m_2$ if $\deg(m_1) < \deg(m_2)$ or if $\deg(m_1) = \deg(m_2)$
and $m_1 < m_2$ in the Reverse Lexicographical Ordering, where
$m_1 < m_2$ if, working in {\it reverse}, or from right-to-left,
the first (say $i$-th) letter on which $m_1$ and $m_2$ differ
is such that the $i$-th letter of $m_1$ is lexicographically
{\it greater} than the $i$-th letter of $m_2$.
\index{degree reverse lex. ordering}

\begin{example}
With $x > y > z$,
consider the noncommutative monomials
$m_1 := zxyx$; $m_2 := xzx$ and $m_3 := y^2zx$.
The following table shows the order placed on
the monomials by various noncommutative monomial orderings. As before,
the final column shows the order induced on the polynomial
$p := m_1 + m_2 + m_3$ by the chosen monomial ordering.
\begin{center}
\begin{tabular}{l|c|c|c|c}
Monomial Ordering $\mathrm{O}$ & $\mathrm{O}(m_1, m_2)$
& $\mathrm{O}(m_1, m_3)$ & $\mathrm{O}(m_2, m_3)$ & $p$ \\
\hline
$\mathrm{Lex}$       & 1 & 1 & 0 & $xzx + y^2zx + zxyx$ \\
$\mathrm{InvLex}$    & 0 & 0 & 1 & $zxyx + y^2zx + xzx$ \\
% $\mathrm{RevLex}$    & 1 & 1 & 1 & $y^2zx + xzx + zxyx$ \\
$\mathrm{DegLex}$    & 0 & 1 & 1 & $y^2zx + zxyx + xzx$ \\
$\mathrm{DegInvLex}$ & 0 & 0 & 1 & $zxyx + y^2zx + xzx$ \\
$\mathrm{DegRevLex}$ & 0 & 1 & 1 & $y^2zx + zxyx + xzx$ \\ \hline
\end{tabular}
\end{center}
\end{example}

% \begin{example}
% Suppose that we wish to order the terms in the polynomial 
% $zxyx + x^2 + zyzx$
% according to Lex, DegLex and DegRevLex (with variable
% order $x > y > z$). Using Lex, our polynomial is
% $p = x^2 + zxyx + zyzx$ because $x^2$ appears in the dictionary
% before $zxyx$ (for the first letter, $x > z$ 
% in the alphabet) and because $zxyx$
% appears in the dictionary before $zyzx$ 
% (for the second letter, $x > y$ in the
% alphabet). Using DegLex, our polynomial is $p = zxyx + zyzx + x^2$ 
% because $x^2$ has the smallest degree so comes last, and then
% $zxyx$ comes before $zyzx$ as it is lexicographically bigger. Finally,
% DegRevLex, our polynomial is $p = zyzx + zxyx + x^2$ where $x^2$ again 
% comes last because it has the smallest
% degree, but this time $zyzx$ comes before $zxyx$ because, working from
% right to left, the first letter on which the two words differ is the
% second-last letter, and as $z < y$ in the alphabet this means that 
% $zyzx > zxyx$ in DegRevLex.
% \end{example}

\subsection{Polynomial Division}
\index{polynomial division}

\begin{defn} \label{polydiv} \index{division!polynomial}
Let $\mathcal{R}$ be a polynomial ring, and let $O$ be an
arbitrary admissible monomial ordering.
Given two nonzero polynomials $p_1$, $p_2 \in \mathcal{R}$,
we say that $p_1$ divides $p_2$ (written $p_1 \mid p_2$) if
the lead monomial of $p_1$
divides some monomial $m$ (with coefficient $c$) in $p_2$. For a
commutative polynomial ring, this means that $m = \LM(p_1)m'$
for some monomial $m'$; for a noncommutative polynomial ring, this
means that $m = m_{\ell}\LM(p_1)m_r$ for some
monomials $m_{\ell}$ and $m_r$
($\LM(p_1)$ is a subword of $m$).

To perform the division, we take away an appropriate multiple
of $p_1$ from $p_2$ in order to cancel off $\LT(p_1)$ with
the term involving $m$ in $p_2$.
In the commutative case, we do
$$p_2 - (c\LC(p_1)^{-1})p_1m';$$
in the noncommutative case, we do
$$p_2 - (c\LC(p_1)^{-1})m_{\ell}p_1m_r.$$
It is clear that the coefficient rings of our polynomial
rings have to be division rings in order for the above
expressions to be valid, and so we make the following
assumption about the polynomial rings we will encounter
in the remainder of this thesis.
\end{defn}

\begin{remark}
From now on, all coefficient rings of polynomial rings
will be fields unless otherwise stated.
\end{remark}

\begin{example}
Let $p_1 := 5z^2x + 2y^2 + x + 4$ and 
$p_2 := 3xyxz^2x^3 + 2x^2$
be two $\mathrm{DegLex}$ ordered polynomials
over the noncommutative polynomial ring
$\mathbb{Q}\langle x, y, z\rangle$.
Because $\LM(p_2) = xyx(z^2x)x^2$, it is
clear that $p_1 \mid p_2$, with the quotient and the
remainder of the division being
$$\begin{array}{c}
q := \left(\frac{3}{5}\right)xyx(5z^2x + 2y^2 + x + 4)x^2
\end{array}
$$
and
$$
\begin{array}{ccl}
r & := & 3xyxz^2x^3 + 2x^2
        - \left(\frac{3}{5}\right)xyx(5z^2x + 2y^2 + x + 4)x^2 \\
  & = & 3xyxz^2x^3 + 2x^2 - 3xyxz^2x^3 - \left(\frac{6}{5}\right)xyxy^2x^2
        - \left(\frac{3}{5}\right)xyx^4 - \left(\frac{12}{5}\right)xyx^3 \\
  & = & - \left(\frac{6}{5}\right)xyxy^2x^2 - \left(\frac{3}{5}\right)xyx^4
        - \left(\frac{12}{5}\right)xyx^3 + 2x^2
\end{array}
$$
respectively.
\end{example}

Now that we know how to divide one polynomial by another,
what does it mean for a polynomial to be divided
by a set of polynomials?

\begin{defn}
Let $\mathcal{R}$ be a polynomial ring, and let $O$ be an
arbitrary admissible monomial ordering.
Given a nonzero polynomial $p \in \mathcal{R}$
and a set of nonzero polynomials $P = \{p_1, p_2,
\hdots, p_m\}$, with $p_i \in \mathcal{R}$ for
all $1 \leqslant i \leqslant m$,
we divide $p$ by $P$ by working through $p$ term by term,
testing to see if each term is divisible by any of the
$p_i$ in turn. We recursively divide
the remainder of each division using the same method
until no more divisions are possible, in which case the 
remainder is either $0$ or is {\it irreducible}.
\end{defn}

% MAKE SURE ALGORITHM 1 APPEARS AFTER THIS PARAGRAPH

Algorithms to divide a polynomial $p$ by a set of polynomials $P$
in the commutative and noncommutative cases are given below
as Algorithms \ref{com-div} and \ref{noncom-div} respectively. 
Note that they take advantage of
the fact that if the first $N$ terms of a polynomial $q$
are irreducible with respect to $P$, then the first $N$ terms
of any reduction of $q$ will also be irreducible with
respect to $P$.

\begin{algorithm}
\setlength{\baselineskip}{3.5ex}
\caption{The Commutative Division Algorithm}
\label{com-div}
\begin{algorithmic}
\vspace*{2mm}
\REQUIRE{A nonzero polynomial $p$ and a set of nonzero polynomials
         $P = \{p_1, \hdots, p_m\}$ over a 
         polynomial ring $R[x_1, \hdots x_n]$;
         an admissible monomial ordering $\mathrm{O}$.}
\ENSURE{$\Rem(p, P) := r$, the remainder of $p$ with respect to $P$.}
\vspace*{1mm}
\STATE
$r = 0$;
\WHILE{($p \neq 0$)}
\STATE
$u = \LM(p)$; $c = \LC(p)$; $j = 1$; found = false; \\
\WHILE{($j \leqslant m$) \textbf{and} (found == false)}
\IF{($\LM(p_j) \mid u$)}
\STATE
found = true; $u' = u/\LM(p_j)$;
$p = p - (c\LC(p_j)^{-1})p_ju'$;
\ELSE
\STATE
$j = j+1$;
\ENDIF
\ENDWHILE
\IF{(found == false)}
\STATE
$r = r+\LT(p)$; $p = p-\LT(p)$;
\ENDIF
\ENDWHILE
\STATE
{\bf return} $r$;
\end{algorithmic}
\vspace*{1mm}
\end{algorithm}

\begin{algorithm}
\setlength{\baselineskip}{3.5ex}
\caption{The Noncommutative Division Algorithm}
\label{noncom-div}
\vspace*{2mm}
To divide a nonzero polynomial $p$ with respect to a
set of nonzero polynomials $P = \{p_1, \hdots, p_m\}$, where $p$
and the $p_i$ are elements of a
noncommutative polynomial ring $R\langle x_1, \hdots, x_n\rangle$,
we apply Algorithm \ref{com-div} with the following changes.
\begin{enumerate}[(a)]
\item
In the inputs, replace the commutative polynomial ring $R[x_1, \hdots x_n]$
by the noncommutative polynomial ring $R\langle x_1, \hdots, x_n\rangle$.
\item
Change the first \textbf{if} condition to read
\begin{algorithmic}
\IF{($\LM(p_j) \mid u$)}
\STATE
found = true; \\
choose $u_{\ell}$ and $u_r$ such that $u = u_{\ell}\LM(p_j)u_r$; \\
$p = p - (c\LC(p_j)^{-1})u_{\ell}p_ju_r$;
\ELSE
\STATE
$j = j+1$;
\ENDIF
\end{algorithmic}
\end{enumerate}
%\vspace*{1mm}
\end{algorithm}

\begin{remark}
All algorithms in this thesis use the conventions
that `=' denotes an assignment and `==' denotes a test.
\end{remark}

\newpage % FINAL VERSION ONLY!!

\begin{remark}
In Algorithm \ref{noncom-div},
if there are several candidates for $u_{\ell}$ (and therefore for $u_r$)
in the line `choose $u_{\ell}$ and $u_r$ such that $u = u_{\ell}\LM(p_j)u_r$',
the convention in this thesis will be
to choose the $u_{\ell}$ with the smallest degree.
\end{remark}

\begin{example} \label{nur}
To demonstrate that the process of dividing a polynomial by a set
of polynomials does not necessarily give a unique result,
consider the polynomial $p := xyz + x$ and the set of polynomials
$P := \{p_1, \: p_2\} = \{xy - z, \: yz + 2x + z\}$,
all polynomials being ordered by
$\mathrm{DegLex}$ and originating from the polynomial ring
$\mathbb{Q}[x, y, z]$. If we choose to divide $p$
by $p_1$ to begin with, we see that $p$ reduces to
$xyz + x - (xy - z)z = z^2 + x$, which is irreducible. But if
we choose to divide $p$ by $p_2$ to begin with, we see that
$p$ reduces to $xyz + x - (yz + 2x + z)x = -2x^2 - xz + x$,
which is again irreducible. This gives rise to the question of which
answer (if any!) is the correct one here? In Chapter \ref{ChCGB}, we will
discover that one way of obtaining a unique answer to this question
will be to calculate a {\it Gr\"obner Basis} for the dividing set $P$.
\end{example}

\begin{defn} \label{arrows}
% In the remainder of this thesis,
In order to describe how one polynomial is obtained from
another through the process of division, we introduce the
following notation.
\vspace*{-2mm}
\begin{enumerate}[(a)]
\item
If the polynomial $r$ is obtained by dividing a polynomial
$p$ by a polynomial $q$, then we will use the notation
% \index{$\rightarrow$}
\index{$<$@\hspace*{-2mm}$\xymatrix{\ar[r] &}$\hspace*{-1mm}}
$p \rightarrow r$ or $p \rightarrow_q r$ (with the latter notation
used if we wish to show how $r$ is obtained from $p$).
\item
If the polynomial $r$ is obtained by dividing a polynomial
$p$ by a sequence of polynomials $q_1, q_2, \hdots, q_{\alpha}$,
then we will use the notation
% \index{$\stackrel{\ast}{\longrightarrow}$}
\index{$<$@\hspace*{-2mm}$\xymatrix{\ar[r]^{*} &}$\hspace*{-1mm}}
% $p \rightarrow_{q_1} \hdots \rightarrow_{q_{\alpha}} r$ or
$p \stackrel{\ast}{\longrightarrow} r$.
\item
If the polynomial $r$ is obtained by dividing a polynomial
$p$ by a set of polynomials $Q$, then we will use the
notation $p \rightarrow_Q r$.
\end{enumerate}
\end{defn}

%
% Chapter 2
% Author: Gareth Evans
% Last Modified: 3rd February 2006
%

\chapter{Commutative Gr\"obner Bases} \label{ChCGB}
% \thispagestyle{empty}

% In this chapter, we will introduce the
% theory and algorithms required in order to
% define and compute commutative Gr\"obner Bases.

% \typeout{INSERT PARAGRAPH ON WHAT I DO IN THIS CHAPTER.}

Given a basis $F$ generating an ideal $J$,
the central idea in Gr\"obner Basis theory is to
use $F$ to find a basis $G$ for $J$
with the property that the remainder of the division
of any polynomial by $G$ is unique. Such a basis
is known as a
\index{Gr\"obner basis} {\it Gr\"obner Basis}.
\index{basis!Gr\"obner}

In particular, if a polynomial $p$ is a member of the ideal $J$,
then the remainder of the division of $p$
by a Gr\"obner Basis $G$ for $J$ is always zero. This gives us a way to
solve the Ideal Membership Problem for $J$ -- if the remainder of the
division of a polynomial $p$ by $G$ is zero,
then $p \in J$ (otherwise $p \notin J$).

\section{S-polynomials} \label{SC3}

How do we determine whether or not an arbitrary basis $F$ generating
an ideal $J$ is a Gr\"obner Basis? Using the informal definition
shown above, in order to show that a basis is {\it not} a Gr\"obner Basis,
it is sufficient to find a polynomial $p$ whose remainder on division
by $F$ is non-unique. Let us now construct an example in which this
is the case, and let us analyse what can to be done to eliminate the
non-uniqueness of the remainder.

Let $p_1 = a_1 + a_2 + \cdots + a_{\alpha}$;
$p_2 = b_1 + b_2 + \cdots + b_{\beta}$ and
$p_3 = c_1 + c_2 + \cdots + c_{\gamma}$ be three
polynomials ordered with respect to some fixed admissible monomial
ordering $O$ (the $a_i$, $b_j$ and $c_k$ are
all nontrivial terms). Assume that $p_1 \mid p_3$ and $p_2 \mid p_3$,
so that we are able to take away from $p_3$ 
multiples $s$ and $t$ of $p_1$ and $p_2$ respectively to obtain
remainders $r_1$ and $r_2$.
\begin{eqnarray*}
r_1 & = & p_3 - sp_1 \\
    & = & c_1 + c_2 + \cdots + c_{\gamma} 
          - s(a_1 + a_2 + \cdots + a_{\alpha}) \\
    & = & c_2 + \cdots + c_{\gamma} - sa_2 - \cdots - sa_{\alpha}; \\
r_2 & = & p_3 - tp_2 \\
    & = & c_2 + \cdots + c_{\gamma} - tb_2 - \cdots - tb_{\beta}.
\end{eqnarray*}
If we assume that $r_1$ and $r_2$ are irreducible
and that $r_1 \neq r_2$, it is clear that the
remainder of the division of the polynomial
$p_3$ by the set of polynomials $P = \{p_1, p_2\}$ is
non-unique, from which we deduce that $P$ is not a
Gr\"obner Basis for the ideal that it generates. We must
therefore change $P$ in some way in order for it to
become a Gr\"obner Basis, but what changes are required
and indeed allowed?

Consider that we want to add a polynomial
to $P$. To avoid changing the ideal that is being
generated by $P$, any polynomial added to $P$ must be a
member of the ideal. It is clear that $r_1$ and $r_2$ are
members of the ideal, as is the polynomial
$p_4 = r_2 - r_1 = -tp_2 + sp_1$. Consider that we
add $p_4$ to $P$, so that $P$ becomes the set
$$
\{a_1 + a_2 + \cdots + a_{\alpha}, \;
b_1 + b_2 + \cdots + b_{\beta}, \;
-tb_2 - tb_3 - \cdots - tb_{\beta} + sa_2 + sa_3 + \cdots + sa_{\alpha}\}.$$
If we now divide the polynomial $p_3$ by the
enlarged set $P$, to begin with (as before)
we can either divide $p_3$ by $p_1$ or $p_2$
to obtain remainders $r_1$ or $r_2$.
Here however, if we assume
(without loss of generality\footnote{The other
possible case is $\LT(p_4) = sa_2$, in which case it is $r_1$
that reduces to $r_2$ and not $r_2$ to $r_1$.})
that $\LT(p_4) = -tb_2$, we can now divide $r_2$
by $p_4$ to obtain a new remainder
\begin{eqnarray*}
r_3 & = & r_2 - p_4 \\
    & = & c_2 + \cdots + c_{\gamma} - tb_2 - \cdots - tb_{\beta}
          -(-tb_2 - tb_3 - \cdots - tb_{\beta}
            + sa_2 + sa_3 + \cdots + sa_{\alpha}) \\
    & = & c_2 + \cdots + c_{\gamma} - sa_2 - \cdots - sa_{\alpha} \\
    & = & r_1.
\end{eqnarray*}
It follows that
by adding $p_4$ to $P$, we have ensured that the
remainder of the division of $p_3$ by $P$ is unique\footnote{This may
not strictly be true if $p_3$ is divisible by $p_4$; for the time
being we shall assume that this is not the case, noting that
the important concept here is of eliminating the non-uniqueness
given by the choice of dividing $p_3$ by $p_1$ or $p_2$ first.}
no matter which of the polynomials $p_1$ and $p_2$ we choose
to divide $p_3$ by first. This solves our original problem of non-unique
remainders in this restricted situation.

At first glance, the polynomial added to $P$ to solve
this problem is dependent upon the polynomial $p_3$.
The reason for saying this is that
the polynomial added to $P$ has the form $p_4 = sp_1 - tp_2$, where
$s$ and $t$ are terms chosen to multiply the
polynomials $p_1$ and $p_2$ so that the lead terms
of $sp_1$ and $tp_2$ equal $\LT(p_3)$
(in fact $s = \frac{\LT(p_3)}{\LT(p_1)}$ and
$t = \frac{\LT(p_3)}{\LT(p_2)}$).

However, by definition, $\LM(p_3)$ is
a common multiple of $\LM(p_1)$ and $\LM(p_2)$.
Because all such common multiples are multiples
of the least common multiple of $\LM(p_1)$ and $\LM(p_2)$
(so that $\LM(p_3) = \mu(\lcm(\LM(p_1), \LM(p_2)))$
for some monomial $\mu$),
it follows that we can rewrite $p_4$ as
$$p_4 = \LC(p_3)\mu\left(\frac{\lcm(\LM(p_1), \LM(p_2))}{\LT(p_1)}p_1 -
\frac{\lcm(\LM(p_1), \LM(p_2))}{\LT(p_2)}p_2\right).$$
Consider now that we add the polynomial 
$p_5 = \frac{p_4}{\LC(p_3)\mu}$ to $P$
instead of adding $p_4$ to $P$. It follows that even though
this polynomial does not depend on the polynomial $p_3$,
we can still obtain a unique remainder when dividing $p_3$ by
$p_1$ and $p_2$, because we can do $r_3 = r_2 - \LC(p_3)\mu p_5$.
Moreover, the polynomial $p_5$ solves the problem of non-unique
remainders for {\it any} polynomial $p_3$ that is divisible by
both $p_1$ and $p_2$ (all that changes is the multiple of $p_5$
used in the reduction of $r_2$); we call such a polynomial
an {\it S-polynomial}\footnote{The
S stands for Syzygy, as in a pair of connected objects.} for
$p_1$ and $p_2$.

\begin{defn}
The \index{S-polynomial} {\it S-polynomial}
of two distinct polynomials $p_1$ and
$p_2$ is given by the expression
$$\mathrm{S\mbox{-}pol}(p_1, p_2) =
\frac{\lcm(\LM(p_1), \LM(p_2))}{\LT(p_1)}p_1 -
\frac{\lcm(\LM(p_1), \LM(p_2))}{\LT(p_2)}p_2.$$
\end{defn}

\begin{remark}
The terms $\frac{\lcm(\LM(p_1), \LM(p_2))}{\LT(p_1)}$
and $\frac{\lcm(\LM(p_1), \LM(p_2))}{\LT(p_2)}$ can be
thought of as the terms used to multiply the
polynomials $p_1$ and $p_2$ so that the lead monomials
of the multiples are equal to the
monomial $\lcm(\LM(p_1), \LM(p_2))$.
\end{remark}

Let us now illustrate how adding an S-polynomial to a basis
solves the problem of non-unique remainders in a
particular example.

\begin{example}
Recall that in Example \ref{nur} we showed how dividing the
polynomial $p := xyz + x$ by the two polynomials in the set
$P := \{p_1, \: p_2\} = \{xy - z, \: yz + 2x + z\}$ gave two different
remainders, $r_1 := z^2 + x$ and $r_2 := -2x^2 - xz + x$ respectively.
Consider now that we add $\mathrm{S\mbox{-}pol}(p_1, p_2)$ to
$P$, where
\begin{eqnarray*}
\mathrm{S\mbox{-}pol}(p_1, p_2) & = & \frac{xyz}{xy}(xy-z)
- \frac{xyz}{yz}(yz+2x+z) \\
& = & (xyz - z^2) - (xyz+2x^2+xz) \\
& = & -2x^2-xz-z^2.
\end{eqnarray*}
Dividing $p$ by the enlarged set, if we choose to divide $p$
by $p_1$ to begin with, we see that $p$ reduces (as before) to give
$xyz + x - (xy - z)z = z^2 + x$, which is irreducible. Similarly,
dividing $p$ by $p_2$ to begin with, we obtain the remainder
$xyz + x - (yz + 2x + z)x = -2x^2 - xz + x$. However, whereas
before this remainder was irreducible, now we can reduce
it by the S-polynomial to give
$-2x^2 - xz + x - (-2x^2-xz-z^2) = z^2+x$, which is equal to
the first remainder.
\end{example}

Let us now formally define a Gr\"obner Basis in terms of
S-polynomials, noting that there are many other
equivalent definitions (see for example \cite{becker93}, page 206).

\begin{defn} \label{GB}
Let $G = \{g_1, \hdots, g_m\}$ be a basis for an ideal $J$
over a commutative polynomial ring $\mathcal{R} = R[x_1, \hdots, x_n]$.
If all the S-polynomials involving members of $G$
reduce to zero using $G$ ($\mathrm{S\mbox{-}pol}(g_i, g_j)
\rightarrow_G 0$ for all $i \neq j$), then $G$ is a
\index{Gr\"obner basis!commutative}
{\it Gr\"obner Basis} for \index{basis!Gr\"obner} $J$.
\index{commutative Gr\"obner basis}
\end{defn}

\begin{thm} \label{URC}
\index{unique remainder}
Given \index{remainder!unique}
any polynomial $p$ over a polynomial
ring $\mathcal{R} = R[x_1, \hdots, x_n]$, the
remainder of the division of $p$ by a basis $G$ for an
ideal $J$ in $\mathcal{R}$ is unique 
if and only if $G$ is a Gr\"obner Basis.
\end{thm}
\begin{pf}
($\Rightarrow$)
By Newman's Lemma (cf. \cite{becker93}, page 176), showing that the
remainder of the division of $p$ by $G$ is unique is equivalent to
showing that the division process is 
{\it locally confluent}, \index{locally confluent} that is
if there are polynomials $f$, $f_1$, $f_2 \in \mathcal{R}$ with
$f_1 = f-t_1g_1$ and $f_2 = f-t_2g_2$ for terms $t_1, t_2$
and $g_1, g_2 \in G$, then there exists a polynomial
$f_3 \in \mathcal{R}$ such that both $f_1$ and $f_2$ reduce to $f_3$.
By the Translation Lemma (cf. \cite{becker93}, page 200),
this in turn is equivalent to showing that the polynomial
$f_2-f_1 = t_1g_1 - t_2g_2$ reduces to zero, which is what we
shall now do.

There are two cases to deal with, $\LT(t_1g_1) \neq \LT(t_2g_2)$
and $\LT(t_1g_1) = \LT(t_2g_2)$. In the first case, notice that
the remainders $f_1$ and $f_2$ are obtained by cancelling off different
terms of the original $f$ (the reductions of $f$ are {\it disjoint}),
so it is possible, assuming (without loss of generality) that
$\LT(t_1g_1) > \LT(t_2g_2)$,
to directly reduce the polynomial $f_2-f_1 = t_1g_1 - t_2g_2$ in the
following manner: $t_1g_1 - t_2g_2 \rightarrow_{g_1}
 -t_2g_2 \rightarrow_{g_2} 0$.
In the second case, the reductions of $f$ are not disjoint
(as the same term $t$ from $f$ is cancelled off during both reductions),
so that the term $t$ does not appear in the polynomial
$t_1g_1 - t_2g_2$. However, the term $t$ is a common multiple of
$\LT(t_1g_1)$ and $\LT(t_2g_2)$, and thus the polynomial
$t_1g_1 - t_2g_2$ is a multiple of the S-polynomial
$\mathrm{S\mbox{-}pol}(g_1, g_2)$, say
$$t_1g_1 - t_2g_2 = \mu(\mathrm{S\mbox{-}pol}(g_1, g_2))$$
for some term $\mu$. Because $G$ is a Gr\"obner Basis, the S-polynomial
$\mathrm{S\mbox{-}pol}(g_1, g_2)$ reduces to zero, and hence by
extension the polynomial $t_1g_1 - t_2g_2$ also reduces to zero.

($\Leftarrow$)
As all S-polynomials are members of the ideal $J$, to 
complete the proof it is sufficient to show that
there is always a reduction path of an arbitrary
member of the ideal that leads to a zero remainder
(the uniqueness of remainders will then imply
that members of the ideal will always reduce to zero).
Let $f \in J = \langle G \rangle$. Then, by definition, there exist
$g_i \in G$ and $f_i \in \mathcal{R}$ 
(where $1 \leqslant i \leqslant j$)
such that $$f = \sum_{i=1}^j f_ig_i.$$
We proceed by induction on $j$. If $j = 1$, then $f = f_1g_1$,
and it is clear that we can use $g_1$ to reduce $f$ to give
a zero remainder ($f \rightarrow f - f_1g_1 = 0$).
Assume that the result is true for $j = k$, and let us
look at the case $j = k+1$, so that
$$f = \left(\sum_{i=1}^{k} f_ig_i \right) + f_{k+1}g_{k+1}.$$
By the inductive hypothesis, $\sum_{i=1}^k f_ig_i$ is a member
of the ideal that reduces to zero. The polynomial $f$ therefore 
reduces to the polynomial $f' := f_{k+1}g_{k+1}$, and we can 
now use $g_{k+1}$ to reduce $f'$ to give a zero remainder 
($f' \rightarrow f' - f_{k+1}g_{k+1} = 0$).
\end{pf}

We are now in a position to be able to define an
algorithm to compute a Gr\"obner Basis. However, to
be able to prove that this algorithm always
terminates, we must first prove a result stating
that all ideals over commutative polynomial rings are
finitely generated. This proof takes place in two stages --
first for monomial ideals (Dickson's Lemma) and then
for polynomial ideals (Hilbert's Basis Theorem).

% In the algorithm used to compute a Gr\"obner Basis, all
% S-polynomials between pairs of basis elements are considered
% in order to deal with the problem of two different
% reductions of a particular polynomial leading to two different
% remainders, the root cause of a non-unique remainder.

\section{Dickson's Lemma and Hilbert's Basis Theorem}

% Our goal in this section is to show that all polynomial
% ideals are finitely generated. To do this, we deal with the
% special case of monomial ideals first.
%
% \subsubsection{Dickson's Lemma}

\begin{defn}
A \index{ideal!monomial}
{\it monomial ideal} 
\index{monomial ideal} is an ideal generated by a
set of monomials.
\end{defn}

\begin{remark}
Any polynomial $p$ that is a member of a monomial ideal
is a sum of terms $p = \sum_i t_i$, where each $t_i$ is
a member of the monomial ideal.
\end{remark}

\begin{lem}[Dickson's Lemma]
\index{Dickson's lemma}
Every monomial ideal over the polynomial ring
$\mathcal{R} = R[x_1, \hdots, x_n]$ is finitely generated.
\end{lem}

\noindent
{\bf Proof (cf. \cite{Froberg98}, page 47):}
Let $J$ be a monomial ideal over $\mathcal{R}$
generated by a set of monomials $S$.
We proceed by induction on $n$, our goal being to show that
$S$ always has a finite subset $T$ generating $J$.
For $n = 1$, notice that all elements of $S$ will be of the
form $x_1^j$ for some $j \geqslant 0$. Let $T$ be the
singleton set containing the member of $S$ with the
lowest degree (that is the $x_1^j$ with the lowest
value of $j$). Clearly $T$ is finite, and because
any element of $S$ is a multiple of the chosen $x_1^j$,
it is also clear that $T$ generates the same ideal as $S$.

For the inductive step, assume that all monomial ideals
over the polynomial ring $\mathcal{R}' = R[x_1, \hdots, x_{n-1}]$
are finitely generated.
% Let $J_1 \subseteq J_2 \subseteq \cdots$ be an ascending chain of
% monomial ideals over $\mathcal{R}'$, each of which is finitely generated
% by the inductive hypothesis.
%
% This enables us to deduce\footnote{Using
% Corollary \ref{ACC} in the case of monomial ideals, replacing
% `Hilbert's Basis Theorem' by `Dickson's Lemma' in the proof.}
% an Ascending Chain Condition for monomial ideals
% over $\mathcal{R}'$, so that any ascending chain of monomial ideals
% $J_1 \subseteq J_2 \subseteq \cdots$ over
% $\mathcal{R}'$ is eventually constant.
Let $C_0 \subseteq C_1 \subseteq C_2
\subseteq \cdots$ be an ascending chain of monomial ideals
over $\mathcal{R}'$, where\footnote{
Think of $C_0$ as the set of monomials $m \in J$
which are also members of $\mathcal{R}'$; think of $C_j$
(for $j \geqslant 1$) as containing
all the elements of $C_{j-1}$ plus the monomials $m \in J$ of
the form $m = m'x_n^j$, $m' \in \mathcal{R'}$.
}
$$C_j = \langle S_j \rangle \cap \mathcal{R}',
\; S_j = \left\{ \frac{s}{\gcd(s, x_n^j)}
\left|\right. s \in S\right\}.$$
Let the monomial $m$ be an arbitrary member
% \footnote{It is sufficient
% to consider monomials only as any polynomial
% which is a member of a monomial ideal is made up of monomials which
% are all members of the monomial ideal.}
of the ideal $J$,
expressed as $m = m'x^k_n$, where $m' \in \mathcal{R}'$
and $k \geqslant 0$. By definition, $m' \in C_k$, and
so $m \in x^k_nC_k$. By the inductive hypothesis, each $C_k$
is finitely generated by a set $T_k$, and so
$m \in x^k_n\langle T_k\rangle$. From this we can deduce that
$$T = T_0 \cup x_nT_1 \cup x_n^2T_2 \cup \cdots$$
is a generating set for $J$.

Consider the ideal $C = \cup C_j$ for $j \geqslant 0$. This is
another monomial ideal over $\mathcal{R}'$, and so by
the inductive hypothesis is finitely generated. It follows that
the chain must stop as soon as the generators of $C$
are contained in some $C_r$, so that $C_r = C_{r+1} = \cdots$
(and hence $T_r = T_{r+1} = \cdots$). It follows that
$T_0 \cup x_nT_1 \cup x_n^2T_2 \cup \cdots \cup x_n^rT_r$
is a finite subset of $S$ generating $J$.
% Note that we will later state this result (and its proof)
% in more generality as Corollary \ref{ACC}.
%
% By the Ascending Chain Condition,
% the chain $C_0 \subseteq C_1 \subseteq \cdots$ is eventually
% constant, and so $T_r = T_{r+1} = \cdots$ for some finite $r$.
% It follows that $T$ is a finite subset of $S$ generating $J$.
{\hfill $\Box$}

\begin{example}
Let $S = \{y^{4}, \: xy^{4}, \: x^2y^{3}, \:
x^3y^{3}, \: x^4y, \: x^k\}$ be an
infinite set of monomials generating an ideal
$J$ over the polynomial ring $\mathbb{Q}[x, y]$, where $k$
is an integer such that $k \geqslant 5$. We can visualise
$J$ by using the following monomial lattice,
where a point $(a, b)$ in the lattice (for non-negative
integers $a,b$) corresponds to the monomial $x^ay^b$, and
the shaded region contains all monomials which are reducible
by some member of $S$ (and hence belong to $J$).
\begin{center}
\begin{picture}(0,0)%
\includegraphics{ch2d1.pstex}%
\end{picture}%
\setlength{\unitlength}{1973sp}%
\begingroup\makeatletter\ifx\SetFigFont\undefined%
\gdef\SetFigFont#1#2#3#4#5{%
  \reset@font\fontsize{#1}{#2pt}%
  \fontfamily{#3}\fontseries{#4}\fontshape{#5}%
  \selectfont}%
\fi\endgroup%
\begin{picture}(5100,4171)(1351,-6530)
\put(1351,-2611){\makebox(0,0)[lb]{\smash{\SetFigFont{12}{14.4}{\familydefault}{\mddefault}{\updefault}{\color[rgb]{0,0,0}$b$}%
}}}
\put(6451,-6436){\makebox(0,0)[lb]{\smash{\SetFigFont{12}{14.4}{\familydefault}{\mddefault}{\updefault}{\color[rgb]{0,0,0}$a$}%
}}}
\end{picture}

\end{center}
To show that $J$ can be finitely generated, we need to
construct the set $T$ as described in the proof of
Dickson's Lemma. The first step in doing this is to
construct the sequence of sets 
$S_j = \left\{ \frac{s}{\gcd(s,\: y^j)} \mid s \in S \right\}$
for all $j \geqslant 0$.
\begin{eqnarray*}
S_0 & = & \{y^{4}, \: xy^{4}, \: x^2y^{3}, \:
x^3y^{3}, \: x^4y, \: x^k\} = S \\
S_1 & = &  \{y^{3}, \: xy^{3}, \: x^2y^{2}, \:
x^3y^{2}, \: x^4, \: x^k\} \\
S_2 & = & \{y^{2}, \: xy^{2}, \: x^2y, \:
x^3y, \: x^4, \: x^k\} \\
S_3 & = & \{y, \: xy, \: x^2, \:
x^3, \: x^4, \: x^k\} \\
S_4 & = & \{y^0 = 1, \: x, \: x^2, \:
x^3, \: x^4, \: x^k\} \\
S_{j+1} & = & S_j ~~ \mbox{for all} ~ j+1 \geqslant 5.
\end{eqnarray*}
Each set $S_j$ gives rise to an ideal $C_j$
consisting of all monomials $m \in \langle S_j\rangle$ of
the form $m = x^i$ for some $i \geqslant 0$.
Because each of these ideals is an ideal over the polynomial
ring $\mathbb{Q}[x]$, we can use an inductive hypothesis
to give us a finite generating set $T_j$ for each $C_j$.
In this case, the first paragraph of the proof of Dickson's
Lemma tells us how to apply the inductive hypothesis ---
each set $T_j$ is formed by choosing the monomial $m \in S_j$
of lowest degree such that $m = x^i$ for some $i \geqslant 0$.
\begin{eqnarray*}
T_0 & = & \{x^5\} \\
T_1 & = & \{x^4\} \\
T_2 & = & \{x^4\} \\
T_3 & = & \{x^2\} \\
T_4 & = & \{x^0 = 1\} \\
T_{j+1} & = & T_j ~~ \mbox{for all} ~ k+1 \geqslant 5.
\end{eqnarray*}
We can now deduce that
$$
T = \{x^5\} \cup \{x^4y\} \cup \{x^4y^2\} \cup
    \{x^2y^3\} \cup \{y^4\} \cup \{y^5\} \cup \cdots
$$
is a generating set for $J$. Further, because
$T_{j+1} = T_j$ for all $k+1 \geqslant 5$, we
can also deduce that the set
$$
% T_0 \cup yT_1 \cup y^2T_2 \cup y^3T_3 \cup y^4T_4
% =
T' = \{x^5, \: x^4y, \: x^4y^2, \: x^2y^3, \: y^4\}
$$
is a finite generating set for $J$ (a fact that can be
verified by drawing a monomial lattice for $T'$
and comparing it with the above monomial lattice for
the set $S$).
\end{example}

\begin{thm}[Hilbert's Basis Theorem]
\index{Hilbert's basis theorem}
Every ideal $J$ over a polynomial ring
$\mathcal{R} = R[x_1, \hdots, x_n]$ is finitely generated.
\end{thm}
\begin{pf}
\index{$LM$@$\LM(J)$}
Let $O$ be a fixed arbitrary admissible monomial ordering, and
define $\LM(J) = \langle \LM(p) \mid p \in J\rangle$.
Because $\LM(J)$ is a monomial ideal, by Dickson's
Lemma it is finitely generated, say by the set of monomials
$M = \{m_1, \hdots, m_r\}$. By definition, each
$m_i \in M$ (for $1 \leqslant i \leqslant r$) has a corresponding
$p_i \in J$ such that $\LM(p_i) = m_i$. We claim that
$P = \{p_1, \hdots, p_r\}$ is a generating set for $J$.
To prove the claim, notice that
$\langle P\rangle \subseteq J$
so that $f \in \langle P\rangle \Rightarrow f \in J$.
Conversely, given a polynomial $f \in J$, we know that
$\LM(f) \in \langle M\rangle$ so that
$\LM(f) = \alpha m_j$ for some monomial $\alpha$ and
some $1 \leqslant j \leqslant r$. From this, if we
define $\alpha' = \frac{\LC(f)}{\LC(p_j)}\alpha$, we can
deduce that $\LM(f - \alpha' p_j) < \LM(f)$. Since
$f - \alpha' p_j \in J$, and because of the
admissibility of $O$, by recursion on 
$f - \alpha' p_j$
(define $f_{k+1} = f_k - \alpha'_kp_{j_k}$
for $k \geqslant 1$,
where $f_1 - \alpha'_1p_{j_1}
:= f - \alpha' p_j$),
we can deduce that $f \in \langle P\rangle$ (in fact
%we will eventually find that $f - \alpha s_j = 0$ so that
$f = \sum_{k=1}^{K} \alpha'_k p_{j_k}$
for some finite $K$).
\end{pf}

\begin{cor}[The Ascending Chain Condition] \label{ACC}
\index{ascending chain condition}
Every ascending sequence of ideals
$J_1 \subseteq J_2 \subseteq \cdots$
over a polynomial ring
$\mathcal{R} = R[x_1, \hdots, x_n]$
is eventually constant, so that there is an
$i$ such that $J_i = J_{i+1} = \cdots$.
\end{cor}
\begin{pf}
By Hilbert's Basis Theorem, each
ideal $J_k$ (for $k \geqslant 1$) is finitely generated.
Consider the ideal $J = \cup J_k$. This is
another ideal over $\mathcal{R}$, and so by Hilbert's Basis Theorem
is also finitely generated. From this we deduce that
the chain must stop as soon as the generators of $J$
are contained in some $J_i$, so that $J_i = J_{i+1} = \cdots$.
\end{pf}

\section{Buchberger's Algorithm}
\index{Buchberger's algorithm}

The algorithm used to compute a Gr\"obner Basis is
known as Buchberger's Algorithm.
\index{Bruno Buchberger}
Bruno Buchberger \index{Buchberger, Bruno}
was a student of Wolfgang Gr\"obner
at the University of Innsbruck, Austria, and the
publication of his PhD thesis in
1965 \cite{buch65} marked the start of Gr\"obner
Basis theory. 

In Buchberger's algorithm, S-polynomials
for pairs of elements from the current basis are
computed and reduced using the current basis. If
the S-polynomial does not reduce to zero, it is
added to the current basis, and this process continues
until all S-polynomials reduce to zero.
The algorithm works on the principle that if an S-polynomial
$\mathrm{S\mbox{-}pol}(g_i, g_j)$ does not reduce to zero using a
set of polynomials $G$, then it will certainly reduce to zero
using the set of polynomials $G \cup \{\mathrm{S\mbox{-}pol}(g_i, g_j)\}$.
% It is startling to think that such a difficult problem possesses
% such a simple solution.

\begin{algorithm}
\setlength{\baselineskip}{3.5ex}
\caption{A Basic Commutative Gr\"obner Basis Algorithm}
\label{com-buch}
\begin{algorithmic}
\vspace*{2mm}
\REQUIRE{A Basis $F = \{f_1, f_2, \hdots, f_m\}$ for an ideal $J$ over a
         commutative polynomial ring $R[x_1, \hdots x_n]$; 
         an admissible monomial ordering $\mathrm{O}$.}
\ENSURE{A Gr\"obner Basis $G = \{g_1, g_2, \hdots, g_p\}$ for $J$.}
\vspace*{1mm}
\STATE
Let $G = F$ and let $A = \emptyset$; \\
For each pair of polynomials $(g_i, g_j)$ in $G$ ($i < j$), \\
add the S-polynomial $\mathrm{S\mbox{-}pol}(g_i, g_j)$ to $A$;
\WHILE{($A$ is not empty)}
\STATE
Remove the first entry $s_1$ from $A$; \\
$s'_1 = \Rem(s_1, G)$; % (Algorithm \ref{com-div}) \\
% Reduce $s_1$ with respect to $G$ (with Algorithm \ref{com-div}); \\
\IF{($s'_1 \neq 0$)}
\STATE
Add $s'_1$ to $G$ and add all the S-polynomials 
$\mathrm{S\mbox{-}pol}(g_i, s'_1)$ to $A$ ($g_i \in G$, $g_i \neq s'_1$);
\ENDIF
\ENDWHILE
\STATE
{\bf return} $G$;
\end{algorithmic}
\vspace*{1mm}
\end{algorithm}

\begin{thm} \label{GBterm}
Algorithm \ref{com-buch} always terminates with a
Gr\"obner Basis for the ideal $J$.
\end{thm}

\noindent
{\bf Proof (cf. \cite{becker93}, page 213):}
{\it Correctness.} If the algorithm terminates, it does
so with a set of polynomials $G$ with the property that
all S-polynomials involving members of $G$
reduce to zero using $G$ ($\mathrm{S\mbox{-}pol}(g_i, g_j)
\rightarrow_G 0$ for all $i \neq j$). $G$ is therefore a 
Gr\"obner Basis by Definition \ref{GB}.
{\it Termination.} If the algorithm does not terminate,
then an endless sequence of polynomials must be added to
the set $G$ so that the set $A$ never becomes empty.
Let $G_0 \subset G_1 \subset G_2 \subset \cdots$
be the successive values of $G$. If we consider the 
corresponding sequence $\LM(G_0) \subset \LM(G_1)
\subset \LM(G_2) \subset \cdots$ of lead monomials,
we note that these sets generate an ascending chain
of ideals $J_0 \subset J_1 \subset J_2 \subset \cdots$
because each time we add a monomial to a particular
set $\LM(G_k)$ to form the set $\LM(G_{k+1})$, the monomial
we choose is irreducible with respect to
$\LM(G_k)$, and hence does not belong to the ideal $J_k$.
However the Ascending Chain Condition tells
us that such a chain of ideals must eventually
become constant, so there must be some $i \geqslant 0$
such that $J_i = J_{i+1} = \cdots$. It follows that the
algorithm will terminate once the set $G_i$ has been
constructed, as all of the S-polynomials left in $A$
will now reduce to zero (if not, some S-polynomial left
in $A$ will reduce to a non-zero polynomial $s'_1$ whose
lead monomial is irreducible with respect to $\LM(G_i)$,
allowing us to construct an ideal $J_{i+1} = 
\langle \LM(G_i) \cup \{\LM(s'_1)\}\rangle \supset
\langle \LM(G_i) \rangle = J_i$, contradicting the fact
that $J_{i+1} = J_i$.)
{\hfill $\Box$}

\begin{example} \label{exBuch}
Let $F := \{f_1, f_2\} = 
\{x^2 - 2xy+3, \: 2xy+y^2+5\}$ generate an ideal over
the commutative polynomial ring $\mathbb{Q}[x, y]$, and let
the monomial ordering be DegLex. Running
Algorithm \ref{com-buch} on $F$, there is only one S-polynomial to
consider initially, namely $\mathrm{S\mbox{-}pol}(f_1, f_2)
= y(f_1) - \frac{1}{2}x(f_2) = -\frac{5}{2}xy^2-\frac{5}{2}x+3y$.
This polynomial reduces (using $f_2$)
to give the irreducible polynomial
$\frac{5}{4}y^3-\frac{5}{2}x+\frac{37}{4}y =: f_3$,
which we add to our current basis. This produces two more
S-polynomials to look at,
$\mathrm{S\mbox{-}pol}(f_1, f_3) = y^3(f_1) - \frac{4}{5}x^2(f_3) =
-2xy^4+2x^3-\frac{37}{5}x^2y+3y^3$ and
$\mathrm{S\mbox{-}pol}(f_2, f_3) = \frac{1}{2}y^2(f_2) - \frac{4}{5}x(f_3) =
\frac{1}{2}y^4+2x^2-\frac{37}{5}xy+\frac{5}{2}y^2$,
both of which reduce to zero.
The algorithm therefore terminates with the set
$\{x^2 - 2xy+3, \: 2xy+y^2+5, \: 
\frac{5}{4}y^3-\frac{5}{2}x+\frac{37}{4}y\}$
as the output Gr\"obner Basis. 

Here is a dry run for Algorithm \ref{com-buch} in this instance.
\begin{center}
\begin{small}
\begin{tabular}{c|c|c|c|c|c}
$G$ & $i$ & $j$ & $A$ & $s_1$ & $s'_1$ \\
\hline
$\{f_1, f_2\}$ & 1 & 2 & $\emptyset$ & & \\
& & & $\{\mathrm{S\mbox{-}pol}(f_1, f_2)\}$ && \\
$\{f_1, f_2, f_3\}$ & 1 && $\emptyset$
& $-\frac{5}{2}xy^2-\frac{5}{2}x+3y$ & $f_3$ \\
&
  2 && $\{\mathrm{S\mbox{-}pol}(f_1, f_3)\}$ && \\
&&& $\{\mathrm{S\mbox{-}pol}(f_2, f_3), \:
          \mathrm{S\mbox{-}pol}(f_1, f_3)\}$ && \\
&&& $\{\mathrm{S\mbox{-}pol}(f_1, f_3)\}$ &
$\frac{1}{2}y^4+2x^2-\frac{37}{5}xy+\frac{5}{2}y^2$ & $0$ \\
&&& $\emptyset$ &
$-2xy^4+2x^3-\frac{37}{5}x^2y+3y^3$ & $0$ \\ \hline
\end{tabular}
\end{small}
\end{center}
\end{example}

\section{Reduced Gr\"obner Bases}

\begin{defn} \label{RGB}
Let $G = \{g_1, \hdots, g_p\}$ be a Gr\"obner Basis
for an ideal
over the polynomial ring $R[x_1, \hdots, x_n]$. $G$ is
a \index{Gr\"obner basis!reduced} {\it reduced}
\index{reduced Gr\"obner basis}
Gr\"obner Basis if the following conditions are satisfied.
\begin{enumerate}[(a)]
\item
$\LC(g_i) = 1_R$ for all $g_i \in G$.
\item
No term in any polynomial $g_i \in G$
is divisible by any $\LT(g_j)$, $j \neq i$.
\end{enumerate}
\end{defn}

\begin{thm} \label{urgb}
Every ideal over a commutative polynomial ring has a
unique reduced Gr\"obner Basis.
\end{thm}
\begin{pf}
{\it Existence.}
By Theorem \ref{GBterm}, there exists a Gr\"obner
Basis $G$ for every ideal over a commutative
polynomial ring. We claim that the following 
procedure transforms $G$ into a reduced Gr\"obner Basis $G'$.
\begin{enumerate}[(i)]
\item
Multiply each $g_i \in G$ by $\LC(g_i)^{-1}$.
\item
Reduce each $g_i \in G$ by $G\setminus \{g_i\}$, removing
from $G$ all polynomials that reduce to zero.
\end{enumerate}
It is clear that $G'$ satisfies the conditions of
Definition \ref{RGB}, so it remains to show that
$G'$ is a Gr\"obner Basis, which we shall do by showing that
the application of each step of instruction (ii) above produces
a basis which is still a Gr\"obner Basis.

Let $G = \{g_1, \hdots, g_p\}$ be a Gr\"obner Basis,
and let $g'_i$ be the reduction of an arbitrary $g_i \in G$
with respect to $G \setminus \{g_i\}$, carried out as follows
(the $t_k$ are terms).
\begin{equation} \label{gdashred}
g'_i = g_i - \sum_{k=1}^{\kappa} t_kg_{j_k}.
\end{equation}
Set $H = (G \setminus \{g_i\}) \cup \{g'_i\}$ if
$g'_i \neq 0$, and set $H = G \setminus \{g_i\}$
if $g'_i = 0$. As $G$ is a Gr\"obner Basis,
all S-polynomials involving elements
of $G$ reduce to zero using $G$, so there are expressions
\begin{equation} \label{sprz}
t_ag_a - t_bg_b - \sum_{u=1}^{\mu}t_ug_{c_u} = 0
\end{equation}
for every S-polynomial $\mathrm{S\mbox{-}pol}(g_a, g_b) =
t_ag_a - t_bg_b$, where $g_a, g_b, g_{c_u} \in G$. To show that $H$
is a Gr\"obner Basis, we must show that all S-polynomials
involving elements of $H$ reduce to zero using $H$. For distinct
polynomials $g_a, g_b \in H$ not equal to
$g'_i$, we can reduce the
S-polynomial $\mathrm{S\mbox{-}pol}(g_a, g_b)$
using the reduction shown in Equation
(\ref{sprz}), substituting for $g_i$ from 
Equation (\ref{gdashred})
if any of the $g_{c_u}$ in 
Equation (\ref{sprz}) are equal to $g_i$.
This gives a reduction to zero of $\mathrm{S\mbox{-}pol}(g_a, g_b)$
in terms of elements of $H$.

If $g'_i = 0$, our proof is complete. Otherwise
consider the S-polynomial $\mathrm{S\mbox{-}pol}(g'_i, g_a)$.
We claim that $\mathrm{S\mbox{-}pol}(g_i, g_a) = t_1g_i - t_2g_a
\Rightarrow \mathrm{S\mbox{-}pol}(g'_i, g_a) = t_1g'_i - t_2g_a$. To prove
this claim, it is sufficient to show that $\LT(g_i) = \LT(g'_i)$.
Assume for a contradiction that $\LT(g_i) \neq \LT(g'_i)$.
It follows that during the reduction of $g_i$ we were able to
reduce its lead term, so that
$\LT(g_i) = t\LT(g_j)$ for some term $t$ and some $g_j \in G$.
By the admissibility of the chosen monomial ordering,
the polynomial $g_i - tg_j$ reduces to zero without
using $g_i$, leading to the conclusion that $g'_i = 0$, 
a contradiction.

It remains to show that $\mathrm{S\mbox{-}pol}(g'_i, g_a)
\rightarrow_H 0$. We know that
$\mathrm{S\mbox{-}pol}(g_i, g_a) = t_1g_i - t_2g_a
\rightarrow_G 0$, and Equation (\ref{sprz}) tells us that
$t_1g_i - t_2g_a - \sum_{u=1}^{\mu}t_ug_{c_u} = 0$.
Substituting for $g_i$ from 
Equation (\ref{gdashred}), we obtain\footnote{Substitutions
for $g_i$ may also occur in the summation
$\sum_{u=1}^{\mu}t_ug_{c_u}$; 
these substitutions have not been considered in the displayed formulae.}
$$
t_1\left(g'_i + \sum_{k=1}^{\kappa} t_kg_{j_k}\right)
- t_2g_a - \sum_{u=1}^{\mu}t_ug_{c_u} = 0 $$
or
$$
t_1g'_i - t_2g_a -
\left(\sum_{u=1}^{\mu}t_ug_{c_u}
- \sum_{k=1}^{\kappa} t_1t_kg_{j_k}\right) = 0,
$$
which implies that $\mathrm{S\mbox{-}pol}(g'_i, g_a)
\rightarrow_H 0$.

{\it Uniqueness.}
Assume for a contradiction that $G = \{g_1, \hdots, g_p\}$
and $H = \{h_1, \hdots, h_q\}$ are two reduced Gr\"obner
Bases for an ideal $J$, with $G \neq H$. Let $g_i$ be an arbitrary
element from $G$ (where $1 \leqslant i \leqslant p$).
Because $g_i$ is a member of the ideal,
then $g_i$ must reduce to zero using $H$
($H$ is a Gr\"obner Basis). This means that there must exist
a polynomial $h_j \in H$ such that $\LT(h_j) \mid \LT(g_i)$. If
$\LT(h_j) \neq \LT(g_i)$, then $\LT(h_j) \times m = \LT(g_i)$ for
some nontrivial monomial $m$. But $h_j$ is also a member of the ideal,
so it must reduce to zero using $G$. Therefore there exists 
a polynomial $g_k \in G$
such that $\LT(g_k) \mid \LT(h_j)$, which implies that
$\LT(g_k) \mid \LT(g_i)$, with $k \neq i$.
This contradicts condition (b) of Definition \ref{RGB}, so that $G$
cannot be a reduced Gr\"obner Basis for $J$ if $\LT(h_j) \neq \LT(g_i)$.
From this we deduce that each $g_i \in G$ has a corresponding
$h_j \in H$ such that $\LT(g_i) = \LT(h_j)$. Further, because $G$ and $H$ are
assumed to be reduced Gr\"obner Bases, this is a one-to-one correspondence.

It remains to show that if $\LT(g_i) = \LT(h_j)$, then $g_i = h_j$.
Assume for a contradiction that $g_i \neq h_j$, and consider the
polynomial $g_i - h_j$. Without loss of generality, assume that
$\LM(g_i - h_j)$ appears in $g_i$. Because $g_i - h_j$ is a member
of the ideal, then there is a polynomial $g_k \in G$ such that
$\LT(g_k) \mid \LT(g_i - h_j)$. But this again contradicts condition
(b) of Definition \ref{RGB}, as we have shown that there is a term in
$g_i$ that is divisible by $\LT(g_k)$ for some $k \neq i$. It follows
that $G$ cannot be a reduced Gr\"obner Basis if $g_i \neq h_j$, 
which means that $G = H$ and therefore reduced Gr\"obner Bases are unique.
\end{pf}

Given a Gr\"obner Basis $G$, we saw in the proof of Theorem \ref{urgb}
that if the lead term of any polynomial $g_i \in G$ is reducible by some
polynomial $g_j \in G$ (where $g_j \neq g_i$), then $g_i$ reduces to zero.
We can use this information to refine the procedure for finding a unique
reduced Gr\"obner Basis (as given in the aforementioned proof)
by allowing the removal of any polynomial $g_i \in G$ whose lead monomial
is a multiple of some other lead monomial $\LM(g_j)$.
This process, which if often referred to as {\it minimising}
\index{Gr\"obner basis!minimal}
a \index{minimal Gr\"obner basis}
Gr\"obner Basis (as in finding a Gr\"obner Basis with the minimal
number of elements), is incorporated into our refined procedure,
which we state as Algorithm \ref{red-com}.

\begin{algorithm}
\setlength{\baselineskip}{3.5ex}
\caption{The Commutative Unique Reduced Gr\"obner Basis Algorithm}
\label{red-com}
\begin{algorithmic}
\vspace*{2mm}
\REQUIRE{A Gr\"obner Basis $G = \{g_1, g_2, \hdots, g_m\}$
         for an ideal $J$ over a
         commutative polynomial ring $R[x_1, \hdots x_n]$;
         an admissible monomial ordering $\mathrm{O}$.}
\ENSURE{The unique reduced Gr\"obner Basis
        $G' = \{g'_1, g'_2, \hdots, g'_p\}$ for $J$.}
\vspace*{1mm}
\STATE
$G' = \emptyset$;
\FOR{\textbf{each} $g_i \in G$}
\STATE
Multiply $g_i$ by $\LC(g_i)^{-1}$; \\
\IF{($\LM(g_i) = u\LM(g_j)$ for some
monomial $u$ and some $g_j \in G$ ($g_j \neq g_i$))}
\STATE
$G = G\setminus\{g_i\}$;
\ENDIF
\ENDFOR
\FOR{\textbf{each} $g_i \in G$}
\STATE
$g'_i = \Rem(g_i, (G \setminus \{g_i\})\cup G')$; \\
$G = G\setminus \{g_i\}$; $G' = G' \cup \{g'_i\}$;
% Add $g'_i$ to $G'$, where $g'_i$ is the reduction of
% $g_i$ with respect to $G\setminus \{g_i\}$;
\ENDFOR
\STATE
{\bf return} $G'$;
\end{algorithmic}
\vspace*{1mm}
\end{algorithm}

\section{Improvements to Buchberger's Algorithm} \label{ImprovC}

Nowadays, most general purpose
symbolic computation systems possess
an implementation of Buchberger's algorithm.
These implementations often take advantage of the
numerous improvements made to
Buchberger's algorithm over the years, some of which
we shall now describe.

% (one of the most efficient can be found in
% Singular \cite{Sing01}).

\subsection{Buchberger's Criteria}

In 1979, Buchberger published a paper \cite{buch79} which
gave criteria that enable the {\it a priori} detection of
S-polynomials that reduce to zero. This speeds up
Algorithm \ref{com-buch} by drastically reducing the number
of S-polynomials that must be reduced with respect to the
current basis.

\begin{prop}[Buchberger's First Criterion]
\index{Buchberger's first criterion}
Let \index{criterion!Buchberger's first}
$f$ and $g$ be two polynomials over a commutative polynomial
ring ordered with respect to some fixed admissible monomial ordering $O$.
If the lead terms of $f$ and $g$ are
disjoint (so that $\lcm(\LM(f), \LM(g)) = \LM(f)\LM(g)$),
then $\mathrm{S\mbox{-}pol}(f, g)$ reduces to zero using
the set $\{f, g\}$.
\end{prop}
\noindent \textbf{Proof (Adapted from \cite{becker93}, Lemma 5.66):}
Assume that $f = \sum_{i=1}^{\alpha} s_i$ and
$g = \sum_{j=1}^{\beta} t_j$, where the $s_i$ and the $t_j$ are terms.
Because $s_1$ and $t_1$ are disjoint, it follows that
\begin{eqnarray}
\nonumber
\mathrm{S\mbox{-}pol}(f, g) & \equiv & t_1f - s_1g \\
\label{B1diff}
 & = & t_1(s_2 + \cdots + s_{\alpha}) - s_1(t_2 + \cdots + t_{\beta}).
\end{eqnarray}
We claim that no two terms in Equation (\ref{B1diff}) are the same.
Assume to the contrary that $t_1s_i = s_1t_j$ for some
$2 \leqslant i \leqslant \alpha$ and $2 \leqslant j \leqslant \beta$.
Then $t_1s_i$ is a multiple of both $t_1$ and $s_1$, which means
that $t_1s_i$ is a multiple of $\lcm(t_1, s_1) = t_1s_1$. But then
we must have $t_1s_i \geqslant t_1s_1$, which gives a contradiction
(by definition $s_1 > s_i$).

As every term in $t_1(s_2 + \cdots + s_{\alpha})$ is a multiple of
$t_1$, we can use $g$ to eliminate each of the terms
$t_1s_{\alpha}$, $t_1s_{\alpha-1}$, $\hdots$, $t_1s_2$ in 
Equation (\ref{B1diff}) in turn:
\begin{eqnarray}
\nonumber
&& t_1(s_2 + \cdots + s_{\alpha}) - s_1(t_2 + \cdots + t_{\beta}) \\
\nonumber
& \rightarrow & t_1(s_2 + \cdots + s_{\alpha}) - s_1(t_2 + \cdots + t_{\beta})
- s_{\alpha}g \\
\nonumber
& = & t_1(s_2 + \cdots + s_{\alpha-1}) - s_1(t_2 + \cdots + t_{\beta})
- s_{\alpha}(t_2 + \cdots + t_{\beta}) \\
% \nonumber
% & = & t_1(s_2 + \cdots + s_{\alpha-1}) -
% (s_1 + s_{\alpha})(t_2 + \cdots + t_{\beta}) \\
\nonumber
& \rightarrow & t_1(s_2 + \cdots + s_{\alpha-2}) -
      (s_1 + s_{\alpha-1} + s_{\alpha})(t_2 + \cdots + t_{\beta}) \\
\nonumber
& \vdots & \\
\nonumber
& \rightarrow & -(s_1 + s_2 + \cdots 
+ s_{\alpha})(t_2 + \cdots + t_{\beta}) \\
\label{B12nd}
& = & -s_1(t_2 + \cdots + t_{\beta}) -
\cdots - s_{\alpha}(t_2 + \cdots + t_{\beta}).
\end{eqnarray}
We do this in reverse order because, having eliminated a term
$t_1s_{\gamma}$ (where $3 \leqslant \gamma \leqslant \alpha$),
to continue
the term $t_1s_{\gamma-1}$ must appear in the reduced polynomial
(which it does because $t_1s_{\gamma-1} > s_{\delta}t_{\eta}$ for
all $\gamma \leqslant \delta \leqslant \alpha$ and
$2 \leqslant \eta \leqslant \beta$).

We now use the same argument on $-s_1(t_2 + \cdots + t_{\beta})$,
using $f$ to eliminate each of its terms in turn, giving the
following reduction sequence.
\begin{eqnarray*}
& & -s_1(t_2 + \cdots + t_{\beta}) -
\cdots - s_{\alpha}(t_2 + \cdots + t_{\beta}) \\
& \rightarrow & -s_1(t_2 + \cdots + t_{\beta}) -
\cdots - s_{\alpha}(t_2 + \cdots + t_{\beta}) + t_2f \\
& = & -s_1(t_2 + \cdots + t_{\beta}) - \cdots
- s_{\alpha}(t_2 + \cdots + t_{\beta}) + t_2(s_1 + \cdots + s_{\alpha}) \\
& = & -s_1(t_3 + \cdots + t_{\beta}) - \cdots
- s_{\alpha}(t_3 + \cdots + t_{\beta}) \\
& \rightarrow & -s_1(t_4 + \cdots + t_{\beta}) - \cdots
- s_{\alpha}(t_4 + \cdots + t_{\beta}) \\
& \vdots & \\
& \rightarrow & 0.
\end{eqnarray*}
{\it Technical point:} If some term $s_it_j$ 
(for $i,j \geqslant 2$) cancels the
term $s_1t_k$ (for $k \geqslant 3$) in Equation (\ref{B12nd}),
then as we must have $j < k$ in order to have $s_it_j = s_1t_k$, the term
$s_1t_k$ will reappear as $s_it_j$ when the term $s_1t_j$ is eliminated,
allowing us to continue the reduction as shown. This argument can be
extended to the case where a combination of terms of the form $s_it_j$
cancel the term $s_1t_k$, as the term $s_1t_k$ will reappear after all
the terms $s_1t_{\kappa}$ (for 
$2 \leqslant \kappa < k$) have been eliminated.
{\hfill $\Box$}

\begin{prop}[Buchberger's Second Criterion]
\index{Buchberger's second criterion}
Let \index{criterion!Buchberger's second}
$f$, $g$ and $h$ be three members of a finite set of
polynomials $P$ over a commutative polynomial ring satisfying
the following conditions.
\begin{enumerate}[(a)]
\item
$\LM(h) \mid \lcm(\LM(f), \LM(g))$.
\item
$\mathrm{S\mbox{-}pol}(f, h) \rightarrow_P 0$ and
$\mathrm{S\mbox{-}pol}(g, h) \rightarrow_P 0$.
\end{enumerate}
Then $\mathrm{S\mbox{-}pol}(f, g) \rightarrow_P 0$.
\end{prop}
\begin{pf}
If $\LM(h) \mid \lcm(\LM(f), \LM(g))$, then $m_h\LM(h) =
\lcm(\LM(f), \LM(g))$ for some monomial $m_h$. Assume that
$\lcm(\LM(f), \LM(g)) = m_f\LM(f) = m_g\LM(g)$ for some monomials
$m_f$ and $m_g$. Then it is clear that $m_f\LM(f) = m_h\LM(h)$ is
a common multiple of $\LM(f)$ and $\LM(h)$, and $m_g\LM(g) =
m_h\LM(h)$ is a common multiple of $\LM(g)$ and $\LM(h)$. It
follows that $\lcm(\LM(f), \LM(g))$ is a multiple of both
$\lcm(\LM(f), \LM(h))$ and $\lcm(\LM(g), \LM(h))$, so that
\begin{equation} \label{lcmjuggle}
\lcm(\LM(f), \LM(g)) = m_{fh}\lcm(\LM(f), \LM(h)) =
m_{gh}\lcm(\LM(g), \LM(h))
\end{equation}
for some monomials $m_{fh}$ and $m_{gh}$.

Because the S-polynomials $\mathrm{S\mbox{-}pol}(f, h)$ and
$\mathrm{S\mbox{-}pol}(g, h)$ both reduce to zero using $P$,
there are expressions
$$\mathrm{S\mbox{-}pol}(f, h) - \sum_{i=1}^{\alpha} s_ip_i = 0$$
and $$\mathrm{S\mbox{-}pol}(g, h) - \sum_{j=1}^{\beta} t_jp_j = 0,$$
where the $s_i$ and the $t_j$ are terms, and $p_i, p_j \in P$ for all
$i$ and $j$. It follows that
\begin{eqnarray*}
m_{fh}\left(\mathrm{S\mbox{-}pol}(f, h) - \sum_{i=1}^{\alpha} s_ip_i\right)
& = & m_{gh}\left(\mathrm{S\mbox{-}pol}(g, h)
- \sum_{j=1}^{\beta} t_jp_j\right);
\end{eqnarray*}
\vspace*{-12mm}
\begin{eqnarray*}
\lefteqn{
m_{fh}\left(\frac{\lcm(\LM(f), \LM(h))}{\LT(f)}f -
\frac{\lcm(\LM(f), \LM(h))}{\LT(h)}h - \sum_{i=1}^{\alpha} s_ip_i\right)
=} \\
&& \hfill
m_{gh}\left(\frac{\lcm(\LM(g), \LM(h))}{\LT(g)}g -
\frac{\lcm(\LM(g), \LM(h))}{\LT(h)}h
- \sum_{j=1}^{\beta} t_jp_j\right); \\
\lefteqn{
m_{fh}\left(\frac{\lcm(\LM(f), \LM(g))}{m_{fh}\LT(f)}f -
\frac{\lcm(\LM(f), \LM(g))}{m_{fh}\LT(h)}h - \sum_{i=1}^{\alpha} s_ip_i\right)
=} \\
&& \hfill
m_{gh}\left(\frac{\lcm(\LM(f), \LM(g))}{m_{gh}\LT(g)}g -
\frac{\lcm(\LM(f), \LM(g))}{m_{gh}\LT(h)}h
- \sum_{j=1}^{\beta} t_jp_j\right);
\end{eqnarray*}
\vspace*{-12mm}
\begin{eqnarray*}
\frac{\lcm(\LM(f), \LM(g))}{\LT(f)}f - m_{fh}\sum_{i=1}^{\alpha} s_ip_i
& = &
\frac{\lcm(\LM(f), \LM(g))}{\LT(g)}g - m_{gh}\sum_{j=1}^{\beta} t_jp_j; \\
\mathrm{S\mbox{-}pol}(f, g) - \sum_{i=1}^{\alpha} m_{fh}s_ip_i
+ \sum_{j=1}^{\beta} m_{gh}t_jp_j & = & 0.
\end{eqnarray*}
To conclude that the S-polynomial $\mathrm{S\mbox{-}pol}(f, g)$ reduces to
zero using $P$, it remains to show that the algebraic expression
$-\sum_{i=1}^{\alpha} m_{fh}s_ip_i
+ \sum_{j=1}^{\beta} m_{gh}t_jp_j$ corresponds to a valid reduction of
$\mathrm{S\mbox{-}pol}(f, g)$. To do this, it is sufficient to show
that no term in either of the summations is
greater than $\lcm(\LM(f), \LM(g))$ (so that $\LM(m_{fh}s_ip_i) <
\lcm(\LM(f), \LM(g))$ and $\LM(m_{gh}t_jp_j) < \lcm(\LM(f), \LM(g))$
for all $i$ and $j$). But this follows from 
Equation (\ref{lcmjuggle}) and from
the fact that the original reductions of $\mathrm{S\mbox{-}pol}(f, h)$ and
$\mathrm{S\mbox{-}pol}(g, h)$ are valid, so that
$\LM(s_ip_i) < \lcm(\LM(f), \LM(h))$ and
$\LM(t_jp_j) < \lcm(\LM(g), \LM(h))$ for all $i$ and $j$.
\end{pf}

\subsection{Homogeneous Gr\"obner Bases} \label{HomSecC}

\begin{defn}
A polynomial is {\it homogeneous} \index{homogeneous} if all its terms
have the same degree. For example, the polynomial $x^2y + 4yz^2 + 3z^3$ is
homogeneous, but the polynomial $x^3y + 4x^2 + 45$ is not homogeneous.
\end{defn}

Of the many systems available for computing commutative Gr\"obner Bases,
some (such as Bergman \cite{Bergman98}) only admit sets 
of homogeneous polynomials
as input. This restriction leads to gains in efficiency 
% in the mechanics of Buchberger's algorithm
as we can take advantage of some of the properties of homogeneous
polynomial arithmetic. For example, the
S-polynomial of two homogeneous polynomials is homogeneous, and
the reduction of a homogeneous polynomial by a set of homogeneous
polynomials yields another homogeneous polynomial. It follows that if $G$ is
a Gr\"obner Basis for a set $F$ of homogeneous polynomials, then $G$
is another set of homogeneous polynomials.

At first glance, it seems that a system accepting only sets of
homogeneous polynomials as input is not able to compute a
Gr\"obner Basis for a set of polynomials containing one or more
non-homogeneous polynomials. However, we can still use the system if
we use an extendible monomial ordering and the processes of
homogenisation and dehomogenisation.

\begin{defn} \label{defhom}
Let $p = p_0 + \cdots + p_m$ be a polynomial over the
polynomial ring $R[x_1, \hdots, x_n]$, where
each $p_i$ is the sum of the degree $i$ terms in $p$
(we assume that $p_m \neq 0$). The
% the $p_i$ are homogeneous polynomials of degree $i$ and $p_m \neq 0$. The
{\it homogenisation} \index{homogenisation} of $p$
with respect to a new (homogenising) variable $y$ is the polynomial
$$h(p) := p_0y^m + p_1y^{m-1} + \cdots +
p_{m-1}y + p_m,$$
where $h(p)$ belongs to a polynomial ring determined
by where $y$ is placed in the lexicographical ordering of
the variables.
\end{defn}

\begin{defn}
The {\it dehomogenisation} \index{dehomogenisation} of a
% homogeneous
polynomial $p$ is the polynomial $d(p)$ given by
substituting $y = 1$ in $p$, where $y$ is the homogenising
variable. For example, the dehomogenisation of the
polynomial $x_1^3 + x_1x_2y + x_1y^2 \in
\mathbb{Q}[x_1, x_2, y]$ is the polynomial $x_1^3 + x_1x_2 + x_1
\in \mathbb{Q}[x_1, x_2]$.
% $p \in \mathbb{Q}[x_1, \hdots, x_n, y]$ is the polynomial
% $d(p) \in \mathbb{Q}[x_1, \hdots, x_n]$ 
% given by substituting $y = 1$ in $p$.
\end{defn}

\begin{defn}
A monomial ordering $O$ is {\it extendible}
\index{monomial ordering!extendible} if, given any polynomial
\index{extendible monomial ordering}
$p = t_1 + \cdots + t_\alpha$ ordered with respect to $O$
(where $t_1 > \cdots > t_\alpha$), 
the homogenisation of $p$ preserves the order
on the terms ($t'_i > t'_{i+1}$ for all
$1 \leqslant i \leqslant \alpha-1$,
where the homogenisation process maps the term $t_i \in p$
to the term $t'_i \in h(p)$).
\end{defn}

Of the monomial orderings defined in Section \ref{CMO},
two of them (Lex and DegRevLex)
are extendible as long as we ensure
that the new variable $y$ is lexicographically {\it less}
than any of the variables $x_1, \hdots, x_n$; another (InvLex)
is extendible as long as we ensure that the new variable
$y$ is lexicographically {\it greater} than any of the 
variables $x_1, \hdots, x_n$.

The other monomial orderings are {\it not} extendible as,
no matter where we place the new variable $y$ in the
ordering of the variables, we can always find two monomials
$m_1$ and $m_2$ such that, if $p = m_1 + m_2$
(with $m_1 > m_2$),
then in $h(p) = m'_1 + m'_2$, we have $m'_1 < m'_2$. For
example, $m_1 := x_1x_2^2$ and $m_2 := x_1^2$ provides a
counterexample for the DegLex monomial ordering.

\begin{defn} \label{homprocC}
Let $F = \{f_1, \hdots, f_m\}$ be a non-homogeneous
set of polynomials. To compute a Gr\"obner Basis for
$F$ using a program that only accepts sets of
homogeneous polynomials as input, we proceed as follows.
\begin{enumerate}[(a)]
\item
Construct a homogeneous set of polynomials
$F' = \{h(f_1), \hdots, h(f_m)\}$.
\item
Compute a Gr\"obner Basis $G'$ for $F'$.
\item
Dehomogenise each polynomial $g' \in G'$ to obtain
a set of polynomials $G$.
\end{enumerate}
\end{defn}
As long as the chosen monomial ordering $O$ is extendible, $G$ will be
a Gr\"obner Basis for $F$ with respect to $O$
\cite[page 113]{Froberg98}. A word of warning
however -- this process is not necessarily more efficient that the
direct computation of a Gr\"obner Basis for $F$ using a program that does
accept non-homogeneous sets of polynomials as input.

\subsection{Selection Strategies} \label{SSC}
\index{selection strategies}

One of the most important factors when considering the efficiency
of Buchberger's algorithm is the order in which S-polynomials
are processed during the algorithm. A particular choice of a
{\it selection strategy} to use can often cut down substantially
the amount of work required in order to obtain a particular Gr\"obner Basis.

In 1979, Buchberger defined the {\it normal strategy} \cite{buch79}
\index{normal strategy} that 
\index{strategy!normal} chooses to process an S-polynomial
$\mathrm{S\mbox{-}pol}(f, g)$ if the monomial
$\lcm(\LM(f), \LM(g))$ is minimal (in the chosen monomial ordering)
amongst all such lowest common multiples. This strategy was
refined in 1991 to give the {\it sugar strategy} \cite{GMNRT91},
\index{sugar strategy} a 
\index{strategy!sugar} strategy that chooses an S-polynomial to process
if the {\it sugar} of the S-polynomial (a value associated to
the S-polynomial) is minimal amongst all such values (the normal
strategy is used in the event of a tie).
% The sugar strategy is a
% popular choice because it usually performs well in most examples.

Motivation for the sugar strategy comes from the
observation that the normal strategy performs well when
used with a degree-based monomial ordering and a homogeneous basis; the
sugar strategy was developed as a way to proceed based on what
% choose which S-polynomial to process next based on what
would happen when using the normal strategy in the computation of
a Gr\"obner Basis for the corresponding homogenised input basis.
We can therefore think of the sugar of an S-polynomial as representing
the degree of the corresponding S-polynomial in the
homogeneous computation.

% Indeed in \cite{GMNRT91}, the sugar of an S-polynomial is
% described as being ``the degree that it would have if computed with the
% homogeneous algorithm''.

% When computing a Gr\"obner Basis for a set of polynomials
% $F = \{f_1, \hdots, f_m\}$, the sugar of an S-polynomial
% in this computation
% represents the degree of the corresponding S-polynomial
% in the computation of a Gr\"obner Basis for the
% homogenised basis $F' = \{h(f_1), \hdots, h(f_m)\}$. This
% definition comes from the observation that the normal strategy
% performs well when used with a degree-based monomial ordering and
% a homogeneous basis.

The sugar of an S-polynomial is computed by
using the following rules on the sugars of polynomials we
encounter during the computation of a Gr\"obner Basis for the
set of polynomials $F = \{f_1, \hdots, f_m\}$.
\begin{enumerate}[(1)]
\item
The sugar $\Sug_{f_i}$ of a polynomial $f_i \in F$ is the total
degree of the polynomial $f_i$ (which is the degree of the term
of maximal degree in $f_i$).
\item
If $p$ is a polynomial and if $t$ is a term, then
$\Sug_{tp} = \deg(t) + \Sug_p$.
\item
If $p = p_1 + p_2$, then
$\Sug_p = \max(\Sug_{p_1}, \Sug_{p_2})$.
\end{enumerate}
It follows that the sugar of the S-polynomial
$\mathrm{S\mbox{-}pol}(g, h) =
\frac{\lcm(\LM(g), \LM(h))}{\LT(g)}g -
\frac{\lcm(\LM(g), \LM(h))}{\LT(h)}h$
is given by the formula
$$\Sug_{\mathrm{S\mbox{-}pol}(g, h)} =
 \max(\Sug_g - \deg(\LM(g)), \Sug_h - \deg(\LM(h)))
 + \deg(\lcm(\LM(g), \LM(h))).$$
\vspace*{-12pt}
\begin{example}
To illustrate how a selection strategy reduces the amount of work
required to compute a Gr\"obner Basis, consider the ideal
generated by the basis $\{x^{31}-x^6-x-y, \, x^8-z, \, x^{10}-t\}$ over the
polynomial ring $\mathbb{Q}[x, y, z, t]$. In our own implementation of
Buchberger's algorithm, here is the number of S-polynomials
processed during the algorithm when different selection
strategies and different monomial orderings
are used (the numbers quoted take into account
the application of both of Buchberger's criteria).
% As you can see, there is a substantial decrease in the number
% of S-polynomials that we are required to look at when a
% selection strategy is used.
\begin{center}
\begin{tabular}{c|c|c|c}
Selection Strategy & Lex & DegLex & DegRevLex \\
\hline
No strategy & 640 & 275 & 320 \\
Normal strategy & 123 & 63 & 61 \\
Sugar strategy & 96 & 55 & 54 \\ \hline
\end{tabular}
\end{center}
\end{example}

% \subsubsection{Lexicographic Orderings}
%
% Optimises the lexicographical ordering according to
% the frequency of the variables in the input basis
% (most frequent = lexicographically smallest).
% Comes from a paper by Freyja Hreinsd\'ottir
% \cite{Freyja94}.

\subsection{Basis Conversion Algorithms} \label{BCAC}

One factor which heavily influences the amount of time
taken to compute a Gr\"obner Basis is the monomial
ordering chosen. It is well known that some monomial
orderings (such as Lex) are characterised as being `slow',
while other monomial orderings (such as DegRevLex) are
said to be `fast'. In practice
what this means is that it usually takes far more time
to calculate (say) a Lex Gr\"obner Basis than it
does to calculate a DegRevLex Gr\"obner Basis for the
same generating set of polynomials.

Because many of the useful applications of Gr\"obner Bases
(such as solving systems of polynomial equations) depend
on using `slow' monomial orderings,
a number of algorithms were developed in the 1990's that
allow us to obtain a Gr\"obner Basis with respect to one
monomial ordering from a Gr\"obner Basis with respect to
another monomial ordering. 

The idea is that the time it takes
to compute a Gr\"obner Basis with respect to a `fast'
monomial ordering and then to convert it to a Gr\"obner
Basis with respect to a `slow' monomial ordering may be significantly
less than the time it takes to compute a Gr\"obner Basis
for the `slow' monomial ordering directly. Although
seemingly counterintuitive, the idea works well in practice.

One of the first conversion methods developed was the {\it FGLM}
method, named after the four authors who published the paper
\index{FGLM}
\cite{FGLM} introducing it. The method relies on linear
algebra to do the conversion, working with coefficient
matrices and irreducible monomials. Its only drawback lies in the fact
that it can only be used with zero-dimensional ideals, which
are the ideals containing only a finite number of irreducible monomials
(for each variable $x_i$ in the polynomial ring, a Gr\"obner Basis for a
zero-dimensional ideal must contain a polynomial which has a power of
$x_i$ as the leading monomial). This restriction does not apply in the
case of the {\it Gr\"obner Walk} \cite{CKMWalk},
\index{Gr\"obner walk} a 
\index{walk!Gr\"obner} basis conversion method
we shall study in further detail in Chapter \ref{ChWalk}.

\subsection{Optimal Variable Orderings}

In many cases, the ordering of the variables in a polynomial
ring can have a significant effect on the time it takes to
compute a Gr\"obner Basis for a particular ideal (an
example can be found in \cite{Carlson05}). This
is worth bearing in mind if we are searching for
{\it any} Gr\"obner Basis with respect to a certain ideal,
so do not mind which variable ordering is being used.
A heuristically optimal variable ordering
is described in \cite{Freyja94} (deriving from 
a discussion in \cite{BGK}), where we order the variables so
that the variable that occurs least often in the
polynomials of the input basis is the largest variable;
the second least common variable is the second
largest variable; and so on (ties are broken randomly).

\begin{example}
Let $F := \{y^2z^2 + x^2y, \; x^2y^4z + xy^2z + y^3, \; y^7 + x^3z\}$
generate an ideal over the polynomial ring $\mathbb{Q}[x, y, z]$.
Because $x$ occurs 8 times in $F$, $y$ occurs 19 times and
$z$ occurs 5 times, the heuristically optimal variable
ordering is $z > x > y$. This is supported by the following
table showing the times taken to compute a Lex Gr\"obner
Basis for $F$ using all six possible variable orderings, where we see
that the time for the heuristically optimal variable ordering is close
to the time for the true optimal variable ordering.
\begin{center}
\begin{tabular}{c|c|c}
Variable Order & Time & Size of Gr\"obner Basis \\ \hline
$x > y > z$ & 1:15.10 & 6 \\ % (49)
$x > z > y$ & 0:02.85 & 7 \\ % (36)
$y > x > z$ & 2:19.45 & 7 \\ % (69)
$y > z > x$ & 2:16.09 & 7 \\ % (67)
$z > x > y$ & 0:05.91 & 8 \\ % (37)
$z > y > x$ & 5:44.38 & 8 \\ \hline % (65)
\end{tabular}
\end{center}
\end{example}

\subsection{Logged Gr\"obner Bases} \label{LGB}

In some situations, such as in the algorithm for the
Gr\"obner Walk, it is desirable to be able to express each member of
a Gr\"obner Basis in terms of members of the original basis from which
the Gr\"obner Basis was computed.
When we have such representations, our Gr\"obner Basis is said to be a
\index{logged Gr\"obner basis}
{\it Logged Gr\"obner Basis}. \index{Gr\"obner basis!logged}

%\newpage
\begin{defn}
Let $G = \{g_1, \hdots, g_p\}$ be
a Gr\"obner Basis computed from
an initial basis $F = \{f_1, \hdots, f_m\}$. We say that $G$ is a
{\it Logged Gr\"obner Basis} if, for each $g_i \in G$,
we have an explicit expression of the form
$$g_i = \sum_{\alpha=1}^{\beta} t_{\alpha}f_{k_{\alpha}},$$
where the $t_{\alpha}$ are terms and
$f_{k_{\alpha}} \in F$ for all $1 \leqslant \alpha \leqslant \beta$.
\end{defn}

\begin{prop} \label{LGBC}
Given a finite basis $F = \{f_1, \hdots, f_m\}$,
it is always possible to compute a Logged Gr\"obner Basis for $F$.
\end{prop}
\begin{pf}
We are required to prove that every polynomial added to the input basis
$F = \{f_1, \hdots, f_m\}$
during Buchberger's algorithm has a representation in terms of
members of $F$. But any such polynomial must be
% any polynomial added to $F$ must be
a reduced S-polynomial, so it follows that the first
polynomial $f_{m+1}$ added to $F$ will always have the form
$$f_{m+1} = \mathrm{S\mbox{-}pol}(f_i, f_j)
  - \sum_{\alpha=1}^{\beta} t_{\alpha}f_{k_{\alpha}},$$
where $f_i, f_j, f_{k_{\alpha}} \in F$ and the $t_{\alpha}$ are terms. This
expression clearly gives a representation of our new polynomial in
terms of members of $F$, and by induction (using substitution)
it is also clear that each subsequent polynomial
added to $F$ will also have a representation in terms of members of $F$.
% Therefore it is always possible to compute a Logged Gr\"obner Basis.
\end{pf}

\begin{example}
% CONTINUATION OF EXAMPLE 3.3.2
Let $F := \{f_1, f_2, f_3\} = \{xy-z, \, 2x+yz+z, \, x+yz\}$
generate an ideal over the polynomial
ring $\mathbb{Q}[x, y, z]$, and let the monomial ordering be Lex.
In obtaining a Gr\"obner Basis for $F$ using Buchberger's
algorithm, three new polynomials are added to $F$, giving a
Gr\"obner Basis $G := \{g_1, g_2, g_3, g_4, g_5, g_6\} =
\{xy-z, \, 2x+yz+z, \, x+yz, \, -\frac{1}{2}yz+\frac{1}{2}z,
\, -2z^2, \, -2z\}$.
These three new polynomials are obtained from the S-polynomials
$\mathrm{S\mbox{-}pol}(2x+yz+z, x+yz)$,
$\mathrm{S\mbox{-}pol}(xy-z, -\frac{1}{2}yz+\frac{1}{2}z)$ and
$\mathrm{S\mbox{-}pol}(xy-z, 2x+yz+z)$ respectively:
\begin{eqnarray*}
\mathrm{S\mbox{-}pol}(2x+yz+z, \: x+yz) & = &
\frac{1}{2}\left(2x+yz+z\right) - (x+yz) \\
& = & -\frac{1}{2}yz+\frac{1}{2}z; \\[4mm]
\mathrm{S\mbox{-}pol}\left(xy-z, \: 
-\frac{1}{2}yz+\frac{1}{2}z\right) & = &
z(xy-z) + 2x\left(-\frac{1}{2}yz+\frac{1}{2}z\right) \\
& = & xz-z^2 \\
& \rightarrow_{f_2} & xz-z^2 - \frac{1}{2}z\left(2x+yz+z\right) \\
& = & -\frac{1}{2}yz^2 - \frac{3}{2}z^2 \\
& \rightarrow_{g_4} & -\frac{1}{2}yz^2
  - \frac{3}{2}z^2 - z\left(-\frac{1}{2}yz+\frac{1}{2}z\right) \\
& = & -2z^2; \\[4mm]
\mathrm{S\mbox{-}pol}(xy-z, \: 2x+yz+z) & = &
(xy-z) - \frac{1}{2}y\left(2x+yz+z\right) \\
& = & -\frac{1}{2}y^2z - \frac{1}{2}yz - z \\
& \rightarrow_{g_4} & -\frac{1}{2}y^2z - \frac{1}{2}yz - z 
- y\left(-\frac{1}{2}yz+\frac{1}{2}z\right) \\
& = & -yz - z \\
& \rightarrow_{g_4} & -yz - z - 
2\left(-\frac{1}{2}yz+\frac{1}{2}z\right) \\
& = & -2z.
\end{eqnarray*}
These reductions enable us to give the following
Logged Gr\"obner Basis for $F$.
% the reductions of which
% (shown below) enable us to give the Logged Gr\"obner Basis for $F$ (in the
% table shown immediately after the S-polynomial reductions).
% These reductions give the logged representations of each
% of our Gr\"obner Basis elements.
\begin{center}
\begin{tabular}{l|l} \label{logtab}
Member of $G$ & Logged Representation \\
\hline
$g_1 = xy-z$    & $f_1$ \\
$g_2 = 2x+yz+z$ & $f_2$ \\
$g_3 = x+yz$    & $f_3$ \\
$g_4 = -\frac{1}{2}yz+\frac{1}{2}z$   & $\frac{1}{2}f_2 - f_3$ \\
$g_5 = -2z^2$   & $zf_1 + (x-z)f_2 + (-2x+z)f_3$ \\
$g_6 = -z$      & $f_1 + (-y-1)f_2 + (y+2)f_3$ \\ \hline
\end{tabular}
\end{center}
\end{example}

%
% Chapter 3
% Author: Gareth Evans
% Last Modified: 2nd February 2006
%

\chapter{Noncommutative Gr\"obner Bases} \label{ChNCGB}
% \thispagestyle{empty}

% \typeout{INSERT PARAGRAPH ON WHAT I DO IN THIS CHAPTER.}

Once the potential of Gr\"obner Basis theory started
to be realised in the 1970's,
it was only natural to try to generalise the theory
to related areas such as noncommutative polynomial rings.
In 1986, \index{Teo Mora} Teo Mora 
\index{Mora, Teo} published a paper
\cite{mora86} giving an algorithm for constructing a
noncommutative Gr\"obner Basis. This work built upon
the work of George Bergman;
\index{George Bergman} in 
\index{Bergman, George} particular his
``diamond lemma for ring theory'' \cite{Bergman78}.

In this chapter, we will describe Mora's algorithm and the
theory behind it, in many ways giving a `noncommutative
version' of the previous chapter. This means that some material
from the previous chapter will be duplicated; this however will be
justified when the subtle differences between the cases
becomes apparent, differences that are all too often
overlooked when an `easy generalisation' is made!
% In order to differentiate
% between the Gr\"obner Basis algorithms in the commutative
% and noncommutative cases, we will refer to the noncommutative
% algorithm as Mora's algorithm (although it should be noted
% that many researchers worked on the theory of noncommutative
% Gr\"obner Bases, most notably George Bergman \cite{Bergman78}).
% We will discover that the main differences between the two theories
% will be in the definition of an S-polynomial and in the termination
% of the algorithm.

As in the previous chapter, we will consider the theory from the point
of view of S-polynomials, in particular defining a noncommutative Gr\"obner
Basis as a set of polynomials for which the S-polynomials all reduce to zero.
At the end of the chapter, in order to give a flavour of a noncommutative
Gr\"obner Basis program, we will give an extended example
of the computation of a noncommutative Gr\"obner Basis, taking
advantage of some of the improvements to Mora's algorithm such as
Buchberger's criteria and selection strategies.

\section{Overlaps}

For a (two-sided) ideal $J$ over a noncommutative polynomial ring,
the concept of a Gr\"obner Basis for $J$ remains the same:
it is a set of polynomials $G$ generating $J$ such that
remainders with respect to $G$ are unique. How we obtain
that Gr\"obner Basis also remains the same (we add S-polynomials
to an initial basis as required); the difference comes in the
definition of an S-polynomial.

Recall (from Section \ref{SC3}) that the purpose
of an S-polynomial $\mathrm{S\mbox{-}pol}(p_1, p_2)$ is to ensure
that any polynomial $p$ reducible by both $p_1$ and $p_2$ has a
unique remainder when divided by a set of polynomials
containing $p_1$ and $p_2$. In the commutative case, 
there is only one way to divide $p$ by $p_1$ or $p_2$
(giving reductions $p - t_1p_1$ or $p - t_2p_2$ respectively,
where $t_1$ and $t_2$ are terms); this
means that there is only one S-polynomial for
each pair of polynomials. In the noncommutative case however, a
polynomial may divide another polynomial in many different ways
(for example the polynomial $xyx - z$ divides the polynomial
$xyxyx + 4x^2$ in two different ways, giving reductions
$zyx + 4x^2$ and $xyz + 4x^2$).
For this reason, we do not have a fixed number of S-polynomials for each
pair $(p_1, p_2)$ of polynomials in the noncommutative case -- that
number will depend on the number of {\it overlaps} between
the lead monomials of $p_1$ and $p_2$.

In order to explain what an overlap is, we first need the following
preliminary definitions allowing us to select a particular part of a
noncommutative monomial. 
% We will then define the three different types of overlap.

\begin{defn}
Consider a monomial $m$ of degree $d$ over a noncommutative
polynomial ring $\mathcal{R}$.
\begin{itemize}
\item
\index{prefix}
Let $\PRE(m, i)$ denote the prefix of $m$ of
degree $i$ (where $1 \leqslant i \leqslant d$). For example,
$\PRE(x^2yz, 3) = x^2y$; $\PRE(zyx^2, 1) = z$ and
$\PRE(y^2zx, 4) = y^2zx$.
\item
\index{suffix}
Let $\SUFF(m, i)$ denote the suffix of $m$ of
degree $i$ (where $1 \leqslant i \leqslant d$). For example,
$\SUFF(x^2yz, 3) = xyz$; $\SUFF(zyx^2, 1) = x$ and
$\SUFF(y^2zx, 4) = y^2zx$.
\item
\index{subword}
Let $\SUB(m, i, j)$ denote the subword of
$m$ starting at position $i$ and finishing at position
$j$ (where $1 \leqslant i \leqslant j \leqslant d$).
For example, $\SUB(zyx^2, 2, 3) = yx$; $\SUB(zyx^2, 3, 3) = x$
and $\SUB(y^2zx, 1, 4) = y^2zx$.
\end{itemize}
\end{defn}

\begin{defn} \label{ov-def} \index{overlap}
Let $m_1$ and $m_2$ be two monomials over a noncommutative
polynomial ring $\mathcal{R}$ with
respective degrees $d_1 \geqslant d_2$. We say that $m_1$ and $m_2$
{\it overlap} if any of the following conditions are satisfied.
\begin{enumerate}[(a)]
\item
$\PRE(m_1, i) = \SUFF(m_2, i)$ ($1 \leqslant i < d_2$);
\item
$\SUB(m_1, i, i+d_2-1) = m_2$ ($1 \leqslant i \leqslant d_1-d_2+1$);
\item
$\SUFF(m_1, i) = \PRE(m_2, i)$ ($1 \leqslant i < d_2$).
\end{enumerate}
We will refer to the above overlap types as being
prefix, subword and suffix overlaps respectively; we
can picture the overlap types as follows.
\begin{center}
\begin{tabular}{ccc}
Prefix & Subword & Suffix \\[1mm]
% $\xymatrix @R=0.5pc{
% & \ar@{-}[rrr]^*+{m_1} &&&& \ar@{-}[rrr]^*+{m_1} &&&&
% \ar@{-}[rrr]^*+{m_1} &&& \\
% \ar@{-}[rrr]_*+{m_2} &&&&&& \ar@{-}[r]_*+{m_2} &&&&
% \ar@{-}[rrr]_*+{m_2} &&& \\
% }$ \\
$\xymatrix @R=0.5pc{
& \ar@{<->}[rrr]^*+{m_1} &&& \\
\ar@{<->}[rrr]_*+{m_2} &&&
}$
&
$\xymatrix @R=0.5pc{
\ar@{<->}[rrr]^*+{m_1} &&& \\
& \ar@{<->}[r]_*+{m_2} &
}$
&
$\xymatrix @R=0.5pc{
\ar@{<->}[rrr]^*+{m_1} &&& \\
& \ar@{<->}[rrr]_*+{m_2} &&&
}$
\end{tabular}
\end{center}
\end{defn}
\begin{remark}
We have defined the cases where $m_2$ is a prefix or a suffix
of $m_1$ to be subword overlaps.
\end{remark}

\begin{prop} \label{B1NC}
Let $p$ be a polynomial over a noncommutative polynomial
ring $\mathcal{R}$ that is divisible by two
polynomials $p_1, p_2 \in \mathcal{R}$, so that $\ell_1\LM(p_1)r_1 =
\LM(p) = \ell_2\LM(p_2)r_2$ for some monomials
$\ell_1, \ell_2, r_1, r_2$. As positioned in $\LM(p)$,
% Let $p$, $p_1$, $p_2$, $\ell_1$, $\ell_2$, $r_1$ and $r_2$
% be as in Remark \ref{TwoTabs}. If
if $\LM(p_1)$ and $\LM(p_2)$ do not overlap,
then no matter which of the two reductions of $p$ we apply first,
we can always obtain a common remainder.
\end{prop}
\begin{pf}
We picture the situation as follows ($u$ is a monomial).
$$\xymatrix @R=2pc @C = 4pc{
&& \ar@{<->}[r]^*+{u} & \\
\ar@{<->}[r]^*+{\ell_1} & \ar@{<->}[r]^*+{\LM(p_1)}
& \ar@{<->}[rrr]^*+{r_1} &&& \\
\ar@{<->}[rrrrr]^*+{\LM(p)} &&&&& \\
\ar@{<->}[rrr]^*+{\ell_2} &&& \ar@{<->}[r]^*+{\LM(p_2)}
& \ar@{<->}[r]^*+{r_2} &
}$$
We construct the common remainder by using
$p_2$ to divide the remainder we obtain by dividing 
$p$ by $p_1$ (and vice versa).
% (for simplicity, we assume that all
% terms have unit lead coefficients):
% $\ell_1(p_1-\LT(p_1))u(p_2-\LT(p_2))r_2$:
\begin{center}
\begin{tabular}{c}
Reduction by $p_1$ first \\ \hline
\begin{tabular}{ccl}
$p$ & $\rightarrow$ & $p - (\LC(p)\LC(p_1)^{-1})\ell_1p_1r_1$ \\
& = & $(p - \LT(p)) - (\LC(p)\LC(p_1)^{-1})\ell_1(p_1 - \LT(p_1))r_1$ \\
& = & $(p - \LT(p)) 
      - (\LC(p)\LC(p_1)^{-1})\ell_1(p_1 - \LT(p_1))u\LM(p_2)r_2$ \\
& $\stackrel{\ast}{\rightarrow}$ & $(p - \LT(p)) -
(\LC(p)\LC(p_1)^{-1}\LC(p_2)^{-1})\ell_1(p_1 
- \LT(p_1))u(p_2 - \LT(p_2))r_2$
\end{tabular} \\ \hline
\end{tabular}

\begin{tabular}{c}
Reduction by $p_2$ first \\ \hline
\begin{tabular}{ccl}
$p$ & $\rightarrow$ & $p - (\LC(p)\LC(p_2)^{-1})\ell_2p_2r_2$ \\
& = & $(p - \LT(p)) - (\LC(p)\LC(p_2)^{-1})\ell_2(p_2 - \LT(p_2))r_2$ \\
& = & $(p - \LT(p)) 
- (\LC(p)\LC(p_2)^{-1})\ell_1\LM(p_1)u(p_2 - \LT(p_2))r_2$ \\
& $\stackrel{\ast}{\rightarrow}$ & $(p - \LT(p)) -
(\LC(p)\LC(p_1)^{-1}\LC(p_2)^{-1})\ell_1(p_1 
- \LT(p_1))u(p_2 - \LT(p_2))r_2$
\end{tabular} \\ \hline
\end{tabular}
\end{center}
\end{pf}

Let $p$, $p_1$, $p_2$, $\ell_1$, $\ell_2$, $r_1$ and $r_2$
be as in Proposition \ref{B1NC}.
% Let $p$ be a polynomial over a noncommutative polynomial
% ring $\mathcal{R}$ that is divisible by two
% polynomials $p_1, p_2 \in \mathcal{R}$, so that $\ell_1\LM(p_1)r_1 =
% \LM(p) = \ell_2\LM(p_2)r_2$ for some monomials
% $\ell_1, \ell_2, r_1, r_2$. 
As positioned in $\LM(p)$, in general
the lead monomials of $p_1$ and $p_2$ may or may not overlap,
giving four different possibilities, each of which is illustrated
by an example in the following table.
\begin{center}
\begin{tabular}{c||c|c|c||c|c|c||c}
$\LM(p)$ & $\ell_1$ & $\LM(p_1)$ & $r_1$ & $\ell_2$ & $\LM(p_2)$ & $r_2$
& Overlap? \\
\hline
$x^2yzxy^3$ & $x^2yz$ & $xy^3$ & $1$ & $x^2y$ & $zx$ & $y^3$
& Prefix overlap \\
$x^2yzxy^3$ & $x$ & $xyzxy$ & $y^2$ & $x^2$ & $yzx$ & $y^3$
& Subword overlap \\
$x^2yzxy^3$ & $x$ & $xyz$ & $xy^3$ & $x^2y$ & $zx$ & $y^3$
& Suffix overlap \\
$x^2yzxy^3$ & $x^2$ & $y$ & $zxy^3$ & $x^2yz$ & $xy^2$ & $y$
& No overlap \\ \hline
\end{tabular}
\end{center}
In the cases that $\LM(p_1)$ and $\LM(p_2)$ do overlap,
we are not guaranteed to be able to obtain a common
remainder when we divide $p$ by both $p_1$ and $p_2$.
To counter this, we introduce (as in the commutative case)
an S-polynomial into our dividing set to ensure a common
remainder, requiring one S-polynomial for every
possible way that $\LM(p_1)$ and $\LM(p_2)$ overlap, including
\index{overlap!self}
{\it self overlaps} \index{self overlap}
(where $p_1 = p_2$, for example $\PRE(xyx, 1)
= \SUFF(xyx, 1)$).

\begin{defn}
Let the lead monomials of two polynomials $p_1$ and $p_2$
overlap in such a way that $\ell_1\LM(p_1)r_1$ $=$
$\ell_2\LM(p_2)r_2$, where $\ell_1, \ell_2, r_1$ and $r_2$
are monomials chosen so that at least one of $\ell_1$ and
$\ell_2$  and at least one of
$r_1$ and $r_2$ is equal to the unit monomial.
The {\it S-polynomial} \index{S-polynomial} associated with this overlap
is given by the expression
$$\mathrm{S\mbox{-}pol}(\ell_1, p_1, \ell_2, p_2) =
c_{1}\ell_{1}p_1r_{1} - c_{2}\ell_{2}p_2r_{2},$$
where $c_1 = \LC(p_2)$ and $c_2 = \LC(p_1)$.
\end{defn}

\begin{remark}
The monomials $\ell_1$ and $\ell_2$ are included in the notation 
$\mathrm{S\mbox{-}pol}(\ell_1, p_1, \ell_2, p_2)$
in order to differentiate between distinct
S-polynomials involving the two polynomials $p_1$ and $p_2$
(there is no need to include $r_1$ and $r_2$ in the
notation because $r_1$ and $r_2$ are uniquely determined by
$\ell_1$ and $\ell_2$ respectively).
\end{remark}

\begin{example}
Consider the polynomial $p := xyz + 2y$ and the set of polynomials
$P := \{p_1, p_2\} = \{xy - z, yz - x\}$,
all polynomials being ordered by DegLex and 
originating from the polynomial ring
$\mathbb{Q}\langle x, y, z\rangle$. We see that $p$
is divisible (in one way) by both of the polynomials in $P$,
giving remainders $z^2 + 2y$ and $x^2 + 2y$ respectively, both of
which are irreducible by $P$. It follows that $p$ does not have a
unique remainder with respect to $P$.

Because there is only one overlap involving the
lead monomials of $p_1$ and $p_2$, namely $\SUFF(xy, 1) = \PRE(yz, 1)$,
there is only one S-polynomial for the set $P$, which is the
polynomial $(xy-z)z - x(yz-x) = x^2-z^2$. When we add this polynomial
to the set $P$, we see that the remainder of $p$ with respect to
the enlarged $P$ is now unique, as the remainder of $p$ with
respect to $p_2$ (the polynomial $x^2+2y$)
is now reducible by our new polynomial,
giving a new remainder $z^2+2y$ which agrees with the remainder of
$p$ with respect to $p_1$.
\end{example}

Let us now give a definition of a noncommutative Gr\"obner Basis
in terms of S-polynomials.

\begin{defn} \label{grob-defn-noncom}
\index{Gr\"obner basis!noncommutative}
Let \index{noncommutative Gr\"obner basis}
$G = \{g_1, \hdots, g_m\}$ 
\index{basis!Gr\"obner} be a basis for an ideal $J$ over
a noncommutative polynomial ring $\mathcal{R} = R\langle x_1,
\hdots, x_n\rangle$. If all the S-polynomials involving members
of $G$ reduce to zero using $G$, then $G$ is a 
{\it noncommutative Gr\"obner Basis} for $J$.
\end{defn}

\begin{thm}
\index{unique remainder}
Given \index{remainder!unique}
any polynomial $p$ over a polynomial ring
$\mathcal{R} = R\langle x_1, \hdots, x_n\rangle$, the remainder
of the division of $p$ by a basis $G$ for an ideal $J$ in $\mathcal{R}$
is unique if and only if $G$ is a Gr\"obner Basis.
\end{thm}
\begin{pf}
($\Rightarrow$)
Following the proof of Theorem \ref{URC}, we need to show
that the division process is
{\it locally confluent}, \index{locally confluent} that is
if there are polynomials $f$, $f_1$, $f_2 \in \mathcal{R}$ with
$f_1 = f-\ell_1g_1r_1$ and $f_2 = f-\ell_2g_2r_2$ for terms
$\ell_1, \ell_2, r_1, r_2$
and $g_1, g_2 \in G$, then there exists a polynomial
$f_3 \in \mathcal{R}$ such that both $f_1$ and $f_2$ reduce to $f_3$.
As before, this is equivalent to showing that the polynomial
$f_2-f_1 = \ell_1g_1r_1 - \ell_2g_2r_2$ reduces to zero.

If $\LT(\ell_1g_1r_1) \neq \LT(\ell_2g_2r_2)$,
% There are two cases to deal with, 
% $\LT(\ell_1g_1r_1) \neq \LT(\ell_2g_2r_2)$
% and $\LT(\ell_1g_1r_1) = \LT(\ell_2g_2r_2)$. 
% In the first case, notice that
then the remainders $f_1$ and $f_2$ are obtained by cancelling off different
terms of the original $f$ (the reductions of $f$ are {\it disjoint}),
so it is possible, assuming (without loss of generality) that
$\LT(\ell_1g_1r_1) > \LT(\ell_2g_2r_2)$,
to directly reduce the polynomial 
$f_2-f_1 = \ell_1g_1r_1 - \ell_2g_2r_2$ in the
following manner: $\ell_1g_1r_1 - \ell_2g_2r_1 \rightarrow_{g_1}
 -\ell_2g_2r_2 \rightarrow_{g_2} 0$.

On the other hand, if $\LT(\ell_1g_1r_1) = \LT(\ell_2g_2r_2)$, then
the reductions of $f$ are not disjoint
(as the same term $t$ from $f$ is cancelled off during both reductions),
so that the term $t$ does not appear in the polynomial
$\ell_1g_1r_1 - \ell_2g_2r_2$. However, the monomial
$\LM(t)$ must contain the monomials $\LM(g_1)$ and $\LM(g_2)$ as subwords
if both $g_1$ and $g_2$ cancel off the term $t$, so it follows
that $\LM(g_1)$ and $\LM(g_2)$ will either overlap or not overlap
in $\LM(t)$. If they do not overlap, then we know from Proposition \ref{B1NC}
that $f_1$ and $f_2$ will have a common remainder
($f_1 \stackrel{\ast}{\longrightarrow} f_3$ and 
$f_2 \stackrel{\ast}{\longrightarrow} f_3$), so that
$f_2 - f_1 \stackrel{\ast}{\longrightarrow} f_3 - f_3 = 0$. 
Otherwise, because of the overlap between $\LM(g_1)$ and $\LM(g_2)$, the 
polynomial $\ell_1g_1r_1 - \ell_2g_2r_2$ 
will be a multiple of an S-polynomial,
say $\ell_1g_1r_1 - \ell_2g_2r_2 =
\ell_3(\mathrm{S\mbox{-}pol}(\ell'_1, g_1, \ell'_2, g_2))r_3$ for some terms
$\ell_3, r_3$ and some monomials $\ell'_1, \ell'_2$.
But $G$ is a Gr\"obner Basis, so the S-polynomial
$\mathrm{S\mbox{-}pol}(\ell'_1, g_1, \ell'_2, g_2)$ will reduce to zero, and
hence by extension the polynomial $\ell_1g_1r_1 - \ell_2g_2r_2$
will also reduce to zero.

($\Leftarrow$)
As all S-polynomials are members of the ideal $J$, to
complete the proof it is sufficient to show that
there is always a reduction path of an arbitrary
member of the ideal that leads to a zero remainder
(the uniqueness of remainders will then imply
that members of the ideal always reduce to zero).
Let $f \in J = \langle G \rangle$. Then, by definition, there exist
$g_i \in G$ (not necessarily all different)
and terms $\ell_i, r_i \in \mathcal{R}$ (where $1 \leqslant i \leqslant j$)
such that $$f = \sum_{i=1}^j \ell_ig_ir_i.$$
We proceed by induction on $j$. If $j = 1$, then $f = \ell_1g_1r_1$,
and it is clear that we can use $g_1$ to reduce $f$ to give
a zero remainder ($f \rightarrow_{g_1} f - \ell_1g_1r_1 = 0$).
Assume that the result is true for $j = k$, and let us
look at the case $j = k+1$, so that
$$f = \left(\sum_{i=1}^{k} \ell_ig_ir_i \right) + \ell_{k+1}g_{k+1}r_{k+1}.$$
By the inductive hypothesis, $\sum_{i=1}^k \ell_ig_ir_i$ is a member
of the ideal that reduces to zero. The polynomial $f$ therefore
reduces to the polynomial $f' := \ell_{k+1}g_{k+1}r_{k+1}$, and we can 
now use $g_{k+1}$ to reduce $f'$
to give a zero remainder ($f' \rightarrow_{g_{k+1}} f'
- \ell_{k+1}g_{k+1}r_{k+1} = 0$).
\end{pf}

\begin{remark}
The above Theorem forms part of Bergman's Diamond Lemma
\cite[Theorem 1.2]{Bergman78}.
\end{remark}

\section{Mora's Algorithm}
\index{Mora's algorithm}

Let us now consider the following pseudo code representing
Mora's algorithm for computing noncommutative Gr\"obner Bases
\cite{mora86}.

\begin{algorithm}
\setlength{\baselineskip}{3.5ex}
\caption{Mora's Noncommutative Gr\"obner Basis Algorithm}
\label{noncom-buch}
\begin{algorithmic}
\vspace*{2mm}
\REQUIRE{A Basis $F = \{f_1, f_2, \hdots, f_m\}$ for an ideal
         $J$ over a noncommutative polynomial ring
         $R\langle x_1, \hdots x_n\rangle$; an admissible
         monomial ordering $\mathrm{O}$.}
\ENSURE{A Gr\"obner Basis $G = \{g_1, g_2, \hdots, g_p\}$ for $J$
        (in the case of termination).}
\vspace*{1mm}
\STATE
Let $G = F$ and let $A = \emptyset$; \\
For each pair of polynomials $(g_i, g_j)$ in $G$ ($i \leqslant j$), add an
S-polynomial $\mathrm{S\mbox{-}pol}(\ell_1, g_i, \ell_2, g_j)$ to $A$
for each overlap $\ell_1\LM(g_i)r_1 = \ell_2\LM(g_j)r_2$
between the lead monomials of $\LM(g_i)$ and $\LM(g_j)$.
% overlaps between the lead monomials 
% $\LM(g_i)$ and $\LM(g_j)$ to a list $A$,
% where the $k$-th entry in the list is a six-tuple
% $T_k = (g_k, g'_k, c_{k}\ell_{k}, r_{k}, c'_{k}\ell'_{k}, r'_{k})$
% such that the following holds: $c_k\ell_{k}(\LT(g_k))r_{k} =
% c'_k\ell'_{k}(\LT(g'_k))r'_{k}$;
\WHILE{($A$ is not empty)}
\STATE
%   Remove the first entry from $A$ and compute the S-polynomial \\
%   $s_k := \:\mathrm{S\mbox{-}pol}(T_k) = c_{k}\ell_{k}g_kr_{k} -
%   c'_{k}\ell'_{k}g'_kr'_{k}$; \\
%   Reduce $s_{k}$ with respect to the current basis
%   (using the noncommutative
%   division algorithm); If $s_{k}$ reduces to zero then do nothing;
%   otherwise if $s_{k}$ reduces to $r_{k} \neq 0$ add $r_{k}$ to
%   $G$ and add all
%   the overlaps between $r_{k}$ and elements of $G$ to $A$;
  Remove the first entry $s_1$ from $A$; \\
  $s'_1 = \Rem(s_1, G)$; % (Algorithm \ref{noncom-div})
  % Reduce $s_1$ with respect to $G$ (with Algorithm \ref{noncom-div}); \\
  \IF{($s'_1 \neq 0$)}
  \STATE
  % If $s_1$ reduces to zero then do nothing; 
  % otherwise if $s_1$ reduces to $r_1 \neq 0$ 
  Add $s'_1$ to $G$ and then (for all $g_i \in G$) add all 
  the S-polynomials of the form
  $\mathrm{S\mbox{-}pol}(\ell_1, g_i, \ell_2, s'_1)$ to $A$; % ($g_i \in G$);
  \ENDIF
\ENDWHILE
\STATE
{\bf return} $G$;
\end{algorithmic}
\vspace*{1mm}
\end{algorithm}

Structurally, Mora's algorithm is virtually identical
to Buchberger's algorithm, in that we compute and reduce each
S-polynomial in turn; we add a reduced S-polynomial to our
basis if it does not reduce to zero; and we continue
until all S-polynomials reduce to zero --- exactly as in
Algorithm \ref{com-buch}. Despite this, there are major
differences from an implementation standpoint, not
least in the fact that noncommutative polynomials are
much more difficult to handle on a computer; and
noncommutative S-polynomials need more complicated data
structures. This may explain why implementations of the
noncommutative Gr\"obner Basis algorithm are currently sparser 
than those for the commutative algorithm; and also why such
implementations often impose restrictions on the problems
that can be handled --- Bergman \cite{Bergman98} for instance
only allows input bases which are homogeneous.

\subsection{Termination}

In the commutative case, Dickson's Lemma and Hilbert's Basis Theorem 
allow us to prove that Buchberger's algorithm always terminates for all
possible inputs. It is a fact however that Mora's algorithm
does not terminate for all possible inputs (so that an ideal may have an
infinite Gr\"obner Basis in general) because there is no analogue of
Dickson's Lemma for noncommutative monomial ideals.

% Because there is no analogue of Dickson's Lemma for
% noncommutative monomial ideals, it is a fact that Mora's algorithm
% does not terminate for all possible inputs, so that an ideal will have an
% infinite Gr\"obner Basis in general.

\begin{prop}
Not all noncommutative monomial ideals are finitely generated.
\end{prop}
\begin{pf}
Assume to the contrary that all noncommutative monomial
ideals are finitely generated, and consider an ascending chain
of such ideals $J_1 \subseteq J_2 \subseteq \cdots$. By our
assumption, the ideal $J = \cup J_i$ (for $i \geqslant 1$) will be
finitely generated, which means that there must be some
$k \geqslant 1$ such that $J_k = J_{k+1} = \cdots$. For a
counterexample, let $\mathcal{R} = \mathbb{Q}\langle x, y \rangle$
be a noncommutative polynomial ring,
and define $J_i$ (for $i \geqslant 1$) to be the ideal in $\mathcal{R}$
generated by the set of monomials $\{xyx, xy^2x, \hdots, xy^ix\}$.
% $s^it^iu$ ($1 \leqslant i \leqslant j$), where
% $s := xy$, $t := yx$ and $u := x^2$.
Because no member of this set is a multiple of any other member
of the set, it is clear that there cannot be a $k \geqslant 1$
such that $J_k = J_{k+1} = \cdots$ because $xy^{k+1}x \in J_{k+1}$
and $xy^{k+1}x \notin J_k$ for all $k \geqslant 1$.
\end{pf}

Another way of explaining why Mora's algorithm does not terminate
comes from considering the link between noncommutative Gr\"obner Bases
and the \index{Knuth-Bendix algorithm} Knuth-Bendix Critical Pairs
Completion Algorithm for monoid rewrite systems \cite{KB}, an algorithm
that attempts to find a complete rewrite system for
any given monoid presentation. Because Mora's algorithm can be
used to emulate the Knuth-Bendix algorithm (for the details,
see for example \cite{Hey00a}), if we assume that Mora's algorithm
always terminates, then we have found a way to solve the
\index{word problem} {\it word problem} for monoids
(so that we can determine whether
any word in a given monoid is equal to the identity word);
this however contradicts the fact that the word problem is actually an
unsolvable problem (so that it is impossible to define
an algorithm that can tell whether two words in a given monoid 
are identical).

\section{Reduced Gr\"obner Bases}

\begin{defn} \label{RGBNC}
Let $G = \{g_1, \hdots, g_p\}$ be a Gr\"obner Basis for an ideal
over a polynomial ring $R\langle x_1, \hdots, x_n\rangle$. $G$ is
a \index{Gr\"obner basis!reduced} {\it reduced}
\index{reduced Gr\"obner basis}
Gr\"obner Basis if the following conditions are satisfied.
\begin{enumerate}[(a)]
\item
$\LC(g_i) = 1_R$ for all $g_i \in G$.
\item
No term in any polynomial $g_i \in G$
is divisible by any $\LT(g_j)$, $j \neq i$.
\end{enumerate}
\end{defn}

\begin{thm} \label{MinNC}
If there exists a Gr\"obner Basis $G$ for an ideal $J$
over a noncommutative polynomial ring, then $J$ has a
unique reduced Gr\"obner Basis.
\end{thm}
\begin{pf}
{\it Existence.}
We claim that the following
procedure transforms $G$ into a reduced Gr\"obner Basis $G'$.
\begin{enumerate}[(i)]
\item
Multiply each $g_i \in G$ by $\LC(g_i)^{-1}$.
\item
Reduce each $g_i \in G$ by $G\setminus \{g_i\}$, removing
from $G$ all polynomials that reduce to zero.
\end{enumerate}
It is clear that $G'$ satisfies the conditions of
Definition \ref{RGBNC}, so it remains to show that
$G'$ is a Gr\"obner Basis, which we shall do by showing that
the application of each step of instruction (ii) above produces
a basis which is still a Gr\"obner Basis.

Let $G = \{g_1, \hdots, g_p\}$ be a Gr\"obner Basis, and
let $g'_i$ be the reduction of an arbitrary $g_i \in G$
with respect to $G \setminus \{g_i\}$, carried out as follows
(the $\ell_k$ and the $r_k$ are terms).
\begin{equation} \label{gdashredNC}
g'_i = g_i - \sum_{k=1}^{\kappa} \ell_kg_{j_k}r_k.
\end{equation}
Set $H = (G \setminus \{g_i\}) \cup \{g'_i\}$ if
$g'_i \neq 0$, and set $H = G \setminus \{g_i\}$
if $g'_i = 0$. As $G$ is a Gr\"obner Basis,
all S-polynomials involving elements
of $G$ reduce to zero using $G$, so there are expressions
\begin{equation} \label{sprzNC}
c_b\ell_ag_ar_a - c_a\ell_bg_br_b - \sum_{u=1}^{\mu} \ell_ug_{c_u}r_u = 0
\end{equation}
for every S-polynomial $\mathrm{S\mbox{-}pol}(\ell_a, g_a, \ell_b, g_b) =
c_b\ell_ag_ar_a - c_a\ell_bg_br_b$, where $c_a = \LC(g_a)$;
$c_b = \LC(g_b)$; the $\ell_u$ and
the $r_u$ are terms (for $1 \leqslant u \leqslant \mu$); and
$g_a, g_b, g_{c_u} \in G$. To show that $H$
is a Gr\"obner Basis, we must show that all S-polynomials
involving elements of $H$ reduce to zero using $H$. 
For polynomials $g_a, g_b \in H$ not equal to $g'_i$,
we can reduce an S-polynomial of the form
$\mathrm{S\mbox{-}pol}(\ell_a, g_a, \ell_b, g_b)$
using the reduction shown in Equation
(\ref{sprzNC}), substituting for $g_i$ from 
Equation (\ref{gdashredNC})
if any of the $g_{c_u}$ in Equation (\ref{sprzNC}) are equal to $g_i$.
This gives a reduction to zero of 
$\mathrm{S\mbox{-}pol}(\ell_a, g_a, \ell_b, g_b)$
in terms of elements of $H$.

If $g'_i = 0$, our proof is complete. Otherwise
consider all S-polynomials 
$\mathrm{S\mbox{-}pol}(\ell'_i, g'_i, \ell_b, g_b)$
involving the pair of polynomials $(g'_i, g_b)$, 
where $g_b \in G\setminus\{g_i\}$.
We claim that there exists an S-polynomial
$\mathrm{S\mbox{-}pol}(\ell_1, g_i, \ell_2, g_b)
= c_b\ell_1g_ir_1 - c_i\ell_2g_br_2$ such that
$\mathrm{S\mbox{-}pol}(\ell'_i, g'_i, \ell_b, g_b)
= c_b\ell_1g'_ir_1 - c_i\ell_2g_br_2$. To prove
this claim, it is sufficient to show that $\LT(g_i) = \LT(g'_i)$.
Assume for a contradiction that $\LT(g_i) \neq \LT(g'_i)$.
It follows that during the reduction of $g_i$ we were able to
reduce its lead term, so that
$\LT(g_i) = \ell\LT(g_j)r$ for some terms 
$\ell$ and $r$ and some $g_j \in G$.
Because $\LM(g_i - \ell g_jr) < \LM(g_i)$, the
polynomial $g_i - \ell g_jr$ must reduce to zero without
using $g_i$, so that $g'_i = 0$, giving a contradiction.

It remains to show that 
$\mathrm{S\mbox{-}pol}(\ell'_i, g'_i, \ell_b, g_b)
\rightarrow_H 0$. We know that
$\mathrm{S\mbox{-}pol}(\ell_1, g_i, \ell_2, g_b)
= c_b\ell_1g_ir_1 - c_i\ell_2g_br_2
\rightarrow_G 0$, and Equation (\ref{sprzNC}) tells us that
$c_b\ell_1g_ir_1 - c_i\ell_2g_br_2 - 
\sum_{u=1}^{\mu}\ell_ug_{c_u}r_u = 0$.
Substituting for $g_i$ from Equation (\ref{gdashredNC}),
we obtain\footnote{Substitutions
for $g_i$ may also occur in the 
summation $\sum_{u=1}^{\mu}\ell_ug_{c_u}r_u$;
these substitutions have not been considered in the
displayed formulae.}
$$
c_b\ell_1\left(g'_i + \sum_{k=1}^{\kappa} \ell_kg_{j_k}r_k\right)r_1
- c_i\ell_2g_br_2 - \sum_{u=1}^{\mu}\ell_ug_{c_u}r_u = 0
$$
or
$$
c_b\ell_1g'_ir_1 - c_i\ell_2g_br_2 -
\left(\sum_{u=1}^{\mu}\ell_ug_{c_u}r_u
- \sum_{k=1}^{\kappa} c_b\ell_1\ell_kg_{j_k}r_kr_1\right) = 0,$$
which implies that $\mathrm{S\mbox{-}pol}(\ell'_i, g'_i, \ell_b, g_b)
\rightarrow_H 0$. The only other case to consider is the case of an
S-polynomial coming from a self overlap involving $\LM(g'_i)$.
But because we now know that $\LT(g'_i) = \LT(g_i)$, we can
use exactly the same argument as above to show that the
S-polynomial $\mathrm{S\mbox{-}pol}(\ell_1, g'_i, \ell_2, g'_i)$
reduces to zero using $H$ because an S-polynomial of the form
$\mathrm{S\mbox{-}pol}(\ell_1, g_i, \ell_2, g_i)$ will exist.

{\it Uniqueness.}
Assume for a contradiction that $G = \{g_1, \hdots, g_p\}$
and $H = \{h_1, \hdots, h_q\}$ are two reduced Gr\"obner
Bases for an ideal $J$, with $G \neq H$. Let $g_i$ be an arbitrary
element from $G$ (where $1 \leqslant i \leqslant p$).
Because $g_i$ is a member of the ideal,
then $g_i$ must reduce to zero using $H$ 
($H$ is a Gr\"obner Basis). This means that there must exist
a polynomial $h_j \in H$ such that $\LT(h_j) \mid \LT(g_i)$. If
$\LT(h_j) \neq \LT(g_i)$, 
then $\ell \times \LT(h_j) \times r = \LT(g_i)$ for
some monomials $\ell$ and $r$, at least one of which is not equal to the 
unit monomial. But $h_j$ is also a member of the ideal,
so it must reduce to zero using $G$. Therefore there
exists a polynomial $g_k \in G$
such that $\LT(g_k) \mid \LT(h_j)$, which implies that
$\LT(g_k) \mid \LT(g_i)$, with $k \neq i$.
This contradicts condition (b) of Definition \ref{RGBNC} so that $G$
cannot be a reduced Gr\"obner Basis for $J$ if $\LT(h_j) \neq \LT(g_i)$.
From this we deduce that each $g_i \in G$ has a corresponding
$h_j \in H$ such that $\LT(g_i) = \LT(h_j)$. Further, because $G$ and $H$ are
assumed to be reduced Gr\"obner Bases, this is a one-to-one correspondence.

It remains to show that if $\LT(g_i) = \LT(h_j)$, then $g_i = h_j$.
Assume for a contradiction that $g_i \neq h_j$ and consider the
polynomial $g_i - h_j$. Without loss of generality, assume that
$\LM(g_i - h_j)$ appears in $g_i$. Because $g_i - h_j$ is a member
of the ideal, then there is a polynomial $g_k \in G$ such that
$\LT(g_k) \mid \LT(g_i - h_j)$. But this again contradicts condition
(b) of Definition \ref{RGBNC}, as we have shown that there is a term in
$g_i$ that is divisible by $\LT(g_k)$ for some $k \neq i$. It follows
that $G$ cannot be a reduced Gr\"obner Basis if $g_i \neq h_j$,
which means that $G = H$ and therefore reduced Gr\"obner Bases are unique.
\end{pf}

As in the commutative case, we may refine the procedure for finding a 
unique reduced Gr\"obner Basis (as given in the proof of Theorem
\ref{MinNC}) by removing from the Gr\"obner Basis all polynomials
whose lead monomials are multiples of the lead monomials of other
Gr\"obner Basis elements. This leads to the definition of 
Algorithm \ref{red-noncom}.

\begin{algorithm}
\setlength{\baselineskip}{3.5ex}
\caption{The Noncommutative Unique Reduced Gr\"obner Basis Algorithm}
\label{red-noncom} \index{Gr\"obner basis!minimal} \index{minimal Gr\"obner basis}
\begin{algorithmic}
\vspace*{2mm}
\REQUIRE{A Gr\"obner Basis $G = \{g_1, g_2, \hdots, g_m\}$ 
         for an ideal $J$ over a
         noncommutative polynomial ring $R\langle x_1, \hdots x_n\rangle$; 
         an admissible monomial ordering $\mathrm{O}$.}
\ENSURE{The unique reduced Gr\"obner Basis
        $G' = \{g'_1, g'_2, \hdots, g'_p\}$ for $J$.}
\vspace*{1mm}
\STATE
$G' = \emptyset$;
\FOR{\textbf{each} $g_i \in G$}
\STATE
Multiply $g_i$ by $\LC(g_i)^{-1}$; \\
\IF{($\LM(g_i) = \ell\LM(g_j)r$ for some
monomials $\ell, r$ and some $g_j \in G$ ($g_j \neq g_i$))}
\STATE
$G = G\setminus\{g_i\}$;
\ENDIF
\ENDFOR
\FOR{\textbf{each} $g_i \in G$}
\STATE
$g'_i = \Rem(g_i, (G \setminus \{g_i\})\cup G')$; \\
$G = G\setminus \{g_i\}$; $G' = G' \cup \{g'_i\}$;
% Add $g'_i$ to $G'$, where $g'_i$ is the reduction of
% $g_i$ with respect to $G\setminus \{g_i\}$;
\ENDFOR
\STATE
{\bf return} $G'$;
\end{algorithmic}
\vspace*{1mm}
\end{algorithm}

\section{Improvements to Mora's Algorithm}

In Section \ref{ImprovC}, we surveyed some of the
numerous improvements of Buchberger's algorithm.
Let us now demonstrate that many of
these improvements can also be applied in the
noncommutative case.

\subsection{Buchberger's Criteria} \label{BCNC}

\index{Buchberger's first criterion}
In \index{criterion!Buchberger's first}
the commutative case, Buchberger's first criterion states
that we can ignore any S-polynomial
$\mathrm{S\mbox{-}pol}(f, g)$ in which
$\lcm(\LM(f), \LM(g)) = \LM(f)\LM(g)$. In the noncommutative
case, this translates as saying that we can ignore any
`S-polynomial' $\mathrm{S\mbox{-}pol}(\ell_1, f, \ell_2, g)
= \LC(g)\ell_1fr_1 - \LC(f)\ell_2gr_2$ such that $\LM(f)$ and
$\LM(g)$ do not overlap in the monomial $\ell_1\LM(f)r_1 =
\ell_2\LM(g)r_2$. We can certainly show that such an `S-polynomial'
will reduce to zero by utilising Proposition \ref{B1NC},
but we will never be able to use this result as,
by definition, an S-polynomial
is only defined when we have an overlap between $\LM(f)$ and
$\LM(g)$. It follows that an `S-polynomial' of the above type
will never occur in Mora's algorithm, and so
Buchberger's first criterion is redundant in the
noncommutative case. The same cannot be said of his second
criterion however, which certainly does improve
the efficiency of Mora's algorithm.

\begin{prop}[Buchberger's Second Criterion] \label{B2NC}
\index{Buchberger's second criterion}
Let \index{criterion!Buchberger's second}
$f$, $g$ and $h$ be three members of a finite set of
polynomials $P$ over a noncommutative polynomial ring,
and consider an S-polynomial of the form
\begin{equation} \label{B2NCE1}
\mathrm{S\mbox{-}pol}(\ell_1, f, \ell_2, g)
= c_2\ell_1fr_1 - c_1\ell_2gr_2.
\end{equation}
% where $c_1 = \LC(f)$ and $c_2 = \LC(g)$.
If $\LM(h) \mid \ell_1\LM(f)r_1$, so that
\begin{equation} \label{B2NCE2}
\ell_1\LM(f)r_1 = \ell_3\LM(h)r_3 = \ell_2\LM(g)r_2
\end{equation}
for some monomials $\ell_3, r_3$, then
$\mathrm{S\mbox{-}pol}(\ell_1, f, \ell_2, g) \rightarrow_P 0$
if all S-polynomials corresponding to overlaps
(as placed in the monomial $\ell_1\LM(f)r_1$) between $\LM(h)$ and
either $\LM(f)$ or $\LM(g)$ reduce to zero using $P$.
% if $\mathrm{S\mbox{-}pol}(\ell_1, f, \ell_3, h) \rightarrow_P 0$
% and $\mathrm{S\mbox{-}pol}(\ell_2, g, \ell_3, h) \rightarrow_P 0$.
\end{prop}

\noindent {\bf Proof (cf. \cite{Keller97}, Appendix A):}
To be able to describe an S-polynomial corresponding to an overlap
(as placed in the monomial $\ell_1\LM(f)r_1$) between $\LM(h)$ and
either $\LM(f)$ or $\LM(g)$, we introduce the following notation.
\begin{itemize}
\item
Let $\ell_{13}$ be the monomial corresponding to the common prefix of
$\ell_1$ and $\ell_3$ of maximal degree,
so that $\ell_1 = \ell_{13}\ell'_1$ and $\ell_3 = \ell_{13}\ell'_3$.
(Here, and similarly below, if there is no common prefix of
$\ell_1$ and $\ell_3$, then $\ell_{13} = 1$, $\ell'_1 = \ell_1$ and
$\ell'_3 = \ell_3$.)
\item
Let $\ell_{23}$ be the monomial corresponding to the common prefix of
$\ell_2$ and $\ell_3$ of maximal degree,
so that $\ell_2 = \ell_{23}\ell''_2$ and $\ell_3 = \ell_{23}\ell''_3$.
\item
Let $r_{13}$ be the monomial corresponding to the common suffix of
$r_1$ and $r_3$ of maximal degree,
so that $r_1 = r'_1r_{13}$ and $r_3 = r'_3r_{13}$.
\item
Let $r_{23}$ be the monomial corresponding to the common suffix of
$r_2$ and $r_3$ of maximal degree,
so that $r_2 = r''_2r_{23}$ and $r_3 = r''_3r_{23}$.
\end{itemize}
We can now manipulate Equation (\ref{B2NCE1}) as follows
(where $c_3 = \LC(h)$).
\begin{eqnarray*}
% \mathrm{S\mbox{-}pol}(\ell_1, f, \ell_2, g)
% & = & c_2\ell_1fr_1 - c_1\ell_2gr_2; \\
c_3(\mathrm{S\mbox{-}pol}(\ell_1, f, \ell_2, g))
& = & c_3c_2\ell_1fr_1 - c_3c_1\ell_2gr_2 \\
& = & c_3c_2\ell_1fr_1
      - c_1c_2\ell_3hr_3 + c_1c_2\ell_3hr_3 - c_3c_1\ell_2gr_2 \\
& = & c_2(c_3\ell_1fr_1 - c_1\ell_3hr_3)
      - c_1(c_3\ell_2gr_2 - c_2\ell_3hr_3) \\
& = & c_2(c_3\ell_{13}\ell'_1fr'_1r_{13}
      - c_1\ell_{13}\ell'_3hr'_3r_{13}) \\
&   & -~c_1(c_3\ell_{23}\ell''_2gr''_2r_{23}
      - c_2\ell_{23}\ell''_3hr''_3r_{23}) \\
& = & c_2\ell_{13}(c_3\ell'_1fr'_1 - c_1\ell'_3hr'_3)r_{13}
- c_1\ell_{23}(c_3\ell''_2gr''_2 - c_2\ell''_3hr''_3)r_{23}.
\end{eqnarray*}
As placed in $\ell_1\LM(f)r_1 = \ell_3\LM(h)r_3$,
if $\LM(f)$ and $\LM(h)$ overlap,
then the S-polynomial corresponding to this overlap is\footnote{For
completeness, we note that the S-polynomial corresponding to
the overlap can also be of the form
$\mathrm{S\mbox{-}pol}(\ell'_3, h, \ell'_1, f)$;
this (inconsequentially) swaps the first two terms
of Equation (\ref{B2NCE3}).}
$\mathrm{S\mbox{-}pol}(\ell'_1, f, \ell'_3, h)$.
Similarly, if $\LM(g)$ and $\LM(h)$ overlap as placed in
$\ell_2\LM(g)r_2 = \ell_3\LM(h)r_3$,
then the S-polynomial corresponding to this
overlap is $\mathrm{S\mbox{-}pol}(\ell''_2, g, \ell''_3, h)$.
By assumption, these S-polynomials reduce to zero using $P$, so there
are expressions
\begin{equation} \label{B2NCE3}
c_3\ell'_1fr'_1 - c_1\ell'_3hr'_3 -
\sum_{i=1}^{\alpha} u_ip_iv_i = 0
\end{equation}
and
\begin{equation} \label{B2NCE4}
c_3\ell''_2gr''_2 - c_2\ell''_3hr''_3 -
\sum_{j=1}^{\beta} u_jp_jv_j = 0,
\end{equation}
where the $u_i$, $v_i$, $u_j$ and $v_j$ are terms;
and $p_i, p_j \in P$ for all $i$ and $j$.
Using Proposition \ref{B1NC}, we can state that these
expressions will still exist even
if $\LM(f)$ and $\LM(h)$ do not overlap as placed in
$\ell_1\LM(f)r_1 = \ell_3\LM(h)r_3$;
and if $\LM(g)$ and $\LM(h)$ do not overlap as placed in
$\ell_2\LM(g)r_2 = \ell_3\LM(h)r_3$. It follows that
\begin{eqnarray*}
c_3(\mathrm{S\mbox{-}pol}(\ell_1, f, \ell_2, g))
& = & c_2\ell_{13}(c_3\ell'_1fr'_1 - c_1\ell'_3hr'_3)r_{13}
- c_1\ell_{23}(c_3\ell''_2gr''_2 - c_2\ell''_3hr''_3)r_{23} \\
& = & c_2\ell_{13}\left(\sum_{i=1}^{\alpha} u_ip_iv_i\right)r_{13}
- c_1\ell_{23}\left(\sum_{j=1}^{\beta} u_jp_jv_j\right)r_{23} \\
& = & \sum_{i=1}^{\alpha} c_2\ell_{13}u_ip_iv_ir_{13}
- \sum_{j=1}^{\beta} c_1\ell_{23}u_jp_jv_jr_{23}; \\
\mathrm{S\mbox{-}pol}(\ell_1, f, \ell_2, g)
& = & \sum_{i=1}^{\alpha} c_3^{-1}c_2\ell_{13}u_ip_iv_ir_{13}
- \sum_{j=1}^{\beta} c_3^{-1}c_1\ell_{23}u_jp_jv_jr_{23}.
\end{eqnarray*}
To conclude that the S-polynomial
$\mathrm{S\mbox{-}pol}(\ell_1, f, \ell_2, g)$ reduces to zero using $P$,
it remains to show that the algebraic expression
$- \sum_{i=1}^{\alpha} c_3^{-1}c_2\ell_{13}u_ip_iv_ir_{13}
+ \sum_{j=1}^{\beta} c_3^{-1}c_1\ell_{23}u_jp_jv_jr_{23}$
corresponds to a valid reduction of
$\mathrm{S\mbox{-}pol}(\ell_1, f, \ell_2, g)$. To do this, it is
sufficient to show that no term in either of the summations is
greater than the term $\ell_1\LM(f)r_1$
(so that $\LM(\ell_{13}u_ip_iv_ir_{13}) < \ell_1\LM(f)r_1$ and
$\LM(\ell_{23}u_jp_jv_jr_{23}) < \ell_1\LM(f)r_1$ for all $i$ and $j$).
But this follows from
Equation (\ref{B2NCE2}) and from the fact that the
reductions of the expressions $c_3\ell'_1fr'_1 - c_1\ell'_3hr'_3$
and $c_3\ell''_2gr''_2 - c_2\ell''_3hr''_3$ in Equations
(\ref{B2NCE3}) and (\ref{B2NCE4}) are valid, so that
$\LM(u_ip_iv_i) < \LM(\ell'_1fr'_1)$ and
$\LM(u_jp_jv_j) < \LM(\ell''_2gr''_2)$ for all $i$ and $j$.
{\hfill $\Box$}

\begin{remark}
The three polynomials $f$, $g$ and $h$ in the above
proposition do not necessarily have to be distinct (indeed,
$f = g = h$ is allowed) --- the only restriction is that the
S-polynomial $\mathrm{S\mbox{-}pol}(\ell_1, f, \ell_2, g)$
has to be different from the S-polynomials
$\mathrm{S\mbox{-}pol}(\ell'_1, f, \ell'_3, h)$ and
$\mathrm{S\mbox{-}pol}(\ell''_2, g, \ell''_3, h)$;
for example, if $f = h$, then we cannot have $\ell'_1 = \ell'_3$.
\end{remark}

\subsection{Homogeneous Gr\"obner Bases}

Because it is computationally more expensive to do
noncommutative polynomial arithmetic than it is to do
commutative polynomial arithmetic, gains in efficiency due
to working with homogeneous bases are even more
significant in the noncommutative case. For this reason, some systems
for computing noncommutative Gr\"obner Bases will only work
with homogeneous input bases, although (as in the commutative case)
it is still sometimes possible to use these systems on non-homogeneous input
bases by using the concepts of homogenisation, dehomogenisation and
extendible monomial orderings.

\begin{defn}
Let $p = p_0 + \cdots + p_m$ be a polynomial over the
polynomial ring $R\langle x_1, \hdots, x_n\rangle$, where
each $p_i$ is the sum of the degree $i$ terms in $p$
(we assume that $p_m \neq 0$). The
% the $p_i$ are homogeneous polynomials of degree $i$ and $p_m \neq 0$. The
{\it left homogenisation} \index{homogenisation} of $p$
with respect to a new (homogenising) variable $y$ is the polynomial
$$h_{\ell}(p) := y^mp_0 + y^{m-1}p_1 + \cdots +
yp_{m-1} + p_m;$$ % \in R\langle x_1, \hdots, x_n, y \rangle;$$
and the {\it right homogenisation} of $p$
with respect to a new (homogenising) variable $y$ is the polynomial
$$h_r(p) := p_0y^m + p_1y^{m-1} + \cdots +
p_{m-1}y + p_m.$$ % \in R\langle x_1, \hdots, x_n, y\rangle.$$
Homogenised polynomials 
belong to polynomial rings determined by where
$y$ is placed in the lexicographical ordering of the variables.
\end{defn}

\begin{defn}
The {\it dehomogenisation} \index{dehomogenisation} of a
% homogeneous
polynomial $p$ is the polynomial $d(p)$ given by
substituting $y = 1$ in $p$, where $y$ is the homogenising variable.
% $p \in R\langle x_1, \hdots, x_n, y\rangle$ is the polynomial
% $d(p) \in R\langle x_1, \hdots, x_n\rangle$ 
% given by substituting $y = 1$ in $p$.
\end{defn}

\begin{defn}
A monomial ordering $O$ is {\it extendible}
\index{monomial ordering!extendible} if, 
\index{extendible monomial ordering} given any polynomial
$p = t_1 + \cdots + t_\alpha$ ordered with respect to $O$
(where $t_1 > \cdots > t_\alpha$),
the homogenisation of $p$ preserves the order
on the terms ($t'_i > t'_{i+1}$ for all $1 \leqslant i \leqslant \alpha-1$,
where the homogenisation process maps the term $t_i \in p$
to the term $t'_i$).
\end{defn}

In the noncommutative case, an extendible monomial ordering must
specify how to homogenise a polynomial (by multiplying with the
homogenising variable on the left or on the right) as well as
stating where the new variable $y$ appears in the ordering of
the variables. Here are the conventions for those monomial orderings
defined in Section \ref{NCMO} that are extendible, assuming that
we start with a polynomial ring $R\langle x_1, \hdots, x_n\rangle$.
\begin{center}
\begin{tabular}{c|c|c}
Monomial Ordering & Type of Homogenisation &
Position of the new variable $y$ \\
& & in the ordering of the variables \\
\hline
% Lex & Right &  $y < x_i$ for all $x_i$ \\
InvLex & Right & $y < x_i$ for all $x_i$ \\
% RevLex & Left & $y < x_i$ for all $x_i$ \\
DegLex & Left & $y < x_i$ for all $x_i$ \\
DegInvLex & Left & $y > x_i$ for all $x_i$ \\
DegRevLex & Right & $y > x_i$ for all $x_i$ \\ \hline
\end{tabular}
\end{center}

Noncommutativity also provides the possibility of the new variable $y$
becoming `trapped' in the middle of some monomial forming part of a
polynomial computed during the course of Mora's algorithm.
For example, working with DegRevLex, consider the homogenised polynomial
$h_r(x_1^2 + x_1) = x_1^2 + x_1y$ and the S-polynomial
$$\mathrm{S\mbox{-}pol}(x_1, x_1^2 + x_1y, 1, x_1^2 + x_1y)
= x_1(x_1^2 + x_1y) - (x_1^2 + x_1y)x_1 = x_1^2y - x_1yx_1.$$
Because $y$ appears in the middle
of the monomial $x_1yx_1$, the S-polynomial does not immediately
reduce to zero as it does in the non-homogenised version of the
S-polynomial,
$$\mathrm{S\mbox{-}pol}(x_1, x_1^2 + x_1, 1, x_1^2 + x_1)
= x_1(x_1^2 + x_1) - (x_1^2 + x_1)x_1 = 0.$$
We must therefore make certain that $y$ only
appears on one side of any given monomial by
% This can either be done
% artificially in software -- in which case Lex {\it is}
% an extendible monomial ordering
% (with right homogenisation and $y < x_i$ for all $x_i$) --
introducing the set of polynomials $H = \{h_1, h_2, \hdots, h_n\} =
\{yx_1 - x_1y, \; yx_2 - x_2y, \; \hdots, \;
yx_n - x_ny\}$ into our initial homogenised basis, ensuring that
$y$ commutes with all the other variables in the polynomial
ring. This way, the first S-polynomial will reduce to zero as follows:
$$x_1^2y - x_1yx_1 \rightarrow_{h_1} x_1^2y - x_1^2y = 0.$$
% so that it always gets pushed to one side of a monomial.
Which side $y$ will appear on will be
determined by whether $\LM(yx_i - x_iy) = yx_i$ or
$\LM(yx_i - x_iy) = x_iy$ in our chosen monomial ordering
(pushing $y$ to the right or to the left respectively). 
This side must match the method of homogenisation, which
explains why Lex is not an extendible monomial ordering --- 
for Lex to be extendible, we must
homogenise on the right and have $y < x_i$ for all $x_i$, but then
because $\LM(yx_i - x_iy) = x_iy$ with respect to Lex, the variable 
$y$ will always in practice appear on the left.
% (so that the order on the terms is not preserved).

% The following table summarises which side the homogenising variable
% gets `pushed' to by virtue of reductions using $H$
% for each noncommutative monomial ordering,
% and also explains why Lex is not an extendible monomial ordering
% % if we use $C$
% -- we require the homogenising variable to appear on
% the left, but it always gets pushed to the right because of $C$.
% \begin{center}
% \begin{tabular}{c|c|c|c}
% Monomial Ordering & $\LM(yx_i - x_iy)$ & 
% Pushes to which side? & Side Wanted \\
% \hline
% Lex & $x_iy$ & Left & Right \\
% InvLex & $yx_i$ & Right & Right \\
% RevLex & $x_iy$ & Left & Left \\
% DegLex & $yx_i$ & Left & Left \\
% DegInvLex & $x_iy$ & Left & Left \\
% DegRevLex & $yx_i$ & Right & Right \\ \hline
% \end{tabular}
% \end{center}

\begin{defn} \label{homprocNC}
Let $F = \{f_1, \hdots, f_m\}$ be a non-homogeneous
set of polynomials over the polynomial ring
$R\langle x_1, \hdots, x_n\rangle$.
To compute a Gr\"obner Basis for $F$
% the ideal $J$ generated by $F$
using a program that
only accepts sets of homogeneous polynomials as input,
we use the following procedure (which will only work in conjunction
with an extendible monomial ordering).
\begin{enumerate}[(a)]
\item
Construct a homogeneous set of polynomials
$F' = \{h_{\ell}(f_1), \hdots, h_{\ell}(f_m)\}$ or
$F' = \{h_r(f_1), \hdots, h_r(f_m)\}$ (dependent on the
monomial ordering used). %, where each $f'_i \in F'$
% is a member of the extended polynomial ring $R\langle
% x_1, \hdots, x_n, y\rangle$.
\item
Compute a Gr\"obner Basis $G'$ for the set $F' \cup H$,
where $H = \{yx_1 - x_1y, \: yx_2 - x_2y, \: \hdots,
\: yx_n - x_ny\}$.
\item
Dehomogenise each polynomial $g' \in G'$ to obtain
a Gr\"obner Basis $G$ for $F$, noting that no polynomial originating
from $H$ will appear in $G$ ($d(h_i) = 0$ for all $h_i \in H$).
\end{enumerate}
\end{defn}

\subsection{Selection Strategies} \label{SSNC}
\index{selection strategies}

As in the commutative case, the order in which S-polynomials
are processed during Mora's algorithm has an important
effect on the efficiency of the algorithm. Let us now
generalise the selection strategies defined in
Section \ref{SSC} for use in the noncommutative setting,
basing our decisions on the {\it overlap words} of
S-polynomials.

\begin{defn}
The {\it overlap word} \index{overlap word} of an S-polynomial
$\mathrm{S\mbox{-}pol}(\ell_1, f, \ell_2, g) =
\LC(g)\ell_1fr_1 - \LC(f)\ell_2gr_2$ is the monomial
$\ell_1\LM(f)r_1$ ($ = \ell_2\LM(g)r_2$).
\end{defn}

\begin{defn}
In the noncommutative {\it normal strategy},
\index{normal strategy} we 
\index{strategy!normal} choose an S-polynomial to process if
its overlap word is minimal in the chosen monomial ordering
amongst all such overlap words.
\end{defn}

\begin{defn}
In the noncommutative {\it sugar strategy},
\index{sugar strategy} we 
\index{strategy!sugar} choose an S-polynomial to process
if its sugar (a value associated to the S-polynomial) is
minimal amongst all such values (we use the normal
strategy in the event of a tie).

The sugar of an S-polynomial is computed by
using the following rules on the sugars of
polynomials we encounter during the computation of a
Gr\"obner Basis for the set of polynomials $F = \{f_1, \hdots, f_m\}$.
\begin{enumerate}[(1)]
\item
The sugar $\Sug_{f_i}$ of a polynomial $f_i \in F$ is the total
degree of the polynomial $f_i$ (which is the degree of the term
of maximal degree in $f_i$).
\item
If $p$ is a polynomial and if $t_1$ and $t_2$ are terms, then
$\Sug_{t_1pt_2} = \deg(t_1) + \Sug_p + \deg(t_2)$.
\item
If $p = p_1 + p_2$, then
$\Sug_p = \max(\Sug_{p_1}, \Sug_{p_2})$.
\end{enumerate}
It follows that the sugar of the S-polynomial
$\mathrm{S\mbox{-}pol}(\ell_1, g, \ell_2, h) =
\LC(h)\ell_1gr_1 - \LC(g)\ell_2hr_2$
is given by the formula
$$\Sug_{\mathrm{S\mbox{-}pol}(\ell_1, g, \ell_2, h)} =
\max(\deg(\ell_1) + \Sug_g + \deg(r_1), \;
     \deg(\ell_2) + \Sug_h + \deg(r_2)).$$
\end{defn}

% \subsection{Basis Conversion Algorithms}
%
% Compared to the commutative case, less progress has been made
% in the area of noncommutative basis conversion algorithms.
% We will explore whether the Gr\"obner Walk can be generalised
% to the noncommutative case in Chapter \ref{ChWalk}, but until then
% we note that the FGLM method (which we briefly mentioned in
% Section \ref{BCAC}) has been generalised to the noncommutative
% case \cite{BBM}.

\subsection{Logged Gr\"obner Bases}

\begin{defn}
Let $G = \{g_1, \hdots, g_p\}$ be
a noncommutative Gr\"obner Basis computed from
an initial basis $F = \{f_1, \hdots, f_m\}$. We say that $G$ is a
\index{Gr\"obner basis!logged}
{\it Logged Gr\"obner Basis} if, 
\index{logged Gr\"obner basis} for each $g_i \in G$,
we have an explicit expression of the form
$$g_i = \sum_{\alpha=1}^{\beta} \ell_{\alpha}f_{k_{\alpha}}r_{\alpha},$$
where the $\ell_{\alpha}$ and the $r_{\alpha}$ are terms and
$f_{k_{\alpha}} \in F$ for all $1 \leqslant \alpha \leqslant \beta$.
\end{defn}

\begin{prop}
Let $F = \{f_1, \hdots, f_m\}$ be a finite basis
over a noncommutative polynomial ring.
If we can compute a Gr\"obner Basis for $F$, then
it is always possible to compute a Logged Gr\"obner Basis for $F$.
\end{prop}
\begin{pf}
We refer to the proof of Proposition \ref{LGBC}, substituting
$$\mathrm{S\mbox{-}pol}(\ell_1, f_i, \ell_2, f_j)
  - \sum_{\alpha=1}^{\beta} \ell_{\alpha}g_{k_{\alpha}}r_{\alpha}$$
for $f_{m+1}$ (the $\ell_{\alpha}$ and the $r_{\alpha}$ are terms).
\end{pf}

\section{A Worked Example}
\label{NCEx}

To demonstrate Mora's algorithm in action,
let us now calculate a Gr\"obner Basis for
the ideal $J$ generated by the set of polynomials
$F := \{f_1, f_2, f_3\} = \{x y - z, y z + 2 x + z, y z + x\}$
over the polynomial ring $\mathbb{Q}\langle x, y, z\rangle$.
We shall use the DegLex monomial ordering (with $x > y > z$);
use the normal selection strategy; calculate a Logged Gr\"obner Basis;
and use Buchberger's criteria.

\subsection{Initialisation}

The first part of Mora's algorithm requires
us to find all the overlaps between the lead
monomials of the three polynomials in the initial basis
$G := \{g_1, g_2, g_3\} = \{x y - z, y z + 2 x + z, y z + x\}$.
There are three overlaps in total, summarised by the following table.
\begin{center}
\begin{tabular}{c|ccc}
& Overlap 1 & Overlap 2 & Overlap 3 \\
\hline
Overlap Word & $yz$ & $xyz$ & $xyz$ \\
Polynomial 1 & $yz+2x+z$ & $xy-z$ & $xy-z$ \\
Polynomial 2 & $yz+x$ & $yz+2x+z$ & $yz+x$ \\
$\ell_1$ & $1$ & $1$ & $1$ \\
$r_1$ & $1$ & $z$ & $z$\\
$\ell_2$ & $1$ & $x$ & $x$\\
$r_2$ & $1$ & $1$ & $1$ \\
Degree of Overlap Word & 2 & 3 & 3 \\ \hline
\end{tabular}
\end{center}

Because we are using the normal selection strategy,
it is clear that Overlap 1 will appear in the list $A$ first, but we are free
to choose the order in which the other two overlaps appear (because
their overlap words are identical). To eliminate this choice, we
will use the following {\it tie-breaking strategy} 
\index{tie-breaking strategy} to 
\index{strategy!tie-breaking} order any two S-polynomials
whose overlap words are identical.

\begin{defn}
Let
$s_1 = \mathrm{S\mbox{-}pol}(\ell_1, g_a, \ell_2, g_b)$ and
$s_2 = \mathrm{S\mbox{-}pol}(\ell_3, g_c, \ell_4, g_d)$ be two
S-polynomials with identical overlap words, where
$g_a, g_b, g_c, g_d \in G = \{g_1, \hdots, g_{\alpha}\}$.
Assuming (without loss of generality) that
$a < b$ and $c < d$, the {\it tie-breaking strategy} places
$s_1$ before $s_2$ in $A$ if $a < c$ or if $a = c$ and $b \leqslant d$;
and later in $A$ otherwise.
\end{defn}

Applying the tie-breaking strategy for Overlaps 2 and 3,
it follows that Overlap 2 = $\mathrm{S\mbox{-}pol}(1, g_1, x, g_2)$
will appear in $A$ before Overlap 3 =
$\mathrm{S\mbox{-}pol}(1, g_1, x, g_3)$.

Before we start the main part of the algorithm, 
let us note that for the Logged
Gr\"obner Basis, we begin the algorithm with trivial 
expressions for each of the
three polynomials in the initial basis $G$ in terms of the polynomials of
the input basis $F$:
$g_1 = x y - z = f_1$; $g_2 = y z + 2 x + z = f_2$;
and $g_3 = y z + x = f_3$.
% Later on when we find new polynomials we will express them as more exotic
% combinations of $f_1$, $f_2$ and $f_3$!

\subsection{Calculating and Reducing S-polynomials}

The first S-polynomial to analyse corresponds to Overlap 1 and is the
polynomial $$1(yz+2x+z)1 - 1(yz+x)1 = 2x+z-x = x+z.$$
This polynomial is irreducible with respect to $G$, and so we add
it to $G$ to obtain a new basis
$G = \{x y - z, y z + 2 x + z, y z + x, x + z\} = \{g_1, g_2, g_3, g_4\}$.
Looking for overlaps between the lead monomial of
$x+z$ and the lead monomials
of the four elements of $G$, we see that there 
is one such overlap (with $g_1$)
whose overlap word has degree 2, so this overlap is added to the
beginning of the list $A$ to obtain $A =
\{\mathrm{S\mbox{-}pol}(1, xy-z, 1, x+z), \:
\mathrm{S\mbox{-}pol}(1, xy-z, x, yz+2x+z), \:
\mathrm{S\mbox{-}pol}(1, xy-z, x, yz+x)\}$.
As far as the Logged Gr\"obner Basis goes,
$g_4 = x+z = 1(yz+2x+z)1 - 1(yz+x)1 = f_2 - f_3$.

The next entry in $A$ produces the polynomial
$$1(xy-z)1 - 1(x+z)y = -zy-z.$$
As before, this polynomial is irreducible with respect to $G$,
so we add it to $G$ as the fifth element.
There are also four overlaps
between the lead monomial of $-zy-z$ and the lead monomials of the five
polynomials in $G$:
\begin{center}
\begin{tabular}{c|cccc}
& Overlap 1 & Overlap 2 & Overlap 3 & Overlap 4 \\
\hline
Overlap Word & $zyz$ & $zyz$ & $yzy$ & $yzy$ \\
Polynomial 1 & $yz+2x+z$ & $yz+x$ & $yz+2x+z$ & $yz+x$ \\
Polynomial 2 & $-zy-z$ & $-zy-z$ & $-zy-z$ & $-zy-z$ \\
$\ell_1$ & $z$ & $z$ & $1$ & $1$ \\
$r_1$ & $1$ & $1$ & $y$ & $y$\\
$\ell_2$ & $1$ & $1$ & $y$ & $y$\\
$r_2$ & $z$ & $z$ & $1$ & $1$ \\
Degree of Overlap Word & 3 & 3 & 3 & 3 \\ \hline
\end{tabular}
\end{center}
Inserting these overlaps into the list $A$, we obtain
\begin{eqnarray*}
A = \{&&\mathrm{S\mbox{-}pol}(z, yz+2x+z, 1, -zy-z), \:
        \mathrm{S\mbox{-}pol}(z, yz+x, 1, -zy-z), \\
      &&\mathrm{S\mbox{-}pol}(1, yz+2x+z, y, -zy-z), \:
        \mathrm{S\mbox{-}pol}(1, yz+x, y, -zy-z), \\
      &&\mathrm{S\mbox{-}pol}(1, xy-z, x, yz+2x+z), \:
        \mathrm{S\mbox{-}pol}(1, xy-z, x, yz+x)
        \; \; \; \; \; \; \; \; \; \; \; \}.
\end{eqnarray*}

The logged representation of the fifth basis element again comes
straight from the S-polynomial (as no reduction was performed),
and is as follows: $g_5 = -zy-z = 1(xy-z)1 - 1(x+z)y  =
1(f_1)1 - 1(f_2 - f_3)y = f_1 - f_2y + f_3y$.

The next entry in $A$ yields the polynomial
$$-z(yz+2x+z)1 - 1(-zy-z)z = -2zx - z^2 + z^2 = -2zx.$$
This time, the fourth polynomial in our basis
reduces the S-polynomial in question, giving a
reduction $-2zx \rightarrow_{g_4} 2z^2$. When we add this
polynomial to $G$ and add all five new overlaps to $A$,
we are left with a six element basis
$G = \{x y - z, y z + 2 x + z, y z + x, x + z, - z y - z, 2 z^2\}$
and a list
\begin{eqnarray*}
A = \{&&\mathrm{S\mbox{-}pol}(1, 2z^2, z, 2z^2), \:
        \mathrm{S\mbox{-}pol}(z, 2z^2, 1, 2z^2), \\
      &&\mathrm{S\mbox{-}pol}(z, -zy-z, 1, 2z^2), \:
        \mathrm{S\mbox{-}pol}(z, yz+x, 1, -zy-z), \\
      &&\mathrm{S\mbox{-}pol}(1, yz+2x+z, y, 2z^2), \:
        \mathrm{S\mbox{-}pol}(1, yz+x, y, 2z^2), \\
      &&\mathrm{S\mbox{-}pol}(1, yz+2x+z, y, -zy-z), \:
        \mathrm{S\mbox{-}pol}(1, yz+x, y, -zy-z), \\
      &&\mathrm{S\mbox{-}pol}(1, xy-z, x, yz+2x+z), \:
        \mathrm{S\mbox{-}pol}(1, xy-z, x, yz+x)
      \; \; \; \; \; \; \; \; \; \; \; \}.
\end{eqnarray*}
We obtain the logged version of the sixth basis element 
by working backwards through our calculations:
\begin{eqnarray*}
g_6 & = & 2z^2 \\
& = & -2zx + 2z(x+z) \\
& = & (-z(yz+2x+z)1 - 1(-zy-z)z) + 2z(x+z) \\
& = & (-z(f_2) - (f_1 - f_2y + f_3y)z) + 2z(f_2 - f_3) \\
& = & -f_1z + zf_2 + f_2yz - 2zf_3 - f_3yz.
% & = & -\frac{1}{2}f_1z + \frac{1}{2}zf_2 + \frac{1}{2}f_2yz
% - zf_3 - \frac{1}{2}f_3yz
\end{eqnarray*}

\subsection{Applying Buchberger's Second Criterion}

The next three entries in $A$ all yield S-polynomials that are either zero
or reduce to zero (for example, the first entry corresponds to the polynomial
$2(2z^2)z - 2z(2z^2)1 = 4z^3 - 4z^3 = 0$). The fourth entry in $A$,
$\mathrm{S\mbox{-}pol}(z, yz+x, 1, -zy-z)$, then enables us
(for the first time) to apply Buchberger's second criterion,
allowing us to move on to look at the fifth entry of $A$.
Before we do this however, let us explain why we
can apply Buchberger's second criterion in this particular case.

Recall (from Proposition \ref{B2NC})
that in order to apply Buchberger's second criterion for the S-polynomial
$\mathrm{S\mbox{-}pol}(z, yz+x, 1, -zy-z)$, we need to find a polynomial
$g_i \in G$ such that $\LM(g_i)$ divides 
the overlap word of our S-polynomial,
and any S-polynomials corresponding to overlaps
(as positioned in the overlap word) between $\LM(g_i)$ and
either $\LM(yz+x)$ or $\LM(-zy-z)$ reduce to zero using $G$ 
(which will be the case if those particular S-polynomials have been 
processed earlier in the algorithm).

Consider the polynomial $g_2 = yz+2x+z$. The lead monomial of this polynomial
divides the overlap word $zyz$ of our S-polynomial, which we illustrate as
follows.
$$\xymatrix @R=1pc{
& \ar@{<->}[rr]^{\LM(g_3)} && \\
\ar@{<->}[rr]^{\LM(g_5)} && \\
\ar@{-}[r]^{\displaystyle z} & \ar@{-}[r]^{\displaystyle y}
& \ar@{-}[r]^{\displaystyle z} & \\
& \ar@{<->}[rr]^{\LM(g_2)} &&
}$$
As positioned in the overlap word, we note that $\LM(g_2)$
overlaps with both $\LM(g_3)$ and $\LM(g_5)$, with the
overlaps corresponding to the S-polynomials
$\mathrm{S\mbox{-}pol}(1, g_2, 1, g_3) =
\mathrm{S\mbox{-}pol}(1, yz+2x+z, 1, yz+x)$ and
$\mathrm{S\mbox{-}pol}(z, g_2, 1, g_5) =
\mathrm{S\mbox{-}pol}(z, yz+2x+z, 1, -zy-z)$ respectively.
But these S-polynomials have been processed earlier
in the algorithm (they were the first and third S-polynomials
to be processed); we can therefore apply
Buchberger's second criterion in this instance.

There are now six S-polynomials left in $A$, all of whom
either reduce to zero or are ignored due to Buchberger's second criterion.
Here is a summary of what happens during the remainder of the algorithm.
\begin{center}
\begin{tabular}{l|l}
S-polynomial & Action \\
\hline
$\mathrm{S\mbox{-}pol}(1, yz+2x+z, y, 2z^2)$
& Reduces to zero using the division algorithm\\[1mm]
$\mathrm{S\mbox{-}pol}(1, yz+x, y, 2z^2)$
& Ignored due to Buchberger's second criterion\\
& (using $y z + 2 x + z$)\\[1mm]
$\mathrm{S\mbox{-}pol}(1, yz+2x+z, y, -zy-z)$
& Reduces to zero using the division algorithm\\[1mm]
$\mathrm{S\mbox{-}pol}(1, yz+x, y, -zy-z)$
& Ignored due to Buchberger's second criterion\\
& (using $y z + 2 x + z$)\\[1mm]
$\mathrm{S\mbox{-}pol}(1, xy-z, x, yz+2x+z)$
& Ignored due to Buchberger's second criterion\\
& (using $x + z$)\\[1mm]
$\mathrm{S\mbox{-}pol}(1, xy-z, x, yz+x)$
& Ignored due to Buchberger's second criterion \\
& (using $y z + 2 x + z$)\\ \hline
\end{tabular}
\end{center}

As the list $A$ is now empty, the algorithm terminates
with the following (Logged) Gr\"obner Basis.
\begin{center}
\begin{tabular}{l|l}
Input Basis $F$ & Gr\"obner Basis $G$ \\
\hline
$f_1 = xy-z$ & $g_1 = xy-z = f_1$ \\
$f_2 = yz+2x+z$ & $g_2 = yz+2x+z = f_2$ \\
$f_3 = yz+x$ & $g_3 = yz+x = f_3$ \\
& $g_4 = x+z = f_2 - f_3$ \\
& $g_5 = -zy-z = f_1 - f_2y + f_3y$ \\
& $g_6 = 2z^2 = -f_1z + zf_2 + f_2yz - 2zf_3 - f_3yz$ \\ \hline
\end{tabular}
\end{center}

\subsection{Reduction}

Now that we have constructed a Gr\"obner Basis for our ideal $J$, 
let us go on to find the unique reduced Gr\"obner Basis for $J$
by applying Algorithm \ref{red-noncom} to $G$.

In the first half of the algorithm, we must multiply each polynomial by
the inverse of its lead coefficient and remove from the basis each
polynomial whose lead monomial is a multiple of the lead
monomial of some other polynomial in the basis. For the Gr\"obner
Basis in question, we multiply $g_5$ by $-1$ and
$g_6$ by $\frac{1}{2}$; and we remove $g_1$ and
$g_2$ from the basis (because $\LM(g_1) = \LM(g_4)\times y$
and $\LM(g_2) = \LM(g_3)$). This leaves us with the following
(minimal) Gr\"obner Basis.
\begin{center}
\begin{tabular}{l|l}
Input Basis $F$ & Gr\"obner Basis $G$ \\
\hline
$f_1 = xy-z$ & $g_3 = yz+x = f_3$ \\
$f_2 = yz+2x+z$ & $g_4 = x+z = f_2 - f_3$ \\
$f_3 = yz+x$ & $g_5 = zy+z = -f_1 + f_2y - f_3y$ \\
& $g_6 = z^2 = -\frac{1}{2}f_1z + \frac{1}{2}zf_2 + \frac{1}{2}f_2yz
- zf_3 - \frac{1}{2}f_3yz$ \\ \hline
\end{tabular}
\end{center}

In the second half of the algorithm, we reduce each $g_i \in G$
with respect to $(G\setminus \{g_i\})\cup G'$, placing the remainder
in the (initially empty) set $G'$ and removing $g_i$ from $G$.
For the Gr\"obner Basis in question, we summarise what happens
in the following table, noting that the only
reduction that takes place is the reduction 
$yz+x \rightarrow_{g_4} yz+x - (x+z) = yz-z$.
\begin{center}
\begin{tabular}{l|l|l|l}
$G$ & $G'$ & $g_i$ & $g'_i$ \\ \hline
$\{yz+x, x+z, zy+z, z^2\}$ & $\emptyset$ & $yz+x$ & $yz-z$ \\
$\{x+z, zy+z, z^2\}$ & $\{yz-z\}$ & $x+z$ & $x+z$ \\
$\{zy+z, z^2\}$ & $\{yz-z, x+z\}$ & $zy+z$ & $zy+z$ \\
$\{z^2\}$ & $\{yz-z, x+z, zy+z\}$ & $z^2$ & $z^2$ \\
$\emptyset$ & $\{yz-z, x+z, zy+z, z^2\}$ && \\ \hline
\end{tabular}
\end{center}

We can now give the unique reduced (Logged) Gr\"obner Basis for $J$.
\begin{center}
\begin{tabular}{l|l}
Input Basis $F$ & Unique Reduced Gr\"obner Basis $G'$ \\
\hline
$f_1 = xy-z$    & $yz-z = -f_2+2f_3$ \\
$f_2 = yz+2x+z$ & $x+z = f_2-f_3$ \\
$f_3 = yz+x$    & $zy+z = -f_1 + f_2y - f_3y$ \\
                & $z^2 = -\frac{1}{2}f_1z 
                  + \frac{1}{2}zf_2 + \frac{1}{2}f_2yz
- zf_3 - \frac{1}{2}f_3yz$ \\ \hline
\end{tabular}
\end{center}

%
% Chapter 4
% Author: Gareth Evans
% Last Modified: 2nd February 2006
%

\chapter{Commutative Involutive Bases} \label{ChCIB}

Given a Gr\"obner Basis $G$ for an ideal $J$ over a polynomial
ring $\mathcal{R}$, we know that the remainder of any
polynomial $p \in \mathcal{R}$ with respect to $G$ is
unique. But although this remainder is unique, there may be
many ways of obtaining the remainder, as it is possible that
several polynomials in $G$ divide our polynomial $p$, giving
several {\it reduction paths} for $p$.

\begin{example} \label{ch4ex1}
Consider the DegLex Gr\"obner Basis 
$G := \{g_1, g_2, g_3\} = \{x^2 - 2xy+3, \: 2xy+y^2+5, \:
\frac{5}{4}y^3-\frac{5}{2}x+\frac{37}{4}y\}$ over the polynomial ring
$\mathcal{R} := \mathbb{Q}[x, y]$ from Example \ref{exBuch},
and consider the polynomial $p := x^2y+y^3+8y \in \mathcal{R}$.
The remainder of $p$ with respect to $G$ is 0 (so that $p$ is a member
of the ideal $J$ generated by $G$), but there are two ways of
obtaining this remainder, as shown in the following diagram.
\begin{equation} \label{ch4ex1f}
\xymatrix @R=1.75pc{
& x^2y+y^3+8y \ar[dl]_{g_1} \ar[dr]^{g_2} \\
2xy^2+y^3+5y \ar[d]_{g_2}
&& -\frac{1}{2}xy^2 + y^3 - \frac{5}{2}x+8y \ar[d]^{g_2} \\
0 && \frac{5}{4}y^3-\frac{5}{2}x+\frac{37}{4}y \ar[d]^{g_3} \\
&& 0
}
\end{equation}
\end{example}

An Involutive Basis is a Gr\"obner Basis $G$ for $J$ such that
there is only one possible reduction path for any
polynomial $p \in \mathcal{R}$. In order to find such a basis,
we must restrict which reductions or divisions may take place by
requiring, for each potential reduction of a polynomial
$p$ by a polynomial $g_i \in G$ (so that $\LM(p) = \LM(g_i)\times u$
for some monomial $u$), some extra conditions on the variables
in $u$ to be satisfied, namely that all variables in $u$ have
to be in a set of {\it multiplicative variables} for $g_i$,
a set that is determined by a particular choice of an
{\it involutive division}.

\section{Involutive Divisions} \label{4point1}
\index{involutive divisions} 

In Definition \ref{polydiv}, we saw that a commutative monomial $u_1$ is
divisible by another monomial $u_2$ if there exists a third monomial
$u_3$ such that $u_1 = u_2u_3$; we also introduced the notation
$u_2 \mid u_1$ to denote that $u_2$ is a divisor of $u_1$, a divisor we
shall now refer to as a {\it conventional} \index{conventional divisor}
divisor \index{divisor!conventional}
of $u_1$. For a particular choice of an involutive division $I$,
we say that $u_2$ is an {\it involutive} \index{involutive divisor}
divisor \index{divisor!involutive}
of $u_1$, written $u_2 \mid_I u_1$, if, given a partitioning
(by $I$) of the variables in the polynomial ring into sets of
{\it multiplicative} and {\it nonmultiplicative} variables for $u_2$,
all variables in $u_3$ are in the set of multiplicative variables for $u_2$.

\begin{example}
Let $u_1 := xy^2z^2$; $u'_1 := x^2yz$ and $u_2 := xz$ be three monomials
over the polynomial ring $\mathcal{R} := \mathbb{Q}[x, y, z]$, and let
an involutive division $I$ partition the variables in $\mathcal{R}$ into
the following two sets of variables for the monomial $u_2$:
multiplicative = $\{y, z\}$; nonmultiplicative = $\{x\}$.
It is true that $u_2$ conventionally divides both monomials $u_1$ and $u'_1$,
but $u_2$ only involutively divides monomial $u_1$ as, defining $u_3 := y^2z$
and $u'_3 := xy$ (so that $u_1 = u_2u_3$ and $u'_1 = u_2u'_3$),
we observe that all variables in $u_3$ 
are in the set of multiplicative variables
for $u_2$, but the variables in $u'_3$ 
(in particular the variable $x$) are not
all in the set of multiplicative variables for $u_2$.
\end{example}

More formally, an involutive division $I$ works with
a set of monomials $U$ over a polynomial ring $R[x_1, \hdots, x_n]$
and assigns a set of multiplicative variables
$\mathcal{M}_I(u, U) \subseteq \{x_1, \hdots, x_n\}$
to each element $u \in U$.
It follows that, {\it with respect to $U$},
a monomial $w$ is divisible by a monomial
$u \in U$ if $w = uv$ for some monomial $v$ and all the
variables that appear in $v$ also
appear in the set $\mathcal{M}_I(u, U)$.
% The set of all monomials involutively divisible
% by $u$ forms the {\it involutive cone}
% of $u$, defined as follows.

\begin{defn}
Let $M$ denote the set of all monomials in the
polynomial ring $\mathcal{R} = R[x_1, \hdots, x_n]$,
and let $U \subset M$. The
{\it involutive \index{cone!involutive}
cone} \index{involutive cone}
$\mathcal{C}_I(u, U)$ \index{$C$@$\mathcal{C}_I(u, U)$}
of any monomial $u \in U$ with respect to some
involutive division $I$ is defined as follows.
$$\mathcal{C}_I(u, U) = \{uv \: \mbox{such that} \:
v \in M \: \mbox{and} \: u \mid_I uv\}.$$
\end{defn}

\begin{remark}
We may think of an involutive cone of a particular
monomial $u$ as containing all monomials that are
involutively divisible by $u$.
\end{remark}

Up to now, we have not mentioned any restriction on how
we may assign multiplicative variables to a particular set
of monomials. Let us now specify the rules that ensure
that a particular scheme of assigning multiplicative
variables may be referred to as an involutive division.

\begin{defn} \label{inv-div-defn}
\index{involutive division} \index{division!involutive}
Let $M$ denote the set of all monomials in the
polynomial ring $\mathcal{R} = R[x_1, \hdots, x_n]$.
An {\it involutive division} $I$ on $M$ is
defined if, given any finite set of monomials $U \subset M$,
we can assign a set of {\it multiplicative
variables} \index{multiplicative variables}
\index{$Ma$@$\mathcal{M}_I(u, U)$}
$\mathcal{M}_I(u, U) \subseteq \{x_1, \hdots, x_n\}$
to any monomial $u \in U$ such that
the following two conditions are satisfied.
\begin{enumerate}[(a)]
\item
If there exist two monomials $u_1, u_2 \in U$
such that $\mathcal{C}_I(u_1, U) \cap
\mathcal{C}_I(u_2, U) \neq \emptyset$, \\
then either $\mathcal{C}_I(u_1, U) \subset
\mathcal{C}_I(u_2, U)$ or
$\mathcal{C}_I(u_2, U) \subset
\mathcal{C}_I(u_1, U)$.
\item
If $V \subset U$, then $\mathcal{M}_I(v, U) \subseteq
\mathcal{M}_I(v, V)$ for all $v \in V$.
\end{enumerate}
\end{defn}

\begin{remark}
Informally, condition (a) above ensures that a
monomial can only appear in two involutive cones
$\mathcal{C}_I(u_1, U)$ and $\mathcal{C}_I(u_2, U)$
if $u_1$ is an involutive divisor of $u_2$ or
vice-versa; while condition (b)
ensures that the multiplicative variables of a polynomial
$v \in V \subset U$ with respect to $U$ all appear in the
set of multiplicative variables of $v$ with respect to $V$.
\end{remark}

\begin{defn}
Given an involutive division $I$, the involutive 
\index{span!involutive} span
\index{involutive span} $\mathcal{C}_I(U)$ \index{$C$@$\mathcal{C}_I(U)$}
of a set of monomials $U$ with respect to $I$ is given
by the expression
$$
\mathcal{C}_I(U) = \bigcup_{u \in U} \mathcal{C}_I(u, U).
$$
\end{defn}

\begin{remark}
The (conventional) \index{conventional span}
span \index{span!conventional} of a set of monomials $U$ is
given by the expression
\index{$C$@$\mathcal{C}(U)$}
$$
\mathcal{C}(U) = \bigcup_{u \in U} \mathcal{C}(u, U),
$$
where
\index{$C$@$\mathcal{C}(u, U)$}
$\mathcal{C}(u, U) = \{uv \mid v \; \mbox{is a monomial}\}$ is the
(conventional) \index{conventional cone} cone 
\index{cone!conventional} of a monomial $u \in U$.
\end{remark}

\begin{defn}
If an involutive division $I$ determines the
multiplicative variables for a monomial $u \in U$
independent of the set $U$, then $I$ is a {\it global} division.
\index{involutive division!global} Otherwise, 
\index{global involutive division} $I$
is a {\it local} \index{local involutive division}
division. \index{involutive division!local}
\end{defn}

\begin{remark}
The multiplicative variables for a set of polynomials
$P$ (whose terms are ordered by a monomial ordering $O$)
are determined by the multiplicative variables for the set
of leading monomials $\LM(P)$.
\end{remark}

\subsection{Involutive Reduction}

In Algorithm \ref{com-inv-div}, we specify how to
involutively divide a polynomial $p$ with respect to a set
of polynomials $P$. 

\begin{algorithm}
\setlength{\baselineskip}{3.5ex}
\caption{The Commutative Involutive Division Algorithm}
\label{com-inv-div}
\begin{algorithmic}
\vspace*{2mm}
\REQUIRE{A nonzero polynomial $p$ and a set of nonzero polynomials
         $P = \{p_1, \hdots, p_m\}$ over a 
         polynomial ring $R[x_1, \hdots x_n]$;
         an admissible monomial ordering $\mathrm{O}$; 
         an involutive division $I$.}
\ENSURE{$\Rem_I(p, P) := r$, the involutive remainder of $p$
        with respect to $P$.}
\vspace*{1mm}
\STATE
$r = 0$;
\WHILE{($p \neq 0$)}
\STATE
$u = \LM(p)$; $c = \LC(p)$; $j = 1$; found = false; \\
\WHILE{($j \leqslant m$) \textbf{and} (found == false)}
\IF{($\LM(p_j) \mid_I u$)}
\STATE
found = true; $u' = u/\LM(p_j)$;
$p = p - (c\LC(p_j)^{-1})p_ju'$;
\ELSE
\STATE
$j = j+1$;
\ENDIF
\ENDWHILE
\IF{(found == false)}
\STATE
$r = r+\LT(p)$; $p = p-\LT(p)$;
\ENDIF
\ENDWHILE
\STATE
{\bf return} $r$;
\end{algorithmic}
\vspace*{1mm}
\end{algorithm}

\begin{remark}
The only difference between
Algorithms \ref{com-div} and \ref{com-inv-div} is that the line
``{\bf if} ($\LM(p_j) \mid u$) {\bf then}''
in Algorithm \ref{com-div} has been changed to the line
``{\bf if} ($\LM(p_j) \mid_I u$) {\bf then}''
in Algorithm \ref{com-inv-div}.
\end{remark}

\begin{defn}
If the polynomial $r$ is obtained by involutively dividing
(with respect to some involutive division $I$) the polynomial
$p$ by one of (a) a polynomial $q$; (b) a sequence of polynomials
$q_1, q_2, \hdots, q_{\alpha}$; or (c) a set of polynomials $Q$,
we will use the notation
% \index{$\rightarrow$@\hspace*{-2mm}$\xymatrix{\ar[r]_{I} &}$\hspace*{-1mm}}
\index{$<$@\hspace*{-2mm}$\xymatrix{\ar[r]_{I} &}$\hspace*{-1mm}}
$p \xymatrix{\ar[r]_{I}_(1){q} &}r$;
% \index{$\stackrel{\ast}{\longrightarrow}$@
% \hspace*{-2mm}$\xymatrix{\ar[r]_{I}^{*} &}$\hspace*{-1mm}}
\index{$<$@\hspace*{-2mm}$\xymatrix{\ar[r]_{I}^{*} &}$\hspace*{-1mm}}
$p \xymatrix{\ar[r]_{I}^{*} &}r$ and
$p \xymatrix{\ar[r]_{I}_(1){Q} &}r$ respectively (matching the
notation introduced in Definition \ref{arrows}).
\end{defn}

\subsection{Thomas, Pommaret and Janet divisions}

Let us now consider three different involutive divisions,
all named after their creators in the theory of
Partial Differential Equations (see \cite{Thomas37},
\cite{Pommaret78} and \cite{Janet29}).

\begin{defn}[Thomas]
Let $U = \{u_1 , \hdots, u_m\}$ be a set of monomials
over a polynomial ring $R[x_1, \hdots, x_n]$, where
the monomial $u_j \in U$ (for $1 \leqslant j \leqslant m$)
has corresponding multidegree $(e_j^1, e_j^2, \hdots, e_j^n)$.
The {\it Thomas} \index{involutive division!Thomas}
involutive \index{$T$@$\mathcal{T}$}
division \index{Thomas division}
$\mathcal{T}$ assigns multiplicative
variables to elements
of $U$ as follows: the variable $x_i$ is multiplicative for monomial
$u_j$ (written $x_i \in \mathcal{M}_{\mathcal{T}}(u_j, U)$)
if $e^i_j = \max_k e^i_k$ for all $1 \leqslant k \leqslant m$.
\end{defn}

\begin{defn}[Pommaret]
Let $u$ be a monomial over a polynomial ring \linebreak
$R[x_1, \hdots, x_n]$ with multidegree $(e^1, e^2, \hdots, e^n)$.
The {\it Pommaret} \index{involutive division!Pommaret}
involutive \index{$P$@$\mathcal{P}$}
division \index{Pommaret division}
$\mathcal{P}$ assigns multiplicative variables to $u$
as follows: if $1 \leqslant i \leqslant n$ is the smallest
integer such that $e^i > 0$, then all variables $x_1, x_2,
\hdots, x_i$ are multiplicative for $u$
(we have $x_j \in \mathcal{M}_{\mathcal{P}}(u)$
for all $1 \leqslant j \leqslant i$).
\end{defn}

\begin{defn}[Janet]
Let $U = \{u_1 , \hdots, u_m\}$ be a set of monomials
over a polynomial ring $R[x_1, \hdots, x_n]$, where
the monomial $u_j \in U$ (for $1 \leqslant j \leqslant m$)
has corresponding multidegree $(e_j^1, e_j^2, \hdots, e_j^n)$.
The {\it Janet} \index{involutive division!Janet}
involutive \index{$J$@$\mathcal{J}$} division 
\index{Janet division} $\mathcal{J}$ assigns
multiplicative variables to elements
of $U$ as follows: the variable $x_n$ is multiplicative for monomial
$u_j$ (written $x_n \in \mathcal{M}_{\mathcal{J}}(u_j, U)$)
if $e^n_j = \max_k e^n_k$ for all $1 \leqslant k \leqslant m$; 
the variable $x_i$ (for $1 \leqslant i < n$) is multiplicative for monomial
$u_j$ (written $x_i \in \mathcal{M}_{\mathcal{J}}(u_j, U)$)
if $e^i_j = \max_k e^i_k$ for all monomials $u_k \in U$ such that
$e^{l}_j = e^{l}_k$ for all $i < l \leqslant n$.
\end{defn}

\begin{remark}
Thomas and Janet are local involutive divisions; Pommaret is a
global involutive division.
\end{remark}

\begin{example}
Let $U := \{x^5y^2z, \, x^4yz^2, \, x^2y^2z,
\, xyz^3, \, xz^3, \, y^2z, \, z\}$ be a set of monomials
over the polynomial ring $\mathbb{Q}[x, y, z]$,
with $x > y > z$. Here are the multiplicative
variables for $U$ according to the three involutive
divisions defined above.
\begin{center}
\begin{tabular}{c|c|c|c}
Monomial & Thomas & Pommaret & Janet \\ \hline
$x^5y^2z$ & $\{x, y\}$  & $\{x\}$       & $\{x, y\}$ \\
$x^4yz^2$ & $\emptyset$ & $\{x\}$       & $\{x, y\}$ \\
$x^2y^2z$ & $\{y\}$     & $\{x\}$       & $\{y\}$ \\
$xyz^3$   & $\{z\}$     & $\{x\}$       & $\{x, y, z\}$ \\
$xz^3$    & $\{z\}$     & $\{x\}$       & $\{x, z\}$ \\
$y^2z$    & $\{y\}$     & $\{x, y\}$    & $\{y\}$ \\
$z$       & $\emptyset$ & $\{x, y, z\}$ & $\{x\}$ \\ \hline
\end{tabular}
\end{center}
\end{example}

\begin{prop}
All three involutive divisions defined above satisfy the conditions
of Definition \ref{inv-div-defn}.
\end{prop}
\begin{pf}
Throughout, let $M$ denote the set of all monomials in the
polynomial ring $\mathcal{R} = R[x_1, \hdots, x_n]$;
let $U = \{u_1 , \hdots, u_m\} \subset M$ be a set of monomials with
corresponding multidegrees $(e_k^1, e_k^2, \hdots, e_k^n)$
(where $1 \leqslant k \leqslant m$);
let $u_i, u_j \in U$ (where $1 \leqslant i, j \leqslant m$, $i \neq j$); and
let $m_1, m_2 \in M$ be two monomials with corresponding
multidegrees $(f_1^1, f_1^2, \hdots, f_1^n)$ and
$(f_2^1, f_2^2, \hdots, f_2^n)$.
For condition (a), we need to show that if there exists
a monomial $m \in M$ such that $m_1u_i = m = m_2u_j$ and all variables
in $m_1$ and $m_2$ are multiplicative for $u_i$ and $u_j$ respectively,
then either $u_i$ is an involutive divisor of $u_j$ or vice-versa.
For condition (b), we need to show that all variables that are
multiplicative for $u_i \in U$ are still multiplicative for
$u_i \in V \subseteq U$.

{\bf Thomas}. (a) It is sufficient to prove that
$u_i = u_j$. Assume to the contrary that $u_i \neq u_j$,
so that there is some $1 \leqslant k \leqslant n$ 
such that $e_i^k \neq e_j^k$.
Without loss of generality, assume that $e_i^k < e_j^k$. Because
$e_i^k + f_1^k = e_j^k + f_2^k$, 
it follows that $f_1^k > 0$ so that the variable
$x_k$ must be multiplicative for the monomial $u_i$.
But this contradicts the fact that
$x_k$ cannot be multiplicative for $u_i$ in the Thomas involutive division
because $e_j^k > e_i^k$. We therefore have $u_i = u_j$.

(b) By definition,
if $x_j \in \mathcal{M}_{\mathcal{T}}(u_i, U)$,
then $e_i^j = \max_k e^j_k$ for all
$u_k \in U$. Given a set $V \subseteq U$, it is clear
that $e_i^j = \max_k e^j_k$ for all $u_k \in V$, so that
$x_j \in \mathcal{M}_{\mathcal{T}}(u_i, V)$ as required.

% Comment: u_i can be a suffix of u_j
{\bf Pommaret}. (a) Let $\alpha$ and $\beta$
($1 \leqslant \alpha, \beta \leqslant n$) be the smallest integers
such that $e^{\alpha}_i > 0$ and $e^{\beta}_j > 0$ respectively, and
assume (without loss of generality) that $\alpha \geqslant \beta$.
By definition, we must have $f_1^k = f_2^k = 0$ for all
$\alpha < k \leqslant n$ because the $x_k$
are all nonmultiplicative for $u_i$ and $u_j$. It follows that
$e_i^k = e_j^k$ for all $\alpha < k \leqslant n$. If $\alpha = \beta$,
then it is clear that $u_i$ is an involutive divisor of $u_j$ if
$e^{\alpha}_i < e^{\alpha}_j$, and $u_j$ is an involutive divisor
of $u_i$ if $e^{\alpha}_i > e^{\alpha}_j$.
If $\alpha > \beta$, then $f_2^{\alpha} = 0$ as variable $x_{\alpha}$ is
nonmultiplicative for $u_j$,
% and $e_i^{\alpha} + f_1^{\alpha} = e_j^{\alpha} + f_2^{\alpha}$,
so it follows that $e^{\alpha}_i \leqslant e^{\alpha}_j$ and 
hence $u_i$ is an involutive divisor of $u_j$.

(b) Follows immediately because
Pommaret is a global involutive division.

{\bf Janet}. (a) We prove that $u_i = u_j$.
Assume to the contrary that $u_i \neq u_j$,
so there exists a maximal $1 \leqslant k \leqslant n$ 
such that $e_i^k \neq e_j^k$.
Without loss of generality, assume that $e_i^k < e_j^k$. 
If $k = n$, we get an
immediate contradiction because Janet is equivalent to Thomas for the
final variable. %(so we must have $e_i^n = e_j^n$).
If $k = n-1$, then because
$e_i^{n-1} + f_1^{n-1} = e_j^{n-1} + f_2^{n-1}$,
it follows that $f_1^{n-1} > 0$ so that the variable
$x_{n-1}$ must be multiplicative for the monomial $u_i$.
But this contradicts the
fact that $x_{n-1}$ cannot be multiplicative 
for $u_i$ in the Janet involutive
division because $e_j^{n-1} > e_i^{n-1}$ 
and $e_j^n = e_i^n$. By induction
on $k$, we can show that $e_i^k = e_j^k$ 
for all $1 \leqslant k \leqslant n$,
so that $u_i = u_j$ as required.

(b) By definition,
if $x_j \in \mathcal{M}_{\mathcal{J}}(u_i, U)$, then
$e_i^j = \max_k e^j_k$ for all monomials $u_k \in U$ such that
$e^{l}_i = e^{l}_k$ for all $i < l \leqslant n$.
Given a set $V \subseteq U$, it is clear
that $e_i^j = \max_k e^j_k$ for all $u_k \in V$ such that
$e^{l}_i = e^{l}_k$ for all $i < l \leqslant n$, so that
$x_j \in \mathcal{M}_{\mathcal{J}}(u_i, V)$ as required.
\end{pf}

The conditions of Definition \ref{inv-div-defn} ensure that
any polynomial is involutively divisible by at most
one polynomial in any Involutive Basis.
One advantage of this important combinatorial property is that
the Hilbert function of an ideal $J$ is easily computable
with respect to an Involutive Basis (see \cite{Apel98b}).

\begin{example}
Returning to Example \ref{ch4ex1}, consider again the
DegLex Gr\"obner Basis $G := \{x^2-2xy+3,
2xy+y^2+5, \frac{5}{4}y^3-\frac{5}{2}x+\frac{37}{4}y\}$ over
the polynomial ring $\mathbb{Q}[x, y]$. A Pommaret Involutive
Basis for $G$ is the set $P :=
G \cup \{g_4 := -5xy^2-5x+6y\}$, with the
variable $x$ being multiplicative for all
polynomials in $P$, and the variable $y$
being multiplicative for just $g_3$. We can
illustrate the difference between the overlapping cones of 
$G$ and the non-overlapping involutive cones of $P$ by the
following diagram.
\begin{center}
\begin{picture}(0,0)%
\includegraphics{ch4d1.pstex}%
\end{picture}%
\setlength{\unitlength}{2633sp}%
\begingroup\makeatletter\ifx\SetFigFont\undefined%
\gdef\SetFigFont#1#2#3#4#5{%
  \reset@font\fontsize{#1}{#2pt}%
  \fontfamily{#3}\fontseries{#4}\fontshape{#5}%
  \selectfont}%
\fi\endgroup%
\begin{picture}(8150,4122)(-1537,-5602)
\put(-554,-1846){\makebox(0,0)[lb]{\smash{\SetFigFont{12}{14.4}{\familydefault}{\mddefault}{\updefault}{\color[rgb]{0,0,0}Gr\"obner Basis $G$}%
}}}
\put(-1304,-1696){\makebox(0,0)[lb]{\smash{\SetFigFont{12}{14.4}{\familydefault}{\mddefault}{\updefault}{\color[rgb]{0,0,0}$y$}%
}}}
\put(1996,-4996){\makebox(0,0)[lb]{\smash{\SetFigFont{12}{14.4}{\familydefault}{\mddefault}{\updefault}{\color[rgb]{0,0,0}$x$}%
}}}
\put(3826,-1861){\makebox(0,0)[lb]{\smash{\SetFigFont{12}{14.4}{\familydefault}{\mddefault}{\updefault}{\color[rgb]{0,0,0}Pommaret Basis $P$}%
}}}
\put(3151,-1711){\makebox(0,0)[lb]{\smash{\SetFigFont{12}{14.4}{\familydefault}{\mddefault}{\updefault}{\color[rgb]{0,0,0}$y$}%
}}}
\put(6451,-5011){\makebox(0,0)[lb]{\smash{\SetFigFont{12}{14.4}{\familydefault}{\mddefault}{\updefault}{\color[rgb]{0,0,0}$x$}%
}}}
\put(2626,-3616){\makebox(0,0)[lb]{\smash{\SetFigFont{12}{14.4}{\rmdefault}{\mddefault}{\updefault}{\color[rgb]{0,0,0}$g_3$}%
}}}
\put(3481,-4306){\makebox(0,0)[lb]{\smash{\SetFigFont{12}{14.4}{\rmdefault}{\mddefault}{\updefault}{\color[rgb]{0,0,0}$g_4$}%
}}}
\put(3481,-4906){\makebox(0,0)[lb]{\smash{\SetFigFont{12}{14.4}{\rmdefault}{\mddefault}{\updefault}{\color[rgb]{0,0,0}$g_2$}%
}}}
\put(4096,-5521){\makebox(0,0)[lb]{\smash{\SetFigFont{12}{14.4}{\rmdefault}{\mddefault}{\updefault}{\color[rgb]{0,0,0}$g_1$}%
}}}
\end{picture}

% Remark: I used 18pt fonts in the diagram, so in order
% to use 12pt fonts in the final article I used an
% export magnification value of 66.7% in XFig
% (exporting as "Combined PS/LaTeX (both parts)")
\end{center}
The diagram also demonstrates that the polynomial
$p := x^2y+y^3+8y$ is initially conventionally
divisible by two members of
the Gr\"obner Basis $G$ (as seen in Equation (\ref{ch4ex1f})),
but is only involutively divisible by one
member of the Involutive Basis $P$,
starting the following unique involutive reduction path for $p$.
$$\xymatrix @R=1.75pc{
x^2y+y^3+8y \ar[d]^{g_2} \\
-\frac{1}{2}xy^2+y^3-\frac{5}{2}x+8y \ar[d]^{g_4} \\
y^3-2x+\frac{37}{5}y \ar[d]^{g_3} \\
0
}$$
\end{example}

\section{Prolongations and Autoreduction}  \label{4point2}

Whereas Buchberger's algorithm constructs a Gr\"obner
Basis by using S-polynomials, the involutive algorithm
will construct an Involutive Basis by using
{\it prolongations} and {\it autoreduction}.

\begin{defn}
\index{prolongation}
Given a set of polynomials $P$, a {\it prolongation} of
a polynomial $p \in P$ is a product $px_i$, where 
$x_i \notin \mathcal{M}_I(\LM(p), \LM(P))$
with respect to some involutive division $I$.
\end{defn}

\begin{defn} \label{com-inv-defn}
\index{autoreduction}
A set of polynomials $P$ is said to be {\it autoreduced} if
no polynomial $p \in P$ exists such that $p$
contains a term which is involutively divisible (with respect
to $P$) by some
polynomial $p' \in P\setminus\{p\}$.
Algorithm \ref{com-auto} provides a way of performing
autoreduction, and introduces the following notation:
Let $\Rem_I(A, B, C)$ denote the involutive
remainder of the polynomial $A$ with respect to the set
of polynomials $B$, where reductions are only to be
performed by elements of the set $C \subseteq B$.
% autoreduction is the process
% of involutively reducing each member of
% a set of polynomials by the rest of the set
% until all members become involutively irreducible
% (an algorithm for performing autoreduction
% is given in Algorithm \ref{com-auto}).
\end{defn}

\begin{remark}
The involutive cones associated to an 
autoreduced set of polynomials are always disjoint,
meaning that a given monomial can only appear in at
most one of the involutive cones.
\end{remark}

\begin{algorithm}
\setlength{\baselineskip}{3.5ex}
\caption{The Commutative Autoreduction Algorithm}
\label{com-auto}
\begin{algorithmic}
\vspace*{2mm}
\REQUIRE{A set of polynomials $P = \{p_1, p_2, \hdots, p_{\alpha}\}$;
         an involutive division $I$.}
\ENSURE{An autoreduced set of polynomials 
        $Q = \{q_1, q_2, \hdots, q_{\beta}\}$.}
\vspace*{1mm}
\WHILE{($\exists \; p_i \in P$ such that $\Rem_I(p_i,
P, P\setminus\{p_i\}) \neq p_i$)}
\STATE
$p'_i = \Rem_I(p_i, P, P\setminus\{p_i\})$; \\
$P = P \setminus\{p_i\}$;
\IF{($p'_i \neq 0$)}
\STATE
$P = P \cup \{p'_i\}$;
\ENDIF
\ENDWHILE
\STATE
$Q = P$; \\
{\bf return} $Q$;
\end{algorithmic}
\vspace*{1mm}
\end{algorithm}

\begin{prop} \label{plusC}
Let $P$ be a set of polynomials over a polynomial ring
$\mathcal{R} = R[x_1, \hdots, x_n]$, and let
$f$ and $g$ be two polynomials also in $\mathcal{R}$.
If $P$ is autoreduced with respect to an
involutive division $I$, then
$\Rem_I(f, P) + \Rem_I(g, P) = \Rem_I(f+g, P)$.
\end{prop}
\begin{pf}
Let $f' := \Rem_I(f, P)$; $g' := \Rem_I(g, P)$ and
$h' := \Rem_I(h, P)$, where $h := f+g$. Then, 
by the respective involutive reductions,
we have expressions
$$
f' = f - \sum_{a=1}^A p_{{\alpha}_a}t_a;
$$
$$
g' = g - \sum_{b=1}^B p_{{\beta}_b}t_b
$$
and
$$
h' = h - \sum_{c=1}^C p_{{\gamma}_c}t_c,
$$
where $p_{{\alpha}_{a}}, \, p_{{\beta}_{b}}, \,
p_{{\gamma}_{c}} \in P$
and $t_{a}, t_{b}, t_{c}$ are terms which
are multiplicative (over $P$) for each $p_{{\alpha}_{a}}$,
$p_{{\beta}_{b}}$ and $p_{{\gamma}_{c}}$ respectively.

Consider the polynomial $h'-f'-g'$. By the above
expressions, we can deduce\footnote{For
$1 \leqslant d \leqslant A$, $p_{{\delta}_d}t_d = p_{{\alpha}_a}t_a$
($1 \leqslant a \leqslant A$); for $A+1 \leqslant d
\leqslant A+B$, $p_{{\delta}_d}t_d = p_{{\beta}_b}t_b$
($1 \leqslant b \leqslant B$); and for $A+B+1 \leqslant d
\leqslant A+B+C =: D$, $p_{{\delta}_d}t_d = p_{{\gamma}_c}t_c$
($1 \leqslant c \leqslant C$).}
that
$$h'-f'-g' = \sum_{a=1}^A p_{{\alpha}_a}t_a
+ \sum_{b=1}^B p_{{\beta}_b}t_b
- \sum_{c=1}^C p_{{\gamma}_c}t_c
=: \sum_{d=1}^D p_{{\delta}_d}t_d.$$
{\bf Claim:} $\Rem_I(h'-f'-g', P) = 0$.

{\bf Proof of Claim:}
Let $t$ denote the leading term of the
polynomial $\sum_{d=1}^{D} p_{{\delta}_d}t_d$.
Then $\LM(t) = \LM(p_{{\delta}_{k}}t_{k})$ for some 
$1 \leqslant k \leqslant D$ since, if not,
there exists a monomial $\LM(p_{{\delta}_{k'}}t_{k'}) =
\LM(p_{{\delta}_{k''}}t_{k''}) =: u$ for some 
$1 \leqslant k', k'' \leqslant D$
(with $p_{{\delta}_{k'}} \neq p_{{\delta}_{k''}}$) 
such that $u$ is involutively
divisible by the two polynomials $p_{{\delta}_{k'}}$ 
and $p_{{\delta}_{k''}}$,
contradicting Definition \ref{inv-div-defn} (recall that our
set $P$ is autoreduced, so that the involutive cones of 
$P$ are disjoint). It follows that we can use
$p_{\delta_k}$ to eliminate $t$ by involutively
reducing $h'-f'-g'$ as shown below.
\begin{equation} \label{GIE1}
\sum_{d=1}^{D} p_{{\delta}_d}t_d
\xymatrix{\ar[r]_I_(1){p_{{\delta}_k}} &}
\sum_{d=1}^{k-1} p_{{\delta}_d}t_d + \sum_{d=k+1}^{D} p_{{\delta}_d}t_d.
\end{equation}
By induction, % and by the admissibility of the chosen monomial ordering,
we can apply a chain of involutive reductions to the
right hand side of Equation (\ref{GIE1})
to obtain a zero remainder, so that
$\Rem_I(h'-f'-g', P) = 0$.
\hfill ${}_{\Box}$

To complete the proof, we note that since 
$f'$, $g'$ and $h'$ are all involutively irreducible, we
must have $\Rem_I(h'-f'-g', P) = h'-f'-g'$. It therefore
follows that $h'-f'-g' = 0$, or $h' = f'+g'$ as required.
\end{pf}

\begin{remark}
The above proof is based on the proofs of Theorem 5.4
and Corollary 5.5 in \cite{Gerdt98a}.
\end{remark}

Let us now give a definition of a Locally Involutive Basis
in terms of prolongations. Later on in this chapter, we will
discover that the Involutive Basis algorithm only constructs
Locally Involutive Bases, and it is the
extra properties of each involutive division used with the
algorithm that ensures that any computed
Locally Involutive Basis is an Involutive Basis.

\begin{defn} \label{LIBC}
\index{involutive basis!locally}
Given \index{locally involutive basis}
an involutive division $I$ and an admissible monomial
ordering $O$, an autoreduced set of polynomials $P$
is a {\it Locally Involutive Basis} with
respect to $I$ and $O$ if any prolongation
of any polynomial $p_i \in P$ involutively reduces to
zero using $P$.
% is involutively
% reducible by some $p_j \in P$. In other words, for all
% $p_i \in P$ and for any variable
% $x_k \notin \mathcal{M}_I(\LM(p_i), \LM(P))$,
% $$\exists \: p_j \in P \mathrm{\; such \; that \;}
% p_j \mid_I p_ix_k.$$
\end{defn}

\begin{defn} \label{IBC}
\index{involutive basis}
Given \index{basis!involutive}
an involutive division $I$ and an admissible monomial
ordering $O$, an autoreduced set of polynomials $P$ is an
{\it Involutive Basis} with respect to $I$ and $O$
if any multiple $p_it$ of any polynomial $p_i \in P$
by any term $t$ involutively reduces to zero using $P$.
% is involutively reducible by some $p_j \in P$. In
% other words, for all $p_i \in P$ and for any monomial $u$,
% $$\exists \: p_j \in P \mathrm{\; such \; that \;}
% p_j \mid_I p_iu.$$
\end{defn}

% \begin{remark}
% By the admissibility of the chosen monomial ordering,
% any prolongation of
% a member of a Locally Involutive Basis involutively
% reduces to zero, and any multiple of a member of an
% Involutive Basis involutively reduces to zero.
% \end{remark}

% \begin{remark}
%  $P$ is an Involutive Basis with respect to an
% involutive division $I$ and an admissible monomial
% ordering $O$, we will often refer to $P$ as being
% ``an $I$-Basis''. For example, a Thomas Basis is
% an Involutive Basis with respect to the
% Thomas involutive division and an admissible monomial
% ordering $O$.
% \end{remark}

\section{Continuity and Constructivity} \label{4point3}

In the theory of commutative Gr\"obner Bases, 
Buchberger's algorithm returns
a Gr\"obner Basis as long as an admissible monomial ordering is
used. In the theory of commutative Involutive Bases however, 
not only must an admissible
monomial ordering be used, but the involutive division chosen must
be {\it continuous} and {\it constructive}.

\begin{defn}[Continuity]
\index{involutive division!continuous}
Let \index{continuous involutive division}
$I$ be an involutive division, and let $U$ be an arbitrary set
of monomials over a polynomial ring $R[x_1, \hdots, x_n]$.
We say that $I$ is {\it continuous} if, given any
sequence of monomials $\{u_1, u_2, \hdots, u_m\}$ from $U$
such that for all $i < m$, we have $u_{i+1} \mid_I u_ix_{j_i}$
for some variable $x_{j_i}$ that is nonmultiplicative for monomial
$m_i$ (or $x_{j_i} \notin \mathcal{M}_I(u_i, U)$), no
two monomials in the sequence are the same
($u_r \neq u_s$ for all $r \neq s$,
where $1 \leqslant r, s \leqslant m$).
\end{defn}

\begin{prop}
The Thomas, Pommaret and Janet involutive divisions are all
continuous.
\end{prop}
\begin{pf}
Throughout, let the sequence of monomials $\{u_1, \hdots, u_i,
\hdots, u_m\}$ have corresponding multidegrees
$(e_i^1, e_i^2, \hdots, e_i^n)$
(where $1 \leqslant i \leqslant m$).

{\bf Thomas}. If the variable $x_{j_i}$ is nonmultiplicative
for monomial $u_i$, % ($1 \leqslant i < m$)
then, by definition,
$e_i^{j_i} \neq \max_{t} e_{t}^{j_i}$
for all $u_{t} \in U$.
Variable $x_{j_i}$ cannot therefore be multiplicative for monomial
$u_{i+1}$ if $e_{i+1}^{j_{i}} \leqslant e_i^{j_i}$, so we must
have $e_{i+1}^{j_i} = e_i^{j_i}+1$ in order to have
$u_{i+1} \mid_{\mathcal{T}} u_{i}x_{j_i}$. Further, for all
$1 \leqslant k \leqslant n$ such that
$k \neq j_{i}$, we must have $e_{i+1}^k = e_{i}^k$ as,
if $e_{i+1}^k < e_i^k$, then $x_k$ cannot be multiplicative
for monomial $u_{i+1}$ (which contradicts $u_{i+1}
\mid_{\mathcal{T}} u_ix_{j_i}$).
Thus $u_{i+1} = u_ix_{j_i}$, and so it is clear that the
monomials in the sequence $\{u_1, u_2, \hdots, u_m\}$ are
all different.

{\bf Pommaret}. Let $\alpha_i$ ($1 \leqslant \alpha_i \leqslant n$)
be the smallest integer such that $e^{\alpha_i}_i > 0$
(where $1 \leqslant i \leqslant m$), so that $e_i^k = 0$ for
all $k < \alpha_i$. Because $u_{i+1} \mid_{\mathcal{P}} u_ix_{j_i}$ for
all $1 \leqslant i < m$, and because (by definition)
% the variable $x_{j_i}$ that is
% nonmultiplicative for monomial $u_i$ is such that
$j_i > \alpha_i$, it follows that we must have $e_{i+1}^k = 0$ for all 
$k < \alpha_i$. Therefore $\alpha_{i+1} \geqslant \alpha_i$ for all
$1 \leqslant i < n$. If $\alpha_{i+1} = \alpha_i$, we note that 
$e_{i+1}^{\alpha_i} \leqslant e_i^{\alpha_i}$ because variable 
$x_{\alpha_i}$ is multiplicative for monomial $u_{i+1}$.
If then we have $e_{i+1}^{\alpha_i} = e_i^{\alpha_i}$, then because the
variable $x_{j_i}$ is also nonmultiplicative for monomial $u_{i+1}$,
we must have $e_{i+1}^{j_i} = e_i^{j_i}+1$.

It is now clear that the monomials in the sequence 
$\{u_1, u_2, \hdots, u_m\}$
are all different because (a) the values in the sequence
$\alpha = \{\alpha_1, \alpha_2, \hdots, \alpha_m\}$ 
monotonically increase; (b)
for consecutive values $\alpha_s, \alpha_{s+1}, \hdots, \alpha_{s+\sigma}$
in $\alpha$ that are identical ($1 \leqslant s < m$, $s+\sigma \leqslant m$),
the values in the corresponding sequence
$E = \{e_s^{\alpha_s}, e_{s+1}^{\alpha_s}, \hdots, e_{s+\sigma}^{\alpha_s}\}$
monotonically decrease; (c) for consecutive values $e_t^{\alpha_s},
e_{t+1}^{\alpha_s}, \hdots, e_{t+\tau}^{\alpha_s}$ in $E$ that are identical
($s \leqslant t < s+\sigma$, $t+\tau \leqslant s+\sigma$),
the degrees of the monomials $u_{t}, u_{t+1}, \hdots, u_{t+\tau}$
strictly increase.

{\bf Janet}. Consider the monomials $u_1$, $u_2$ and the variable
$x_{j_1}$ that is nonmultiplicative for $u_1$. We will first prove
(by induction) that $e_2^i = e_1^i$ for all $j_1 < i \leqslant n$.
For the case $i = n$, we must have $e_2^n = e_1^n$ otherwise (by definition)
variable $x_n$ is nonmultiplicative for monomial $u_2$ 
(we have $e_2^n < e_1^n$),
contradicting that fact that $u_2 \mid_{\mathcal{J}} u_1x_{j_1}$. 
For the inductive
step, assume that $e_2^i = e_1^i$ for all $k \leqslant i \leqslant n$,
and let us look at the case $i = k-1$. If
$e_2^{k-1} < e_1^{k-1}$, then (by definition) variable
$x_{k-1}$ is nonmultiplicative for monomial $u_2$,
again contradicting the fact that $u_2 \mid_{\mathcal{J}} u_1x_{j_1}$. 
It follows that we must have $e_2^{k-1} = e_1^{k-1}$.

Let us now prove that $e_2^{j_1} = e_1^{j_1}+1$. We can rule out
the case $e_2^{j_1} < e_1^{j_1}$ immediately because this implies
that the variable $x_{j_1}$ is nonmultiplicative for monomial $u_2$
(by definition), contradicting the fact that 
$u_2 \mid_{\mathcal{J}} u_1x_{j_1}$.
The case $e_2^{j_1} = e_1^{j_1}$ can also be ruled out
because we cannot have $e_2^i = e_1^i$ for all $j_1 \leqslant i \leqslant n$
and variable $x_{j_1}$ being simultaneously nonmultiplicative
for monomial $u_1$ and multiplicative for monomial $u_2$.
Thus $e_2^{j_1} = e_1^{j_1}+1$. It follows
that $u_1 < u_2$ in the InvLex monomial ordering (see Section
\ref{CMO}) and so, by induction, $u_1 < u_2 < \cdots < u_m$ in the InvLex
monomial ordering. The monomials in the sequence
$\{u_1, u_2, \hdots, u_m\}$ are therefore all different.
\end{pf}

\begin{prop} \label{LocalToGlobal}
If an involutive division $I$ is continuous,
and a given set of polynomials $P$ is a Locally Involutive
Basis with respect to $I$ and some admissible monomial
ordering $O$, then $P$ is an Involutive Basis with respect
to $I$ and $O$.
\end{prop}
\begin{pf}
Let $I$ be a continuous involutive
division; let $O$ be an admissible monomial ordering; and let
$P$ be a Locally Involutive Basis with respect to $I$ and $O$.
Given any polynomial $p \in P$ and any term $t$,
in order to show that $P$ is an Involutive Basis with
respect to $I$ and $O$, we
must show that
$pt \xymatrix{\ar[r]_I_(1){P} &} 0$.
% $\exists \, p' \in P$ such that
% $p' \mid_I pu$.

If $p \mid_I pt$ we are done, as we can use $p$ to
involutively reduce $pt$ to obtain a zero remainder.
Otherwise, $\exists \, y_1
\notin \mathcal{M}_I(\LM(p), \LM(P))$ such that $t$ contains
$y_1$. By Local Involutivity, the prolongation $py_1$
involutively reduces to zero using $P$. Assuming that
the first step of this involutive reduction
involves the polynomial $p_1 \in P$,
we can write
\begin{equation} \label{ctsE1C}
py_1 = p_1t_1 + \sum_{a=1}^A p_{\alpha_a}t_{\alpha_a},
\end{equation}
where $p_{\alpha_a} \in P$ and $t_1, t_{\alpha_a}$ are
terms which are multiplicative (over $P$) for $p_1$
and each $p_{\alpha_a}$ respectively. Multiplying
both sides of Equation (\ref{ctsE1C}) by $\frac{t}{y_1}$,
we obtain the equation
\begin{equation} \label{ctsE2C}
pt = p_1t_1\frac{t}{y_1} +
\sum_{a=1}^A p_{\alpha_a}t_{\alpha_a}\frac{t}{y_1}.
\end{equation}
If $p_1 \mid_I pt$, it is clear that
we can use $p_1$ to involutively
reduce the polynomial $pt$ to obtain the polynomial
$\sum_{a=1}^A p_{\alpha_a}t_{\alpha_a}\frac{t}{y_1}$.
By Proposition \ref{plusC}, we can then continue
to involutively reduce $pt$ by repeating this proof
on each polynomial $p_{\alpha_a}t_{\alpha_a}\frac{t}{y_1}$
individually (where $1 \leqslant a \leqslant A$),
noting that this process will terminate because of the
admissibility of $O$ (we have
$\LM(p_{\alpha_a}t_{\alpha_a}\frac{t}{y_1}) < \LM(pt)$ for
all $1 \leqslant a \leqslant A$).

Otherwise, if $p_1$ does not involutively divide $pt$,
there exists a variable
$y_2 \in \frac{t}{y_1}$ such that
$y_2 \notin \mathcal{M}_I(\LM(p_1), \LM(P))$.
By Local Involutivity, the prolongation $p_1y_2$
involutively reduces to zero using $P$. Assuming that
the first step of this involutive reduction
involves the polynomial $p_2 \in P$,
we can write
\begin{equation} \label{ctsE3C}
p_1y_2 = p_2t_2 + \sum_{b=1}^B p_{\beta_b}t_{\beta_b},
\end{equation}
where $p_{\beta_b} \in P$ and $t_2, t_{\beta_b}$ are
terms which are multiplicative (over $P$) for $p_2$
and each $p_{\beta_b}$ respectively. Multiplying
both sides of Equation (\ref{ctsE3C}) by $\frac{t_1t}{y_1y_2}$,
we obtain the equation
\begin{equation} \label{ctsE4C}
p_1t_1\frac{t}{y_1} = p_2t_2\frac{t_1t}{y_1y_2} +
\sum_{b=1}^B p_{\beta_b}t_{\beta_b}\frac{t_1t}{y_1y_2}.
\end{equation}
Substituting for $p_1t_1\frac{t}{y_1}$ from Equation
(\ref{ctsE4C}) into Equation (\ref{ctsE2C}), we obtain
the equation
\begin{equation} \label{ctsE5C}
pt = p_2t_2\frac{t_1t}{y_1y_2} +
\sum_{a=1}^A p_{\alpha_a}t_{\alpha_a}\frac{t}{y_1} +
\sum_{b=1}^B p_{\beta_b}t_{\beta_b}\frac{t_1t}{y_1y_2}.
\end{equation}
If $p_2 \mid_I pt$, it is clear that
we can use $p_2$ to involutively
reduce the polynomial $pt$ to obtain the polynomial
$\sum_{a=1}^A p_{\alpha_a}t_{\alpha_a}\frac{t}{y_1} +
\sum_{b=1}^B p_{\beta_b}t_{\beta_b}\frac{t_1t}{y_1y_2}$.
As before, we can then use Proposition \ref{plusC}
to continue the involutive reduction of $pt$ by
repeating this proof on each summand individually.

Otherwise, if $p_2$ does not involutively divide $pt$,
we continue by induction, obtaining a sequence
$p, p_1, p_2, p_3, \hdots$
of elements in $P$. By construction, each element in the sequence
divides $pt$. By continuity, each element in the
sequence is different. Because $P$ is finite and
because $pt$ has a finite number of distinct
divisors, the sequence must be finite,
terminating with an involutive divisor $p' \in P$ of $pt$,
which then allows us to finish the proof through use of
Proposition \ref{plusC} and the admissibility of $O$.
\end{pf}

\begin{remark}
The above proof is a slightly clarified version of
the proof of Theorem 6.5 in \cite{Gerdt98a}.
\end{remark}

\begin{defn}[Constructivity] \label{ConstC}
\index{involutive division!constructive}
Let \index{constructive involutive division}
$I$ be an involutive division, and let $U$ be an arbitrary set
of monomials over a polynomial ring $R[x_1, \hdots, x_n]$.
We say that $I$ is {\it constructive} if, given any
monomial $u \in U$ and any nonmultiplicative variable
$x_i \notin \mathcal{M}_I(u, U)$ satisfying the following
two conditions, no monomial
$w \in \mathcal{C}_I(U)$ exists such that $ux_{i} \in
\mathcal{C}_I(w, U\cup\{w\})$.
\begin{enumerate}[(a)]
\item
$ux_{i} \notin \mathcal{C}_I(U)$.
\item
If there exists a monomial $v \in U$ and
a nonmultiplicative variable $x_{j} \notin
\mathcal{M}_I(v, U)$ such that $vx_{j} \mid ux_{i}$
but $vx_{j} \neq ux_{i}$, then $vx_{j} \in \mathcal{C}_I(U)$.
\end{enumerate}
\end{defn}

\begin{remark}
Constructivity allows us to consider only polynomials
whose lead monomials lie {\it outside} the current involutive
span as potential new Involutive Basis elements.
\end{remark}

\begin{prop}
The Thomas, Pommaret and Janet involutive divisions are all
constructive.
\end{prop}
\begin{pf}
Throughout, let the monomials $u$, $v$ and $w$ that appear
in Definition \ref{ConstC} have
corresponding multidegrees $(e_u^1, e_u^2, \hdots, e_u^n)$,
$(e_v^1, e_v^2, \hdots, e_v^n)$ and $(e_w^1, e_w^2, \hdots, e_w^n)$;
and let the monomials $w_1$, $w_2$, $w_3$ and $\mu$ that
appear in this proof have corresponding multidegrees
$(e_{w_1}^1, e_{w_1}^2, \hdots, e_{w_1}^n)$,
$(e_{w_2}^1, e_{w_2}^2, \hdots, e_{w_2}^n)$,
$(e_{w_3}^1, e_{w_3}^2, \hdots, e_{w_3}^n)$ and
$(e_{\mu}^1, e_{\mu}^2, \hdots, e_{\mu}^n)$.

To prove that a particular involutive division $I$
is constructive, we will
assume that a monomial $w \in \mathcal{C}_I(U)$
exists such that $ux_i \in \mathcal{C}_I(w, U \cup \{w\})$.
Then $w = \mu w_1$ for some monomial $\mu \in U$ and some monomial
$w_1$ that is multiplicative for $\mu$ over the set $U$
($e_{w_1}^k > 0 \Rightarrow x_k \in \mathcal{M}_I(\mu, U)$
for all $1 \leqslant k \leqslant n$); and
$ux_i = ww_2$ for some monomial $w_2$ that is multiplicative
for $w$ over the set $U \cup \{w\}$
($e_{w_2}^k > 0 \Rightarrow x_k \in \mathcal{M}_I(w, U \cup \{w\})$
for all $1 \leqslant k \leqslant n$). It follows that
$ux_i = \mu w_1w_2$. If we can show that all variables
appearing in $w_2$ are multiplicative for $\mu$ over the set $U$
($e_{w_2}^k > 0 \Rightarrow x_k \in \mathcal{M}_I(\mu, U)$
for all $1 \leqslant k \leqslant n$),
then $\mu$ is an involutive divisor of $ux_i$,
contradicting the assumption $ux_i \notin \mathcal{C}_I(U)$.

{\bf Thomas}.
Let $x_k$ be an arbitrary variable 
($1 \leqslant k \leqslant n$) such that $e_{w_2}^k > 0$.
If $e_{w_1}^k > 0$, then it is clear that $x_k$
is multiplicative for $\mu$. Otherwise $e_{w_1}^k = 0$
so that $e_w^k = e_{\mu}^k$. By definition, this implies
that $x_k \in \mathcal{M}_{\mathcal{T}}(\mu, U)$ as
$x_k \in \mathcal{M}_{\mathcal{T}}(w, U\cup\{w\})$.
Thus $x_k \in \mathcal{M}_{\mathcal{T}}(\mu, U)$.

{\bf Pommaret}.
Let $\alpha$ and $\beta$ ($1 \leqslant \alpha, \beta \leqslant n$)
be the smallest integers such that $e^{\alpha}_{\mu} > 0$
and $e^{\beta}_{w} > 0$ respectively. By definition,
$\beta \leqslant \alpha$ (because $w = \mu w_1$), so
for an arbitrary $1 \leqslant k \leqslant n$, it follows that
$e_{w_2}^k > 0 \Rightarrow k \leqslant \beta \leqslant \alpha
\Rightarrow x_k \in \mathcal{M}_{\mathcal{P}}(\mu, U)$ as required.

{\bf Janet}.
Here we proceed by searching for a 
monomial $\nu \in U$ such that
$ux_i \in \mathcal{C}_{\mathcal{J}}(\nu, U)$, contradicting the
assumption $ux_i \notin \mathcal{C}_{\mathcal{J}}(U)$.
Let $\alpha$ and $\beta$ ($1 \leqslant \alpha, \beta \leqslant n$)
be the largest integers such that
$e_{w_1}^{\alpha} > 0$ and $e_{w_2}^{\beta} > 0$  respectively
(such integers will exist because if $\deg(w_1) = 0$ or
$\deg(w_2) = 0$, we obtain an immediate contradiction
$ux_i \in \mathcal{C}_{\mathcal{J}}(U)$).
We claim that $i > \max\{\alpha, \beta\}$.
\begin{itemize}
\item
If $i < \beta$, then $e_{w}^{\beta} < e_{u}^{\beta}$ which
contradicts $x_{\beta} \in \mathcal{M}_{\mathcal{J}}(w, U \cup\{w\})$ as
$e_{w}^{\gamma} = e_u^{\gamma}$ for all $\gamma > \beta$.
Thus $i \geqslant \beta$. 
\item
If $i < \alpha$, then as
$\beta \leqslant i$ we must have $e_{\mu}^{\gamma} = e_u^{\gamma}$
for all $\alpha < \gamma \leqslant n$. 
Therefore $e_{\mu}^{\alpha} < e_u^{\alpha}
\Rightarrow x_{\alpha} \notin \mathcal{M}_{\mathcal{J}}(\mu, U)$, 
a contradiction; it follows that $i \geqslant \alpha$.
\item
If $i = \alpha$, then either $\beta < \alpha$
or $\beta = \alpha$. If $\beta = \alpha$, then as
$e_{w_1}^i > 0$; $e_{w_2}^i > 0$ and
$e_u^i + 1 = e_{\mu}^i + e_{w_1}^i + e_{w_2}^i$,
we have $e_u^i > e_{\mu}^i
\Rightarrow x_{\alpha} \notin \mathcal{M}_{\mathcal{J}}(\mu, U)$,
a contradiction. If $\beta < \alpha$, then
$e_u^i + 1 = e_{\mu}^i + e_{w_1}^i$. If $e_{w_1}^i \geqslant 2$,
we get the same contradiction as before
($x_{\alpha} \notin \mathcal{M}_{\mathcal{J}}(\mu, U)$). Otherwise
$e_{w_1}^i = 1$ so that $e_u^{\gamma} = e_{\mu}^{\gamma}$
for all $\alpha \leqslant \gamma \leqslant n$. If $w = \mu x_i$, then as
$e_w^{\beta} < e_u^{\beta}$ we have $x_{\beta} \notin
\mathcal{M}_{\mathcal{J}}(w, U \cup \{w\})$, a contradiction.
Else let $\delta$ (where $1 \leqslant \delta < \alpha$)
be the second greatest integer such that
$e_{w_1}^{\delta} > 0$. Then, as $e_{\mu}^{\delta} <
e_u^{\delta}$ and $e_{\mu}^{\gamma} = e_u^{\gamma}$ for
all $\delta < \gamma \leqslant n$,
we have $x_{\delta} \notin \mathcal{M}_{\mathcal{J}}(\mu, U)$,
another contradiction. It follows that $i > \max\{\alpha, \beta\}$,
so that $e_u^{\gamma} = e_{\mu}^{\gamma}$ for all $i < \gamma \leqslant n$
and $e_u^i + 1 = e_{\mu}^i$.
\end{itemize}

If $ux_i \notin \mathcal{C}_{\mathcal{J}}(U)$, then there must exist a
variable $x_k$ (where $1 \leqslant k < i$) such that $e_{w_2}^k > 0$ and
$x_k \notin \mathcal{M}_{\mathcal{J}}(\mu, U)$. 
Because $e_{w_1}^{\alpha} > 0$,
we can use condition (b) of Definition \ref{ConstC} to give
us a monomial $\mu_1 \in U$ and a monomial $w_3$ multiplicative
for $\mu_1$ over $U$ ($e_{w_3}^{\gamma} > 0 \Rightarrow
x_{\gamma} \in \mathcal{M}_{\mathcal{J}}(\mu_1, U)$ for all
$1 \leqslant \gamma \leqslant n$) such that
\begin{eqnarray*}
ux_i & = & \mu w_1w_2 \\
& = & \mu x_kw_1\left(\frac{w_2}{x_k}\right) \\
& = & \mu_1 w_3w_1\left(\frac{w_2}{x_k}\right).
\end{eqnarray*}
If $\mu_1 \mid_{\mathcal{J}} ux_i$, then
the proof is complete, with $\nu = \mu_1$.
Otherwise there must be a variable $x_{k'}$
appearing in the monomial $w_1(\frac{w_2}{x_k})$ such that
$x_{k'} \notin \mathcal{M}_{\mathcal{J}}(\mu_1, U)$. To use 
condition (b) of Definition \ref{ConstC} to yield
a monomial $\mu_2 \in U$ and a monomial $w_4$ multiplicative
for $\mu_2$ over $U$ such that
$$
\mu_1 w_3w_1\left(\frac{w_2}{x_k}\right)
= \mu_2 w_4 \left(\frac{w_1w_2}{x_kx_{k'}}\right)w_3,
$$
it is sufficient to demonstrate that at least one variable appearing in
the monomial $w_3w_1(\frac{w_2}{x_k})$ is
multiplicative for $\mu_1$ over the set $U$. We will
do this by showing that $x_{\alpha} \in 
\mathcal{M}_{\mathcal{J}}(\mu_1, U)$
(recall that $e_{w_1}^{\alpha} > 0$).

By the definition of the Janet involutive division,
\begin{equation} \label{cstyE1}
e_{\mu_1}^{\gamma} = e_{\mu}^{\gamma} \; \: \mbox{for all} \;
k < \gamma \leqslant n
\end{equation}
and
\begin{equation} \label{cstyE2}
e_{\mu_1}^k = e_{\mu}^k + 1,
\end{equation}
so that $\mu < \mu_1$ in the InvLex monomial ordering.
If we can show that $\alpha > k$, then it is clear from 
Equation (\ref{cstyE1})
and $x_{\alpha} \in \mathcal{M}_{\mathcal{J}}(\mu, U)$ that
$x_{\alpha} \in \mathcal{M}_{\mathcal{J}}(\mu_1, U)$.
\begin{itemize}
\item
If $\alpha > \beta$, then $\alpha > k$ because $k \leqslant \beta$
by definition.
\item
If $\alpha = \beta$, then $\alpha > k$ if $k < \beta$;
otherwise $k = \beta$ in which case
% but this gives a contradiction because
$x_{\alpha} \in \mathcal{M}_{\mathcal{J}}(\mu, U)$ is contradicted by
Equations (\ref{cstyE1}) and (\ref{cstyE2}).
\item
If $\alpha < \beta$, then $e_{\mu}^{\gamma} =
e_w^{\gamma}$ for all $\alpha < \gamma \leqslant n$.
Thus $k \leqslant \alpha$
otherwise $x_{k} \in \mathcal{M}_{\mathcal{J}}(w, U\cup\{w\}) \Rightarrow
x_{k} \in \mathcal{M}_{\mathcal{J}}(\mu, U)$, a contradiction. Further,
$k = \alpha$ is not allowed because
$x_{\alpha} \in \mathcal{M}_{\mathcal{J}}(\mu, U)$ and
$x_k \notin \mathcal{M}_{\mathcal{J}}(\mu, U)$ cannot both be true;
therefore $\alpha > k$ again.
\end{itemize}

If $\mu_2 \mid_{\mathcal{J}} ux_i$, then
the proof is complete, with $\nu = \mu_2$.
Otherwise we proceed by induction to obtain
the sequence shown below (Equation (\ref{seqC})), which is valid because
$\mu_{\sigma-1} < \mu_{\sigma}$ (for $\sigma \geqslant 2$)
in the InvLex monomial ordering allows us to prove that
the variable $x_{\alpha}$ (that appears in the monomial $w_1$)
is multiplicative (over $U$) for the monomial $\mu_{\sigma}$;
this in turn enables us to construct the next entry in the
sequence by using condition (b) of Definition \ref{ConstC}.
\begin{equation} \label{seqC}
\mu w_1w_2 = \mu_1 w_3w_1\left(\frac{w_2}{x_k}\right)
= \mu_2w_4\left(\frac{w_1w_2}{x_kx_{k'}}\right)w_3 =
\mu_3w_5\left(\frac{w_1w_2w_3}{x_kx_{k'}x_{k''}}\right)w_4 = \cdots
\end{equation}
Because $\mu < \mu_1 < \mu_2 < \cdots$ in the InvLex monomial ordering,
elements of the sequence $\mu, \mu_1, \mu_2, \hdots$
are distinct. It follows that the sequence in Equation
(\ref{seqC}) is finite (terminating with the required $\nu$)
because $\mu$ and the $\mu_{\sigma}$ (for $\sigma \geqslant 1$) are
all divisors of the monomial $ux_i$, of which there are only a finite
number of.
\end{pf}

\begin{remark}
The above proof that Janet is a constructive involutive
division does not use the property of Janet being a 
continuous involutive division, unlike the proofs
found in \cite{Gerdt98a} and \cite{Seiler02b}.
\end{remark}

\section{The Involutive Basis Algorithm} \label{4point4}

To compute an Involutive Basis for an ideal
$J$ with respect to some admissible monomial ordering $O$
and some involutive division $I$, it is sufficient
to compute a Locally Involutive Basis for $J$ with
respect to $I$ and $O$ if $I$ is continuous; and we can
compute this Locally Involutive Basis by considering
only prolongations whose lead monomials lie outside
the current involutive span if $I$ is constructive.
Let us now consider Algorithm \ref{com-inv}, an algorithm to
construct an Involutive Basis for $J$ (with respect
to $I$ and $O$) in exactly this manner.

\begin{algorithm}
\setlength{\baselineskip}{3.5ex}
\caption{The Commutative Involutive Basis Algorithm}
\label{com-inv}
\begin{algorithmic}
\vspace*{2mm}
\REQUIRE{A Basis $F = \{f_1, f_2, \hdots, f_m\}$ for an ideal $J$
         over a commutative polynomial ring $R[x_1, \hdots x_n]$;
         an admissible monomial ordering $O$; a continuous and
	 constructive involutive division $I$.}
\ENSURE{An Involutive Basis $G = \{g_1, g_2, \hdots, g_p\}$ for $J$ (in the
        case of termination).}
\vspace*{1mm}
\STATE
$G  = \emptyset$; \\
$F = \mathrm{Autoreduce}(F)$;
\WHILE{($G == \emptyset$)}
\STATE
$S = \{x_if \mid f \in F, \: x_i \notin \mathcal{M}_I(f, F)\}$; \\
$s' = 0$;
\WHILE{($S \neq \emptyset$) {\bf and} ($s' == 0$)}
\STATE
Let $s$ be a polynomial in $S$ whose lead monomial is minimal 
with respect to $O$; \\
$S = S\setminus \{s\}$; \\
$s' = \Rem_I(s,  F)$;
\ENDWHILE
\IF{($s' \neq 0$)}
\STATE
$F = \mathrm{Autoreduce}(F\cup\{s'\})$;
\ELSE
\STATE
$G = F$;
\ENDIF
\ENDWHILE
\STATE
{\bf return} $G$;
\end{algorithmic}
\vspace*{1mm}
\end{algorithm}

The algorithm starts by autoreducing the input
basis $F$ using Algorithm \ref{com-auto}. We then
construct a set $S$ containing all the
possible prolongations of elements of $F$, before
recursively (a) picking a polynomial $s$ from $S$ such that
$\LM(s)$ is minimal in the chosen monomial ordering; (b)
removing $s$ from $S$; and (c) finding the involutive
remainder $s'$ of $s$ with respect to $F$. 

If during this loop a remainder $s'$ is found that 
is nonzero, we exit the loop and autoreduce
the set $F\cup\{s'\}$, continuing thereafter to
construct a new set $S$ and repeating the above
process on this new set. If however all the prolongations
in $S$ involutively reduce to zero, then by definition
$F$ is a Locally Involutive Basis, and so we can exit
the algorithm with this basis. The correctness of 
Algorithm \ref{com-inv} is therefore clear;
termination however requires us to show that each
% The correctness of Algorithm \ref{com-inv} is clear
% from the above discussion; the question of termination
% still remains however, a question we shall
% answer by showing that the algorithm always terminates
involutive division used
with the algorithm is {\it Noetherian} and {\it stable}.

\begin{defn}
\index{involutive division!Noetherian}
An \index{Noetherian involutive division}
involutive division $I$ is {\it Noetherian} if,
given any finite set of monomials $U$, there is a finite
Involutive Basis $V \supseteq U$ with respect to $I$
and some arbitrary admissible monomial ordering $O$.
\end{defn}

\begin{prop}
The Thomas and Janet divisions are Noetherian.
\end{prop}
\begin{pf}
Let $U = \{u_1, \hdots, u_m\}$ be an arbitrary set of
monomials over a polynomial ring 
$\mathcal{R} = R[x_1, \hdots, x_n]$
generating an ideal $J$. We will explicitly construct
an Involutive Basis $V$ for $U$ with respect to some
arbitrary admissible monomial ordering $O$.

{\bf Janet (Adapted from \cite{Seiler02b}, Lemma 2.13).}
Let $\mu \in \mathcal{R}$ be the monomial with multidegree
$(e_{\mu}^1, e_{\mu}^2, \hdots, e_{\mu}^n)$ defined as
follows: $e_{\mu}^i = \max_{u \in U} e_u^i$
($1 \leqslant i \leqslant n$). We claim that the set $V$
containing all monomials $v \in J$ such that $v \mid \mu$
is an Involutive Basis for $U$ with respect to the
Janet involutive division and $O$. To prove the claim, first note
that $V$ is a basis for $J$ because
$U \subseteq V$ and $V \subset J$; to prove that $V$ is a
Janet Involutive Basis for $J$ we have to show that all multiples
of elements of $V$ involutively reduce to zero using $V$,
which we shall do by showing that all members of the ideal
involutively reduce to zero using $V$.
 
Let $p$ be an arbitrary element of $J$. 
If $p \in V$, then trivially $p \in \mathcal{C}_{\mathcal{J}}(V)$
and so $p$ involutively reduces to zero using $V$.
Otherwise set $X = \{x_i \; \mbox{such that} \;
e_{\LM(p)}^i > e_{\mu}^i\}$,
and define the monomial $p'$ by $e_{p'}^i = e_{\LM(p)}^i$ for
$x_i \notin X$; and $e_{p'}^i = e_{\mu}^i$ for
$x_i \in X$ (so that $e_{p'}^i = \min\{e_{\LM(p)}^i, e_{\mu}^i\}$).
By construction of the set $V$ and by the definition of
$\mu$, it follows that $v' \in V$ and
$X \subseteq \mathcal{M}_{\mathcal{J}}(p', V)$. But this implies that
$\LM(p) \in \mathcal{C}_{\mathcal{J}}(p', V)$, and thus
$p \xymatrix{\ar[r]_{\mathcal{J}}_(1){p'} &}(p-\LM(p))$.
By induction and by the admissibility of $O$, 
$p \xymatrix{\ar[r]_{\mathcal{J}}_(1){V} &}0$ and thus
$V$ is a finite Janet Involutive Basis for $J$.

{\bf Thomas.} We use the same proof as for Janet above,
replacing ``Janet'' by ``Thomas'' and ``$\mathcal{J}$''
by ``$\mathcal{T}$''.
\end{pf}

\begin{prop}
The Pommaret division is not Noetherian.
\end{prop}
\begin{pf}
Let $J$ be the ideal generated by the monomial
$u := xy$ over the polynomial ring
$\mathbb{Q}[x, y]$. For the Pommaret division, 
$\mathcal{M}_{\mathcal{P}}(u) = \{x\}$,
and it is clear that $\mathcal{M}_{\mathcal{P}}(v) = \{x\}$ for all
$v \in J$ as $v \in J \Rightarrow v = (xy)p$ for some
polynomial $p$. It follows that no finite Pommaret Involutive
Basis exists for $J$ as no prolongation by the variable $y$
of any polynomial $p \in J$ is involutively divisible by
some other polynomial $p' \in J$; the Pommaret Involutive Basis
for $J$ is therefore the infinite basis
$\{xy, \; xy^2, \; xy^3, \; \hdots\}$.
\end{pf}

\begin{defn}
\index{involutive division!stable}
Let \index{stable involutive division}
$u$ and $v$ be two distinct monomials such that $u \mid v$.
An involutive division $I$ is {\it stable}
if $\Rem_I(v, \{u,v\}, \{u\}) = v$. In other words,
$u$ is not an involutive divisor of $v$ with
respect to $I$ when multiplicative variables are taken
over the set $\{u, v\}$.
\end{defn}

\begin{prop}
The Thomas and Janet divisions are stable.
\end{prop}
\begin{pf}
Let $u$ and $v$ have corresponding
multidegrees $(e_u^1, \hdots, e_u^n)$ and
$(e_v^1, \hdots, e_v^n)$. If $u \mid v$ and if $u$ and $v$ are
different, then we must have
$e_u^i < e_v^i$ for at least one $1 \leqslant i \leqslant n$.

{\bf Thomas}. By definition,
$x_i \notin \mathcal{M}_{\mathcal{T}}(u, \{u,v\})$,
so that $\Rem_{\mathcal{T}}(v, \{u,v\}, \{u\}) = v$.
% $u$ is not a Thomas involutive divisor of $v$
% when multiplicative variables are taken
% over the set $\{u, v\}$.

{\bf Janet}. Let $j$ be the greatest integer
such that $e_u^j < e_v^j$. Then, as $e_u^k = e_v^k$ for
all $j < k \leqslant n$, it follows that
$x_j \notin \mathcal{M}_{\mathcal{J}}(u, \{u,v\})$,
and so $\Rem_{\mathcal{J}}(v, \{u,v\}, \{u\}) = v$.
% $u$ is not a Janet involutive divisor of $v$
% when multiplicative variables are taken
% over the set $\{u, v\}$.
\end{pf}

\begin{prop}
The Pommaret division is not stable.
\end{prop}
\begin{pf}
Consider the two monomials $u := x$ and $v := x^2$ over the polynomial
ring $\mathbb{Q}[x]$. Because $\mathcal{M}_{\mathcal{P}}(u, \{u, v\})
= \{x\}$, it is clear that $u \mid_{\mathcal{P}} v$, and so
the Pommaret involutive division is not stable.
\end{pf}

\begin{remark}
Stability ensures that any set of distinct monomials is
autoreduced. In particular, if a set $U$ of monomials is
autoreduced, and if we add a monomial $u \notin U$ to $U$,
then the resultant set $U\cup\{u\}$ is also autoreduced.
This contradicts a statement made on page 24 of
\cite{Seiler02b}, where it is claimed that if we add
an involutively irreducible prolongation $ux_i$ of a monomial $u$
from an autoreduced set of monomials $U$
to that set, then the resultant set is also autoreduced
regardless of whether or not the involutive division used is
stable\footnote{This claim is integral to the proof of
Theorem 6.4 in \cite{Seiler02b}, a theorem that states than an
algorithm corresponding to 
Algorithm \ref{com-inv} in this thesis terminates.}.
For a counterexample, consider the set
of monomials $U := \{u_1, u_2\} = \{xy, x^2y^2\}$ over
the polynomial ring $\mathbb{Q}[x,y]$,
% let the monomial ordering be DegLex; 
and let the involutive division be Pommaret. 
\begin{center}
\begin{tabular}{c|c}
$u$ & $\mathcal{M}_{\mathcal{P}}(u, U)$ \\ \hline
$xy$ & $\{x\}$ \\
$x^2y^2$ & $\{x\}$ \\ \hline
\end{tabular}
\end{center}
Because the variable $y$ is
nonmultiplicative for the monomial $xy$, it is clear
that the set $U$ is autoreduced.
Consider the prolongation $xy^2$ of the monomial $u_1$ by
the variable $y$. This prolongation is
involutively irreducible with respect to $U$, but if
we add the prolongation to $U$ to obtain the set
$V := \{v_1, v_2, v_3\} = \{xy, x^2y^2, xy^2\}$, then
$v_3$ will involutively reduce $v_2$, contradicting the
claim that the set $V$ is autoreduced.
\begin{center}
\begin{tabular}{c|c}
$v$ & $\mathcal{M}_{\mathcal{P}}(v, V)$ \\ \hline
$xy$ & $\{x\}$ \\
$x^2y^2$ & $\{x\}$ \\
$xy^2$ & $\{x\}$ \\ \hline
\end{tabular}
\end{center}
\end{remark}

\begin{prop} \label{InvTermNoeth}
Algorithm \ref{com-inv} always terminates when used with
a Noetherian and stable involutive division.
\end{prop}
\begin{pf}
Let $I$ be a Noetherian and stable involutive division, and
consider the computation (using Algorithm \ref{com-inv})
of an Involutive Basis for a set of polynomials $F$
with respect to $I$ and some admissible monomial ordering $O$.
The algorithm begins by autoreducing $F$ to give a basis
(which we shall denote by $F_1$) generating the same
ideal $J$ as $F$. Each pass of the algorithm then produces
a basis $F_{i+1} = \mathrm{Autoreduce}(F_i\cup\{s'_i\})$ generating $J$
($i \geqslant 1$), where each $s'_i \neq 0$ is an involutively
reduced prolongation. Consider the monomial ideal
$\langle \LM(F_i)\rangle$ generated by the lead
monomials of the set $F_i$. {\bf Claim:}
\begin{equation} \label{AMon}
\langle \LM(F_1)\rangle \subseteq \langle \LM(F_2)\rangle
\subseteq \langle \LM(F_3)\rangle \subseteq \cdots
\end{equation}
is an ascending chain of monomial ideals. 

{\bf Proof of Claim:} It
is sufficient to show that if an arbitrary polynomial $f \in F_i$ does
not appear in $F_{i+1}$, then there must be a polynomial
$f' \in F_{i+1}$ such that $\LM(f') \mid \LM(f)$. It is clear that
such an $f'$ will exist if the lead monomial of $f$
is not reduced during autoreduction; otherwise a polynomial $p$
reduces the lead monomial of $f$ during autoreduction, so that
$\LM(p) \mid_I \LM(f)$. If there exists a 
polynomial $p' \in F_{i+1}$ such that
$\LM(p') = \LM(p)$, we are done; otherwise we proceed by induction on
$p$ to obtain a polynomial $q$ such that $\LM(q) \mid_I \LM(p)$.
Because $\deg(\LM(f)) > \deg(\LM(p)) > \deg(\LM(q)) > \cdots$, this
process is guaranteed to terminate with the required $f'$.
\hfill ${}_{\Box}$

By the Ascending Chain Condition (Corollary \ref{ACC}),
the chain in Equation
(\ref{AMon}) must eventually become constant, so
there must be an integer $N$ ($N \geqslant 1$) such that
$$
\langle \LM(F_N)\rangle = \langle \LM(F_{N+1})\rangle = \cdots.
$$
{\bf Claim:} If $F_{k+1} = \mathrm{Autoreduce}(F_k\cup\{s'_k\})$
for some $k \geqslant N$, then $\LM(s'_k) = \LM(fx_i)$ for some
polynomial $f \in F_k$ and some variable
$x_i \notin \mathcal{M}_I(f, F_k)$ such that
$s'_k = \Rem_I(fx_i, F_k)$. 

{\bf Proof of Claim:} Assume to the contrary
that $\LM(s'_k) \neq \LM(fx_i)$. Then because $s'_k = \Rem_I(fx_i, F_k)$,
it follows that $\LM(s'_k) < \LM(fx_i)$. But
$\langle \LM(F_k)\rangle = \langle \LM(F_{k+1})\rangle$,
so that $\LM(s'_k) = \LM(f'u)$ for some $f' \in F_k$
and some monomial $u$ containing at least one variable
$x_j \notin \mathcal{M}_I(f', F_k)$
(otherwise $s'_k$ can be involutively reduced with respect
to $F_k$, a contradiction). 

Because $O$ is admissible, $1 \leqslant \frac{u}{x_j}$ and therefore
$x_j \leqslant u$, so that
$\LM(f'x_j) \leqslant \LM(f'u) < \LM(fx_i)$. But the prolongation
$fx_i$ was chosen so that its lead monomial is minimal amongst the
lead monomials of all prolongations of elements of $F_k$ that
do not involutively reduce to zero; the prolongation
$f'x_k$ must therefore involutively reduce to zero, so that
$\LM(f'x_j) = \LM(f''u')$ for some polynomial $f'' \in F_k$ and some
monomial $u'$ that is multiplicative for $f''$ over $F_k$.
But $s'_k$ is involutively irreducible with respect to $F_k$,
so a variable $x'_j \notin \mathcal{M}_I(f'', F_k)$ must
appear in the monomial $\frac{u}{x_j}$.

It is now clear that
we can construct a sequence $f'x_j, \, f''x'_j, \hdots$ of
prolongations. But $I$ is continuous, so all elements in the
corresponding sequence $\LM(f'), \, \LM(f''), \hdots$ of
monomials must be distinct. Because $F_k$ is finite, 
it follows that the sequence of prolongations will terminate with
a prolongation that does not involutively reduce to zero and
whose lead monomial is less than the monomial $\LM(fx_i)$,
contradicting our assumptions.
% , which contradicts
% the fact that the prolongation $fx_i$ was chosen so that
% its lead monomial is minimal amongst the lead monomials of
% all prolongations of elements of $F_k$.
Thus $\LM(s'_k)$ for
$k \geqslant N$ is always equal to the lead monomial of
some prolongation of some polynomial $f \in F_k$.
\hfill ${}_{\Box}$

% We can use the above argument to show that the sets
% $\LM(F_k)$ and $\LM(F_{k+1})$ are identical up until
% the monomial $\LM(s'_k)$:
% \begin{itemize}
% \item
% For all polynomials $f' \in F_k$ such that
% $\LM(f') < \LM(s'_k)$, a polynomial $g'$ will appear in
% $F_{k+1}$ such that $\LM(g') = \LM(f')$;
% \item
% No polynomial $g'$ will appear in $F_{k+1}$ such that
% $\LM(g') < \LM(s'_k)$ and no polynomial
% $f' \in F_k$ exists such that $\LM(f') = \LM(g')$;
% \item
% A polynomial $g'$ appears in $F_{k+1}$ such that
% $\LM(g') = \LM(s'_k)$.
% \end{itemize}
%
% Prove that $\LM(F_{k+1}) = \LM(F_k) \cup \LM(fx_i)$ for
% all $k \geqslant N$.

Consider now the set of monomials $\LM(F_{k+1})$. 
{\bf Claim:} $\LM(F_{k+1}) = \LM(F_k) \cup \{\LM(s'_k)$\} for
all $k \geqslant N$, so that when
autoreducing the set $F_k \cup \{s'_k\}$, no leading
monomial is involutively reducible.

{\bf Proof of Claim:} Consider an arbitrary polynomial
$p \in F_k \cup \{s'_k\}$. If $p = s'_k$, then (by
definition) $p$ is irreducible with respect to the
set $F_k$, and so (by condition (b) of Definition
\ref{inv-div-defn}) $p$ will also be irreducible with
respect to the set $F_k \cup \{s'_k\}$.
If $p \neq s'_k$, then $p$ is irreducible with respect
to the set $F_k$ (as the set $F_k$ is autoreduced),
and so (again by condition (b) of Definition \ref{inv-div-defn})
the only polynomial in the set $F_k \cup \{s'_k\}$
that can involutively reduce the polynomial $p$ is
the polynomial $s'_k$. But $I$ is stable, so that
$s'_k$ cannot involutively reduce $\LM(p)$. It follows that
a polynomial $p'$ will appear in the autoreduced set
$F_{k+1}$ such that $\LM(p') = \LM(p)$, and thus
$\LM(F_{k+1}) = \LM(F_k) \cup \{\LM(s'_k)$\} as required.
\hfill ${}_{\Box}$

For the final part of the proof, consider the basis $F_N$.
Because $I$ is Noetherian, there exists a finite Involutive
Basis $G_N$ for the ideal generated by the set of lead
monomials $\LM(F_N)$, where $G_N \supseteq \LM(F_N)$.
Let $fx_i$ be the prolongation chosen during the $N$-th
iteration of Algorithm \ref{com-inv}, so that $\LM(fx_i)
\notin \mathcal{C}_I(F_N)$. Because $G_N$ is an
Involutive Basis for $\LM(F_N)$, there must be a monomial
$g \in G_N$ such that $g \mid_I \LM(fx_i)$. % with respect to $G_N$.
{\bf Claim:} $g = \LM(fx_i)$.

{\bf Proof of Claim:} We proceed by
showing that if $g \neq \LM(fx_i)$,
then $g \in \mathcal{C}_I(\LM(F_N))$ so that (because
of condition (b) of Definition \ref{inv-div-defn})
$\LM(fx_i) \in \mathcal{C}_I(G_N) \Rightarrow
\LM(fx_i) \in \mathcal{C}_I(g, \LM(F_N)\cup \{g\})$,
contradicting the constructivity of $I$
(Definition \ref{ConstC}).

Assume that $g \neq \LM(fx_i)$.
Because $\langle G_N\rangle = \langle\LM(F_N)\rangle$,
there exists a polynomial $f_1 \in F_N$ such that
$\LM(f_1) \mid g$. If $\LM(f_1) \mid_I g$ with respect
to $F_N$, then we are done. Otherwise $\LM(g) = \LM(f_1)u_1$
for some monomial $u_1 \neq 1$ containing at least one
variable $x_{j_1} \notin \mathcal{M}_I(f_1, F_N)$.
Because $\deg(g) < \deg(\LM(fx_i))$ and $\LM(f_1)x_{j_1} \mid
\LM(fx_i)$, we must have $\LM(f_1)x_{j_1} < \LM(fx_i)$
with respect to our chosen monomial ordering, so that
$\LM(f_1)x_{j_1} \in \mathcal{C}_I(F_N)$ by definition of
how the prolongation $fx_i$ was chosen. It follows that
there exists a polynomial $f_2 \in F_N$ such that
$\LM(f_2) \mid_I \LM(f_1)x_{j_1}$ with respect
to $F_N$. If $\LM(f_2) \mid_I g$ with
respect to $F_N$, then we are done. Otherwise we
iterate ($\LM(f_1)x_{j_1} = \LM(f_2)u_2$ for some
monomial $u_2$ containing at least one variable
$x_{j_2} \notin \mathcal{M}_I(f_2, F_N)$\ldots) to obtain
the sequence $(f_1, f_2, f_3, \hdots)$ of polynomials, 
where the lead monomial
of each element in the sequence divides $g$ and
$\LM(f_{k+1}) \mid_I \LM(f_k)x_{j_k}$ with respect to $F_N$ for
all $k \geqslant 1$. Because $I$ is continuous, this sequence
must be finite, terminating with a polynomial $f_k \in F_N$
(for some $k \geqslant 1$) such that $f_k \mid_I g$ with
respect to $F_N$.  \hfill ${}_{\Box}$

It follows that during the $N$-th iteration of the algorithm, a polynomial
is added to the current basis $F_N$ whose lead monomial
is a member of the Involutive Basis $G_N$. By induction,
every step of the algorithm after the $N$-th step also
adds a polynomial to the current basis whose lead monomial
is a member of $G_N$. Because $G_N$ is a finite set, after
a finite number of steps the basis $\LM(F_k)$ (for some
$k \geqslant N$) will contain all the elements of $G_N$.
We can therefore deduce that $\LM(F_k) = G_N$; it follows 
that $\LM(F_k)$ is an Involutive Basis, and so $F_k$ is also an
Involutive Basis.
% Claim: $\LM(F_k)$ is an Involutive Basis. To prove the claim,
% we have to show that all prolongations $fx_i$ ($f \in F_k$,
% $x_i \notin \mathcal{M}_I(f, F_k)$) involutively
% reduce to zero. If this is not the case, then one such
% prolongation involutively reduces to give a polynomial
% $s'_k \neq 0$ such that $\LM(s'_k) = \LM(fx_i)$. But we know
% that a polynomial $g \in G_N$ exists such that
% $\LM(g) = \LM(fx_i)$ and, because $G_N \subseteq \LM(F_k)$,
% it follows that there exists a polynomial $f' \in F_k$
% such that $\LM(f') = g$, contradicting the fact that
% $s'_k$ is involutively irreducible with respect to $F_k$.
% It follows that $\LM(F_k)$ is an Involutive Basis, and
% therefore $F_k$ is also an Involutive Basis.
\end{pf}

\begin{thm} \label{IisG}
Every Involutive Basis is a Gr\"obner Basis.
\end{thm}
\begin{pf}
Let $G = \{g_1, \hdots, g_m\}$ be an Involutive
Basis with respect to some involutive division $I$ and
some admissible monomial ordering $O$, where each
$g_i \in G$ (for all $1 \leqslant i \leqslant m$) is a
member of the polynomial ring $R[x_1, \hdots, x_n]$.
To prove that
$G$ is a Gr\"obner Basis, we must show that all
S-polynomials
$$\mathrm{S\mbox{-}pol}(g_i, g_j) =
\frac{\lcm(\LM(g_i), \LM(g_j))}{\LT(g_i)}g_i -
\frac{\lcm(\LM(g_i), \LM(g_j))}{\LT(g_j)}g_j$$
conventionally reduce to zero using $G$ ($1 \leqslant i,j \leqslant m$,
$i \neq j$). Because $G$ is an Involutive Basis, it is clear that
$\frac{\lcm(\LM(g_i), \LM(g_j))}{\LT(g_i)}g_i
\xymatrix{\ar[r]_{I}_(1){G} &}0$ and
$\frac{\lcm(\LM(g_i), \LM(g_j))}{\LT(g_j)}g_j
\xymatrix{\ar[r]_{I}_(1){G} &}0$. By Proposition
\ref{plusC}, it follows that
$\mathrm{S\mbox{-}pol}(g_i, g_j)
\xymatrix{\ar[r]_{I}_(1){G} &}0$. But every involutive
reduction is a conventional reduction, so we can deduce that
$\mathrm{S\mbox{-}pol}(g_i, g_j) \rightarrow_G 0$ as required.
\end{pf}

\begin{lem}
\index{unique remainder}
Remainders \index{remainder!unique}
are involutively unique with respect to Involutive Bases.
\end{lem}
\begin{pf}
Given an Involutive Basis $G$ with respect to some
involutive division $I$ and some admissible monomial
ordering $O$, Theorem \ref{IisG} tells us that $G$ is
a Gr\"obner Basis with respect to $O$ and thus
remainders are conventionally unique with respect to
$G$. To prove that remainders are involutively unique
with respect to $G$, we must show that the conventional and
involutive remainders of an arbitrary polynomial $p$ with respect
to $G$ are identical. For this it is sufficient to show
that a polynomial $p$ is conventionally reducible by $G$ if
and only if it is involutively reducible by $G$.
($\Rightarrow$) Trivial as every involutive reduction is a
conventional reduction. ($\Leftarrow$)
If a polynomial $p$ is conventionally reducible by a
polynomial $g \in G$, it follows that $\LM(p) = \LM(g)u$ for
some monomial $u$. But $G$ is an Involutive Basis, so there
must exist a polynomial $g' \in G$ such that
$\LM(g)u = \LM(g')u'$ for some monomial $u'$ that is
multiplicative (over $G$) for $g'$. Thus $p$ is also
involutively reducible by $G$.
\end{pf}

\begin{example}
Let us return to our favourite example of an ideal $J$
generated by the set of polynomials $F := \{f_1, f_2\} =
\{x^2-2xy+3, \; 2xy+y^2+5\}$ over the polynomial ring
$\mathbb{Q}[x, y, z]$. To compute an Involutive Basis for $F$ 
with respect to the DegLex monomial ordering and
the Janet involutive division $\mathcal{J}$,
we apply Algorithm \ref{com-inv} to $F$, in which the
first task is to autoreduce $F$. This produces the set
$F = \{f_2, \; f_3\} = \{2xy+y^2+5, \; x^2+y^2+8\}$ as output
(because $f_1 = x^2-2xy+3 \xymatrix{\ar[r]_{\mathcal{J}}_(1){f_2} &}
x^2+y^2+8 =: f_3$ and $f_2$ is involutively irreducible
with respect to $f_3$), with multiplicative variables
as shown below.
\begin{center}
\begin{tabular}{c|c}
Polynomial & $\mathcal{M}_{\mathcal{J}}(f_i, F)$ \\ \hline
$f_2 = 2xy+y^2+5$ & $\{x, y\}$ \\
$f_3 = x^2+y^2+8$ & $\{x\}$ \\ \hline
\end{tabular}
\end{center}
The first set of prolongations of elements of $F$
is the set $S = \{f_3y\} = \{x^2y+y^3+8y\}$. As this set
only has one element, it is clear that on entering the
second while loop of the algorithm, we must remove the
polynomial $s = x^2y+y^3+8y$ from $S$ and involutively
reduce $s$ with respect to $F$ to give the polynomial
$s' = \frac{5}{4}y^3-\frac{5}{2}x+\frac{37}{4}y$ as follows.
\begin{eqnarray*}
s = x^2y + y^3 + 8y & \xymatrix{\ar[r]_{\mathcal{J}}_(1){f_2} &} &
x^2y + y^3 + 8y - \frac{1}{2}x(2xy+y^2+5) \\
& = & -\frac{1}{2}xy^2 + y^3 - \frac{5}{2}x + 8y \\
& \xymatrix{\ar[r]_{\mathcal{J}}_(1){f_2} &} &
-\frac{1}{2}xy^2 + y^3 - \frac{5}{2}x + 8y  + \frac{1}{4}y(2xy+y^2+5) \\
& = & \frac{5}{4}y^3 - \frac{5}{2}x + \frac{37}{4}y = s' =: f_4.
\end{eqnarray*}
As the prolongation did not involutively reduce to zero,
we exit from the second while loop of the algorithm and proceed by
autoreducing the set $F\cup\{f_4\} =
\{2xy+y^2+5, \; x^2+y^2+8, \; \frac{5}{4}y^3 - \frac{5}{2}x
+ \frac{37}{4}y\}$. This process does not alter the set, so now
we consider the prolongations of the three element set
$F = \{f_2, \; f_3, \; f_4\}$.
\begin{center}
\begin{tabular}{c|c}
Polynomial & $\mathcal{M}_{\mathcal{J}}(f_i, F)$ \\ \hline
$f_2 = 2xy+y^2+5$ & $\{x\}$ \\
$f_3 = x^2+y^2+8$ & $\{x\}$ \\
$f_4 = \frac{5}{4}y^3 - \frac{5}{2}x + \frac{37}{4}y$ & $\{x, y\}$ \\ \hline
\end{tabular}
\end{center}
We see that there are 2 prolongations to consider,
so that $S = \{f_2y, \; f_3y\} = 
\{2xy^2+y^3+5y, \; x^2y+y^3+8y\}$. As $xy^2 < x^2y$ in the
DegLex monomial ordering, we must consider 
the prolongation $f_2y$ first.
\begin{eqnarray*}
f_2y = 2xy^2 + y^3 + 5y
& \xymatrix{\ar[r]_{\mathcal{J}}_(1){f_4} &} &
2xy^2 + y^3 + 5y -
\frac{4}{5}\left(\frac{5}{4}y^3 - \frac{5}{2}x + \frac{37}{4}y\right) \\
& = & 2xy^2 + 2x - \frac{12}{5}y =: f_5.
\end{eqnarray*}
As before, the prolongation did not involutively reduce
to zero, so now we autoreduce the set $F \cup \{f_5\} = 
\{2xy+y^2+5, \; x^2+y^2+8, \; \frac{5}{4}y^3 - \frac{5}{2}x
+ \frac{37}{4}y, \; 2xy^2 + 2x - \frac{12}{5}y\}$. Again
this leaves the set unchanged, so we proceed with the
set $F = \{f_2, \; f_3, \; f_4, \; f_5\}$.
\begin{center}
\begin{tabular}{c|c}
Polynomial & $\mathcal{M}_{\mathcal{J}}(f_i, F)$ \\ \hline
$f_2 = 2xy+y^2+5$ & $\{x\}$ \\
$f_3 = x^2+y^2+8$ & $\{x\}$ \\
$f_4 = \frac{5}{4}y^3 - \frac{5}{2}x + \frac{37}{4}y$ & $\{x, y\}$ \\
$f_5 = 2xy^2 + 2x - \frac{12}{5}y$ & $\{x\}$ \\ \hline
\end{tabular}
\end{center}
This time, $S = \{f_2y, \; f_3y, \; f_5y\} = 
\{2xy^2+y^3+5y, \; x^2y+y^3+8y, \; 2xy^3 + 2xy - \frac{12}{5}y^2\}$,
and we must consider the prolongation $f_2y$ first.
\begin{eqnarray*}
f_2y = 2xy^2 + y^3 + 5y
& \xymatrix{\ar[r]_{\mathcal{J}}_(1){f_5} &} &
2xy^2 + y^3 + 5y - \left(2xy^2 + 2x - \frac{12}{5}y\right) \\
& = & y^3-2x+\frac{37}{5}y \\
& \xymatrix{\ar[r]_{\mathcal{J}}_(1){f_4} &} &
y^3-2x+\frac{37}{5}y - 
\frac{4}{5}\left(\frac{5}{4}y^3 - \frac{5}{2}x + \frac{37}{4}y\right) \\
& = & 0.
\end{eqnarray*}
Because the prolongation involutively reduced to zero,
we move on to look at the next prolongation $f_3y$ (which comes from the
revised set $S = \{f_3y, \; f_5y\} =
\{x^2y+y^3+8y, \; 2xy^3 + 2xy - \frac{12}{5}y^2\}$).
\begin{eqnarray*}
f_3y = x^2y + y^3 + 8y & \xymatrix{\ar[r]_{\mathcal{J}}_(1){f_2} &} &
x^2y + y^3 + 8y - \frac{1}{2}x(2xy+y^2+5) \\
& = & -\frac{1}{2}xy^2 + y^3 - \frac{5}{2}x + 8y \\
& \xymatrix{\ar[r]_{\mathcal{J}}_(1){f_5} &} &
-\frac{1}{2}xy^2 + y^3 - \frac{5}{2}x + 8y 
+ \frac{1}{4}\left(2xy^2 + 2x - \frac{12}{5}y\right) \\
& = & y^3-2x+\frac{37}{5}y \\
& \xymatrix{\ar[r]_{\mathcal{J}}_(1){f_4} &} &
y^3-2x+\frac{37}{5}y -
\frac{4}{5}\left(\frac{5}{4}y^3 - \frac{5}{2}x + \frac{37}{4}y\right) \\
& = & 0. 
\end{eqnarray*}
Finally, we look at the prolongation $f_5y$ from the
set $S = \{2xy^3 + 2xy - \frac{12}{5}y^2\}$.
\begin{eqnarray*}
f_5y = 2xy^3 + 2xy - \frac{12}{5}y^2 &
\xymatrix{\ar[r]_{\mathcal{J}}_(1){f_4} &} &
2xy^3 + 2xy - \frac{12}{5}y^2 -
\frac{8}{5}x\left(\frac{5}{4}y^3 - \frac{5}{2}x + \frac{37}{4}y\right) \\
& = & 4x^2 - \frac{64}{5}xy - \frac{12}{5}y^2 \\
& \xymatrix{\ar[r]_{\mathcal{J}}_(1){f_3} &} &
4x^2 - \frac{64}{5}xy - \frac{12}{5}y^2 - 4(x^2+y^2+8) \\
& = & -\frac{64}{5}xy - \frac{32}{5}y^2 - 32 \\
& \xymatrix{\ar[r]_{\mathcal{J}}_(1){f_2} &} &
-\frac{64}{5}xy - \frac{32}{5}y^2 - 32 +
\frac{32}{5}(2xy+y^2+5) \\
& = & 0.
\end{eqnarray*}
Because this prolongation also involutively reduced to zero using $F$,
we are left with $S = \emptyset$, which means that the algorithm now
terminates with the Janet Involutive
Basis $G = \{2xy+y^2+5, \; x^2+y^2+8, \;
\frac{5}{4}y^3 - \frac{5}{2}x + \frac{37}{4}y, \;
2xy^2 + 2x - \frac{12}{5}y\}$ as output.
\end{example}

\section{Improvements to the Involutive Basis Algorithm} \label{4point5}

\subsection{Improved Algorithms}

In \cite{ZharBlink93}, Zharkov and Blinkov introduced
an algorithm for computing an Involutive Basis
and proved its termination for zero-dimensional
ideals. This work led other researchers to produce
improved versions of the algorithm
(see for example \cite{Apel98b}, \cite{CHS00}, \cite{Gerdt02}, \cite{Gerdt97},
\cite{Gerdt01a} and \cite{Gerdt01b}); improvements made to
the algorithm include the introduction of selection strategies
\index{selection strategies}
(which, as we have seen in the proof of Proposition \ref{InvTermNoeth},
are crucial for proving the termination of the algorithm in general),
and the introduction of criteria (analogous to
Buchberger's criteria) allowing the {\it a priori} detection
of prolongations that involutively reduce to zero.

\subsection{Homogeneous Involutive Bases}

When computing an Involutive Basis,
a prolongation of a homogeneous polynomial
is another homogeneous polynomial, and the involutive
reduction of a homogeneous polynomial by a set of
homogeneous polynomials yields another homogeneous
polynomial. It would therefore be entirely feasible
for a program computing Involutive Bases for
homogeneous input bases to take advantage of the
properties of homogeneous polynomial arithmetic.

It would also be desirable to be able to use such a
program on input bases containing non-homogeneous
polynomials. The natural way to do this would be to
modify the procedure outlined in Definition \ref{homprocC}
by replacing every occurrence of the phrase
``a Gr\"obner Basis'' by the phrase ``an Involutive Basis'',
thus creating the following definition.

\begin{defn} \label{homprocIC}
Let $F = \{f_1, \hdots, f_m\}$ be a non-homogeneous
set of polynomials. To compute an Involutive Basis for
$F$ using a program that only accepts sets of
homogeneous polynomials as input, we proceed as follows.
\begin{enumerate}[(a)]
\item
Construct a homogeneous set of polynomials
$F' = \{h(f_1), \hdots, h(f_m)\}$.
\index{homogenisation}
\item
Compute an Involutive Basis $G'$ for $F'$.
\item
Dehomogenise each polynomial $g' \in G'$ to obtain
a set of polynomials $G$.
\index{dehomogenisation}
\end{enumerate}
\end{defn}
Ideally, we would like to say that $G$ is always an
Involutive Basis for $F$ as long as the monomial
ordering used is extendible, mirroring the
conclusion reached in Definition \ref{homprocC}.
However, we will only prove the validity of this statement
in the case that the set $G$ is autoreduced, and also only
for certain combinations of monomial orderings and
involutive divisions --- all combinations
will not work, as the following example
demonstrates.
% However, for $G$ to be an Involutive Basis for $F$,
% not only must the monomial ordering used be extendible
% (as in Definition \ref{homprocC}),
% but the involutive division chosen must also satisfy some
% extra condition. This is because $G$ will not always be an
% Involutive Basis for $F$ with respect to some
% involutive division $I$ (assuming of course
% that $G'$ is an Involutive Basis for $F'$ with
% respect to $I$), as the following example demonstrates.

\begin{example}
Let $F := \{x_1^2+x_2^3, \, x_1+x_3^3\}$ be a basis
generating an ideal $J$ over the polynomial ring
$\mathbb{Q}[x_1, x_2, x_3]$, and let the monomial
ordering be Lex. Computing an Involutive Basis
for $F$ with respect to the Janet involutive
division using Algorithm \ref{com-inv},
we obtain the set $G := \{x_2^3+x_3^6, \, x_1x_2^2+x_2^2x_3^3, \,
x_1x_2+x_2x_3^3, \, x_1^2-x_3^6, \, x_1+x_3^3\}$.

Taking the homogeneous route,
we can homogenise $F$ (with respect to Lex) to obtain the set
$F' := \{x_1^2y+x_2^3, \, x_1y^2+x_3^3\}$ over the polynomial
ring $\mathbb{Q}[x_1, x_2, x_3, y]$.
Computing an Involutive Basis for $F'$
with respect to the Janet involutive division,
we obtain the set
$G' := \{x_2^3y^3+x_3^6, \, x_1x_2^2y^3+x_2^2x_3^3y, \,
x_1x_2y^3+x_2x_3^3y, \, x_1y^3+x_3^3y, x_1y^2+x_3^3, x_1x_3^3y-x_2^3y^2, \,
x_1^2x_3^2y+x_2^3x_3^2, \, x_1^2x_3y+x_2^3x_3,
x_1^2y+x_2^3, \, x_1x_3^3-x_2^3y\}$.
Finally, if we dehomogenise $G'$, we obtain the set
$H := \{x_2^3+x_3^6, \, x_1x_2^2+x_2^2x_3^3, \, x_1x_2+x_2x_3^3, \,
x_1+x_3^3, \, x_1x_3^3-x_2^3, \, x_1^2x_3^2+x_2^3x_3^2, \,
x_1^2x_3+x_2^3x_3, \, x_1^2+x_2^3\}$;
however this set is not a Janet Involutive Basis
for $F$, as can be verified by checking that
(with respect to $H$) the variable $x_3$ is
nonmultiplicative for the polynomial $x_2^3+x_3^6$, and
the prolongation of the polynomial $x_2^3+x_3^6$ by the
variable $x_3$ is involutively irreducible with
respect to $H$.
\end{example}

The reason why $H$ is not an Involutive Basis
for $J$ in the above example is that the Janet
multiplicative variables for the set $G'$ do not
correspond to the Janet multiplicative variables
for the set $H = d(G')$. This means that we
cannot use the fact that all prolongations
of elements of $G'$ involutively reduce to zero
using $G'$ to deduce that all prolongations of
elements of $H$ involutively reduce to zero
using $H$. To do this, our involutive division
must satisfy the following additional property, which 
ensures that the multiplicative variables of $G'$ and 
$d(G')$ do correspond to each other.

\begin{defn}
Let $O$ be a fixed extendible monomial ordering.
An involutive division $I$ is {\it extendible
with respect to} $O$
\index{involutive division!extendible} if,
\index{extendible involutive division}
given any set of polynomials $P$, we have
$$\mathcal{M}_I(p, P) \setminus \{y\} =
\mathcal{M}_I(d(p), d(P))$$
for all $p \in P$, where $y$ is the homogenising variable.
\end{defn}

In Section \ref{HomSecC}, we saw that of the
monomial orderings defined in Section \ref{CMO}, only
Lex, InvLex and DegRevLex are extendible.
Let us now consider which involutive divisions
are extendible with respect to these three
monomial orderings.

\begin{prop}
The Thomas involutive division is extendible
with respect to Lex, \linebreak 
InvLex and DegRevLex.
\end{prop}
\begin{pf}
Let $P$ be an arbitrary set of polynomials over
a polynomial ring containing variables
$x_1, \hdots, x_n$ and a homogenising variable $y$.
Because the Thomas involutive division decides whether a
variable $x_i$ (for $1 \leqslant i \leqslant n$)
is multiplicative for a polynomial
$p \in P$ independent of the variable $y$, it is
clear that $x_i$ is multiplicative for $p$ if and
only if $x_i$ is multiplicative for $d(p)$ with respect to any
of the monomial orderings Lex, InvLex and DegRevLex. It follows
that $\mathcal{M}_{\mathcal{T}}(p, P)
\setminus \{y\} = \mathcal{M}_{\mathcal{T}}(d(p), d(P))$
as required.
\end{pf}

\begin{prop}
The Pommaret involutive division is extendible
with respect to Lex and DegRevLex.
\end{prop}
\begin{pf}
Let $p$ be an arbitrary polynomial over
a polynomial ring containing variables
$x_1, \hdots, x_n$ and a homogenising variable $y$.
Because we are using either the Lex or the DegRevLex
monomial orderings, the variable $y$ must be
lexicographically less than any of the variables
$x_1, \hdots, x_n$, and so we can state (without
loss of generality) that $p$ belongs to the polynomial ring
$R[x_1, \hdots, x_n, y]$. Let
$(e^1, e^2, \hdots, e^n, e^{n+1})$ be
the multidegree corresponding to the monomial
$\LM(p)$, and let $1 \leqslant i \leqslant n+1$ be the
smallest integer such that $e^i > 0$.

If $i = n+1$, then the variables $x_1, \hdots, x_n$
will all be multiplicative for $p$. But then $d(p)$ will be
a constant, so that the variables $x_1, \hdots, x_n$ will also
all be multiplicative for $d(p)$.

If $i \leqslant n$, then the variables
$x_1, \hdots, x_i$ will all be multiplicative for $p$.
But because $y$ is the smallest variable, it is clear
that $i$ will also be the smallest integer such
that $f^i > 0$, where $(f^1, f^2, \hdots, f^n)$ is
the multidegree corresponding to the monomial
$\LM(d(p))$. It follows that the variables
$x_1, \hdots, x_i$ will also all be multiplicative
for $d(p)$, and so we can conclude that
$\mathcal{M}_{\mathcal{P}}(p, P)
\setminus \{y\} = \mathcal{M}_{\mathcal{P}}(d(p), d(P))$
as required.
\end{pf}

\begin{prop}
The Pommaret involutive division is not
extendible with respect to InvLex.
\end{prop}
\begin{pf}
Let $p := yx_2+x_1^2$ be a polynomial
over the polynomial ring $\mathbb{Q}[y, x_1, x_2]$,
where $y$ is the homogenising variable (which
must be greater than all other variables in
order for InvLex to be extendible).
As $\LM(p) = yx_2$ with respect to InvLex, it follows
that $\mathcal{M}_{\mathcal{P}}(p) = \{y\}$.
Further, as $\LM(d(p)) = \LM(x_2+x_1^2) = x_2$ with
respect to InvLex, it follows that
$\mathcal{M}_{\mathcal{P}}(d(p)) = \{x_1, x_2\}$. We can
now deduce that the Pommaret involutive division is not
extendible with respect to InvLex, as
$\mathcal{M}_{\mathcal{P}}(p) \setminus \{y\}
\neq \mathcal{M}_{\mathcal{P}}(d(p))$, or
$\emptyset \neq \{x_1, x_2\}$.
\end{pf}

\begin{prop}
The Janet involutive division is extendible with
respect to InvLex.
\end{prop}
\begin{pf}
Let $P$ be an arbitrary set of polynomials over
a polynomial ring containing variables
$x_1, \hdots, x_n$ and a homogenising variable $y$.
Because we are using the InvLex monomial ordering,
the variable $y$ must be lexicographically greater
than any of the variables $x_1, \hdots, x_n$, and so we
can state (without loss of generality) that $p$ belongs to
the polynomial ring $R[y, x_1, \hdots, x_n]$.
But the Janet involutive division will then decide whether a
variable $x_i$ (for $1 \leqslant i \leqslant n$)
is multiplicative for a polynomial
$p \in P$ independent of the variable $y$, so it is
clear that $x_i$ is multiplicative for $p$ if and
only if $x_i$ is multiplicative for $d(p)$, and so
(with respect to InvLex) $\mathcal{M}_{\mathcal{J}}(p, P)
\setminus \{y\} = \mathcal{M}_{\mathcal{J}}(d(p), d(P))$
as required. ${}$
\end{pf}

\begin{prop}
The Janet involutive division is not extendible with respect
to Lex or DegRevLex.
\end{prop}
\begin{pf}
Let $U := \{x_1^2y, \, x_1y^2\}$ be a set of monomials
over the polynomial ring $\mathbb{Q}[x_1, y]$, where
$y$ is the homogenising variable (which must be less than
$x_1$ in order for Lex and DegRevLex to be extendible).
The Janet multiplicative variables for $U$ (with respect
to Lex and DegRevLex) are shown in the table below.
\begin{center}
\begin{tabular}{c|c}
$u$ & $\mathcal{M}_{\mathcal{J}}(u, U)$ \\ \hline
$x_1^2y$ & $\{x_1\}$ \\
$x_1y^2$ & $\{x_1, y\}$ \\ \hline
\end{tabular}
\end{center}
When we dehomogenise $U$ with respect to $y$, we
obtain the set $d(U) := \{x_1^2, \, x_1\}$ with
multiplicative variables as follows.
\begin{center}
\begin{tabular}{c|c}
$d(u)$ & $\mathcal{M}_{\mathcal{J}}(d(u), d(U))$ \\ \hline
$x_1^2$ & $\{x_1\}$ \\
$x_1$   & $\emptyset$ \\ \hline
\end{tabular}
\end{center}
It is now clear that Janet is not an extendible
involutive division with respect to Lex or DegRevLex,
as $\mathcal{M}_{\mathcal{J}}(x_1y^2, U)
\setminus \{y\} \neq \mathcal{M}_{\mathcal{J}}(x_1, d(U))$,
or $\{x_1\} \neq \emptyset$.
\end{pf}

\begin{prop}
Let $G'$ be a set of polynomials over
a polynomial ring containing variables
$x_1, \hdots, x_n$ and a homogenising variable $y$.
If (i) $G'$ is an Involutive Basis with respect to some
extendible monomial ordering $O$ and some involutive
division $I$ that is extendible with respect to $O$;
and (ii) $d(G')$ is an autoreduced set, then $d(G')$
is an Involutive Basis with respect to $O$ and $I$.
\end{prop}
\begin{pf}
By Definition \ref{IBC}, we can show that $d(G')$ is
an Involutive Basis with respect to $O$ and $I$ by showing
that any multiple $d(g')t$ of any polynomial
$d(g') \in d(G')$ by any term $t$ involutively reduces
to zero using $d(G')$. Because $G'$ is an Involutive Basis
with respect to $O$ and $I$, the polynomial $g't$
involutively reduces to zero using $G'$ by the series
of involutive reductions
$$g't \xymatrix{\ar[r]_{I}_(1){g'_{\alpha_1}} &}
h_1 \xymatrix{\ar[r]_{I}_(1){g'_{\alpha_2}} &}
h_2 \xymatrix{\ar[r]_{I}_(1){g'_{\alpha_3}} &}
\hdots \xymatrix{\ar[r]_{I}_(1){g'_{\alpha_A}} &} 0,$$
where $g'_{\alpha_i} \in G'$ for all
$1 \leqslant i \leqslant A$.

{\bf Claim:} The polynomial $d(g')t$ involutively
reduces to zero using $d(G')$ % (and hence $d(G')$ is an Involutive
% Basis with respect to $I$ and $O$) 
by the series of involutive reductions
$$d(g')t \xymatrix{\ar[r]_(0.4){I}_(1){d(g'_{\alpha_1})} &}
d(h_1) \xymatrix{\ar[r]_(0.4){I}_(1){d(g'_{\alpha_2})} &}
d(h_2) \xymatrix{\ar[r]_(0.4){I}_(1){d(g'_{\alpha_3})} &}
\hdots \xymatrix{\ar[r]_(0.4){I}_(1){d(g'_{\alpha_A})} &} 0,$$
where $d(g'_{\alpha_i}) \in d(G')$ for all
$1 \leqslant i \leqslant A$.

{\bf Proof of Claim:} It is clear that if a polynomial
$g'_j \in G'$ involutively reduces a polynomial $h$, then
the polynomial $d(g'_j) \in d(G')$ will always conventionally
reduce the polynomial $d(h)$. Further, knowing that $I$ is
extendible with respect to $O$, we can state that
$d(g'_j)$ will also always involutively reduce
$d(h)$. The result now follows by
noticing that $d(G')$ is autoreduced, so that $d(g'_j)$ is
the only possible involutive divisor of $d(h)$, and hence
the above series of involutive reductions is the only possible
way of involutively reducing the polynomial $d(g')t$.
\end{pf}

\begin{openquestion} \label{oq1}
If the set $G$ returned by the procedure outlined
in Definition \ref{homprocIC} is not autoreduced,
under what circumstances does autoreducing 
$G$ result in obtaining a set that is an
Involutive Basis for the ideal generated by $F$?
\end{openquestion}

Let us now consider two examples illustrating that
the set $G$ returned by the procedure outlined
in Definition \ref{homprocIC} may or may not be
autoreduced.

\begin{example}
Let $F := \{2x_1x_2+x_1^2+5, \, x_2^2+x_1+8\}$ be a basis
generating an ideal $J$ over the polynomial ring
$\mathbb{Q}[x_1, x_2]$, and let the monomial
ordering be InvLex. Ordinarily, we can compute
an Involutive Basis $G := \{x_2^2 + x_1 + 8, \,
2x_1x_2 + x_1^2 + 5, \, 10x_2 - x_1^3 - 4x_1^2 - 37x_1, \,
x_1^4 + 4x_1^3 + 42x_1^2 +25\}$
for $F$ with respect to the
Janet involutive division by using Algorithm \ref{com-inv}.
%a basis which we know to be a Gr\"obner Basis for $J$.

Taking the homogeneous route (using Definition
\ref{homprocIC}), we can homogenise
$F$ to obtain the basis $F' := \{2x_1x_2+x_1^2+5y^2, \,
x_2^2+yx_1+8y^2\}$ over the polynomial ring
$\mathbb{Q}[y, x_1, x_2]$, where $y$ is
the homogenising variable (which must be greater than all
other variables). Computing an
Involutive Basis for the set $F'$ with respect
to the Janet involutive division using Algorithm \ref{com-inv},
we obtain the basis $G' := \{x_2^2 + yx_1 + 8y^2, \,
2x_1x_2 + x_1^2 + 5y^2, \, 10y^2x_2 - x_1^3 - 4yx_1^2 - 37y^2x_1, \,
x_1^4 + 4yx_1^3 + 42y^2x_1^2 + 25y^4\}$. When we dehomogenise
this basis, we obtain the set $d(G') := \{x_2^2 + x_1 + 8, \,
2x_1x_2 + x_1^2 + 5, \, 10x_2 - x_1^3 - 4x_1^2 - 37x_1, \,
x_1^4 + 4x_1^3 + 42x_1^2 + 25\}$. It is now clear that the set 
$d(G')$ is autoreduced (and hence $d(G')$ is an Involutive
Basis for $J$) because $d(G') = G$.
\end{example}

\begin{example}
Let $F := \{x_2^2 + 2x_1x_2+5, \, x_2+x_1^2+8\}$ be a basis
generating an ideal $J$ over the polynomial ring
$\mathbb{Q}[x_1, x_2]$, and let the monomial
ordering be InvLex. Ordinarily, we can compute
an Involutive Basis $G := \{x_2^2-2x_1^3-16x_1+5, \,
x_2+x_1^2+8, \, x_1^4-2x_1^3+16x_1^2-16x_1+69\}$
for $F$ with respect to the
Janet involutive division by using Algorithm \ref{com-inv}.
%a basis which we know to be a Gr\"obner Basis for $J$.

Taking the homogeneous route (using Definition
\ref{homprocIC}), we can homogenise
$F$ to obtain the basis $F' := \{x_2^2+2x_1x_2+5y^2, \,
yx_2+x_1^2+8y^2\}$ over the polynomial ring
$\mathbb{Q}[y, x_1, x_2]$, where $y$ is
the homogenising variable (which must be greater than all
other variables). Computing an
Involutive Basis for the set $F'$ with respect
to the Janet involutive division using Algorithm \ref{com-inv},
we obtain the basis $G' := \{x_2^2+2x_1x_2+5y^2, \,
x_1^2x_2 + 2x_1^3 - 8yx_1^2 + 16y^2x_1 - 69y^3, \,
yx_1x_2 + x_1^3 + 8y^2x_1, \, yx_2+x_1^2+8y^2, \,
x_1^4-2yx_1^3+16y^2x_1^2-16y^3x_1+69y^4\}$. When we dehomogenise
this basis, we obtain the set $d(G') := \{x_2^2+2x_1x_2+5, \,
x_1^2x_2 + 2x_1^3 - 8x_1^2 + 16x_1 - 69, \,
x_1x_2 + x_1^3 + 8x_1, \, x_2+x_1^2+8, \,
x_1^4-2x_1^3+16x_1^2-16x_1+69\}$. This time however,
because the set $d(G')$ is not autoreduced 
(the polynomial $x_1x_2+x_1^3+8x_1 \in d(G')$
can involutively reduce the second term of the polynomial
$x_2^2+2x_1x_2+5 \in d(G')$), we cannot deduce that
$d(G')$ is an Involutive Basis for $J$.
\end{example}

\begin{remark}
Although the set $G$ returned by the procedure outlined
in Definition \ref{homprocIC} may not always be an Involutive
Basis for the ideal generated by $F$, because the set
$G'$ will always be an Involutive Basis (and hence also
a Gr\"obner Basis), we can state that $G$ will always be
a Gr\"obner Basis for the ideal generated by $F$
(cf. Definition \ref{homprocC}).
\end{remark}

% \begin{remark}
% In the above example, when we dehomogenise the
% set $G'$, we obtain the Involutive Basis $G$
% for $F$ that was computed directly. In general
% however, the set obtained by dehomogenising the
% set $G'$ will not be an Involutive Basis for $F$,
% as the following example demonstrates.
% \end{remark}
% Given a basis $F'$ which is obtained by homogenising a basis
% $F$ with respect to some extendible monomial ordering, it
% would be interesting to investigate whether, if we compute
% an Involutive Basis $G'$ for $F'$ with
% respect to some involutive division $I$, the basis $G$
% obtained by dehomogenising each element of $G'$ is always
% an Involutive Basis for $F$ with respect to $I$ (as well
% as being a Gr\"obner Basis for $F$) --- perhaps
% $I$ would have to satisfy some extra property
% (extendibility?) for this to be the case?

\subsection{Logged Involutive Bases}

Just as a Logged Gr\"obner Basis expresses each member of
the Gr\"obner Basis in terms of members of the original
basis from which the Gr\"obner Basis was computed, a
Logged Involutive Basis expresses each member of
the Involutive Basis in terms of members of the original
basis from which the Involutive Basis was computed.

\begin{defn}
Let $G = \{g_1, \hdots, g_p\}$ be
an Involutive Basis computed from
an initial basis $F = \{f_1, \hdots, f_m\}$. We say that $G$ is a
{\it Logged Involutive Basis}
\index{involutive basis!logged}
if, \index{logged involutive basis}
for each $g_i \in G$, we have an explicit expression of the form
$$g_i = \sum_{\alpha=1}^{\beta} t_{\alpha}f_{k_{\alpha}},$$
where the $t_{\alpha}$ are terms and
$f_{k_{\alpha}} \in F$ for all $1 \leqslant \alpha \leqslant \beta$.
\end{defn}

\begin{prop} \label{LogIBC}
Given a finite basis $F = \{f_1, \hdots, f_m\}$,
it is always possible to compute a Logged Involutive Basis for $F$.
\end{prop}
\begin{pf}
Let $G = \{g_1, \hdots, g_p\}$ be an Involutive Basis
computed from the initial basis $F = \{f_1, \hdots, f_m\}$
using Algorithm \ref{com-inv}
(where $f_i \in R[x_1, \hdots, x_n]$ for all $f_i \in F$).
If an arbitrary $g_i \in G$
is not a member of the original basis $F$, then
either $g_i$ is an involutively reduced prolongation,
or $g_i$ is obtained through the process of autoreduction.
In the former case, we can express $g_i$ in terms of
members of $F$ by substitution because
$$g_i = hx_j - \sum_{\alpha=1}^{\beta} t_{\alpha}h_{k_{\alpha}}$$
for a variable $x_j$; terms $t_{\alpha}$ and
polynomials $h$ and $h_{k_{\alpha}}$ which we already know
how to express in terms of members of $F$.
In the latter case,
$$g_i = h - \sum_{\alpha=1}^{\beta} t_{\alpha}h_{k_{\alpha}}$$
for terms $t_{\alpha}$ and polynomials
$h$ and $h_{k_{\alpha}}$ which we already know how
to express in terms of members of $F$, so it follows that we can again
express $g_i$ in terms of members of $F$.
\end{pf}

%
% Chapter 5
% Author: Gareth Evans
% Last Modified: 2nd February 2006
%

\chapter{Noncommutative Involutive Bases} \label{ChNCIB}

In the previous chapter, we introduced the theory of
commutative Involutive Bases and saw that such bases
are always commutative Gr\"obner Bases with extra
structure. In this chapter, we will follow a similar
path, in that we will define an algorithm to compute a
{\it noncommutative Involutive Basis} that will serve
as an alternative method of obtaining a noncommutative
Gr\"obner Basis, and the noncommutative Gr\"obner Bases
we will obtain will also have some extra structure.

As illustrated by the diagram shown below, the theory
of noncommutative Involutive Bases will draw upon all the
theory that has come before in this thesis, and as a
consequence will inherit
many of the restrictions imposed by this theory.
For example, our noncommutative Involutive Basis algorithm
will not be guaranteed to terminate precisely
because we are working in a noncommutative
setting, and noncommutative involutive divisions will have
properties that will influence the correctness and
termination of the algorithm.

$$\xymatrix{
\mbox{Commutative Gr\"obner Bases} \ar[r] \ar[dd]
& \mbox{Commutative Involutive Bases} \ar[dd] \\ \\
\mbox{Noncommutative Gr\"obner Bases} \ar[r]
& \mbox{Noncommutative Involutive Bases}
}$$

\section{Noncommutative Involutive Reduction} \label{5point1}

Recall that in a commutative polynomial ring, a monomial $u_2$ is an
involutive divisor of a monomial $u_1$ if $u_1 = u_2u_3$ for some
monomial $u_3$ % ($u_2$ is a conventional divisor of $u_1$)
and all variables in $u_3$ are multiplicative for $u_2$. %, where the
%multiplicative variables for $u_2$ are determined by a set of
%monomials $U$ such that $u_2 \in U$.
In other words, we are able to form $u_1$ from $u_2$ by multiplying
$u_2$ with multiplicative variables.

In a noncommutative polynomial ring, an involutive division will again
induce a restricted form of division. However, because left and right
%restrict conventional division. However, because left and right
multiplication are separate processes in noncommutative polynomial rings,
we will require the notion of {\it left} and {\it right multiplicative}
variables in order to determine whether a conventional divisor is an
involutive divisor, so that (intuitively) a monomial
$u_2$ will involutively divide a monomial $u_1$ if we are able to form
$u_1$ from $u_2$ by multiplying $u_2$ on the left with left multiplicative
variables and on the right by right multiplicative variables.

More formally, let $u_1$ and $u_2$ be two monomials over a noncommutative
polynomial ring, and assume that $u_1$ is a conventional divisor of $u_2$,
\index{conventional divisor}
so \index{divisor!conventional}
that $u_1 = u_3u_2u_4$ for some monomials $u_3$ and $u_4$.
Assume that an arbitrary noncommutative involutive division $I$ partitions
the variables in the polynomial ring into sets of {\it left multiplicative}
and {\it left nonmultiplicative} variables for $u_2$, and also partitions
the variables in the polynomial ring into sets of {\it right multiplicative}
and {\it right nonmultiplicative} variables for $u_2$.
Let us now define two methods of deciding
whether $u_2$ is an {\it involutive} \index{involutive divisor} divisor 
\index{divisor!involutive} of
$u_1$ (written $u_2 \mid_I u_1$), the first of which will depend only on the
{\it first} variable we multiply $u_2$ 
with on the left and on the right in order
to form $u_1$, and the second of which will 
depend on {\it all} the variables
we multiply $u_2$ with in order to form $u_1$.

\begin{defn} \label{NCID}
Let $u_1 = u_3u_2u_4$, and let $I$ be
defined as in the previous paragraph.
\begin{itemize}
\item
{\bf (Thin Divisor)}
\index{thin divisor}
$u_2 \mid_I u_1$ if \index{divisor!thin}
the variable $\SUFF(u_3, 1)$ (if it exists)
is in the set of left multiplicative variables for $u_2$,
and the variable $\PRE(u_4, 1)$ (again if it exists) is in the
set of right multiplicative variables for $u_2$.
% In other words, $u_2$ is an involutive divisor of $u_1$ if the
% two variables `closest to $u_2$' are (respectively left or right)
% multiplicative for $u_2$.
\item
{\bf (Thick Divisor)}
\index{thick divisor}
$u_2 \mid_I u_1$ if \index{divisor!thick}
all the variables in $u_3$ are in the set of left
multiplicative variables for $u_2$, and all
the variables in $u_4$ are in the
set of right multiplicative variables for $u_2$.
\end{itemize}
\end{defn}

\begin{remark}
We introduce two methods for determining whether a conventional
divisor is an involutive divisor because each of the
methods has its own advantages and disadvantages.
From a theoretical standpoint, using thin divisors
enables us to follow
the path laid down in Chapter \ref{ChCIB}, in that we are able
to show that a Locally Involutive Basis is an Involutive Basis
by proving that the involutive division used is continuous,
something that we cannot do if thick divisors are being used.
On the other hand, once we have obtained our Locally Involutive Basis,
involutive reduction with respect to
thick divisors is more efficient
than it is with respect to thin divisors, as less work is required
in order to determine whether a monomial is involutively
divisible by a set of monomials. For these reasons, we will use
thin divisors when presenting the theory in this chapter
(hence the following definition), and will only use thick divisors
when, by doing so, we are able to gain some advantage.
\end{remark}

% It is tempting to make an obvious generalisation of the theory
% from the commutative to the noncommutative case by requiring {\it all}
% variables in $u_3$ to be in the set of left multiplicative
% variables for $u_2$ and requiring {\it all} variables in $u_4$ to be in
% the set of right multiplicative variables for $u_2$.
% Indeed, this exact definition appears in \cite{Evans04}.
% There are two reasons why this definition is not used in this thesis:
% first, the definition is more restrictive than the one that
% appears in Definition \ref{NCID}, in that less conventional
% divisors will be involutive divisors; second (and crucially),
% the definition does not allow us to
% define a property of an involutive division that leads to the
% statement that a Locally Involutive Basis is an Involutive Basis,
% creating the need to find an alternative method of showing that
% a Locally Involutive Basis is a Gr\"obner Basis
% (in \cite{Evans04}, the concept of a `Gr\"obner' involutive
% division is used).

\begin{remark}
Unless otherwise stated, from now on we will use thin divisors
to determine whether a conventional divisor is an
involutive divisor.
\end{remark}

\begin{example}
Let $u_1 := xyz^2x$; $u'_1 := yz^2y$ and $u_2 := z^2$ be three monomials
over the polynomial ring
$\mathcal{R} = \mathbb{Q}\langle x, y, z\rangle$, and let
an involutive division $I$ partition the variables in $\mathcal{R}$ into
the following sets of variables for the monomial $u_2$:
left multiplicative = $\{x, y\}$; left nonmultiplicative = $\{z\}$;
right multiplicative = $\{x, z\}$; right nonmultiplicative = $\{y\}$.
It is true that $u_2$ conventionally divides both monomials $u_1$ and $u'_1$,
but $u_2$ only involutively divides monomial $u_1$ as, defining $u_3 := xy$;
$u_4 := x$; $u'_3 = y$ and $u'_4 = y$
(so that $u_1 = u_3u_2u_4$ and $u'_1 = u'_3u_2u'_4$),
we observe that the variable $\SUFF(u_3, 1) = y$ is in the set of left
multiplicative variables for $u_2$;
the variable $\PRE(u_4, 1) = x$ is in the
set of right multiplicative variables for $u_2$; but the variable
$\PRE(u'_4, 1) = y$ is not in the set of right multiplicative variables
for $u_2$.
% all variables in $u_3$ are in the set of left multiplicative variables
% for $u_2$, all variables in $u_4$
% are in the set of right multiplicative variables
% for $u_2$, but all variables in $u'_4$
% (in particular the variable $y$) are not
% all in the set of right multiplicative variables for $u_2$.
\end{example}

Let us now formally define what is meant by a (noncommutative)
involutive division.

\begin{defn} \label{noncom-cones}
Let $M$ denote the set of all monomials in a
noncommutative polynomial ring
$\mathcal{R} = R\langle x_1, \hdots, x_n\rangle$,
and let $U \subset M$. The involutive cone $\mathcal{C}_I(u, U)$
\index{$C$@$\mathcal{C}_I(u, U)$}
of \index{involutive cone}
any \index{cone!involutive}
monomial $u \in U$ with respect to some involutive
division $I$ is defined as follows.
%\begin{eqnarray*}
%\mathcal{C}^L_I(u, U) & = & \{v \in M \: \mbox{such that} \:
%u \mid_I vu\}; \\
%\mathcal{C}^R_I(u, U) & = & \{v \in M \: \mbox{such that} \:
%u \mid_I uv\}; \\
%\mathcal{C}_I(u, U) & = & \{v_1uv_2 \: \mbox{such that} \:
%v_1 \in \mathcal{C}^L_I(u, U), v_2 \in \mathcal{C}^R_I(u, U)\}.
%\end{eqnarray*}
$$
\mathcal{C}_I(u, U)  =  \{v_1uv_2 \: \mbox{such that} \:
v_1, v_2 \in M \: \mbox{and} \: u \mid_I v_1uv_2\}.
$$
\end{defn}

\begin{defn} \label{noncom-div-defn}
Let $M$ denote the set of all monomials in a
noncommutative polynomial ring
$\mathcal{R} = R\langle x_1, \hdots, x_n\rangle$.
\index{involutive division!strong}
A {\it strong} \index{strong involutive division}
involutive division \index{involutive division} $I$ 
\index{division!involutive} is
defined on $M$ if, given any finite
set of monomials $U \subset M$,
we can assign a set of left multiplicative
\index{$Mb$@$\mathcal{M}^L_I(u, U)$}
variables $\mathcal{M}^L_I(u, U) \subseteq \{x_1, \hdots, x_n\}$
and a set of right multiplicative variables
\index{$Mb$@$\mathcal{M}^R_I(u, U)$}
$\mathcal{M}^R_I(u, U) \subseteq \{x_1, \hdots, x_n\}$
to any monomial $u \in U$ such that the
following three conditions are satisfied.
% Further, if $\mathcal{C}_I(u, U)$
% denotes the \index{involutive cone}
% {\it involutive cone} of the monomial
% $u \in U$ (defined as
% \index{$C$@$\mathcal{C}_I(u, U)$}
% $$
% \mathcal{C}_I(u, U) = \{v_1uv_2 \mid
% v_1 \in \mathcal{C}^L_I(u, U), v_2 \in \mathcal{C}^R_I(u, U)\},
% $$
% % where $\mathcal{C}^L_I(u, U) = \{v \in M \mid$
% % all variables in $v$ are left multiplicative for $u\}$ and
% % $\mathcal{C}^R_I(u, U) = \{v \in M \mid$
% % all variables in $v$ are right multiplicative for $u\}$.
% where $\mathcal{C}^L_I(u, U)$ is the set of all monomials
% $v \in M$ such that $u \mid_I vu$, and
% $\mathcal{C}^R_I(u, U)$ is the set of all monomials
% $v \in M$ such that $u \mid_I uv$),
% then $I$ is said to be
%
% {\it strong} if
% the following three conditions are all satisfied;
% and \index{involutive division!weak}
% {\it weak} if any of the conditions are not
% satisfied.
\begin{itemize}
\item
If there exist two elements $u_1, u_2 \in U$
such that $\mathcal{C}_I(u_1, U) \cap
\mathcal{C}_I(u_2, U) \neq \emptyset$, then
either $\mathcal{C}_I(u_1, U) \subset
\mathcal{C}_I(u_2, U)$ or
$\mathcal{C}_I(u_2, U) \subset
\mathcal{C}_I(u_1, U)$.
\item
Any monomial $v \in \mathcal{C}_I(u, U)$ is
involutively divisible by $u$ in one way only, so that if
$u$ appears as a subword of $v$ in more than one
way, then only one of these ways allows us to
deduce that $u$ is an involutive divisor of $v$.
\item
If $V \subset U$, then $\mathcal{M}^L_I(v, U) \subseteq
\mathcal{M}^L_I(v, V)$ and $\mathcal{M}^R_I(v, U) \subseteq
\mathcal{M}^R_I(v, V)$ for all $v \in V$.
\end{itemize}
If any of the above conditions are not satisfied, the
involutive division is called a {\it weak}
\index{involutive division!weak}
involutive \index{weak involutive division} division.
\end{defn}

\begin{remark}
We shall refer to the three conditions of
Definition \ref{noncom-div-defn} as (respectively)
the Disjoint Cones condition, the Unique Divisor
condition and the Subset condition.
\end{remark}

\begin{defn}
Given an involutive division $I$, the involutive 
\index{span!involutive} span
\index{involutive span}
$\mathcal{C}_I(U)$ \index{$C$@$\mathcal{C}_I(U)$}
of a set of noncommutative monomials $U$ with respect to $I$
is given by the expression
$$
\mathcal{C}_I(U) = \bigcup_{u \in U} \mathcal{C}_I(u, U).
$$
\end{defn}

\begin{remark}
The (conventional) \index{span!conventional}
span \index{conventional span}
of a set of noncommutative monomials
$U$ is given by the expression
\index{$C$@$\mathcal{C}(U)$}
$$
\mathcal{C}(U) = \bigcup_{u \in U} \mathcal{C}(u, U),
$$
where \index{$C$@$\mathcal{C}(u, U)$}
$\mathcal{C}(u, U) = \{v_1uv_2 \: \mbox{such that} \:
v_1, v_2 \; \mbox{are monomials}\}$
is the (conventional) \index{cone!conventional}
cone \index{conventional cone} of a monomial $u \in U$.
\end{remark}

\begin{defn}
If an involutive division $I$ determines the left and right
multiplicative variables for a monomial $u \in U$
independent of the set $U$, then $I$ is a {\it global} 
\index{global involutive division} division.
\index{involutive division!global} Otherwise, $I$
is a {\it local} \index{local involutive division} division.
\index{involutive division!local}
\end{defn}

\begin{remark}
The multiplicative variables for a set of polynomials
$P$ (whose terms are ordered by a monomial ordering $O$)
are determined by the multiplicative variables for the set
of leading monomials $\LM(P)$.
\end{remark}

In Algorithm \ref{noncom-inv-div}, we specify how to
involutively divide a polynomial $p$ with respect to a set
of polynomials $P$
using thin divisors. Note that this algorithm
combines the modifications made to Algorithm
\ref{com-div} in Algorithms
\ref{noncom-div} and \ref{com-inv-div}.

\begin{algorithm} 
\setlength{\baselineskip}{3.5ex}
\caption{The Noncommutative Involutive Division Algorithm}
\label{noncom-inv-div}
\begin{algorithmic}
\vspace*{2mm}
\REQUIRE{A nonzero polynomial $p$ and a set of nonzero polynomials
         $P = \{p_1, \hdots, p_m\}$ over a polynomial ring
         $R\langle x_1, \hdots x_n\rangle$;
         an admissible monomial ordering $\mathrm{O}$; 
         an involutive division $I$.}
\ENSURE{$\Rem_I(p, P) := r$, the involutive remainder of $p$ 
        with respect to $P$.}
\vspace*{1mm}
\STATE
$r = 0$;
\WHILE{($p \neq 0$)}
\STATE
$u = \LM(p)$; $c = \LC(p)$; $j = 1$; found = false; \\
\WHILE{($j \leqslant m$) \textbf{and} (found == false)}
\IF{($\LM(p_j) \mid_I u$)}
\STATE
found = true; \\
choose $u_{\ell}$ and $u_r$ such that
$u = u_{\ell}\LM(p_j)u_r$, the variable
$\SUFF(u_{\ell}, 1)$ (if it exists) is left multiplicative for $p_j$,
and the variable $\PRE(u_r, 1)$ (again if it exists) is right
multiplicative for $p_j$; \\
% (If there are several candidates for $u_{\ell}$
% (and therefore for $u_r$), by convention choose the one
% with the smallest degree); \\
% choose\footnotemark \, $u_{\ell}$ and $u_r$ 
% such that all variables in $u_{\ell}$
% are left multiplicative for $p_j$, all variables in $u_r$ are right
% multiplicative for $p_j$, and $u = u_{\ell}\LM(p_j)u_r$; \\
$p = p - (c\LC(p_j)^{-1})u_{\ell}p_ju_r$;
\ELSE
\STATE
$j = j+1$;
\ENDIF
\ENDWHILE
\IF{(found == false)}
\STATE
$r = r+\LT(p)$; $p = p-\LT(p)$;
\ENDIF
\ENDWHILE
\STATE
{\bf return} $r$;
\end{algorithmic}
\vspace*{1mm}
\end{algorithm}

\begin{remark}
Continuing the convention from Algorithm \ref{noncom-div},
we will always choose the $u_{\ell}$ with the smallest degree
in the line `choose $u_{\ell}$ and $u_r$ such that$\hdots$'
in Algorithm \ref{noncom-inv-div}.
\end{remark}

\begin{example}
Let $P := \{x^2-2y, \: xy-x, \: y^3+3\}$ be a set
of polynomials over the polynomial ring
$\mathbb{Q}\langle x, y\rangle$ ordered with respect to
the DegLex monomial ordering, and assume
that an involutive division $I$ assigns multiplicative
variables to $P$ as follows.
\begin{center}
\begin{tabular}{c|c|c}
$p$ & $\mathcal{M}_I^L(\LM(p), \LM(P))$ &
$\mathcal{M}_I^R(\LM(p), \LM(P))$ \\ \hline
$x^2 - 2y$ & $\{x, y\}$ & $\{x\}$ \\
$xy - x$   & $\{y\}$    & $\{x, y\}$ \\
$y^3 + 3$  & $\{x\}$    & $\emptyset$ \\ \hline
\end{tabular}
\end{center}
Here is a dry run for Algorithm \ref{noncom-inv-div} when we
involutively divide the polynomial $p := 2x^2y^3 + yxy$
with respect to $P$ to obtain the polynomial
$yx-12y$, where A; B; C and D refer to the tests
($p \neq 0$)?; ($(j \leqslant 3)$ and (found $==$ false))?;
($\LM(p_j) \mid_I u$)? and (found $==$ false)? respectively.
\begin{small}
\begin{center}
\begin{tabular}{c|c|c|c|c|c|c|c||c|c|c|c}
$p$           & $r$      & $u$         & $c$   & $j$ & found & $u_{\ell}$
& $u_r$ & A     & B     & C     & D     \\ \hline
$2x^2y^3+yxy$ & $0$      &             &       &     &       &
&       & true  &       &       &       \\
              &          & $x^2y^3$    & $2$   & $1$ & false &
              &       &       & true  & false &       \\
              &          &             &       & $2$ &       &
              &       &       & true  & false &       \\
              &          &             &       & $3$ &       &
              &       &       & true  & true  &       \\
$yxy-6x^2$    &          &             &       &     & true  & $x^2$
& $1$   &       & false &       & false \\
              &          &             &       &     &       &
              &       & true  &       &       &       \\
              &          & $yxy$       & $1$   & $1$ & false &
              &       &       & true  & false &       \\
              &          &             &       & $2$ &       &
              &       &       & true  & true  &       \\
$-6x^2+yx$    &          &             &       &     & true  & $y$
& $1$   &       & false &       & false \\
              &          &             &       &     &       &
              &       & true  &       &       &       \\
              &          & $x^2$       & $-6$  & $1$ & false &
              &       &       & true  & true  &       \\
$yx-12y$      &          &             &       &     & true  & $1$
& $1$   &       & false &       & false \\
              &          &             &       &     &       &
              &       & true  &       &       &       \\
              &          & $yx$        & $1$   & $1$ & false &
              &       &       & true  & false &       \\
              &          &             &       & $2$ &       &
              &       &       & true  & false &       \\
              &          &             &       & $3$ &       &
              &       &       & true  & false &       \\
              &          &             &       & $4$ &       &
              &       &       & false &       & true  \\
$-12y$        & $yx$     &             &       &     &       &
&       & true  &       &       &       \\
              &          & $y$         & $-12$ & $1$ & false &
              &       &       & true  & false &       \\
              &          &             &       & $2$ &       &
              &       &       & true  & false &       \\
              &          &             &       & $3$ &       &
              &       &       & true  & false &       \\
              &          &             &       & $4$ &       &
              &       &       & false &       & true  \\
$0$           & $yx-12y$ &             &       &     &       &
&       & false &       &       &       \\ \hline
\end{tabular}
\end{center}
\end{small}
\end{example}

\section{Prolongations and Autoreduction} \label{5point2}

Just as in the commutative case, we will compute
a (noncommutative) Locally Involutive Basis by using
{\it prolongations} and {\it autoreduction}, but
here we have to distinguish between
{\it left prolongations} and {\it right prolongations}.

\begin{defn}
\index{prolongation!left}
Given a set of polynomials $P$, a {\it left prolongation} of
\index{left prolongation}
a polynomial $p \in P$ is a product $x_ip$, where
$x_i \notin \mathcal{M}^L_I(\LM(p), \LM(P))$
with respect to some involutive division $I$;
\index{prolongation!right}
and a {\it right prolongation} of
\index{right prolongation}
a polynomial $p \in P$ is a product $px_i$, where
$x_i \notin \mathcal{M}^R_I(\LM(p), \LM(P))$
with respect to some involutive division $I$.
\end{defn}

\begin{defn} \label{noncom-inv-defn}
\index{autoreduction}
A set of polynomials $P$ is said to be {\it autoreduced} if
no polynomial $p \in P$ exists such that $p$
contains a term which is involutively divisible (with respect
to $P$) by some polynomial $p' \in P\setminus\{p\}$.
\end{defn}

\begin{algorithm}
\setlength{\baselineskip}{3.5ex}
\caption{The Noncommutative Autoreduction Algorithm}
\label{noncom-auto}
\begin{algorithmic}
\vspace*{2mm}
\REQUIRE{A set of polynomials $P = \{p_1, p_2, \hdots, p_{\alpha}\}$;
         an involutive division $I$.}
\ENSURE{An autoreduced set of polynomials 
        $Q = \{q_1, q_2, \hdots, q_{\beta}\}$.}
\vspace*{1mm}
\WHILE{($\exists \; p_i \in P$ such that $\Rem_I(p_i,
P, P\setminus\{p_i\}) \neq p_i$)}
\STATE
$p'_i = \Rem_I(p_i, P, P\setminus\{p_i\})$; \\
$P = P \setminus\{p_i\}$;
\IF{($p'_i \neq 0$)}
\STATE
$P = P \cup \{p'_i\}$;
\ENDIF
\ENDWHILE
\STATE
$Q = P$; \\
{\bf return} $Q$;
\end{algorithmic}
\vspace*{1mm}
\end{algorithm}

\begin{remark}
With respect to a strong involutive division, the
involutive cones of an autoreduced set of polynomials
are always disjoint.
\end{remark}

\begin{remark}
The notation $\Rem_I(p_i, P, P\setminus\{p_i\})$
used in Algorithm \ref{noncom-auto} has the same meaning
as in Definition \ref{com-inv-defn}.
% autoreduction is the process
% of involutively reducing each member of
% a set of polynomials by the rest of the set
% until all members become involutively irreducible
% (an algorithm for performing autoreduction
% is given in Algorithm \ref{com-auto}).
\end{remark}

\begin{prop} \label{plusNC}
Let $P$ be a set of polynomials over a 
noncommutative polynomial ring
$\mathcal{R} = R\langle x_1, \hdots, x_n\rangle$, and let
$f$ and $g$ be two polynomials also in $\mathcal{R}$.
If $P$ is autoreduced with respect to a strong
involutive division $I$, then
$\Rem_I(f, P) + \Rem_I(g, P) = \Rem_I(f+g, P)$.
\end{prop}
\begin{pf}
Let $f' := \Rem_I(f, P)$; $g' := \Rem_I(g, P)$ and
$h' := \Rem_I(h, P)$, where $h := f+g$. Then,
by the respective involutive reductions,
we have expressions
$$
f' = f - \sum_{a=1}^A u_ap_{{\alpha}_a}v_a;
$$
$$
g' = g - \sum_{b=1}^B u_bp_{{\beta}_b}v_b
$$
and
$$
h' = h - \sum_{c=1}^C u_cp_{{\gamma}_c}v_c,
$$
where $p_{{\alpha}_{a}}, \, p_{{\beta}_{b}}, \,
p_{{\gamma}_{c}} \in P$ and
$u_{a}, v_a, u_{b}, v_b, u_{c}, v_c$ are terms
such that each $p_{{\alpha}_{a}}$, $p_{{\beta}_{b}}$
and $p_{{\gamma}_{c}}$ involutively divides each
$u_ap_{{\alpha}_a}v_a$, $u_bp_{{\beta}_b}v_b$ and
$u_cp_{{\gamma}_c}v_c$ respectively.

Consider the polynomial $h'-f'-g'$. By the above
expressions, we can deduce\footnote{For
$1 \leqslant d \leqslant A$,
$u_dp_{{\delta}_d}v_d = u_ap_{{\alpha}_a}v_a$
($1 \leqslant a \leqslant A$); for $A+1 \leqslant d
\leqslant A+B$, $u_dp_{{\delta}_d}v_d = u_bp_{{\beta}_b}v_b$
($1 \leqslant b \leqslant B$); and for $A+B+1 \leqslant d
\leqslant A+B+C =: D$, $u_dp_{{\delta}_d}v_d = u_cp_{{\gamma}_c}v_c$
($1 \leqslant c \leqslant C$).}
that
$$h'-f'-g' = \sum_{a=1}^A u_ap_{{\alpha}_a}v_a
+ \sum_{b=1}^B u_bp_{{\beta}_b}v_b
- \sum_{c=1}^C u_cp_{{\gamma}_c}v_c
=: \sum_{d=1}^D u_dp_{{\delta}_d}v_d.$$
{\bf Claim:} $\Rem_I(h'-f'-g', P) = 0$.

{\bf Proof of Claim:}
Let $t$ denote the leading term of the
polynomial $\sum_{d=1}^{D} u_dp_{{\delta}_d}v_d$.
Then $\LM(t) = \LM(u_kp_{{\delta}_{k}}v_{k})$ for some
$1 \leqslant k \leqslant D$ since, if not,
there exists a monomial
$$\LM(u_{k'}p_{{\delta}_{k'}}v_{k'}) =
\LM(u_{k''}p_{{\delta}_{k''}}v_{k''}) =: w$$ for some
$1 \leqslant k', k'' \leqslant D$
(with $p_{{\delta}_{k'}} \neq p_{{\delta}_{k''}}$) 
such that $w$ is involutively
divisible by the two polynomials $p_{{\delta}_{k'}}$ 
and $p_{{\delta}_{k''}}$,
contradicting Definition \ref{noncom-div-defn} 
(recall that $I$ is strong and
$P$ is autoreduced, so that the involutive cones of $P$
are disjoint). It follows that we can use
$p_{\delta_k}$ to eliminate $t$ by involutively reducing
$h'-f'-g'$ as shown below.
\begin{equation} \label{GIE1N}
\sum_{d=1}^{D} u_dp_{{\delta}_d}v_d
\xymatrix{\ar[r]_I_(1){p_{{\delta}_k}} &}
\sum_{d=1}^{k-1} u_dp_{{\delta}_d}v_d
+ \sum_{d=k+1}^{D} u_dp_{{\delta}_d}v_d.
\end{equation}
By induction, % and by the admissibility of the chosen monomial ordering,
we can apply a chain of involutive reductions to the
right hand side of Equation (\ref{GIE1N})
to obtain a zero remainder, so that
$\Rem_I(h'-f'-g', P) = 0$.
\hfill ${}_{\Box}$

To complete the proof, we note that since
$f'$, $g'$ and $h'$ are all involutively irreducible, we
must have $\Rem_I(h'-f'-g', P) = h'-f'-g'$. It therefore
follows that $h'-f'-g' = 0$, or $h' = f'+g'$ as required.
\end{pf}

\begin{defn} \label{LIBNC}
\index{involutive basis!locally}
Given an involutive division $I$ and an admissible monomial
ordering $O$, an autoreduced set of noncommutative polynomials $P$
is a {\it Locally Involutive Basis} with
\index{locally involutive basis}
respect to $I$ and $O$ if any (left or right) prolongation
of any polynomial $p_i \in P$ involutively reduces to zero
using $P$.
% is involutively
% reducible by some $p_j \in P$. In other words,
% \begin{enumerate}[(a)]
% \item
% for all $p_i \in P$ and for any variable
% $x_k \notin \mathcal{M}^L_I(\LM(p_i), \LM(P))$,
% $$\exists \: p_j \in P \mathrm{\; such \; that \;}
% p_j \mid_I x_kp_i;$$
% \item
% for all $p_i \in P$ and for any variable
% $x_k \notin \mathcal{M}^R_I(\LM(p_i), \LM(P))$,
% $$\exists \: p_j \in P \mathrm{\; such \; that \;}
% p_j \mid_I p_ix_k.$$
% \end{enumerate}
\end{defn}

\begin{defn} \label{IBNC}
\index{involutive basis}
Given an involutive division $I$ and an admissible monomial
ordering $O$, an autoreduced set of noncommutative polynomials $P$ is an
{\it Involutive Basis} \index{basis!involutive}
with respect to $I$ and $O$
if any multiple $up_iv$ of any polynomial $p_i \in P$
by any terms $u$ and $v$ involutively reduces to zero using $P$.
% is involutively reducible by some $p_j \in P$. In
% other words, for all $p_i \in P$ and for any monomials $u$
% and $v$, $$\exists \: p_j \in P \mathrm{\; such \; that \;}
% p_j \mid_I up_iv.$$
\end{defn}

% \begin{remark}
% By the admissibility of the chosen monomial ordering,
% any prolongation of
% a member of a Locally Involutive Basis involutively
% reduces to zero, and any multiple of a member of an
% Involutive Basis involutively reduces to zero.
% \end{remark}

\section{The Noncommutative Involutive Basis Algorithm}
\label{5point4}

To compute a (noncommutative) Locally Involutive Basis,
we use Algorithm \ref{noncom-inv}, an algorithm that
is virtually identical to
Algorithm \ref{com-inv}, apart from the fact that
at the beginning of the first {\bf while} loop, the set $S$ is
constructed in different ways.

\begin{algorithm}
\setlength{\baselineskip}{3.5ex}
\caption{The Noncommutative Involutive Basis Algorithm}
\label{noncom-inv}
\begin{algorithmic}
\vspace*{2mm}
\REQUIRE{A Basis $F = \{f_1, f_2, \hdots, f_m\}$ for an ideal $J$
         over a noncommutative polynomial ring
         $R\langle x_1, \hdots x_n\rangle$;
         an admissible monomial ordering $O$;
         an involutive division $I$.}
\ENSURE{A Locally Involutive Basis $G = \{g_1, g_2, \hdots, g_p\}$ for $J$ (in the
        case of termination).}
\vspace*{1mm}
\STATE
$G  = \emptyset$; \\
$F = \mathrm{Autoreduce}(F)$;
\WHILE{($G == \emptyset$)}
\STATE
$S = \{x_if \mid f \in F, \: x_i \notin \mathcal{M}^L_I(f, F)\}
     \cup \{fx_i \mid f \in F, \: x_i \notin \mathcal{M}^R_I(f, F)\}$; \\
$s' = 0$;
\WHILE{($S \neq \emptyset$) {\bf and} ($s' == 0$)}
\STATE
Let $s$ be a polynomial in $S$ whose lead monomial is minimal
with respect to $O$; \\
$S = S\setminus \{s\}$; \\
$s' = \Rem_I(s,  F)$;
\ENDWHILE
\IF{($s' \neq 0$)}
\STATE
$F = \mathrm{Autoreduce}(F\cup\{s'\})$;
\ELSE
\STATE
$G = F$;
\ENDIF
\ENDWHILE
\STATE
{\bf return} $G$;
\end{algorithmic}
\vspace*{1mm}
\end{algorithm}

\section{Continuity and Conclusivity} \label{CoCo}

In the commutative case, when we construct a Locally Involutive
Basis using Algorithm \ref{com-inv}, we know that
the algorithm will always return a commutative Gr\"obner Basis
as long as we use an admissible monomial ordering and the
chosen involutive division possesses certain properties.
In summary,
\begin{enumerate}[(a)]
\item
Any Locally Involutive Basis returned by Algorithm \ref{com-inv}
is an Involutive Basis if the involutive division used
is continuous (Proposition \ref{LocalToGlobal});
\item
Algorithm \ref{com-inv} always terminates if (in addition)
the involutive division used is constructive, Noetherian and stable
(Proposition \ref{InvTermNoeth});
\item
Every Involutive Basis is a Gr\"obner Basis
(Theorem \ref{IisG}).
\end{enumerate}
In the noncommutative case, we cannot hope to
produce a carbon copy of the above results
because a finitely generated basis may have an
infinite Gr\"obner Basis, leading to the
conclusion that Algorithm \ref{noncom-inv} does not
always terminate.
The best we can therefore hope for is
if an ideal generated by a set of polynomials $F$ possesses
a finite Gr\"obner Basis with respect to some admissible
monomial ordering $O$, then $F$ also possesses a
finite Involutive Basis with respect to $O$ and some
involutive division $I$. We shall call any involutive
division that possesses this property {\it conclusive}.

\begin{defn}
Let $F$ be an arbitrary basis generating an ideal over a
noncommutative polynomial ring,
and let $O$ be an arbitrary admissible monomial ordering.
An involutive division $I$ is {\it conclusive}
\index{involutive division!conclusive} 
\index{conclusive involutive division} if Algorithm
\ref{noncom-inv} terminates with $F$, $I$ and $O$ as
input whenever Algorithm \ref{noncom-buch} terminates with
$F$ and $O$ as input.
\end{defn}

Of course it is easy enough to define the above property,
but much harder to prove that a particular involutive
division is conclusive. In fact, no involutive division
defined in this thesis will be shown to be conclusive,
and the existence of such divisions will be left as an
open question.

\subsection{Properties for Strong Involutive Divisions}

Here is a summary of facts that can be deduced
when using a strong involutive division.
\begin{enumerate}[(a)]
\item
Any Locally Involutive Basis returned by Algorithm \ref{noncom-inv}
is an Involutive Basis if the involutive division used
is strong and continuous (Proposition \ref{LocalToGlobalNC});
\item
Algorithm \ref{noncom-inv} always terminates
whenever Algorithm \ref{noncom-buch} terminates if
(in addition) the involutive division used is conclusive;
\item
Every Involutive Basis with respect to a strong
involutive division is a Gr\"obner Basis
(Theorem \ref{IisGNC}).
\end{enumerate}

Let us now prove the assertions made in parts (a) and (c)
of the above list, beginning by defining what is meant
by a continuous involutive division in the noncommutative case.

\begin{defn} \label{ncc}
Let $I$ be a fixed involutive division;
let $w$ be a fixed monomial;
let $U$ be any set of monomials;
and consider any sequence
$(u_1, \; u_2, \; \hdots, \; u_k)$ of monomials from $U$
($u_i \in U$ for all $1 \leqslant i \leqslant k$),
each of which is a conventional divisor of $w$
(so that $w = \ell_iu_ir_i$ for all $1 \leqslant i \leqslant k$,
where the $\ell_i$ and the $r_i$ are monomials).
For all $1 \leqslant i < k$, suppose that the monomial $u_{i+1}$
satisfies exactly one of the following conditions.
\begin{enumerate}[(a)]
\item
$u_{i+1}$ involutively divides a left prolongation of $u_i$,
so that $\deg(\ell_i) \geqslant 1$; $\SUFF(\ell_i, 1) \notin
\mathcal{M}^L_I(u_i, U)$; and
$u_{i+1} \mid_I (\SUFF(\ell_i, 1))u_i$.
\item
$u_{i+1}$ involutively divides a right prolongation of $u_i$,
so that $\deg(r_i) \geqslant 1$; $\PRE(r_i, 1) \notin
\mathcal{M}^R_I(u_i, U)$; and
$u_{i+1} \mid_I u_i(\PRE(r_i, 1))$.
\end{enumerate}
Then $I$ is {\it continuous at} $w$ if all the pairs $(\ell_i, r_i)$
are distinct ($(\ell_i, r_i) \neq (\ell_j, r_j)$ for all $i \neq j$);
$I$ is a {\it continuous}
\index{involutive division!continuous}
involutive division if $I$ is continuous 
\index{continuous involutive division} for all possible $w$.
\end{defn}

\begin{prop} \label{LocalToGlobalNC}
If an involutive division $I$ is strong and continuous,
and a given set of polynomials $P$ is a Locally Involutive
Basis with respect to $I$ and some admissible monomial
ordering $O$, then $P$ is an Involutive Basis with respect
to $I$ and $O$.
\end{prop}
\begin{pf}
Let $I$ be a strong and continuous involutive
division; let $O$ be an admissible monomial ordering; and let
$P$ be a Locally Involutive Basis with respect to $I$ and $O$.
Given any polynomial $p \in P$ and any terms $u$ and $v$,
in order to show that $P$ is an Involutive Basis with
respect to $I$ and $O$, we
must show that
$upv \xymatrix{\ar[r]_I_(1){P} &} 0$.
% $\exists \, p' \in P$ such that
% $p' \mid_I pu$.

If $p \mid_I upv$ we are done, as we can use $p$ to
involutively reduce $upv$ to obtain a zero remainder.
Otherwise, either $\exists \, y_1
\notin \mathcal{M}^L_I(\LM(p), \LM(P))$ such that
$y_1 = \SUFF(u, 1)$, or $\exists \, y_1 \notin
\mathcal{M}^R_I(\LM(p), \LM(P))$ such that
$y_1 = \PRE(v, 1)$. Without loss of generality,
assume that the first case applies.
By Local Involutivity, the prolongation $y_1p$
involutively reduces to zero using $P$. Assuming that
the first step of this involutive reduction
involves the polynomial $p_1 \in P$,
we can write
\begin{equation} \label{ctsE1NC}
y_1p = u_1p_1v_1 + \sum_{a=1}^A u_{\alpha_a}p_{\alpha_a}v_{\alpha_a},
\end{equation}
where $p_{\alpha_a} \in P$ and $u_1, v_1,
u_{\alpha_a}, v_{\alpha_a}$ are terms
such that $p_1$ and each $p_{{\alpha}_{a}}$
involutively divide $u_1p_1v_1$ and each
$u_{\alpha_a}p_{\alpha_a}v_{\alpha_a}$ respectively.
Multiplying both sides of Equation (\ref{ctsE1NC}) on the left
by $u' := \PRE(u, \deg(u)-1)$ and on the right by $v$,
we obtain the equation
\begin{equation} \label{ctsE2NC}
upv = u'u_1p_1v_1v +
\sum_{a=1}^A u'u_{\alpha_a}p_{\alpha_a}v_{\alpha_a}v.
\end{equation}
If $p_1 \mid_I upv$, it is clear that
we can use $p_1$ to involutively
reduce the polynomial $upv$ to obtain the polynomial
$\sum_{a=1}^A u'u_{\alpha_a}p_{\alpha_a}v_{\alpha_a}v$.
By Proposition \ref{plusNC}, we can then continue
to involutively reduce $upv$ by repeating this proof
on each polynomial $u'u_{\alpha_a}p_{\alpha_a}v_{\alpha_a}v$
individually (where $1 \leqslant a \leqslant A$),
noting that this process will terminate because of the
admissibility of $O$ (we have
$\LM(u'u_{\alpha_a}p_{\alpha_a}v_{\alpha_a}v) < \LM(upv)$ for
all $1 \leqslant a \leqslant A$).

Otherwise, if $p_1$ does not involutively divide $upv$,
either $\exists \, y_2
\notin \mathcal{M}^L_I(\LM(p_1), \LM(P))$ such that
$y_2 = \SUFF(u'u_1, 1)$,
or $\exists \, y_2 \notin
\mathcal{M}^R_I(\LM(p_1), \LM(P))$ such that
$y_2 = \PRE(v_1v, 1)$. This time
(again without loss of generality), assume that the
second case applies.
By Local Involutivity, the prolongation $p_1y_2$
involutively reduces to zero using $P$. Assuming that
the first step of this involutive reduction
involves the polynomial $p_2 \in P$,
we can write
\begin{equation} \label{ctsE3NC}
p_1y_2 = u_2p_2v_2
+ \sum_{b=1}^B u_{\beta_b}p_{\beta_b}v_{\beta_b},
\end{equation}
where $p_{\beta_b} \in P$ and $u_2, v_2,
u_{\beta_b}, v_{\beta_b}$ are
terms such that $p_2$ and each $p_{\beta_b}$
involutively divide $u_2p_2v_2$ and each
$u_{\beta_b}p_{\beta_b}v_{\beta_b}$ respectively.
Multiplying both sides of Equation (\ref{ctsE3NC}) on
the left by $u'u_1$ and on the right by % $v'v$, where
% $v' := \SUFF(v_1, \deg(v_1)-1)$,
$v' := \SUFF(v_1v, \deg(v_1v)-1)$,
we obtain the equation
\begin{equation} \label{ctsE4NC}
% u'u_1p_1v_1v = u'u_1u_2p_2v_2v'v
% + \sum_{b=1}^B u'u_1u_{\beta_b}p_{\beta_b}v_{\beta_b}v'v.
u'u_1p_1v_1v = u'u_1u_2p_2v_2v'
+ \sum_{b=1}^B u'u_1u_{\beta_b}p_{\beta_b}v_{\beta_b}v'.
\end{equation}
Substituting for $u'u_1p_1v_1v$ from Equation
(\ref{ctsE4NC}) into Equation (\ref{ctsE2NC}), we obtain
the equation
\begin{equation} \label{ctsE5NC}
% upv = u'u_1u_2p_2v_2v'v +
upv = u'u_1u_2p_2v_2v' +
\sum_{a=1}^A u'u_{\alpha_a}p_{\alpha_a}v_{\alpha_a}v +
% \sum_{b=1}^B u'u_1u_{\beta_b}p_{\beta_b}v_{\beta_b}v'v.
\sum_{b=1}^B u'u_1u_{\beta_b}p_{\beta_b}v_{\beta_b}v'.
\end{equation}
If $p_2 \mid_I upv$, it is clear that
we can use $p_2$ to involutively
reduce the polynomial $upv$ to obtain the polynomial
$\sum_{a=1}^A u'u_{\alpha_a}p_{\alpha_a}v_{\alpha_a}v +
\sum_{b=1}^B u'u_1u_{\beta_b}p_{\beta_b}v_{\beta_b}v'$.
As before, we can then use Proposition \ref{plusNC}
to continue the involutive reduction of $upv$ by
repeating this proof on each summand individually.

Otherwise, if $p_2$ does not involutively divide $upv$,
we continue by induction, obtaining a sequence
$p, p_1, p_2, p_3, \hdots$
of elements in $P$. By construction, each element in the sequence
divides $upv$. By continuity (at $\LM(upv)$), no two elements in the
sequence divide $upv$ in the same way.
Because $upv$ has a finite number of
subwords, the sequence must be finite,
terminating with an involutive divisor $p' \in P$ of $upv$,
which then allows us to finish the proof through use of
Proposition \ref{plusNC} and the admissibility of $O$.
\end{pf}

\begin{thm} \label{IisGNC}
An Involutive Basis with respect to a strong
involutive division is a Gr\"obner Basis.
\end{thm}
\begin{pf}
Let $G = \{g_1, \hdots, g_m\}$ be an Involutive
Basis with respect to some strong involutive division $I$ and
some admissible monomial ordering $O$, where each
$g_i \in G$ (for all $1 \leqslant i \leqslant m$) is a
member of the polynomial ring $R\langle x_1, \hdots, x_n\rangle$.
To prove that
$G$ is a Gr\"obner Basis, we must show that all
S-polynomials involving elements of $G$ conventionally
reduce to zero using $G$. Recall that each
S-polynomial corresponds to an overlap between
the lead monomials of two (not necessarily distinct)
elements $g_i, g_j \in G$. Consider such an
arbitrary overlap, with corresponding S-polynomial
$$\mathrm{S\mbox{-}pol}(\ell_i, g_i, \ell_j, g_j) =
c_2\ell_ig_ir_i - c_1\ell_jg_jr_j.$$
Because $G$ is an Involutive Basis, it is clear that
$c_2\ell_ig_ir_i \xymatrix{\ar[r]_{I}_(1){G} &}0$
and
$c_1\ell_jg_jr_j \xymatrix{\ar[r]_{I}_(1){G} &}0$.
By Proposition \ref{plusNC}, it follows that
$\mathrm{S\mbox{-}pol}(\ell_i, g_i, \ell_j, g_j)
\xymatrix{\ar[r]_{I}_(1){G} &}0$. But every involutive
reduction is a conventional reduction, so we can deduce that
$\mathrm{S\mbox{-}pol}(\ell_i, g_i, \ell_j, g_j)
\rightarrow_G 0$ as required.
\end{pf}

\begin{lem} \label{urnc}
\index{unique remainder}
Given an Involutive Basis $G$ with respect to a strong involutive
division, remainders are involutively unique with respect to $G$.
\index{remainder!unique}
\end{lem}
\begin{pf}
Let $G$ be an Involutive Basis with respect to some
strong involutive division $I$ and some admissible monomial
ordering $O$. Theorem \ref{IisGNC} tells us that $G$ is
a Gr\"obner Basis with respect to $O$ and thus
remainders are conventionally unique with respect to
$G$. To prove that remainders are involutively unique
with respect to $G$, we must show that the conventional and
involutive remainders of an arbitrary polynomial $p$ with respect
to $G$ are identical. For this it is sufficient to show
that a polynomial $p$ is conventionally reducible by $G$ if
and only if it is involutively reducible by $G$.
($\Rightarrow$) Trivial as every involutive reduction is a
conventional reduction. ($\Leftarrow$)
If a polynomial $p$ is conventionally reducible by a
polynomial $g \in G$, it follows that $\LM(p) = u\LM(g)v$ for
some monomials $u$ and $v$. But $G$ is an Involutive Basis, so there
must exist a polynomial $g' \in G$ such that
$\LM(g') \mid_I u\LM(g)v$. Thus $p$ is also
involutively reducible by $G$.
% $u\LM(g)v = u'\LM(g')v'$ for some monomials $u'$ and
% $v'$ that are (respectively) left and right
% multiplicative (over $G$) for $g'$. Thus $p$ is also
% involutively reducible by $G$.
\end{pf}

\subsection{Properties for Weak Involutive Divisions}

While it is true that the previous three results
(Proposition \ref{LocalToGlobalNC}, Theorem \ref{IisGNC}
and Lemma \ref{urnc}) do not apply if a
weak involutive division has been chosen, we will now
show that corresponding results can be obtained for
weak involutive divisions that are also {\it Gr\"obner}
involutive divisions.

% Despite this, it is still possible
% to prove that an Involutive Basis with respect to a
% weak involutive division is a Gr\"obner Basis on a
% case-by-case basis, and we shall call such involutive
% divisions {\it Gr\"obner divisions}.

\begin{defn} \label{weakGrob}
A weak involutive division $I$ is
a {\it Gr\"obner} involutive division
\index{involutive division!Gr\"obner} if
every Locally Involutive Basis with respect to $I$
is a Gr\"obner Basis.
\index{Gr\"obner involutive division}
\end{defn}

It is an easy consequence of Definition
\ref{weakGrob} that any Involutive Basis with
respect to a weak and Gr\"obner involutive division
is a Gr\"obner Basis; it therefore follows that we
can also prove an analog of Lemma \ref{urnc} for
such divisions. To complete the mirroring of the results of
Proposition \ref{LocalToGlobalNC}, Theorem \ref{IisGNC}
and Lemma \ref{urnc} for weak and Gr\"obner involutive
divisions, it remains to show that a
Locally Involutive Basis with
respect to a weak; continuous and Gr\"obner
involutive division is an Involutive Basis.

\begin{prop} \label{LocalToGlobalNCweak}
If an involutive division $I$ is weak; continuous and Gr\"obner,
and if a given set of polynomials $P$ is a Locally Involutive
Basis with respect to $I$ and some admissible monomial
ordering $O$, then $P$ is an Involutive Basis with respect
to $I$ and $O$.
\end{prop}
\begin{pf}
Let $I$ be a weak; continuous and Gr\"obner involutive
division; let $O$ be an admissible monomial ordering; and let
$P$ be a Locally Involutive Basis with respect to $I$ and $O$.
Given any polynomial $p \in P$ and any terms $u$ and $v$,
in order to show that $P$ is an Involutive Basis with
respect to $I$ and $O$, we
must show that
$upv \xymatrix{\ar[r]_I_(1){P} &} 0$.

For the first part of the proof, we proceed as in the
proof of Proposition \ref{LocalToGlobalNC} to find an
involutive divisor $p' \in P$ of $upv$ using the
continuity of $I$ at $\LM(upv)$. This then allows us to
involutive reduce $upv$ using $p'$ to obtain a
polynomial $q$ of the form
\begin{equation} \label{weakLTGE}
q = \sum_{a=1}^A u_{a}p_{\alpha_a}v_{_a},
\end{equation}
where $p_{\alpha_a} \in P$ and the $u_a$ 
and the $v_a$ are terms.
% such that each $p_{\alpha_a}$ involutively divides
% each $u_{a}p_{\alpha_a}v_{_a}$.

For the second part of the proof, we now use the fact
that $P$ is a Gr\"obner Basis to find a polynomial
$q' \in P$ such that $q'$ conventionally divides $q$
(such a polynomial will always exist because $q$ is clearly
a member of the ideal generated by $P$). If $q'$
is an involutive divisor of $q$, then we can use
$q'$ to involutively reduce $q$ to obtain a
polynomial $r$ of the form shown
in Equation (\ref{weakLTGE}). Otherwise, if $q'$ is not 
an involutive divisor of $q$, we can use the fact that
$I$ is continuous at $\LM(q)$ to find such an involutive 
divisor, which we can then use to involutive reduce $q$
to obtain a polynomial $r$, again of the form shown
in Equation (\ref{weakLTGE}). In both cases, we now
proceed by induction on $r$, noting that this process will
terminate because of the admissibility of $O$
(we have $\LM(r) < \LM(q)$).
\end{pf}

To summarise, here is the situation for weak and
Gr\"obner involutive divisions.
\begin{enumerate}[(a)]
\item
Any Locally Involutive Basis returned by Algorithm \ref{noncom-inv}
is an Involutive Basis if the involutive division used
is weak; continuous and Gr\"obner (Proposition \ref{LocalToGlobalNCweak});
\item
Algorithm \ref{noncom-inv} always terminates
whenever Algorithm \ref{noncom-buch} terminates if
(in addition) the involutive division used is conclusive;
\item
Every Involutive Basis with respect to a weak and Gr\"obner
involutive division is a Gr\"obner Basis.
\end{enumerate}

% \begin{remark}
% By Theorem \ref{IisGNC}, every strong involutive division
% is a Gr\"obner division.
% \end{remark}

\section{Noncommutative Involutive Divisions} \label{5point5}
\index{involutive divisions}

Before we consider some examples of useful
noncommutative involutive divisions,
let us remark that it is possible
to categorise any noncommutative involutive division
somewhere between the following two
{\it extreme} global divisions.

\begin{defn}[The Empty Division]
\index{involutive division!empty}
Given any monomial $u$, let $u$ have no
(left or right) multiplicative variables.
\index{empty division}
\end{defn}

\begin{defn}[The Full Division]
\index{involutive division!full}
Given any monomial $u$, let $u$ have no
(left or right) nonmultiplicative variables
(in other words, all variables are left and
right multiplicative for $u$).
\index{full division}
\end{defn}

\begin{remark}
It is clear that any set of polynomials $G$
will be an Involutive Basis with respect to the (weak)
full division as any multiple of a polynomial
$g \in G$ will be involutively reducible by $g$
(all conventional divisors are involutive
divisors); in contrast it is impossible to find a
finite Locally Involutive Basis for $G$
with respect to the (strong) empty
division as there will always be a prolongation
of an element of the current basis that is
involutively irreducible.
\end{remark}

\subsection{Two Global Divisions} \label{TGD}

Whereas most of the theory seen so far in this chapter has
closely mirrored the corresponding commutative theory from
Chapter \ref{ChCIB}, the commutative involutive divisions
(Thomas, Janet and Pommaret)
seen in the previous chapter do not generalise to the
noncommutative case, or at the very least do not yield
noncommutative involutive divisions of any value. Despite this,
an essential property of these divisions is that they
ensure that the least common multiple
$\lcm(\LM(p_1), \LM(p_2))$ associated with an
S-polynomial $\mathrm{S\mbox{-}pol}(p_1, p_2)$ is
involutively irreducible by at least one of $p_1$ and $p_2$,
ensuring that the
S-polynomial $\mathrm{S\mbox{-}pol}(p_1, p_2)$ 
is constructed and involutively reduced
during the course of the Involutive Basis algorithm.

To ensure that the corresponding process occurs in the
noncommutative Involutive Basis algorithm, we must ensure
that all overlap words associated to the S-polynomials
of a particular basis are involutively irreducible (as placed
in the overlap word) by at least one of the polynomials
associated to each overlap word. This obviously holds
true for the empty division, but it will also hold
true for the following two global involutive divisions,
where all variables are either assigned to be
left multiplicative and right nonmultiplicative, or
left nonmultiplicative and right multiplicative.

\begin{defn}[The Left Division]
\index{involutive division!left}
\index{$<$@$\lhd$}
Given any monomial $u$, the left division
$\lhd$ assigns no left nonmultiplicative
variables to $u$, and assigns no right
multiplicative variables to $u$
(in other words, all variables are
left multiplicative and right
nonmultiplicative for $u$).
\index{left division}
\end{defn}

\begin{defn}[The Right Division]
\index{involutive division!right}
\index{$<$@$\rhd$}
Given any monomial $u$, the right division
$\rhd$ assigns no left multiplicative
variables to $u$, and assigns no right
nonmultiplicative variables to $u$
(in other words, all variables are
left nonmultiplicative and right
multiplicative for $u$).
\index{right division}
\end{defn}

\begin{prop}
The left and right divisions are strong
involutive divisions.
\end{prop}
\begin{pf}
We will only give the proof for the left
division -- the proof for the right division
will follow by symmetry (replacing `left' by `right',
and so on).

To prove that the left division is a strong
involutive division, we need to show that
the three conditions of Definition
\ref{noncom-div-defn} hold.
\begin{itemize}
\item
{\bf Disjoint Cones Condition} \\
Consider two involutive cones
$\mathcal{C}_{\lhd}(u_1)$ and
$\mathcal{C}_{\lhd}(u_2)$ associated to
two monomials $u_1, u_2$ over some noncommutative
polynomial ring $\mathcal{R}$. If
$\mathcal{C}_{\lhd}(u_1) \cap
\mathcal{C}_{\lhd}(u_2) \neq \emptyset$, then there must
be some monomial $v \in \mathcal{R}$ such that
$v$ contains both monomials $u_1$ and $u_2$
as subwords, and (as placed in $v$) both
$u_1$ and $u_2$ must be involutive divisors of
$v$. By definition of $\lhd$, both
$u_1$ and $u_2$ must be suffices of $v$.
Thus, assuming (without loss of generality) that
$\deg(u_1) > \deg(u_2)$, we are able to draw
the following diagram summarising the
situation.
$$
\xymatrix @R=1.5pc{
\ar@{<->}[rrrrrr]^*+{v} &&&&&& \\
&& \ar@{<->}[rrrr]^*+{u_1} &&&& \\
&&&& \ar@{<->}[rr]^*+{u_2} &&
}
$$
But now, assuming that $u_1 = u_3u_2$ for
some monomial $u_3$, it is clear that
$\mathcal{C}_{\lhd}(u_1) \subset \mathcal{C}_{\lhd}(u_2)$
because any monomial $w \in \mathcal{C}_{\lhd}(u_1)$
must be of the form $w = w'u_1$ for some monomial
$w'$; this means that $w = w'u_3u_2 \in
\mathcal{C}_{\lhd}(u_2)$.
\item
{\bf Unique Divisor Condition} \\
As a monomial $v$ is only involutively divisible
by a monomial $u$ with respect to the left division
if $u$ is a suffix of $v$, it is clear that
$u$ can only involutively divide $v$ in
at most one way.
\item
{\bf Subset Condition} \\
Follows immediately due to the left division
being a global division.
\end{itemize}
\end{pf}

\begin{prop}
The left and right divisions are continuous.
\end{prop}
\begin{pf}
Again we will only treat the case of the left division.
Let $w$ be an arbitrary fixed monomial;
let $U$ be any set of monomials; and consider any sequence
$(u_1, \; u_2, \; \hdots, \; u_k)$ of monomials from $U$
($u_i \in U$ for all $1 \leqslant i \leqslant k$),
each of which is a conventional divisor of $w$
(so that $w = \ell_iu_ir_i$ for all $1 \leqslant i \leqslant k$,
where the $\ell_i$ and the $r_i$ are monomials).
For all $1 \leqslant i < k$, suppose that the monomial $u_{i+1}$
satisfies condition (b) of Definition \ref{ncc} (condition (a)
can never be satisfied because $\lhd$ never assigns any left
nonmultiplicative variables). To show that $\lhd$ is continuous, we
must show that no two pairs $(\ell_i, r_i)$ and $(\ell_j, r_j)$ are
the same, where $i \neq j$.

Consider an arbitrary monomial $u_i$ from the sequence,
where $1 \leqslant i < k$. Because $\lhd$ assigns
no right multiplicative variables, the next monomial
$u_{i+1}$ in the sequence must be a suffix of the
prolongation $u_i(\PRE(r_i, 1))$ of $u_i$, so that
$\deg(r_{i+1}) = \deg(r_i) - 1$. It is therefore
clear that no two identical $(\ell, r)$ pairs can be found in the
sequence, as $\deg(r_1) > \deg(r_2) > \cdots > \deg(r_k)$.
\end{pf}

To illustrate the difference between the overlapping
cones of a noncommutative Gr\"obner Basis and the disjoint
cones of a noncommutative Involutive Basis with respect
to the left division, consider the following example.

\begin{example} \label{infloop}
Let $F := \{2xy+y^2+5, \, x^2+y^2+8\}$ be a basis
over the polynomial ring 
$\mathbb{Q}\langle x, y\rangle$,
and let the monomial ordering be DegLex. Applying
Algorithm \ref{noncom-buch} to $F$, we obtain the
Gr\"obner Basis $G := \{2xy+y^2+5, \, x^2+y^2+8, \,
5y^3-10x+37y, \, 2yx+y^2+5\}$. Applying Algorithm
\ref{noncom-inv} to $F$ with respect to the left
involutive division, we obtain the Involutive
Basis $H := \{2xy+y^2+5, \, x^2+y^2+8, \,
5y^3-10x+37y, \, 5xy^2+5x-6y, \, 2yx+y^2+5\}$.

To illustrate which monomials are reducible with respect
to the Gr\"obner Basis, we can draw a monomial lattice,
part of which is shown below. In the lattice, we
draw a path from the (circled) lead monomial
of any Gr\"obner Basis element to any multiple of that 
lead monomial, so that any monomial which lies on some
path in the lattice is reducible 
by one or more Gr\"obner Basis elements.
To distinguish between different Gr\"obner Basis
elements we use different arrow types; we also arrange
the lattice so that monomials of the same degree lie
on the same level.
% Arrows: ~, -, =, --, ., :
$$\xymatrix @C=0.04pc @R=2pc{
&&&&&&&& 1 \\
&&&& x &&&&&&&& y \\
&& *+[o][F-]{x^2} \ar@{~}@<0.75ex>[ddl] \ar@{~}[ddl]
\ar@{~}[ddr] \ar@{~}[ddrrrrr]
&&&& *+[o][F-]{xy} \ar@{-}[llldd]
\ar@{-}[ldd] \ar@{-}[rrrdd] \ar@{-}[rrrrrdd]
&&&& *+[o][F-]{yx} \ar@{--}[ddlll] \ar@{--}[ddlllll]
\ar@{--}[ddr] \ar@{--}[ddrrr]
&&&& y^2 \\ \\
& x^3 \ar@{~}[dddl] \ar@{~}@<-0.75ex>[dddl]
\ar@{~}[ddd] \ar@{~}[dddrrr]
&& x^2y \ar@{-}[dddll] \ar@{-}[dddrr]
\ar@{~}@<0.75ex>[dddl] \ar@{~}@<0.75ex>[rrrrrrrddd]
\ar@{~}@<-0.75ex>[dddll] \ar@{~}@<0.75ex>[dddrr]
\ar@{-}[dddl] \ar@{-}[rrrrrrrddd]
&& xyx \ar@{-}[dddlll] \ar@{-}[dddrr]
\ar@{-}[dddll] \ar@{-}[rrrrddd]
\ar@{--}@<-0.75ex>[dddlll] \ar@{--}@<-0.75ex>[dddrr]
\ar@{--}@<0.75ex>[dddll] \ar@{--}@<0.75ex>[rrrrddd]
&& yx^2 \ar@{--}[dddlll] \ar@{--}[dddllll]
\ar@{--}[dddrrrr] \ar@{--}[dddrrr]
\ar@<0.75ex>@{~}[dddlll] \ar@{~}@<-0.75ex>[dddllll]
\ar@{~}@<0.75ex>[dddrrrr] \ar@{~}@<-0.75ex>[dddrrr]
&& xy^2 \ar@{-}[llllddd] \ar@{-}[rrrddd]
\ar@{-}[dddlll] \ar@{-}[rrrrddd]
&& yxy \ar@{-}[dddllll] \ar@{-}[dddrr]
\ar@{-}[dddrrr] \ar@{-}[llddd]
\ar@{--}@<-0.75ex>[dddllll] \ar@{--}@<-0.75ex>[dddrr]
\ar@{--}@<0.75ex>[dddrrr] \ar@{--}@<0.75ex>[llddd]
&& y^2x \ar@{--}[dddr] \ar@{--}[dddrr]
\ar@{--}[dddll] \ar@{--}[dddlllllll]
&& *+[o][F-]{y^3} \ar@{.}[rddd] \ar@{.}@<0.75ex>[rddd]
\ar@{.}[ddd] \ar@{.}[dddlll]
\\ \\ \\
x^4 & x^3y & x^2yx & xyx^2 &
yx^3 & x^2y^2 & xy^2x & xyxy &&
yxyx & yx^2y & y^2x^2 & xy^3 &
yxy^2 & y^2xy & y^3x & y^4
}$$
Notice that many of the monomials in the lattice
are reducible by several of the Gr\"obner Basis elements.
For example, the monomial $x^2yx$ is reducible
by the Gr\"obner Basis elements $2xy+y^2+5$; $x^2+y^2+8$
and $2yx+y^2+5$. In contrast, any monomial in the
corresponding lattice for the Involutive Basis
may only be involutively reducible by at most one
element in the Involutive Basis. We illustrate this by
the following diagram, where we note that in the 
involutive lattice, a monomial only lies on a particular
path if a member of the Involutive Basis is an
involutive divisor of that monomial.
$$\xymatrix @C=0.04pc @R=2pc{
&&&&&&&& 1 \\
&&&& x &&&&&&&& y \\
&& *+[o][F-]{x^2}  \ar@{~}[ddl] \ar@{~}[ddrrrrr]
&&&& *+[o][F-]{xy} \ar@{-}[llldd]\ar@{-}[rrrrrdd]
&&&& *+[o][F-]{yx} \ar@{--}[ddlllll] \ar@{--}[ddrrr]
&&&& y^2 \\ \\
& x^3 \ar@{~}[dddl] \ar@{~}[dddrrr]
&& x^2y \ar@{-}[dddll] \ar@{-}[rrrrrrrddd]
&& xyx \ar@{--}[dddlll] \ar@{--}[rrrrddd]
&& yx^2 \ar@{~}[dddllll] \ar@{~}[dddrrrr] 
&& *+[o][F-]{xy^2} \ar@{{}{o}{}}[llllddd] \ar@{{}{o}{}}[rrrrddd]
&& yxy \ar@{-}[dddllll] \ar@{-}[dddrrr] 
&& y^2x \ar@{--}[dddrr] \ar@{--}[dddlllllll]
&& *+[o][F-]{y^3} \ar@{.}[rddd] \ar@{.}[dddlll]
\\ \\ \\
x^4 & x^3y & x^2yx & xyx^2 &
yx^3 & x^2y^2 & xy^2x & xyxy &&
yxyx & yx^2y & y^2x^2 & xy^3 &
yxy^2 & y^2xy & y^3x & y^4
}$$
Comparing the two monomial lattices, we
see that any monomial that is conventionally
divisible by the Gr\"obner Basis is
uniquely involutively divisible by the
Involutive Basis. In other words, the
involutive cones of the Involutive Basis
form a disjoint cover of the conventional
cones of the Gr\"obner Basis.
\end{example}

\subsubsection{Fast Reduction}

In the commutative case, we can sometimes use
the properties of an involutive division to
speed up the process of involutively reducing
a polynomial with respect to a set of
polynomials. For example, the Janet tree
\cite{Gerdt01a, Gerdt01b} enables us to quickly
determine whether a polynomial is involutively 
reducible by a set of polynomials with respect to
the Janet involutive division.

In the noncommutative case, 
we usually use Algorithm \ref{noncom-inv-div}
to involutively reduce a polynomial $p$
with respect to a set of polynomials $P$.
When this is done with respect to the left or right
divisions however, we can improve
Algorithm \ref{noncom-inv-div}
by taking advantage of the fact that a monomial
$u_1$ only involutively divides another
monomial $u_2$ with respect to the left (right)
division if $u_1$ is a suffix (prefix) of $u_2$.

For the left division, we can replace the code
found in the first {\bf if} loop of Algorithm
\ref{noncom-inv-div} with the following code in
order to obtain an improved algorithm.
\begin{algorithmic}
\vspace*{2mm}
\IF{($\LM(p_j)$ is a suffix of $u$)}
\STATE
found = true; \\
$p = p - (c\LC(p_j)^{-1})u_{\ell}p_j$, where
$u_{\ell} = \PRE(p, \deg(p)-\deg(p_j))$;
\ELSE
\STATE
$j = j+1$;
\ENDIF
\end{algorithmic}
We note that only one operation is required to
determine whether the monomial $\LM(p_j)$
involutively divides the monomial $u$ here
(test to see if $\LM(p_j)$ is a suffix of $u$);
whereas in general there are many ways
that $\LM(p_j)$ can conventionally divide $u$,
each of which has to be tested to see whether it is
an involutive reduction. This means that, with respect
to the left or right divisions, we can
determine whether a monomial $u$ is involutively
irreducible with respect to a set of polynomials $P$
in linear time (linear in the number of elements
in $P$); whereas in general we
can only do this in quadratic time.

\subsection{An Overlap-Based Local Division}

Even though the left and right involutive
divisions are strong and continuous (so that
any Locally Involutive Basis returned by
Algorithm \ref{noncom-inv} is a noncommutative
Gr\"obner Basis), these divisions are not conclusive
as the following example demonstrates.

\begin{example}  \label{fiveA}
Let $F := \{xy-z, \, x+z, \, yz-z, \, xz, \,
zy+z, \, z^2\}$ be a basis over the polynomial
ring $\mathbb{Q}\langle x, y, z\rangle$,
and let the monomial ordering be DegLex. Applying
Algorithm \ref{noncom-buch} to $F$, we discover
that $F$ is a noncommutative Gr\"obner Basis
($F$ is returned to us as the output of Algorithm
\ref{noncom-buch}).
When we apply Algorithm \ref{noncom-inv} to
$F$ with respect to the left involutive division however,
we notice that the algorithm goes into an infinite
loop, constructing the infinite basis
$G := F \cup \{zy^n-z, \, xy^n+z, zy^m+z, xy^m-z\}$,
where $n \geqslant 2$, $n$ even and
$m \geqslant 3$, $m$ odd.
\end{example}

The reason why Algorithm \ref{noncom-inv} goes into
an infinite loop in the above example is that
the right prolongations of the polynomials
$xy-z$ and $zy+z$ by the variable $y$ do not
involutively reduce to zero (they reduce to the
polynomials $xy^2+z$ and $zy^2-z$ respectively).
These prolongations are the only prolongations
of elements of $F$ that do not involutively
reduce to zero, and this is also true for all
polynomials we subsequently add to $F$, thus
allowing Algorithm \ref{noncom-inv} to construct
the infinite set $G$.

Consider a modification of the left division 
where we assign the variable $y$ to be
right multiplicative for the (lead) monomials
$xy$ and $zy$. Then it is clear that $F$ will
be a Locally Involutive Basis with respect to
this modified division, but will it also be true
that $F$ is an Involutive Basis and
(had we not known so already) a Gr\"obner Basis?

Intuitively, for this particular example,
it would seem that the answer to both of the
above questions should be affirmative,
because the modified division still ensures that
all the overlap words associated with the S-polynomials
of $F$ are involutively irreducible (as placed in the overlap word)
by at least one of the polynomials associated
to each S-polynomial. This leads to the following
% The variable need not be nonmultiplicative because,
% considering just the set of overlap words
% associated with the set of S-polynomials of $F$,
% allowing $y$ to be right multiplicative for
% $xy$ and $yz$ does not change the fact that all
% the overlap words are involutively irreducible
% by at least one of the polynomials associated
% with that S-polynomial. This leads to the following
idea for a {\it local} involutive division, where we
refine the left division by choosing right
nonmultiplicative variables based on the
overlap words of S-polynomials associated to a set
of polynomials only (note that there will also be
a similar local involutive division refining the
right division called the right overlap division).

\begin{defn}[The Left Overlap Division $\mathcal{O}$]
\index{$O$@$\mathcal{O}$}
\index{involutive division!left overlap}
Let $U = \{u_1, \hdots, u_m\}$ be a set of monomials,
and assume that all variables are left and right
multiplicative for all elements of $U$ to begin with.
\index{left overlap division}
\begin{enumerate}[(a)]
\item
For all possible ways that
a monomial $u_j \in U$ is a subword of a (different)
monomial $u_i \in U$, so that
$$\SUB(u_i, k, k+\deg(u_j)-1) = u_j$$
for some integer $k$, if $u_j$ is not a suffix of $u_i$,
assign the variable
$\SUB(u_i, k+\deg(u_j), k+\deg(u_j))$ to be right
nonmultiplicative for $u_j$.
\item
For all possible ways that
a proper prefix of a monomial $u_i \in U$ is equal to
a proper suffix of a (not necessarily different)
monomial $u_j \in U$, so that
$$\PRE(u_i, k) = \SUFF(u_j, k)$$ for some integer $k$
and $u_i$ is not a subword of $u_j$ or vice-versa, assign
the variable $\SUB(u_i, k+1, k+1)$ to be
right nonmultiplicative for $u_j$.
\end{enumerate}
\end{defn}

\begin{remark}
One possible algorithm for the
left overlap division is presented 
in Algorithm \ref{inv-div}, where the reason
for insisting that the input set of monomials is ordered
with respect to DegRevLex is in order to minimise the number
of operations needed to discover all the subword
overlaps (a monomial of degree $d_1$ can never be
a subword of a different monomial of degree 
$d_2 \leqslant d_1$).
\end{remark}

\begin{algorithm}
\setlength{\baselineskip}{3.5ex}
\caption{The Left Overlap Division $\mathcal{O}$}
\label{inv-div}
\begin{algorithmic}
\vspace*{2mm}
% p -> m, m_i -> u_i
\REQUIRE{A set of % distinct 
         monomials $U = \{u_1, \hdots, u_m\}$
         ordered by DegRevLex ($u_1 \geqslant u_2 
         \geqslant \cdots \geqslant u_m$),
         where $u_i \in R\langle x_1, \hdots, x_n\rangle$.}
\ENSURE{A table $T$ of left and right multiplicative variables
        for all $u_i \in U$, where each entry of $T$ is either 1
        (multiplicative) or 0 (nonmultiplicative).}
\vspace*{1mm}
\STATE
Create a table $T$ of multiplicative 
variables as shown below: \\ \vspace*{1mm}
\begin{tabular}{c|ccccccc}
& $x_{1}^L$ & $x_{1}^R$ & $x_{2}^L$ & $x_{2}^R$
& $\cdots$ & $x_{n}^L$ & $x_{n}^R$ \\
\hline
$u_1$ & 1 & 1 & 1 & 1 & $\cdots$ & 1 & 1 \\
$u_2$ & 1 & 1 & 1 & 1 & $\cdots$ & 1 & 1 \\
$\vdots$ & $\vdots$ & $\vdots$ & $\vdots$ & $\vdots$
& $\ddots$ & $\vdots$ & $\vdots$ \\
$u_m$ & 1 & 1 & 1 & 1 & $\cdots$ & 1 & 1
\end{tabular} \\[1.5mm]
\FOR{{\bf each} monomial $u_i \in U$ ($1 \leqslant i \leqslant m$)}
\FOR{{\bf each} monomial $u_j \in U$ ($i \leqslant j \leqslant m$)}
\STATE
Let $u_i = x_{i_1}x_{i_2}\hdots x_{i_{\alpha}}$ and
$u_j = x_{j_1}x_{j_2}\hdots x_{j_{\beta}}$; \\[1mm]
%// Middle Overlaps \\
\IF{($i \neq j$)}
\FOR{{\bf each} $k$ ($1 \leqslant k < \alpha - \beta + 1$)}
  \IF{($\SUB(u_i, k, k+\beta-1) == u_j$)}
  \STATE
%    {\bf if} ($k < \alpha-\beta+1$) {\bf then}
    $T(u_j, x^R_{i_{k+\beta}}) = 0$; \\
%    {\bf end if}
  \ENDIF
\ENDFOR
\ENDIF
%\STATE
\vspace*{1mm}
%// Left \& Right Overlaps
\FOR{{\bf each} $k$ ($1 \leqslant k \leqslant \beta - 1$)}
 \IF{($\PRE(u_i, k) == \SUFF(u_j, k)$)}
 \STATE
 $T(u_j, x^R_{i_{k+1}}) = 0$; 
 \ENDIF
 \IF{($\SUFF(u_i, k) == \PRE(u_j, k)$)}
 \STATE
 $T(u_i, x^R_{j_{k+1}}) = 0$; 
 \ENDIF
%  \STATE
%  {\bf if} ($\PRE(u_i, k) == \SUFF(u_j, k)$) {\bf then}
%  $T(u_j, x^R_{i_{k+1}}) = 0$; \\
%  {\bf end if} \\
%  {\bf if} ($\SUFF(u_i, k) == \PRE(u_j, k)$) {\bf then}
%  $T(u_i, x^R_{j_{k+1}}) = 0$; \\
%  {\bf end if}
\ENDFOR \vspace*{1mm}
\ENDFOR
\ENDFOR
\STATE
{\bf return} $T$;
\end{algorithmic}
\vspace*{1mm}
\end{algorithm}

\begin{example} \label{fiveB}
Consider again the set of polynomials
$F := \{xy-z, \, x+z, \, yz-z, \, xz, \,
zy+z, \, z^2\}$ from Example \ref{fiveA}.
Here are the left and right multiplicative
variables for $\LM(F)$ with respect to the left overlap
division $\mathcal{O}$.
\begin{center}
\begin{tabular}{c|c|c}
$u$ & $\mathcal{M}_{\mathcal{O}}^L(u, \LM(F))$ &
$\mathcal{M}_{\mathcal{O}}^R(u, \LM(F))$ \\ \hline
$xy$  & $\{x, y, z\}$  & $\{x, y\}$ \\
$x$   & $\{x, y, z\}$  & $\{x\}$ \\
$yz$  & $\{x, y, z\}$  & $\{x\}$ \\
$xz$  & $\{x, y, z\}$  & $\{x\}$ \\
$zy$  & $\{x, y, z\}$  & $\{x, y\}$ \\
$z^2$ & $\{x, y, z\}$  & $\{x\}$ \\ \hline
\end{tabular}
\end{center}
When we apply Algorithm \ref{noncom-inv} to
$F$ with respect to the DegLex monomial ordering
and the left overlap division, $F$ is returned to us
as the output, an assertion that is easily
verified by showing that the 10 right prolongations
of elements of $F$ all involutively reduce to zero using $F$.
This means that $F$ is a Locally Involutive
Basis with respect to the left overlap division;
to show that $F$ (and indeed any Locally Involutive
Basis returned by Algorithm \ref{noncom-inv} with
respect to the left overlap division) is also an Involutive
Basis with respect to the left overlap division,
we need to show that the left overlap division is 
continuous and either strong or Gr\"obner; we begin
(after the following remark) by showing that the left
overlap division is continuous.
\end{example}

\begin{remark}
In the above example, the table of multiplicative
variables can be constructed from the table $T$ shown below,
a table that is obtained by applying Algorithm \ref{inv-div} to $\LM(F)$.
\begin{center}
\begin{tabular}{c|cccccc}
Monomial & $x^L$ & $x^R$ & $y^L$ & $y^R$ & $z^L$ & $z^R$ \\ \hline
$xy$  & 1 & 1 & 1 & 1 & 1 & 0 \\
$x$   & 1 & 1 & 1 & 0 & 1 & 0 \\
$yz$  & 1 & 1 & 1 & 0 & 1 & 0 \\
$xz$  & 1 & 1 & 1 & 0 & 1 & 0 \\
$zy$  & 1 & 1 & 1 & 1 & 1 & 0 \\
$z^2$ & 1 & 1 & 1 & 0 & 1 & 0 \\ \hline
\end{tabular}
\end{center}
The zero entries in $T$ correspond to the following
overlaps between the elements of $\LM(F)$.
\begin{center}
\begin{tabular}{c|c}
Table Entry & Overlap \\ \hline
$T(xy, z^R)$
& $\SUFF(xy, 1) = \PRE(yz, 1)$ \\
$T(x, y^R)$
& $\SUB(xy, 1, 1) = x$ \\
$T(x, z^R)$
& $\SUB(xz, 1, 1) = x$ \\
$T(yz, y^R)$
& $\SUFF(yz, 1) = \PRE(zy, 1)$ \\
$T(yz, z^R)$
& $\SUFF(yz, 1) = \PRE(z^2, 1)$ \\
$T(xz, y^R)$
& $\SUFF(xz, 1) = \PRE(zy, 1)$ \\
$T(xz, z^R)$
& $\SUFF(xz, 1) = \PRE(z^2, 1)$ \\
$T(zy, z^R)$
& $\SUFF(zy, 1) = \PRE(yz, 1)$ \\
$T(z^2, y^R)$
& $\SUFF(z^2, 1) = \PRE(zy, 1)$ \\
$T(z^2, z^R)$
& $\SUFF(z^2, 1) = \PRE(z^2, 1)$ \\ \hline
\end{tabular}
\end{center}
\end{remark}

\vspace*{5mm}
\begin{prop} \label{ov-cts-ch6}
The left overlap division $\mathcal{O}$ is continuous.
\end{prop}
\begin{pf}
Let $w$ be an arbitrary fixed monomial;
let $U$ be any set of monomials; and consider any sequence
$(u_1, \; u_2, \; \hdots, \; u_k)$ of monomials from $U$
($u_i \in U$ for all $1 \leqslant i \leqslant k$),
each of which is a conventional divisor of $w$
(so that $w = \ell_iu_ir_i$ for all $1 \leqslant i \leqslant k$,
where the $\ell_i$ and the $r_i$ are monomials).
For all $1 \leqslant i < k$, suppose that the monomial $u_{i+1}$
satisfies condition (b) of Definition \ref{ncc} (condition (a)
can never be satisfied because $\mathcal{O}$ never assigns any left
nonmultiplicative variables). To show that $\mathcal{O}$ is continuous, we
must show that no two pairs $(\ell_i, r_i)$ and $(\ell_j, r_j)$ are
the same, where $i \neq j$.

Consider an arbitrary monomial $u_i$ from the sequence,
where $1 \leqslant i < k$. By definition of $\mathcal{O}$,
the next monomial $u_{i+1}$ in the sequence cannot be
either a prefix or a proper subword of $u_i$. This leaves
two possibilities: (i) $u_{i+1}$ is a suffix of $u_i$
(in which case $\deg(u_{i+1}) < \deg(u_i)$);
or (ii) $u_{i+1}$ is a suffix of the
prolongation $u_iv_i$ of $u_i$, where
$v_i := \PRE(r_i, 1)$.
\begin{center}
\begin{tabular}{ccc}
Example of possibility (i) & \vspace*{5mm} &
Example of possibility (ii) \\
$\xymatrix @R=0.5pc{
\ar@{<->}[rrrr]^*+{u_{i}} &&&&
\ar@{-}[r]^*+{v_i} & \\
&& \ar@{<->}[rr]_*+{u_{i+1}} &&
}$
&&
$\xymatrix @R=0.5pc{
\ar@{<->}[rrrr]^*+{u_{i}} &&&&
\ar@{-}[r]^*+{v_i} & \\
&& \ar@{<->}[rrr]_*+{u_{i+1}} &&&
}$
\end{tabular}
\end{center}
In both cases, it is clear that we have
$\deg(r_{i+1}) \leqslant \deg(r_i)$,
so that $\deg(r_1) \geqslant \deg(r_2) \geqslant \cdots
\geqslant \deg(r_k)$. It follows that no two $(\ell, r)$ pairs
in the sequence can be the same, because for each subsequence
$u_a, u_{a+1}, \hdots, u_b$ such that
$\deg(r_a) = \deg(r_{a+1}) = \cdots = \deg(r_b)$,
we must have $\deg(\ell_a) < \deg(\ell_{a+1})
< \cdots < \deg(\ell_b)$.
\end{pf}

Having shown that the left overlap division is
continuous, one way of showing that every Locally
Involutive Basis with respect to the left overlap 
division is an Involutive Basis would be to show that
the left overlap division is a
strong involutive division. However, the left
overlap division is only a weak involutive
division, as the following counterexample
demonstrates.

\begin{prop} \label {ov-weak-ch6}
The left overlap division is a weak involutive division.
\end{prop}
\begin{pf}
Let $U := \{yz, xy\}$ be a set of monomials
over the polynomial ring
$\mathbb{Q}\langle x, y, z\rangle$. Here are the
multiplicative variables for $U$ with respect to
the left overlap division $\mathcal{O}$.
\begin{center}
\begin{tabular}{c|c|c}
$u$ & $\mathcal{M}_{\mathcal{O}}^L(u, U)$
& $\mathcal{M}_{\mathcal{O}}^R(u, U)$ \\ \hline
$yz$ & $\{x, y, z\}$ & $\{x, y, z\}$ \\
$xy$ & $\{x, y, z\}$ & $\{x, y\}$ \\ \hline
\end{tabular}
\end{center}
Because $yzxy \in \mathcal{C}_{\mathcal{O}}(yz, U)$ and
$yzxy \in \mathcal{C}_{\mathcal{O}}(xy, U)$, one of the
conditions $\mathcal{C}_{\mathcal{O}}(yz, U) \subset
\mathcal{C}_{\mathcal{O}}(xy, U)$ or 
$\mathcal{C}_{\mathcal{O}}(xy, U) \subset
\mathcal{C}_{\mathcal{O}}(yz, U)$ must be satisfied in order for
$\mathcal{O}$ to be a strong involutive division (this is the
Disjoint Cones condition of
Definition \ref{noncom-div-defn}). But neither of these
conditions can be satisfied when we consider that
$xy \notin \mathcal{C}_{\mathcal{O}}(yz, U)$ and
$yz \notin \mathcal{C}_{\mathcal{O}}(xy, U)$, so
$\mathcal{O}$ must be a weak involutive division.
\end{pf}

The weakness of the left overlap division is the
price we pay for refining the left division by
allowing more right multiplicative variables.
All is not lost however, as we can still show
that every Locally Involutive Basis with respect
to the left overlap division is an Involutive Basis
by showing that the left overlap division is a
Gr\"obner involutive division.

\begin{prop} \label{ov-grob-ch6}
The left overlap division $\mathcal{O}$ is a
Gr\"obner involutive division.
\end{prop}
\begin{pf}
We are required to show that if Algorithm
\ref{noncom-inv} terminates with $\mathcal{O}$ and some
arbitrary admissible monomial ordering $O$ as input,
then the Locally Involutive Basis $G$ it returns is a
noncommutative Gr\"obner Basis. By Definition \ref{grob-defn-noncom},
we can do this by showing that all S-polynomials involving
elements of $G$ conventionally reduce to zero using $G$.

Assume that $G = \{g_1, \hdots, g_p\}$ is sorted (by lead
monomial) with respect to
the DegRevLex monomial ordering (greatest first),
and let $U = \{u_1, \hdots, u_p\}
:= \{\LM(g_1), \hdots, \LM(g_p)\}$ be the set of leading monomials.
Let $T$ be the
table obtained by applying Algorithm \ref{inv-div} to $U$.
Because $G$ is a Locally Involutive Basis, every zero
entry $T(u_i, x_j^{\Gamma})$ ($\Gamma \in \{L, R\}$) in the table
corresponds to a prolongation $g_ix_j$ or $x_jg_i$ that involutively
reduces to zero.

Let $S$ be the set of S-polynomials involving elements of $G$, where
the $t$-th entry of $S$ ($1 \leqslant t \leqslant |S|$)
is the S-polynomial
$$s_t = c_{t}\ell_{t}g_ir_{t} - c'_{t}\ell'_{t}g_jr'_{t},$$
with $\ell_{t}u_ir_{t} = \ell'_{t}u_jr'_{t}$ being the
overlap word of the S-polynomial.
We will prove that every S-polynomial in $S$
conventionally reduces to zero using $G$.

Recall (from Definition \ref{ov-def})
that each S-polynomial in $S$ corresponds to a particular
type of overlap --- `prefix', `subword' or `suffix'. For the
purposes of this proof, let us now split the subword
overlaps into three further types  --- `left', `middle' and `right',
corresponding to the cases where a monomial $m_2$ is
a prefix, proper subword and suffix of a monomial $m_1$.
\begin{center}
\begin{tabular}{ccc}
Left & Middle & Right \\[1mm]
$\xymatrix @R=0.5pc{
\ar@{<->}[rrr]^*+{m_1} &&& \\
\ar@{<->}[rr]_*+{m_2} &&
}$
&
$\xymatrix @R=0.5pc{
\ar@{<->}[rrr]^*+{m_1} &&& \\
& \ar@{<->}[r]_*+{m_2} &
}$
&
$\xymatrix @R=0.5pc{
\ar@{<->}[rrr]^*+{m_1} &&& \\
& \ar@{<->}[rr]_*+{m_2} &&
}$
\end{tabular}
\end{center}
This classification provides
us with five cases to deal with in total,
which we shall process in the following order:
right, middle, left, prefix, suffix.

{\bf (1)} Consider an arbitrary entry $s_{t} \in S$
($1 \leqslant t \leqslant |S|$)
corresponding to a right overlap where the monomial
$u_j$ is a suffix of the monomial $u_i$.
Because $\mathcal{O}$ never assigns any left
nonmultiplicative variables, $u_j$ must be an involutive
divisor of $u_i$. But this contradicts
the fact that the set $G$ is autoreduced; it follows
that no S-polynomials corresponding to right overlaps
can appear in $S$.

{\bf (2)} Consider an arbitrary entry $s_{t} \in S$
($1 \leqslant t \leqslant |S|$)
corresponding to a middle overlap where the monomial
$u_j$ is a proper subword of the monomial $u_i$. This means that
$s_{t} = c_tg_i - c'_t\ell'_tg_jr'_t$ for some $g_i, g_j \in G$,
with overlap word $u_i = \ell'_tu_jr'_t$. 
Let $u_i = x_{i_1}\hdots x_{i_{\alpha}}$;
let $u_j = x_{j_1}\hdots x_{j_{\beta}}$; and choose $D$ such that
$x_{i_D} = x_{j_{\beta}}$.
$$\xymatrix @R=0.5pc @C=2.4pc {
u_i = & \ar@{-}[r]_{x_{i_1}} & \ar@{--}[r]
& \ar@{-}[r]_{x_{i_{D-\beta}}} & \ar@{-}[r]_{x_{i_{D-\beta+1}}}
& \ar@{-}[r]_{x_{i_{D-\beta+2}}} & \ar@{--}[r]
& \ar@{-}[r]_{x_{i_{D-1}}} 
& \ar@{-}[r]_{x_{i_{D}}} & \ar@{-}[r]_{x_{i_{D+1}}} &
\ar@{--}[r] & \ar@{-}[r]_{x_{i_{\alpha}}} & \\
u_j = & &&& \ar@{-}[r]_{x_{j_1}} & \ar@{-}[r]_{x_{j_2}} & \ar@{--}[r]
& \ar@{-}[r]_{x_{j_{\beta - 1}}} & \ar@{-}[r]_{x_{j_{\beta}}} &
}$$
Because $u_j$ is a proper subword of $u_i$, it follows that
$T(u_j, x^R_{i_{D+1}}) = 0$.
This gives rise to the prolongation
$g_jx_{i_{D+1}}$ of $g_j$. But we know that all prolongations
involutively reduce to zero ($G$ is a Locally Involutive
Basis), so Algorithm \ref{noncom-inv-div} must find a monomial
$u_{k} = x_{k_1}\hdots x_{k_{\gamma}} \in U$ such that
$u_{k}$ involutively divides $u_jx_{i_{D+1}}$.
Assuming that $x_{k_{\gamma}} = x_{i_{\kappa}}$,
we can deduce that
any candidate for $u_k$ must be a suffix of $u_jx_{i_{D+1}}$
(otherwise $T(u_k, x^R_{i_{\kappa + 1}}) = 0$ because of the
overlap between $u_i$ and $u_k$). This means that the degree of
$u_k$ is in the range
$1 \leqslant \gamma \leqslant \beta+1$; we shall illustrate
this in the following diagram by using a squiggly line to
indicate that the monomial $u_k$ can begin anywhere
(or nowhere if $u_k = x_{i_{D+1}}$) on the squiggly line.
$$\xymatrix @R=0.5pc @C=2.4pc {
u_i = & \ar@{-}[r]_{x_{i_1}} & \ar@{--}[r]
& \ar@{-}[r]_{x_{i_{D-\beta}}} & \ar@{-}[r]_{x_{i_{D-\beta+1}}}
& \ar@{-}[r]_{x_{i_{D-\beta+2}}} & \ar@{--}[r]
& \ar@{-}[r]_{x_{i_{D-1}}}
& \ar@{-}[r]_{x_{i_{D}}} & \ar@{-}[r]_{x_{i_{D+1}}} &
\ar@{--}[r] & \ar@{-}[r]_{x_{i_{\alpha}}} & \\
u_j = & &&& \ar@{-}[r]_{x_{j_1}} & \ar@{-}[r]_{x_{j_2}} & \ar@{--}[r]
& \ar@{-}[r]_{x_{j_{\beta - 1}}} & \ar@{-}[r]_{x_{j_{\beta}}} & \\
u_k = & &&& \ar@{~}[rrrrr] &&&&& \ar@{-}[r]_{x_{k_{\gamma}}} &
}$$
We can now use the monomial $u_k$ 
together with Buchberger's Second Criterion
to simplify our goal of showing that the S-polynomial $s_t$ reduces
to zero. Notice that the monomial $u_k$
is a subword of the overlap word $u_i$ associated to $s_t$,
and so in order to show that $s_t$ reduces to zero,
all we have to do is to show that the two S-polynomials
$$s_u = c_ug_i - c'_{u}(x_{i_1}x_{i_2}\hdots
x_{i_{D+1-\gamma}})g_k(x_{i_{D+2}}\hdots x_{i_{\alpha}})$$
and\footnote{Technical point: if $\gamma \neq \beta+1$,
the S-polynomial $s_v$ could in fact appear as
$s_v = c_{v}g_jx_{i_{D+1}} - c'_v(x_{j_1}\hdots x_{i_{D+1-\gamma}})g_k$
and not as
$s_v = c_v(x_{j_1}\hdots x_{i_{D+1-\gamma}})g_k - c'_{v}g_jx_{i_{D+1}}$;
for simplicity we will treat both cases the same in the
proof as all that changes is the notation and the signs.}
$$s_v = c_v(x_{j_1}\hdots x_{i_{D+1-\gamma}})g_k - c'_{v}g_jx_{i_{D+1}}$$
reduce to zero ($1 \leqslant u,v \leqslant |S|$).

For the S-polynomial $s_v$, there are two cases to consider:
$\gamma = 1$, and $\gamma > 1$. In the former case,
because (as placed in $u_i$) the monomials
$u_j$ and $u_k$ do not overlap, we can use Buchberger's
First Criterion to say that the `S-polynomial' $s_v$ reduces to
zero (for further explanation, see the paragraph at the beginning of
Section \ref{BCNC}). In the latter case, we know that
the first step of the involutive reduction of the
prolongation $g_jx_{i_{D+1}}$ is to take away the multiple
$(\frac{c_v}{c'_v})(x_{j_1}\hdots x_{i_{D+1-\gamma}})g_k$ of $g_k$
from $g_jx_{i_{D+1}}$ to leave the polynomial $g_jx_{i_{D+1}}
- (\frac{c_v}{c'_v})(x_{j_1}\hdots x_{i_{D+1-\gamma}})g_k
= -(\frac{1}{c'_v})s_v$.
But as we know that all prolongations involutively reduce
to zero, we can conclude that the S-polynomial
$s_v$ conventionally reduces to zero.

For the S-polynomial $s_u$, we note that if $D = \alpha-1$, then
$s_u$ corresponds to a right overlap. But we know from part (1)
that right overlaps cannot appear in $S$, and so $s_t$ also
cannot appear in $S$. Otherwise, we proceed by
induction on the S-polynomial $s_u$ to produce a sequence
$\{u_{q_{D+1}}, u_{q_{D+2}}, \hdots, u_{q_{\alpha}}\}$ of monomials, so that
$s_u$ (and hence $s_t$) reduces to zero if the S-polynomial
$$s_{\eta} = c_{\eta}g_i - c'_{\eta}(x_{i_1}\hdots
x_{i_{\alpha-\mu}})g_{q_{\alpha}}$$
reduces to zero ($1 \leqslant \eta \leqslant |S|$),
where $\mu = \deg(u_{q_{\alpha}})$.
$$\xymatrix @R=0.5pc @C=2.0pc {
u_i = & \ar@{-}[r]_{x_{i_1}} & \ar@{--}[r]
& \ar@{-}[r]_{x_{i_{D-\beta}}} & \ar@{-}[r]_{x_{i_{D-\beta+1}}}
& \ar@{--}[r] & \ar@{-}[r]_{x_{i_{D}}} & \ar@{-}[r]_{x_{i_{D+1}}} &
\ar@{-}[r]_{x_{i_{D+2}}} & \ar@{--}[r] & \ar@{-}[r]_{x_{i_{\alpha-1}}} &
\ar@{-}[r]_{x_{i_{\alpha}}} & \\
u_j = & &&& \ar@{-}[r]_{x_{j_1}} & \ar@{--}[r]
& \ar@{-}[r]_{x_{j_{\beta}}} & \\
u_{q_{D+1}} = u_k = & &&& \ar@{~}[rrr] &&& \ar@{-}[r]_{x_{k_{\gamma}}} & \\
u_{q_{D+2}} = & &&& \ar@{~}[rrrr] &&&& \ar@{-}[r] & \\
& &&&&&&& \ddots \\
u_{q_{\alpha}} = & &&& \ar@{~}[rrrrrrr] &&&&&&& \ar@{-}[r] &
}$$
But $s_{\eta}$ always corresponds to a right overlap,
so we must conclude that middle overlaps (as well
as right overlaps) cannot appear in $S$.

{\bf (3)} Consider an arbitrary entry $s_{t} \in S$
($1 \leqslant t \leqslant |S|$)
corresponding to a left overlap where the monomial
$u_j$ is a prefix of the monomial $u_i$. This means that
$s_{t} = c_tg_i - c'_tg_jr'_t$ for some $g_i, g_j \in G$,
with overlap word $u_i = u_jr'_t$. Let $u_i = x_{i_1}\hdots x_{i_{\alpha}}$
and let $u_j = x_{j_1}\hdots x_{j_{\beta}}$.
$$\xymatrix @R=0.5pc {
u_i = & \ar@{-}[r]_{x_{i_1}} & \ar@{-}[r]_{x_{i_2}} & \ar@{--}[r] &
\ar@{-}[r]_{x_{i_{\beta-1}}} & \ar@{-}[r]_{x_{i_{\beta}}} &
\ar@{-}[r]_{x_{i_{\beta+1}}} & \ar@{--}[r] &
\ar@{-}[r]_{x_{i_{\alpha-1}}} & \ar@{-}[r]_{x_{i_{\alpha}}} & \\
u_j = & \ar@{-}[r]_{x_{j_1}} & \ar@{-}[r]_{x_{j_2}} & \ar@{--}[r] &
\ar@{-}[r]_{x_{j_{\beta-1}}} & \ar@{-}[r]_{x_{j_{\beta}}} &
}$$
Because $u_j$ is a prefix of $u_i$, it follows that
$T(u_j, x^R_{i_{\beta+1}}) = 0$.
This gives rise to the prolongation
$g_jx_{i_{\beta+1}}$ of $g_j$. But we know that all prolongations
involutively reduce to zero, so there must exist a monomial
$u_{k} = x_{k_1}\hdots x_{k_{\gamma}} \in U$ such that
$u_{k}$ involutively divides $u_jx_{i_{\beta+1}}$.
Assuming that $x_{k_{\gamma}} = x_{i_{\kappa}}$,
any candidate for $u_k$ must be a suffix of $u_jx_{i_{\beta+1}}$
(otherwise $T(u_k, x^R_{i_{\kappa + 1}}) = 0$ because of the
overlap between $u_i$ and $u_k$). Further, any candidate for $u_k$
cannot be either a suffix or a proper subword of $u_i$ (because
of parts (1) and (2) of this proof). This leaves only
one possibility for $u_k$, namely $u_k = u_jx_{i_{\beta+1}}$.
$$\xymatrix @R=0.5pc {
u_i = & \ar@{-}[r]_{x_{i_1}} & \ar@{-}[r]_{x_{i_2}} & \ar@{--}[r] &
\ar@{-}[r]_{x_{i_{\beta-1}}} & \ar@{-}[r]_{x_{i_{\beta}}} &
\ar@{-}[r]_{x_{i_{\beta+1}}} & \ar@{--}[r] &
\ar@{-}[r]_{x_{i_{\alpha-1}}} & \ar@{-}[r]_{x_{i_{\alpha}}} & \\
u_j = & \ar@{-}[r]_{x_{j_1}} & \ar@{-}[r]_{x_{j_2}} & \ar@{--}[r] &
\ar@{-}[r]_{x_{j_{\beta-1}}} & \ar@{-}[r]_{x_{j_{\beta}}} & \\
u_k = & \ar@{-}[r]_{x_{k_1}} & \ar@{-}[r]_{x_{k_2}} & \ar@{--}[r] &
\ar@{-}[r]_{x_{k_{\gamma-2}}} & \ar@{-}[r]_{x_{k_{\gamma-1}}} &
\ar@{-}[r]_{x_{k_{\gamma}}} &
}$$
If $\alpha = \beta+1$, then it is clear that $u_k = u_i$,
and so the first step in the involutive reduction of the
prolongation $g_jx_{i_{\alpha}}$
is to take away the multiple
$(\frac{c_t}{c'_t})g_i$ of $g_i$
from $g_jx_{i_{\alpha}}$ to leave the polynomial $g_jx_{i_{\alpha}}
- (\frac{c_t}{c'_t})g_i = -(\frac{1}{c'_t})s_t$.
But as we know that all prolongations involutively reduce
to zero, we can conclude that the S-polynomial
$s_t$ conventionally reduces to zero.

Otherwise, if $\alpha > \beta+1$,
we can now use the monomial $u_k$ together with Buchberger's Second Criterion
to simplify our goal of showing that the S-polynomial $s_t$ reduces
to zero. Notice that the monomial $u_k$
is a subword of the overlap word $u_i$ associated to $s_t$,
and so in order to show that $s_t$ reduces to zero,
all we have to do is to show that the two S-polynomials
$$s_u = c_ug_i - c'_{u}g_k(x_{i_{\beta+2}}\hdots x_{i_{\alpha}})$$
and
$$s_v = c_vg_k - c'_{v}g_jx_{i_{\beta+1}}$$
reduce to zero ($1 \leqslant u,v \leqslant |S|$).

The S-polynomial $s_v$ reduces to zero by comparison with
part (2). For the S-polynomial $s_u$, we proceed by
induction (we have another left overlap), eventually
coming across a left overlap of `type
$\alpha = \beta+1$' because we move one letter at a time
to the right after each inductive step.
$$\xymatrix @R=0.5pc {
u_i = & \ar@{-}[r]_{x_{i_1}} & \ar@{-}[r]_{x_{i_2}} & \ar@{--}[r] &
\ar@{-}[r]_{x_{i_{\beta-1}}} & \ar@{-}[r]_{x_{i_{\beta}}} &
\ar@{-}[r]_{x_{i_{\beta+1}}} & \ar@{-}[r]_{x_{i_{\beta+2}}} & \ar@{--}[r] &
\ar@{-}[r]_{x_{i_{\alpha-1}}} & \ar@{-}[r]_{x_{i_{\alpha}}} & \\
u_j = & \ar@{-}[r]_{x_{j_1}} & \ar@{-}[r]_{x_{j_2}} & \ar@{--}[r] &
\ar@{-}[r]_{x_{j_{\beta-1}}} & \ar@{-}[r]_{x_{j_{\beta}}} & \\
u_k = & \ar@{-}[r]_{x_{k_1}} & \ar@{-}[r]_{x_{k_2}} & \ar@{--}[r] &
        \ar@{-}[r]_{x_{k_{\gamma-2}}} & \ar@{-}[r]_{x_{k_{\gamma-1}}} &
        \ar@{-}[r]_{x_{k_{\gamma}}} & \\
      & \ar@{-}[r] & \ar@{-}[r] & \ar@{--}[r]& \ar@{-}[r]& \ar@{-}[r]
      & \ar@{-}[r] & \ar@{-}[r] & \\
      &&&&&& \ddots \\
      & \ar@{-}[r] & \ar@{-}[r] & \ar@{--}[r] & \ar@{-}[r] &
      \ar@{-}[r] & \ar@{-}[r]& \ar@{-}[r] & \ar@{--}[r] & \ar@{-}[r] & \\
}$$

{\bf (4 and 5)} In Definition \ref{ov-def}, we defined a prefix
overlap to be an overlap where, given two monomials $m_1$ and
$m_2$ such that $\deg(m_1) \geqslant \deg(m_2)$, a prefix of
$m_1$ is equal to a suffix of $m_2$; suffix overlaps were defined
similarly. If we drop the condition on the degrees
of the monomials, it is clear that every suffix
overlap can be treated as a prefix overlap (by swapping the
roles of $m_1$ and $m_2$); this allows us to deal with the
case of a prefix overlap only.

Consider an arbitrary entry $s_{t} \in S$
($1 \leqslant t \leqslant |S|$)
corresponding to a prefix overlap where a prefix of the monomial
$u_i$ is equal to a suffix of the monomial $u_j$. This means that
$s_{t} = c_t\ell_tg_i - c'_tg_jr'_t$ for some $g_i, g_j \in G$,
with overlap word $\ell_tu_i = u_jr'_t$.
Let $u_i = x_{i_1}\hdots x_{i_{\alpha}}$;
let $u_j = x_{j_1}\hdots x_{j_{\beta}}$; and choose $D$ such that
$x_{i_D} = x_{j_{\beta}}$.
$$\xymatrix @R=0.5pc @C=2.4pc{
u_i = &&&& \ar@{-}[r]_{x_{i_1}} & \ar@{--}[r]
& \ar@{-}[r]_{x_{i_{D}}}
& \ar@{-}[r]_{x_{i_{D+1}}} & \ar@{--}[r]
& \ar@{-}[r]_{x_{i_{\alpha - 1}}} & \ar@{-}[r]_{x_{i_{\alpha}}} & \\
u_j = & \ar@{-}[r]_{x_{j_1}} & \ar@{--}[r] & \ar@{-}[r]_{x_{j_{\beta-D}}}
& \ar@{-}[r]_{x_{j_{\beta-D+1}}} & \ar@{--}[r]
& \ar@{-}[r]_{x_{j_{\beta}}} &
}$$
By definition of $\mathcal{O}$, we must have
$T(u_j, x^R_{i_{D+1}}) = 0$.

Because we know that the prolongation
$g_jx_{i_{D+1}}$ involutively reduces to zero, there must exist a monomial
$u_{k} = x_{k_1}\hdots x_{k_{\gamma}} \in U$ such that
$u_{k}$ involutively divides $u_jx_{i_{D+1}}$.
This $u_k$ must be a suffix of $u_jx_{i_{D+1}}$
(otherwise, assuming that $x_{k_{\gamma}} = x_{j_{\kappa}}$,
we have $T(u_k, x^R_{i_{D+1}}) = 0$ if 
$\kappa = \beta$ (because of the overlap between
$u_i$ and $u_k$); and $T(u_k, x^R_{j_{\kappa + 1}}) = 0$
if $\kappa < \beta$ (because of the overlap between
$u_j$ and $u_k$)).
$$\xymatrix @R=0.5pc @C=2.4pc{
u_i = &&&& \ar@{-}[r]_{x_{i_1}} & \ar@{--}[r]
& \ar@{-}[r]_{x_{i_{D}}}
& \ar@{-}[r]_{x_{i_{D+1}}} & \ar@{--}[r]
& \ar@{-}[r]_{x_{i_{\alpha - 1}}} & \ar@{-}[r]_{x_{i_{\alpha}}} & \\
u_j = & \ar@{-}[r]_{x_{j_1}} & \ar@{--}[r] & \ar@{-}[r]_{x_{j_{\beta-D}}}
& \ar@{-}[r]_{x_{j_{\beta-D+1}}} & \ar@{--}[r]
& \ar@{-}[r]_{x_{j_{\beta}}} & \\
u_k = & \ar@{~}[rrrrrr] &&&&&& \ar@{-}[r]_{x_{k_{\gamma}}} &
}$$
Let us now
use the monomial $u_k$ together with Buchberger's Second Criterion
to simplify our goal of showing that the S-polynomial $s_t$ reduces
to zero. Because $u_k$ is a subword of the overlap word $\ell_tu_i$
associated to $s_t$, in order to show that $s_t$ reduces to zero,
all we have to do is to show that the two S-polynomials
$$s_u =
\begin{cases}
c_{u}(x_{k_{1}}\hdots x_{j_{\beta-D}})g_i
- c'_{u}g_k(x_{i_{D+2}}\hdots x_{i_{\alpha}})
& \text{if $\gamma > D+1$} \\
c_{u}g_i - c'_{u}\ell'_ug_k(x_{i_{D+2}}\hdots x_{i_{\alpha}})
& \text{if $\gamma \leqslant D+1$}
\end{cases}
$$
and
$$s_v = c_{v}g_jx_{i_{D+1}} 
- c'_{v}(x_{j_1}\hdots x_{j_{\beta+1-\gamma}})g_k$$
reduce to zero ($1 \leqslant u,v \leqslant |S|$).

The S-polynomial $s_v$ reduces to zero by comparison with
part (2). For the S-polynomial $s_u$, first note that if $\alpha = D+1$,
then either $u_k$ is a suffix of $u_i$,
$u_i$ is a suffix of $u_k$, or $u_k = u_i$; it follows that
$s_u$ reduces to zero trivially if $u_k = u_i$, and (by
part (1)) $s_u$ (and hence $s_t$) cannot appear in $S$
in the other two cases.

If however $\alpha \neq D+1$,
then either $s_u$ is a middle overlap (if $\gamma < D+1$), a left overlap
(if $\gamma = D+1$), or another prefix overlap.
The first case leads us to conclude that $s_t$ cannot
appear in $S$; the second case is handled by part (3)
of this proof; and the final case is handled by induction,
where we note that after each step of the induction, the value
$\alpha+\beta-2D$ strictly decreases, so we are
guaranteed at some stage to find an overlap that is not
a prefix overlap, enabling us either to verify that the
S-polynomial $s_t$ conventionally reduces to zero, or
to conclude that $s_t$ can not in fact appear in $S$.
\end{pf}

\subsection{A Strong Local Division}

Thus far, we have encountered two global divisions
that are strong and continuous, and one local
division that is weak, continuous and Gr\"obner.
Our next division can be considered to be a hybrid
of these previous divisions, as it will be
a local division that is continuous and (as long as
thick divisors are being used) strong.

\begin{defn}[The Strong Left Overlap Division $\mathcal{S}$]
\label{SODiv}
\index{$S$@$\mathcal{S}$}
\index{involutive division!strong left overlap}
Let $U = \{u_1, \hdots, u_m\}$ be a set of monomials.
Assign multiplicative variables to $U$ according to
Algorithm \ref{inv-div2}, which (in words) performs the
following two tasks.
\index{strong left overlap division}
\begin{enumerate}[(a)]
\item
Assign multiplicative variables to $U$ according to
the left overlap division.
\item
Using the recipe provided in Algorithm \ref{DJ},
ensure that at least one variable in every monomial 
$u_j \in U$ is right nonmultiplicative for each
monomial $u_i \in U$.
\end{enumerate}
\end{defn}

\begin{remark}
As Algorithm \ref{inv-div2} expects any input set
to be ordered with respect to DegRevLex, we may
sometimes have to reorder a set of monomials $U$ 
to satisfy this condition before we can assign
multiplicative variables to $U$ according to the
strong left overlap division.
\end{remark}

\begin{algorithm}
\setlength{\baselineskip}{3.5ex}
\caption{`DisjointCones' Function for Algorithm 15}
\label{DJ}
\begin{algorithmic}
\vspace*{2mm}
\REQUIRE{A set of % distinct 
         monomials $U = \{u_1, \hdots, u_m\}$
         ordered by DegRevLex ($u_1 \geqslant u_2 \geqslant
         \cdots \geqslant u_m$),
         where $u_i \in R\langle x_1, \hdots, x_n\rangle$;
         a table $T$ of left and right multiplicative variables
         for all $u_i \in U$, where each entry of $T$ is either 1
         (multiplicative) or 0 (nonmultiplicative).}
\ENSURE{$T$.}
\vspace*{1mm}
\FOR{{\bf each} monomial $u_i \in U$ ($m \geqslant i \geqslant 1$)}
\FOR{{\bf each} monomial $u_j \in U$ ($m \geqslant j \geqslant 1$)}
\STATE
Let $u_i = x_{i_1}x_{i_2}\hdots x_{i_{\alpha}}$ and
$u_j = x_{j_1}x_{j_2}\hdots x_{j_{\beta}}$; \\
% foundleft = false;
found = false; \\
$k = 1$; \\
\WHILE{($k \leqslant \beta$)}
% \FOR{each $k$ ($1 \leqslant k \leqslant \beta$)}
%    {\bf if} ($T(u_i, x_{j_k}^L) = 0$) {\bf then}
%    foundleft = true; \\
%    {\bf if} ($T(u_i, x_{j_k}^R) = 0$) {\bf then}
\IF{($T(u_i, x_{j_k}^R) = 0$)}
\STATE
found = true; \\
$k = \beta+1$;
\ELSE
\STATE
$k = k+1$;
\ENDIF
\ENDWHILE
\IF{(found == false)}
\STATE
%  {\bf if} (foundleft == false) {\bf then}
%  $T(u_i, x_{j_{\beta}}^L) = 0$; \\
%    {\bf if} (found == false) {\bf then}
$T(u_i, x_{j_{1}}^R) = 0$; \\
%  {\bf end if}
\ENDIF
\ENDFOR
\ENDFOR \vspace*{1mm}
\STATE
{\bf return} $T$;
\end{algorithmic}
\vspace*{1mm}
\end{algorithm}

\begin{algorithm}
\setlength{\baselineskip}{3.5ex}
% \index{$S$@$\mathcal{S}$}
\caption{The Strong Left Overlap Division $\mathcal{S}$}
\label{inv-div2}
\begin{algorithmic}
\vspace*{2mm}
% p -> m, m_i -> u_i
\REQUIRE{A set of % distinct 
         monomials $U = \{u_1, \hdots, u_m\}$
         ordered by DegRevLex ($u_1 \geqslant u_2 \geqslant
         \cdots \geqslant u_m$),
         where $u_i \in R\langle x_1, \hdots, x_n\rangle$.}
\ENSURE{A table $T$ of left and right multiplicative variables
        for all $u_i \in U$, where each entry of $T$ is either 1
        (multiplicative) or 0 (nonmultiplicative).}
\vspace*{1mm}
\STATE
Create a table $T$ of multiplicative 
variables as shown below: \\ \vspace*{1mm}
\begin{tabular}{c|ccccccc}
& $x_{1}^L$ & $x_{1}^R$ & $x_{2}^L$ & $x_{2}^R$
& $\cdots$ & $x_{n}^L$ & $x_{n}^R$ \\
\hline
$u_1$ & 1 & 1 & 1 & 1 & $\cdots$ & 1 & 1 \\
$u_2$ & 1 & 1 & 1 & 1 & $\cdots$ & 1 & 1 \\
$\vdots$ & $\vdots$ & $\vdots$ & $\vdots$ & $\vdots$
& $\ddots$ & $\vdots$ & $\vdots$ \\
$u_m$ & 1 & 1 & 1 & 1 & $\cdots$ & 1 & 1
\end{tabular} \\[2mm]
\FOR{{\bf each} monomial $u_i \in U$ ($1 \leqslant i \leqslant m$)}
\FOR{{\bf each} monomial $u_j \in U$ ($i \leqslant j \leqslant m$)}
\STATE
Let $u_i = x_{i_1}x_{i_2}\hdots x_{i_{\alpha}}$ and
$u_j = x_{j_1}x_{j_2}\hdots x_{j_{\beta}}$; \\[1mm]
%// Middle Overlaps \\
\IF{($i \neq j$)}
\FOR{{\bf each} $k$ ($1 \leqslant k < \alpha - \beta + 1$)}
  \IF{($\SUB(u_i, k, k+\beta-1) == u_j$)}
  \STATE
%    {\bf if} ($k < \alpha-\beta+1$) {\bf then}
    $T(u_j, x^R_{i_{k+\beta}}) = 0$; \\
%    {\bf end if}
  \ENDIF
\ENDFOR
\ENDIF
%\STATE
\vspace*{1mm}
%// Left \& Right Overlaps
\FOR{{\bf each} $k$ ($1 \leqslant k \leqslant \beta - 1$)}
 \IF{($\PRE(u_i, k) == \SUFF(u_j, k)$)}
 \STATE
 $T(u_j, x^R_{i_{k+1}}) = 0$; 
 \ENDIF
 \IF{($\SUFF(u_i, k) == \PRE(u_j, k)$)}
 \STATE
 $T(u_i, x^R_{j_{k+1}}) = 0$; 
 \ENDIF
%  \STATE
%  {\bf if} ($\PRE(u_i, k) == \SUFF(u_j, k)$) {\bf then}
%  $T(u_j, x^R_{i_{k+1}}) = 0$; \\
%  {\bf end if} \\
%  {\bf if} ($\SUFF(u_i, k) == \PRE(u_j, k)$) {\bf then}
%  $T(u_i, x^R_{j_{k+1}}) = 0$; \\
%  {\bf end if}
\ENDFOR \vspace*{1mm}
\ENDFOR
\ENDFOR
\STATE
$T$ = DisjointCones$( U, T )$; \; {\it (Algorithm \ref{DJ})} \\
{\bf return} $T$;
\end{algorithmic}
\vspace*{1mm}
\end{algorithm}

\newpage
\begin{prop} \label{str-cts-ch6}
The strong left overlap division is continuous.
\end{prop}
\begin{pf}
We refer to the proof of Proposition \ref{ov-cts-ch6},
replacing $\mathcal{O}$ by $\mathcal{S}$.
\end{pf}

\begin{prop} \label{str-grob-ch6}
The strong left overlap division is a Gr\"obner
involutive division.
\end{prop}
\begin{pf}
We refer to the proof of Proposition \ref{ov-grob-ch6},
replacing $\mathcal{O}$ by $\mathcal{S}$.
\end{pf}

\begin{remark}
Propositions \ref{str-cts-ch6} and \ref{str-grob-ch6} apply
either when using thin divisors or when using thick divisors.
\end{remark}

\begin{prop}
With respect to thick divisors,
the strong left overlap division is a strong
involutive division.
\end{prop}
\begin{pf}
To prove that the strong left overlap division is a strong
involutive division, we need to show that
the three conditions of Definition
\ref{noncom-div-defn} hold. 
\begin{itemize}
\item
{\bf Disjoint Cones Condition} \\
Let $\mathcal{C}_{\mathcal{S}}(u_1, U)$ and
$\mathcal{C}_{\mathcal{S}}(u_2, U)$ be the 
involutive cones associated to
the monomials $u_1$ and $u_2$ over some noncommutative
polynomial ring $\mathcal{R}$, where
$\{u_1, u_2\} \subset U \subset \mathcal{R}$. If
$\mathcal{C}_{\mathcal{S}}(u_1, U) \cap
\mathcal{C}_{\mathcal{S}}(u_2, U) \neq \emptyset$, 
then there must
be some monomial $v \in \mathcal{R}$ such that
$v$ contains both monomials $u_1$ and $u_2$
as subwords, and (as placed in $v$) both
$u_1$ and $u_2$ must be involutive divisors of
$v$. By definition of $\mathcal{S}$, % and as long
% as thick divisors are being used\footnote{
% If thin divisors are being used, the monomials
% $u_1$ and $u_2$ may not overlap as placed in $v$,
% allowing us to prove that $\mathcal{S}$ is not
% a strong involutive division.},
both $u_1$ and $u_2$ must be suffices of $v$.
Thus, assuming (without loss of generality) that
$\deg(u_1) > \deg(u_2)$, we are able to draw
the following diagram summarising the situation.
$$
\xymatrix @R=1.5pc{
\ar@{<->}[rrrrrr]^*+{v} &&&&&& \\
&& \ar@{<->}[rrrr]^*+{u_1} &&&& \\
&&&& \ar@{<->}[rr]^*+{u_2} &&
}
$$
For $\mathcal{S}$ to be strong, we must have
$\mathcal{C}_{\mathcal{S}}(u_1, U) \subset 
\mathcal{C}_{\mathcal{S}}(u_2, U)$ (it is
clear that $\mathcal{C}_{\mathcal{S}}(u_2, U) \not\subset 
\mathcal{C}_{\mathcal{S}}(u_1, U)$ because
$u_2 \notin \mathcal{C}_{\mathcal{S}}(u_1, U)$).
This can be verified by proving that a variable
is right nonmultiplicative for $u_1$ if and
only if it is right nonmultiplicative for $u_2$.

($\Rightarrow$) If an arbitrary variable $x$ is
right nonmultiplicative for $u_2$, then
either some monomial $u \in U$ overlaps with
$u_2$ in one of the ways shown below (where
the variable immediately to the right of $u_2$ is
the variable $x$), or
$x$ was assigned right nonmultiplicative for $u_2$
in order to ensure that some variable in some monomial
$u \in U$ is right nonmultiplicative for $u_2$.
\begin{center}
\begin{tabular}{ccc}
Overlap (i) & \vspace*{5mm} & Overlap (ii) \\
$\xymatrix @R=1.5pc{
\ar@{<->}[rrrr]^*+{u} &&&& \\
& \ar@{<->}[rr]^*+{u_2} &&
}
$
&&
$\xymatrix @R=1.5pc{
&& \ar@{<->}[rrrr]^*+{u} &&&& \\
\ar@{<->}[rrrr]^*+{u_2} &&&&
}
$
\end{tabular}
\end{center}
If the former case applies, then it is clear that
for both overlap types there will be another
overlap between $u_1$ and $u$ that will lead
$\mathcal{S}$ to assign $x$ to be right
nonmultiplicative for $u_1$. It follows that after
we have assigned multiplicative variables
to $U$ according to the left overlap division
(which we recall is the first step of assigning
multiplicative variables to $U$ according to $\mathcal{S}$),
the right multiplicative variables of $u_1$
and $u_2$ will be identical.
It therefore remains to show that
if $x$ is assigned right nonmultiplicative for $u_2$
in the latter case (which will happen during the final
step of assigning multiplicative variables to $U$
according to $\mathcal{S}$), then $x$ is also assigned
right nonmultiplicative for $u_1$. But this is clear
when we consider that
Algorithm \ref{DJ} is used to perform this final
step, because for $u_1$ and $u_2$ in Algorithm \ref{DJ}, we will
always analyse each monomial in $U$ in the same order.

($\Leftarrow$) Use the same argument as above,
replacing $u_1$ by $u_2$ and vice-versa.
\item
{\bf Unique Divisor Condition} \\
Given a monomial $u$ belonging to a set of monomials $U$,
$u$ may not involutively divide an arbitrary monomial $v$
in more than one way (and hence the Unique Divisor condition
is satisfied) because (i) $\mathcal{S}$ ensures that no
overlap word involving only $u$ is involutively divisible
in more than one way by $u$; and (ii) $\mathcal{S}$ ensures
that at least one variable in $u$ is right nonmultiplicative
for $u$, so that if $u$ appears twice in $v$ as subwords
that are disjoint from one another, then only the `right-most'
subword can potentially be an involutive divisor of $v$.
\item
{\bf Subset Condition} \\
Let $v$ be a monomial belonging to a set $V$
of monomials, where $V$ itself is a subset of a larger set
$U$ of monomials. Because $\mathcal{S}$ assigns
no left nonmultiplicative variables, it is clear
that $\mathcal{M}^L_{\mathcal{S}}(v, U) \subseteq
\mathcal{M}^L_{\mathcal{S}}(v, V)$. To prove that
$\mathcal{M}^R_{\mathcal{S}}(v, U) \subseteq
\mathcal{M}^R_{\mathcal{S}}(v, V)$, note that if
a variable $x$ is right nonmultiplicative for $v$
with respect to $U$ and $\mathcal{S}$
(so that $x \notin \mathcal{M}^R_{\mathcal{S}}(v, U)$),
then (as in the proof for the Disjoint Cones Condition)
either some monomial $u \in U$ overlaps with
$v$ in one of the ways shown below (where
the variable immediately to the right of $v$ is
the variable $x$), or
$x$ was assigned right nonmultiplicative for $v$
in order to ensure that some variable in some monomial
$u \in U$ is right nonmultiplicative for $v$.
\begin{center}
\begin{tabular}{ccc}
Overlap (i) & \vspace*{5mm} & Overlap (ii) \\
$\xymatrix @R=1.5pc{
\ar@{<->}[rrrr]^*+{u} &&&& \\
& \ar@{<->}[rr]^*+{v} &&
}
$
&&
$\xymatrix @R=1.5pc{
&& \ar@{<->}[rrrr]^*+{u} &&&& \\
\ar@{<->}[rrrr]^*+{v} &&&&
}
$
\end{tabular}
\end{center}
In both cases, it is clear that, with respect to the
set $V$, the variable $x$ may not
be assigned right nonmultiplicative for $v$
if $u \notin V$, so that
$\mathcal{M}^R_{\mathcal{S}}(v, U) \subseteq
\mathcal{M}^R_{\mathcal{S}}(v, V)$ as required.
\end{itemize}
\end{pf}

\begin{prop}
With respect to thin divisors,
the strong left overlap division is a weak
involutive division.
\end{prop}
\begin{pf}
Let $U := \{xy\}$ be a set of monomials
over the polynomial ring
$\mathbb{Q}\langle x, y\rangle$. Here are the
multiplicative variables for $U$ with respect to
the strong left overlap division $\mathcal{S}$.
\begin{center}
\begin{tabular}{c|c|c}
$u$ & $\mathcal{M}_{\mathcal{S}}^L(u, U)$
& $\mathcal{M}_{\mathcal{S}}^R(u, U)$ \\ \hline
$xy$ & $\{x, y\}$ & $\{y\}$ \\ \hline
\end{tabular}
\end{center}
For $\mathcal{S}$ to be strong
with respect to thin divisors, the monomial
$xy^2xy$, which is conventionally divisible by
$xy$ in two ways, must only be involutively
divisible by $xy$ in one way (this is the
Unique Divisor condition of
Definition \ref{noncom-div-defn}). 
However it is clear that $xy^2xy$ is involutively
divisible by $xy$ in two ways with respect to 
thin divisors, so $\mathcal{S}$ must be a weak
involutive division with respect to thin divisors.
\end{pf}

\begin{example} \label{fiveC}
Continuing Examples \ref{fiveA} and \ref{fiveB},
here are the multiplicative variables for the
set $\LM(F)$ of monomials with respect to the
strong left overlap division $\mathcal{S}$, where
we recall that $F := \{xy-z, \, x+z, \, yz-z, \,
xz, \, zy+z, \, z^2\}$.
\begin{center}
\begin{tabular}{c|c|c}
$u$ & $\mathcal{M}_{\mathcal{S}}^L(u, \LM(F))$ &
$\mathcal{M}_{\mathcal{S}}^R(u, \LM(F))$ \\ \hline
$xy$  & $\{x, y, z\}$  & $\{y\}$ \\
$x$   & $\{x, y, z\}$  & $\emptyset$ \\
$yz$  & $\{x, y, z\}$  & $\emptyset$ \\
$xz$  & $\{x, y, z\}$  & $\emptyset$ \\
$zy$  & $\{x, y, z\}$  & $\{y\}$ \\
$z^2$ & $\{x, y, z\}$  & $\emptyset$ \\ \hline
\end{tabular}
\end{center}
When we apply Algorithm \ref{noncom-inv} to
$F$ with respect to the DegLex monomial ordering,
thick divisors and the strong left overlap division, $F$
(as in Example \ref{fiveB}) is returned to us
as the output Locally Involutive Basis. 
\end{example}

\begin{remark}
In the above example, even though we know
that $\mathcal{S}$ is continuous,
we cannot deduce that the Locally Involutive Basis
$F$ is an Involutive Basis
because we are using thick divisors
(Proposition \ref{LocalToGlobalNC} does not apply
in the case of using thick divisors).

What this means is that the involutive cones of $F$
(and in general any Locally Involutive Basis with
respect to $\mathcal{S}$ and thick divisors)
will be disjoint (because $\mathcal{S}$ is strong),
but will not necessarily completely cover the conventional
cones of $F$, so that some monomials that are
conventionally reducible by $F$ may not be involutively
reducible by $F$. It follows that when involutively
reducing a polynomial with respect to $F$, the
reduction path will be unique but the correct remainder
may not always be obtained (in the sense that some of the
terms in our `remainder' may still be conventionally
reducible by members of $F$). One remedy to this problem 
would be to involutively reduce a polynomial $p$ with respect 
to $F$ to obtain a remainder $r$, and then to
conventionally reduce $r$ with respect to $F$ to
obtain a remainder $r'$ which we can be sure
contains no term that is conventionally reducible by $F$.
\end{remark}

Let us now summarise (with respect to thin divisors)
the properties of the involutive
divisions we have encountered so far,
where we note that any strong and continuous involutive
division is by default a Gr\"obner involutive division.
\begin{center}
\begin{tabular}{l|c|c|c}
Division & Continuous & Strong & Gr\"obner \\ \hline
Left                 & Yes & Yes & Yes \\
Right                & Yes & Yes & Yes \\
Left Overlap         & Yes & No  & Yes \\
Right Overlap        & Yes & No  & Yes \\
Strong Left Overlap  & Yes & No  & Yes \\
Strong Right Overlap & Yes & No  & Yes \\ \hline
\end{tabular}
\end{center}
There is a balance to be struck between choosing
an involutive division with nice theoretical
properties and an involutive division which
is of practical use, which is to say that it
is more likely to terminate compared to other divisions.
To this end, one suggestion would be to try
to compute an Involutive Basis with respect
to the left or right divisions to begin with
(as they are easily defined and involutive
reduction with respect to these divisions is
very efficient); otherwise to try one of the
`overlap' divisions, choosing a strong overlap
division if it is important to obtain disjoint
involutive cones.

It is also worth mentioning that for all the
divisions we have encountered so far, if
Algorithm \ref{noncom-inv} terminates then it does so
with a noncommutative Gr\"obner Basis, which
means that Algorithm \ref{noncom-inv} can be
thought of as an alternative algorithm for computing
noncommutative Gr\"obner Bases. Whether this
method is more or less efficient than computing
noncommutative Gr\"obner Bases using
Algorithm \ref{noncom-buch} is a matter
for further discussion.

\subsection{Alternative Divisions}

Having encountered three different types of
involutive division so far (each of which has
two variants -- left and right), let us now
consider if there are
any other involutive divisions with some
useful properties, starting by thinking of
global divisions.

\subsubsection{Alternative Global Divisions}

\begin{openquestion} \label{oq2}
Apart from the empty, left and right divisions,
are there any other global involutive divisions
of the following types:
\begin{enumerate}[(a)]
\item
strong and continuous;
\item
weak, continuous and Gr\"obner?
\end{enumerate}
% If there are any such divisions, are any of them
% conclusive? (recall that all global
% divisions thus far encountered in this
% thesis have been shown to be not conclusive).
\end{openquestion}

\begin{remark}
It seems unlikely that a global division will
exist that affirmatively answers Open Question
\ref{oq2} and does not either assign all variables
to be left nonmultiplicative or all right
nonmultiplicative (thus refining the right or
left divisions respectively). The reason for
saying this is because the moment you have one variable
being left multiplicative and another variable being
right multiplicative for the same monomial
{\it globally}, then you risk not being able
to prove that your division is strong;
similarly the moment you have one variable
being left nonmultiplicative and another
variable being right nonmultiplicative for the same
monomial globally, then you risk not being
able to prove that your division is
continuous. 
% However, at the moment, all of this is just conjecture.
\end{remark}

\subsubsection{Alternative Local Divisions}

So far, all the local divisions we have considered
have assigned all variables to be multiplicative
on one side, and have chosen certain
variables to be nonmultiplicative on the
other side. Let us now consider
a local division that modifies the left overlap
division by assigning some variables
to be nonmultiplicative on both left and
right hand sides.

\begin{defn}[The Two-Sided Left Overlap Division $\mathcal{W}$]
\index{$W$@$\mathcal{W}$}
\index{involutive division!two sided left overlap}
Consider a set $U = \{u_1, \hdots, u_m\}$
of monomials, where all variables are assumed to be left and right
multiplicative for all elements of $U$ to begin with.
Assign multiplicative variables to $U$ according to
Algorithm \ref{inv-div-alt}, which (in words) performs the
following tasks.
\index{two sided left overlap division}
\begin{enumerate}[(a)]
\item
For all possible ways that
a monomial $u_j \in U$ is a subword of a (different)
monomial $u_i \in U$, so that
$$\SUB(u_i, k, k+\deg(u_j)-1) = u_j$$
for some integer $k$, assign
the variable $\SUB(u_i, k-1, k-1)$ to be left
nonmultiplicative for $u_j$ if $u_j$ is a suffix
of $u_i$; and assign the variable
$\SUB(u_i, k+\deg(u_j), k+\deg(u_j))$ to be right
nonmultiplicative for $u_j$ if $u_j$ is not a suffix of $u_i$.
\item
For all possible ways that
a proper prefix of a monomial $u_i \in U$ is equal to
a proper suffix of a (not necessarily different)
monomial $u_j \in U$, so that
$$\PRE(u_i, k) = \SUFF(u_j, k)$$ for some integer $k$
and $u_i$ is not a subword of $u_j$ or vice-versa,
use the recipe provided in the second half
of Algorithm \ref{inv-div-alt}
to ensure that at least one of the following conditions are
satisfied: (i) the variable $\SUB(u_i, k+1, k+1)$ is
right nonmultiplicative for $u_j$; (ii) the variable
$\SUB(u_j, \deg(u_j)-k, \deg(u_j)-k)$ is left
nonmultiplicative for $u_i$.
\end{enumerate}
\end{defn}

\begin{remark}
For task (b) above, Algorithm \ref{inv-div-alt}
gives preference to monomials which are greater in the
DegRevLex monomial ordering (given the choice, it always assigns a
nonmultiplicative variable to whichever monomial out of $u_i$
and $u_j$ is the smallest); it also attempts to
minimise the number of variables made nonmultiplicative
by only assigning a variable to be nonmultiplicative
if both  the variables $\SUB(u_i, k+1, k+1)$ and
$\SUB(u_j, \deg(u_j)-k, \deg(u_j)-k)$ % referred to above
are respectively right multiplicative and left
multiplicative. These refinements will become crucial when
proving the continuity of the division.
\end{remark}

\begin{algorithm}
\setlength{\baselineskip}{3.5ex}
\caption{The Two-Sided Left Overlap Division $\mathcal{W}$}
\label{inv-div-alt}
\begin{algorithmic}
\vspace*{1.75mm}
% p -> m, m_i -> u_i
\REQUIRE{A set of % distinct 
         monomials $U = \{u_1, \hdots, u_m\}$
         ordered by DegRevLex ($u_1 \geqslant u_2 \geqslant 
         \cdots \geqslant u_m$),
         where $u_i \in R\langle x_1, \hdots, x_n\rangle$.}
\ENSURE{A table $T$ of left and right multiplicative variables
        for all $u_i \in U$, where each entry of $T$ is either 1
        (multiplicative) or 0 (nonmultiplicative).}
\vspace*{1mm}
\STATE
Create a table $T$ of multiplicative 
variables as shown below: \\ \vspace*{0.75mm}
\begin{tabular}{c|ccccccc}
& $x_{1}^L$ & $x_{1}^R$ & $x_{2}^L$ & $x_{2}^R$
& $\cdots$ & $x_{n}^L$ & $x_{n}^R$ \\
\hline
$u_1$ & 1 & 1 & 1 & 1 & $\cdots$ & 1 & 1 \\
% $u_2$ & 1 & 1 & 1 & 1 & $\cdots$ & 1 & 1 \\
$\vdots$ & $\vdots$ & $\vdots$ & $\vdots$ & $\vdots$
& $\ddots$ & $\vdots$ & $\vdots$ \\
$u_m$ & 1 & 1 & 1 & 1 & $\cdots$ & 1 & 1
\end{tabular} \\[1.5mm]
\FOR{{\bf each} monomial $u_i \in U$ ($1 \leqslant i \leqslant m$)}
\FOR{{\bf each} monomial $u_j \in U$ ($i \leqslant j \leqslant m$)}
\STATE
Let $u_i = x_{i_1}x_{i_2}\hdots x_{i_{\alpha}}$ and
$u_j = x_{j_1}x_{j_2}\hdots x_{j_{\beta}}$; \\[0.5mm]
%// Middle Overlaps \\
\IF{($i \neq j$)}
\FOR{{\bf each} $k$ ($1 \leqslant k \leqslant \alpha - \beta + 1$)}
  \IF{($\SUB(u_i, k, k+\beta-1) == u_j$)}
  \STATE
    {\bf if} ($k < \alpha-\beta+1$) {\bf then}
    $T(u_j, x^R_{i_{k+\beta}}) = 0$; \\
    {\bf else} $T(u_j, x^L_{i_{k-1}}) = 0$; \\
    {\bf end if}
  \ENDIF
\ENDFOR
\ENDIF
%\STATE
\vspace*{0.5mm}
%// Left \& Right Overlaps
\FOR{{\bf each} $k$ ($1 \leqslant k \leqslant \beta - 1$)}
  \IF{($\PRE(u_i, k) == \SUFF(u_j, k)$)}
    \STATE
    {\bf if} ($T(u_i, x^L_{j_{\beta - k}}) + T(u_j, x^R_{i_{k+1}}) == 2$)
    {\bf then} $T(u_j, x^R_{i_{k+1}}) = 0$; \\
    {\bf end if}
  \ENDIF
  \IF{($\SUFF(u_i, k) == \PRE(u_j, k)$)}
    \STATE
    {\bf if} ($T(u_i, x^R_{j_{k+1}}) + T(u_j, x^L_{i_{\alpha - k}}) == 2$)
    {\bf then} $T(u_j, x^L_{i_{\alpha - k}}) = 0$; \\
    {\bf end if}
  \ENDIF
\ENDFOR \vspace*{0.5mm}
\ENDFOR
\ENDFOR
\STATE
{\bf return} $T$;
\end{algorithmic}
\vspace*{0.5mm}
\end{algorithm}

\begin{example}
Consider the set of monomials $U := \{zx^2yxy, \, yzx, \, xy\}$
over the polynomial ring $\mathbb{Q}\langle x, y, z\rangle$.
Here are the left and right multiplicative
variables for $U$ with respect to the two-sided left overlap
division $\mathcal{W}$. %(as defined in Algorithm \ref{inv-div-alt}).
\begin{center}
\begin{tabular}{c|c|c}
$u$ & $\mathcal{M}_{\mathcal{W}}^L(u, U)$ &
$\mathcal{M}_{\mathcal{W}}^R(u, U)$ \\ \hline
$zx^2yxy$ & $\{x, y, z\}$  & $\{x, y, z\}$ \\
$yzx$     & $\{y, z\}$     & $\{y, z\}$ \\
$xy$      & $\{x\}$        & $\{y, z\}$ \\ \hline
\end{tabular}
\end{center}
The above table is constructed from the table $T$ shown below,
a table which is obtained by applying Algorithm \ref{inv-div-alt} to $U$.
\begin{center}
\begin{tabular}{c|cccccc}
Monomial & $x^L$ & $x^R$ & $y^L$ & $y^R$ & $z^L$ & $z^R$ \\ \hline
$zx^2yxy$ & 1 & 1 & 1 & 1 & 1 & 1 \\
$yzx$     & 0 & 0 & 1 & 1 & 1 & 1 \\
$xy$      & 1 & 0 & 0 & 1 & 0 & 1 \\ \hline
\end{tabular}
\end{center}
The zero entries in $T$ correspond to the following
overlaps between the elements of $U$ (presented in the
order in which Algorithm \ref{inv-div-alt} encounters
them).
\begin{center}
\begin{tabular}{c|c}
Table Entry & Overlap \\ \hline
$T(yzx, x^R)$
& $\PRE(zx^2yxy, 2) = \SUFF(yzx, 2)$ \\
$T(yzx, x^L)$
& $\SUFF(zx^2yxy, 1) = \PRE(yzx, 1)$ \\
$T(xy, x^R)$
& $\SUB(zx^2yxy, 3, 4) = xy$ \\
$T(xy, y^L)$
& $\SUB(zx^2yxy, 5, 6) = xy$ \\
$T(xy, z^L)$
& $\SUFF(yzx, 1) = \PRE(xy, 1)$ \\ \hline
\end{tabular}
\end{center}
Notice that the overlap $\PRE(yzx, 1) =
\SUFF(xy, 1)$ does not produce a zero entry
for $T(xy, z^R)$, as by the time that we encounter this overlap
in the algorithm, we have already assigned
$T(yzx, x^L) = 0$.
\end{example}

\begin{prop}
The two-sided left overlap division $\mathcal{W}$
is a weak involutive division.
\end{prop}
\begin{pf}
We refer to the proof of Proposition \ref{ov-weak-ch6},
making the obvious changes (for example
replacing $\mathcal{O}$ by $\mathcal{W}$).
\end{pf}

For the following two propositions, we defer their
proofs to Appendix \ref{appA} due to their length
and technical nature.

\begin{prop} \label{alt-cts}
The two-sided left overlap division $\mathcal{W}$
is continuous.
\end{prop}
\begin{pf}
We refer to Appendix \ref{appA}.
\end{pf}

\begin{prop} \label{equiv} 
The two-sided left overlap division $\mathcal{W}$
is a Gr\"obner involutive division.
\end{prop}
\begin{pf}
We refer to Appendix \ref{appA}, noting that the
proof is similar to the proof of
Proposition \ref{ov-grob-ch6}.
\end{pf}

\begin{remark}
Because a variable is sometimes only assigned
nonmultiplicative if two other variables are
multiplicative in Algorithm \ref{inv-div-alt},
the subset condition of Definition
\ref{noncom-div-defn} will not always be satisfied with
respect to the two-sided left overlap division.
This will still hold true even if we apply Algorithm
\ref{DJ} at the end of Algorithm \ref{inv-div-alt},
which means that the two-sided left overlap division
cannot be converted to give a strong involutive division
in the same way that we converted the left overlap
division to give the strong left overlap division.
\end{remark}

To finish this section, let us now consider some further
variations of the left overlap division,
variations that will allow us to assign more
multiplicative variables than the left overlap
division (and hence potentially have to deal with fewer
prolongations when using Algorithm \ref{noncom-inv}),
but variations that cannot
be modified to give strong involutive divisions
in the same way that the left overlap division was
modified to give the strong left overlap division
(this is because there are other ways beside a monomial
being a suffix of another monomial that
two involutive cones can be non-disjoint with respect
to these modified divisions).

\begin{defn}[The Prefix-Only Left Overlap Division]
% \index{$\mathcal{O}'$}
\index{involutive division!prefix-only left overlap}
Let $U = \{u_1, \hdots, u_m\}$ be a set of monomials,
and assume that all variables are left and right
multiplicative for all elements of $U$ to begin with.
\index{prefix-only left overlap division}
\begin{enumerate}[(a)]
\item
For all possible ways that
a monomial $u_j \in U$ is a proper prefix of a % (different)
monomial $u_i \in U$,
% so that $$\PRE(u_i, k, k+\deg(u_j)-1) = u_j$$
% for some integer $k$, if $u_j$ is a proper prefix $u_i$,
assign the variable
$\SUB(u_i, \deg(u_j)+1, \deg(u_j)+1)$ to be right
nonmultiplicative for $u_j$.
\item
For all possible ways that
a proper prefix of a monomial $u_i \in U$ is equal to
a proper suffix of a (not necessarily different)
monomial $u_j \in U$, so that
$$\PRE(u_i, k) = \SUFF(u_j, k)$$ for some integer $k$
and $u_i$ is not a subword of $u_j$ or vice-versa, assign
the variable $\SUB(u_i, k+1, k+1)$ to be
right nonmultiplicative for $u_j$.
\end{enumerate}
\end{defn}

\begin{defn}[The Subword-Free Left Overlap Division]
% \index{$\mathcal{O}''$}
\index{involutive division!subword-free left overlap}
Consider a set $U = \{u_1, \hdots, u_m\}$ of monomials,
where all variables are assumed to be left and right
multiplicative for all elements of $U$ to begin with.
\index{subword-free left overlap division}

For all possible ways that
a proper prefix of a monomial $u_i \in U$ is equal to
a proper suffix of a (not necessarily different)
monomial $u_j \in U$, so that
$$\PRE(u_i, k) = \SUFF(u_j, k)$$ for some integer $k$
and $u_i$ is not a subword of $u_j$ or vice-versa, assign
the variable $\SUB(u_i, k+1, k+1)$ to be
right nonmultiplicative for $u_j$.
\end{defn}

\begin{prop}
Both the prefix-only left overlap and the
subword-free left overlap divisions are
continuous, weak and Gr\"obner.
\end{prop}
\begin{pf}
We leave these proofs as exercises for the interested
reader, noting that the proofs will be based on (and in
some cases will be identical to) the proofs
of Propositions \ref{ov-cts-ch6}, \ref{ov-weak-ch6}
and \ref{ov-grob-ch6} respectively.
\end{pf}

\begin{remark}
To help distinguish between the different types
of overlap division we have encountered in this chapter,
let us now give the following table showing
which types of overlap each overlap division
considers.
\begin{center}
\begin{tabular}{cccc}
Type A & Type B & Type C & Type D \\[1mm]
$\xymatrix @R=0.5pc{
& \ar@{<->}[rrr]^*+{} &&& \\
\ar@{<->}[rrr]_*+{} &&&
}$
&
$\xymatrix @R=0.5pc{
\ar@{<->}[rrr]^*+{} &&& \\
\ar@{<->}[rr]_*+{} &&
}$
&
$\xymatrix @R=0.5pc{
\ar@{<->}[rrr]^*+{} &&& \\
& \ar@{<->}[r]_*+{} &
}$
&
$\xymatrix @R=0.5pc{
\ar@{<->}[rrr]^*+{} &&& \\
& \ar@{<->}[rr]_*+{} &&
}$
\end{tabular} \\[5mm]
\begin{tabular}{l|ccccccccc}
Overlap Division Type & \multicolumn{9}{c}{Overlap Type} \\
& \hspace*{2mm} & A & \hspace*{5mm} & B
& \hspace*{5mm} & C & \hspace*{5mm} & D & \hspace*{2mm} \\ \hline
Left               && $\checkmark$ && $\checkmark$ 
&& $\checkmark$ && $\times$ & \\
Right              && $\checkmark$ && $\times$     
&& $\checkmark$ && $\checkmark$ & \\
Strong Left        && $\checkmark$ && $\checkmark$ 
&& $\checkmark$ && $\times$ & \\
Strong Right       && $\checkmark$ && $\times$     
&& $\checkmark$ && $\checkmark$ & \\
Two-Sided Left     && $\checkmark$ && $\checkmark$ 
&& $\checkmark$ && $\checkmark$ & \\
Two-Sided Right    && $\checkmark$ && $\checkmark$ 
&& $\checkmark$ && $\checkmark$ & \\
Prefix-Only Left   && $\checkmark$ && $\checkmark$ 
&& $\times$     && $\times$ & \\
Suffix-Only Right  && $\checkmark$ && $\times$     
&& $\times$     && $\checkmark$ & \\
Subword-Free Left  && $\checkmark$ && $\times$     
&& $\times$     && $\times$ & \\
Subword-Free Right && $\checkmark$ && $\times$     
&& $\times$     && $\times$ & \\ \hline
\end{tabular}
\end{center}
\end{remark}

\section{Termination} \label{5point6}

Given a basis $F$ generating an ideal over a
noncommutative polynomial ring $\mathcal{R}$,
does there exist a finite Involutive Basis for
$F$ with respect to some admissible monomial
ordering $O$ and some
involutive division $I$? Unlike the
commutative case, where the answer to the
corresponding question (for certain divisions)
is always `Yes', the answer to this question can
potentially be `No', as
if the noncommutative Gr\"obner Basis for $F$ with
respect to $O$ is infinite, then the noncommutative
Involutive Basis algorithm will not find a
finite Involutive Basis for $F$
with respect to $I$ and $O$,
as it will in effect be trying to compute the
same infinite Gr\"obner Basis. % as the noncommutative
% Gr\"obner Basis algorithm.

However, a valid follow-up question would be to ask
whether the noncommutative Involutive Basis algorithm
will terminate in the case that the noncommutative Gr\"obner
Basis algorithm terminates.
In Section \ref{CoCo}, we defined a property
of noncommutative involutive divisions (conclusivity)
that ensures, when satisfied, that the answer to this
secondary question is always `Yes'. Despite this,
we will not prove in this thesis that any of the
divisions we have defined are conclusive. Instead,
we leave the following open question for
further investigation. % , and will then prove a result
% that provides some hope that a noncommutative
% involutive division may be able to be shown to be
% conclusive.

\begin{openquestion} \label{oq3}
Are there any conclusive noncommutative involutive divisions
that are also continuous and either strong or Gr\"obner?
\end{openquestion}

To obtain an affirmative answer to the above question,
one approach may be to start by finding a proof for the
following conjecture.

\begin{conj}
Let $O$ be an arbitrary admissible monomial ordering,
and let $I$ be an arbitrary involutive
division that is continuous and either strong or
Gr\"obner. When computing an Involutive Basis for
some basis $F$ with respect to $O$ and $I$
using Algorithm \ref{noncom-inv}, if $F$ possesses
a finite unique reduced Gr\"obner Basis $G$ with
respect to $O$, then after a finite number of steps 
of Algorithm \ref{noncom-inv}, $\LM(G)$ appears as a 
subset of the set of leading monomials of the current basis.
\end{conj}

To prove that a particular involutive division is conclusive, we
would then need to show that once $\LM(G)$ appears as a
subset of the set of leading monomials of the current basis,
% a noncommutative Gr\"obner Basis 
% appears as a subset of the current basis,
then the noncommutative Involutive Basis algorithm
terminates (either immediately or in a finite number
of steps), thus providing the required finite noncommutative
Involutive Basis for $F$.

\section{Examples} \label{Ch6Ex}

\subsection{A Worked Example}

\begin{example} \label{appCex1}
Let $F := \{f_1, f_2\} =
\{x^2y^2 - 2xy^2 + x^2, \, x^2y - 2xy\}$
be a basis for an ideal $J$ over the
polynomial ring $\mathbb{Q}\langle x, y\rangle$,
and let the monomial ordering be DegLex.
Let us now compute a Locally Involutive Basis for $F$
with respect to the strong left overlap division
$\mathcal{S}$ and thick divisors using
Algorithm \ref{noncom-inv}.

To begin with, we must
autoreduce the input set $F$. This leaves the set
unchanged, as we can verify by using the
following table of multiplicative variables
(obtained by using Algorithm \ref{inv-div2}), where $y$
is right nonmultiplicative for $f_2$ because of the
overlap $\LM(f_2) = \SUB(\LM(f_1), 1, 3)$; and
$x$ is right nonmultiplicative for $f_1$ because we
need to have a variable in $\LM(f_2)$ being right
nonmultiplicative for $f_1$.
\begin{center}
\begin{tabular}{c|c|c}
Polynomial & $\mathcal{M}_{\mathcal{S}}^L(f_i, F)$ &
$\mathcal{M}_{\mathcal{S}}^R(f_i, F)$ \\ \hline
$f_1 = x^2y^2 - 2xy^2 + x^2$  & $\{x, y\}$ & $\{y\}$ \\
$f_2 = x^2y-2xy$              & $\{x, y\}$ & $\{x\}$ \\ \hline
\end{tabular}
\end{center}
The above table also provides us with the set
$S = \{f_1x, f_2y\} =
\{x^2y^2x-2xy^2x+x^3, \, x^2y^2-2xy^2\}$ of
prolongations that is required
for the next step of the algorithm.
As $x^2y^2 < x^2y^2x$ in the DegLex monomial ordering, we
involutively reduce the element $f_2y \in S$ first.
\begin{eqnarray*}
f_2y = x^2y^2-2xy^2 & \xymatrix{\ar[r]_{\mathcal{S}}_(1){f_1} &} &
x^2y^2-2xy^2 - (x^2y^2-2xy^2+x^2) \\
& = & -x^2.
\end{eqnarray*}
As the prolongation did not involutively reduce to zero, we now
exit from the second while loop of Algorithm \ref{noncom-inv} and
proceed by autoreducing the set $F\cup\{f_3 := -x^2\} =
\{x^2y^2-2xy^2+x^2, \, x^2y-2xy, \, -x^2\}$.
\begin{center}
\begin{tabular}{c|c|c}
Polynomial & $\mathcal{M}_{\mathcal{S}}^L(f_i, F)$ &
$\mathcal{M}_{\mathcal{S}}^R(f_i, F)$ \\ \hline
$f_1 = x^2y^2 - 2xy^2 + x^2$  & $\{x, y\}$ & $\{y\}$ \\
$f_2 = x^2y-2xy$              & $\{x, y\}$ & $\emptyset$ \\
$f_3 = -x^2$                  & $\{x, y\}$ & $\emptyset$ \\ \hline
\end{tabular}
\end{center}
This process involutively reduces the third term of
$f_1$ using $f_3$, leaving the new
set $\{f_4 := x^2y^2-2xy^2, \, f_2, \, f_3\}$ whose
multiplicative variables are identical to the
multiplicative variables of the set
$\{f_1, \, f_2, \, f_3\}$ shown above.

Next, we construct the set $S = \{f_4x, \, f_2x, \, f_2y, \,
f_3x, \, f_3y\}$ of prolongations, processing the element
$f_3y$ first.
\begin{eqnarray*}
f_3y = -x^2y & \xymatrix{\ar[r]_{\mathcal{S}}_(1){f_2} &} &
-x^2y + (x^2y-2xy) \\
& = & -2xy.
\end{eqnarray*}
Again the prolongation did not involutively reduce to zero,
so we add the involutively reduced prolongation to
our basis to obtain the set $\{f_4, \, f_2, \, f_3,
\, f_5 := -2xy\}$.
\begin{center}
\begin{tabular}{c|c|c}
Polynomial & $\mathcal{M}_{\mathcal{S}}^L(f_i, F)$ &
$\mathcal{M}_{\mathcal{S}}^R(f_i, F)$ \\ \hline
$f_4 = x^2y^2 - 2xy^2$  & $\{x, y\}$ & $\{y\}$ \\
$f_2 = x^2y-2xy$        & $\{x, y\}$ & $\emptyset$ \\
$f_3 = -x^2$            & $\{x, y\}$ & $\emptyset$ \\
$f_5 = -2xy$            & $\{x, y\}$ & $\emptyset$ \\ \hline
\end{tabular}
\end{center}
This time during autoreduction, the polynomial
$f_2$ involutively reduces to zero with respect to the
set $\{f_4, \, f_3, \, f_5\}$:
\begin{eqnarray*}
f_2 = x^2y-2xy & \xymatrix{\ar[r]_{\mathcal{S}}_(1){f_5} &} &
x^2y-2xy + \frac{1}{2}x(-2xy) \\
& = & -2xy \\
& \xymatrix{\ar[r]_{\mathcal{S}}_(1){f_5} &} &
-2xy - (-2xy) \\
& = & 0.
\end{eqnarray*}
This leaves us with the set $\{f_4, \, f_3, \, f_5\}$
after autoreduction is complete.
\begin{center}
\begin{tabular}{c|c|c}
Polynomial & $\mathcal{M}_{\mathcal{S}}^L(f_i, F)$ &
$\mathcal{M}_{\mathcal{S}}^R(f_i, F)$ \\ \hline
$f_4 = x^2y^2 - 2xy^2$  & $\{x, y\}$ & $\{y\}$ \\
$f_3 = -x^2$            & $\{x, y\}$ & $\emptyset$ \\
$f_5 = -2xy$            & $\{x, y\}$ & $\emptyset$ \\ \hline
\end{tabular}
\end{center}
The next step is to construct the set
$S = \{f_4x, f_3x, f_3y, f_5x, f_5y\}$ of prolongations,
from which the element $f_5y$ is processed first.
% (its lead monomial is minimal with respect to DegLex).
$$f_5y = -2xy^2 =: f_6.$$
When the set $\{f_4, \, f_3, \, f_5, \, f_6\}$ is
autoreduced, the polynomial $f_4$ now involutively
reduces to zero, leaving us with the autoreduced set
$\{f_3, \, f_5, \, f_6\} = \{-x^2, \, -2xy, \, -2xy^2\}$.
\begin{center}
\begin{tabular}{c|c|c}
Polynomial & $\mathcal{M}_{\mathcal{S}}^L(f_i, F)$ &
$\mathcal{M}_{\mathcal{S}}^R(f_i, F)$ \\ \hline
$f_3 = -x^2$   & $\{x, y\}$ & $\emptyset$ \\
$f_5 = -2xy$   & $\{x, y\}$ & $\emptyset$ \\
$f_6 = -2xy^2$ & $\{x, y\}$ & $\{y\}$ \\ \hline
\end{tabular}
\end{center}
Our next task is to process the elements of the set
$S = \{f_3x, \, f_3y, \, f_5x, \, f_5y, \, f_6x\}$
of prolongations. The first element $f_5y$ we pick from
$S$ involutively reduces to zero, but the second
element $f_5x$ does not:
\begin{eqnarray*}
f_5y = -2xy^2 & \xymatrix{\ar[r]_{\mathcal{S}}_(1){f_6} &} &
-2xy^2 - (-2xy^2) \\
& = & 0; \\[5mm]
f_5x = -2xyx & =: & f_7.
\end{eqnarray*}
After constructing the set $\{f_3, \, f_5, \, f_6, \, f_7\}$,
autoreduction does not alter the contents of the set, leaving us to
construct our next set of prolongations from the following
table of multiplicative variables.
\begin{center}
\begin{tabular}{c|c|c}
Polynomial & $\mathcal{M}_{\mathcal{S}}^L(f_i, F)$ &
$\mathcal{M}_{\mathcal{S}}^R(f_i, F)$ \\ \hline
$f_3 = -x^2$   & $\{x, y\}$ & $\emptyset$ \\
$f_5 = -2xy$   & $\{x, y\}$ & $\emptyset$ \\
$f_6 = -2xy^2$ & $\{x, y\}$ & $\{y\}$ \\
$f_7 = -2xyx$  & $\{x, y\}$ & $\emptyset$ \\ \hline
\end{tabular}
\end{center}
Whilst processing this (7 element)
set of prolongations, we add
the involutively irreducible prolongation $f_6x
= -2xy^2x =: f_8$ to our basis to give a five
element set which in unaffected by autoreduction.
\begin{center}
\begin{tabular}{c|c|c}
Polynomial & $\mathcal{M}_{\mathcal{S}}^L(f_i, F)$ &
$\mathcal{M}_{\mathcal{S}}^R(f_i, F)$ \\ \hline
$f_3 = -x^2$    & $\{x, y\}$ & $\emptyset$ \\
$f_5 = -2xy$    & $\{x, y\}$ & $\emptyset$ \\
$f_6 = -2xy^2$  & $\{x, y\}$ & $\{y\}$ \\
$f_7 = -2xyx$   & $\{x, y\}$ & $\emptyset$ \\
$f_8 = -2xy^2x$ & $\{x, y\}$ & $\emptyset$ \\ \hline
\end{tabular}
\end{center}
To finish, we analyse the elements of the set
$$S = \{f_3x, \, f_3y, \, f_5x, \, f_5y, \,
f_6x, \, f_7x, \, f_7y, \, f_8x, \, f_8y\}$$
of prolongations in
the order $f_5y, \, f_5x, \, f_3y, \, f_3x, \,
f_6x, \, f_7y, \, f_7x, \, f_8y, \, f_8x$.
\begin{eqnarray*}
f_5y = -2xy^2 & \xymatrix{\ar[r]_{\mathcal{O}}_(1){f_6} &} &
-2xy^2 - (-2xy^2) \\
& = & 0; \\
& \vdots & \\
f_8x = -2xy^2x^2 & \xymatrix{\ar[r]_{\mathcal{O}}_(1){f_3} &} &
-2xy^2x^2 - 2xy^2(-x^2) \\
& = & 0.
\end{eqnarray*}
Because all prolongations involutively reduce to zero
(and hence $S = \emptyset$), the algorithm now terminates with the
Involutive Basis $$G := \{-x^2, \, -2xy, \, -2xy^2, \,
-2xyx, \, -2xy^2x\}$$ as output, a basis which can be visualised by
looking at the following (partial) involutive monomial lattice
for $G$.
$$\xymatrix @C=0.04pc @R=2pc{
&&&&&&&& 1 \\
&&&& x &&&&&&&& y \\
&& *+[o][F-]{x^2} \ar@{~}[ddl] \ar@{~}[ddrrrrr]
&&&& *+[o][F-]{xy} \ar@{-}[llldd] \ar@{-}[rrrrrdd]
&&&& yx &&&& y^2 \\ \\
& x^3 \ar@{~}[dddl] \ar@{~}[dddrrr]
&& x^2y \ar@{-}[dddll] \ar@{-}[rrrrrrrddd]
&& *+[o][F-]{xyx} \ar@{--}[dddlll] \ar@{--}[rrrrddd]
&& yx^2 \ar@{~}[dddllll] \ar@{~}[dddrrrr]
&& *+[o][F-]{xy^2} \ar@{.}[llllddd] \ar@{.}[rrrddd]
\ar@{.}[rrrrddd]
&& yxy \ar@{-}[dddllll] \ar@{-}[dddrrr]
&& y^2x
&& y^3 \\ \\ \\
x^4 & x^3y & x^2yx & xyx^2 &
yx^3 & x^2y^2 & *+[o][F-]{xy^2x} & xyxy &&
yxyx & yx^2y & y^2x^2 & xy^3 &
yxy^2 & y^2xy & y^3x & y^4
}$$
For comparison, the (partial)
monomial lattice of the reduced
DegLex Gr\"obner Basis $H$ for $F$ is shown
below, where $H := \{x^2, \, xy\}$ is
obtained by applying Algorithm \ref{red-noncom}
to $G$.
$$\xymatrix @C=0.04pc @R=2pc{
&&&&&&&& 1 \\
&&&& x &&&&&&&& y \\
&& *+[o][F-]{x^2} \ar@{~}@<0.75ex>[ddl] \ar@{~}[ddl]
\ar@{~}[ddr] \ar@{~}[ddrrrrr]
&&&& *+[o][F-]{xy} \ar@{-}[llldd]
\ar@{-}[ldd] \ar@{-}[rrrdd] \ar@{-}[rrrrrdd]
&&&& yx &&&& y^2 \\ \\
& x^3 \ar@{~}[dddl] \ar@{~}@<-0.75ex>[dddl]
\ar@{~}[ddd] \ar@{~}[dddrrr]
&& x^2y \ar@{-}[dddll] \ar@{-}[dddrr]
\ar@{~}@<0.75ex>[dddl] \ar@{~}@<0.75ex>[rrrrrrrddd]
\ar@{~}@<-0.75ex>[dddll] \ar@{~}@<0.75ex>[dddrr]
\ar@{-}[dddl] \ar@{-}[rrrrrrrddd]
&& xyx \ar@{-}[dddlll] \ar@{-}[dddrr]
\ar@{-}[dddll] \ar@{-}[rrrrddd]
&& yx^2 \ar@{~}[dddlll] \ar@{~}[dddllll]
\ar@{~}[dddrrrr] \ar@{~}[dddrrr]
&& xy^2 \ar@{-}[llllddd] \ar@{-}[rrrddd]
\ar@{-}[dddlll] \ar@{-}[rrrrddd]
&& yxy \ar@{-}[dddllll] \ar@{-}[dddrr]
\ar@{-}[dddrrr] \ar@{-}[llddd]
&& y^2x
&& y^3 \\ \\ \\
x^4 & x^3y & x^2yx & xyx^2 &
yx^3 & x^2y^2 & xy^2x & xyxy &&
yxyx & yx^2y & y^2x^2 & xy^3 &
yxy^2 & y^2xy & y^3x & y^4
}$$
Looking at the lattices, we can verify that the
involutive cones give a disjoint cover of the
conventional cones {\it up to monomials of degree 4}.
However, if we were to draw the next part of
the lattices (monomials of degree 5), we would notice
that the monomial $xy^3x$ is conventionally
reducible by the Gr\"obner Basis, but is not
involutively reducible by the Involutive Basis.
This fact verifies that when thick divisors are
being used, a Locally Involutive Basis is not
necessarily an Involutive Basis, as for $G$ to be
an Involutive Basis with respect to $\mathcal{S}$
and thick divisors, the monomial $xy^3x$
has to be involutively reducible with respect to $G$.
\end{example}

\subsection{Involutive Rewrite Systems}

\begin{remark}
In this section, we use terminology from the theory
of term rewriting that has not yet been introduced in this thesis.
For an elementary introduction to this theory,
we refer to \cite{Baader98}, \cite{Der93} and 
\cite{Jou95}.
\end{remark}

Let $C = \langle A \mid B \rangle$ be a monoid
rewrite system, where $A = \{a_1, \hdots, a_n\}$
is an alphabet and $B = \{b_1, \hdots, b_m\}$ is
a set of rewrite rules of the form
$b_i = \ell_i \rightarrow r_i$
($1 \leqslant i \leqslant m$;
$\ell_i, r_i \in A^{\ast}$). Given a fixed
admissible well-order on the words in $A$
compatible with $C$,
the Knuth-Bendix critical pairs completion
algorithm \cite{KB} attempts to find
a complete rewrite system 
\index{complete rewrite system}
$C'$ for $C$ that is
Noetherian and confluent, so that
any word over the alphabet $A$ has a unique normal form
with respect to $C'$. The algorithm proceeds by
considering overlaps of left hand sides of rules,
forming new rules when two reductions of an overlap
word result in two distinct normal forms.

It is well known (see for example \cite{Hey00a})
that the Knuth-Bendix critical pairs completion
algorithm is a special case of the noncommutative
Gr\"obner Basis algorithm. To find
a complete rewrite system for $C$ using Algorithm
\ref{noncom-buch}, we treat $C$ as a set of polynomials
$F = \{\ell_1 - r_1, \: \ell_2 - r_2, \: \hdots, \:
\ell_m - r_m\}$
generating a two-sided ideal
over the noncommutative polynomial ring
$\mathbb{Z}\langle a_1, \hdots, a_n\rangle$, and we compute a
noncommutative Gr\"obner Basis $G$ for $F$ using a
monomial ordering induced from the fixed admissible
well-order on the words in $A$. 

Because every noncommutative Involutive Basis (with respect
to a strong or Gr\"obner involutive division) is a noncommutative
Gr\"obner Basis, it is clear that 
a complete rewrite system for $C$
can now also be obtained by computing an Involutive
Basis for $F$, a complete rewrite system we shall 
call an \index{complete rewrite system!involutive}
{\it involutive complete rewrite system}. 
\index{involutive complete rewrite system}

The advantage of involutive complete rewrite systems
over conventional complete rewrite systems is
that the unique normal form of any word over the
alphabet $A$ can be obtained uniquely with
respect to an involutive complete rewrite system
(subject of course to certain conditions (such as
working with a strong involutive division) being
satisfied), a fact that will be illustrated in the
following example.

\begin{example} \label{S3}
Let $C := \langle Y, X, y, x
\mid
x^3 \rightarrow \varepsilon, \,
y^2 \rightarrow \varepsilon, \,
(xy)^2 \rightarrow \varepsilon, \,
Xx \rightarrow \varepsilon, \,
xX \rightarrow \varepsilon, \,
Yy \rightarrow \varepsilon, \,
yY \rightarrow \varepsilon\rangle$ be a monoid
rewrite system for the group $S_3$,
where $\varepsilon$ denotes the empty word,
and $Y > X > y > x$ is the alphabet ordering.
If we apply the Knuth-Bendix algorithm to $C$
with respect to the DegLex (word) ordering,
we obtain the complete rewrite system
\begin{center}
$C' := \langle Y, X, y, x
\mid
xyx \rightarrow y, \;
yxy \rightarrow X, \;
x^2 \rightarrow X, \;
Xx \rightarrow \varepsilon, \;
y^2 \rightarrow \varepsilon, \;
Xy \rightarrow yx, \;
xX \rightarrow \varepsilon, \;
yX \rightarrow xy, \;
X^2 \rightarrow x, \;
Y \rightarrow y\rangle.$
\end{center}
With respect to the DegLex monomial ordering
and the left division, if we apply Algorithm 
\ref{noncom-inv} to the basis
$F := \{x^3-1, \, y^2-1, \, (xy)^2-1, \, 
Xx-1, \, xX-1, \, Yy-1, \, yY-1\}$ corresponding
to $C$, we obtain the following Involutive Basis  
for $F$ (which we have converted back to a 
rewrite system to give an involutive complete
rewrite system $C''$ for $C$).
\begin{center}
$C'' := \langle Y, X, y, x \mid y^2 \rightarrow \varepsilon, \, 
Xx \rightarrow \varepsilon, \, xX \rightarrow \varepsilon, \, 
Yy \rightarrow \varepsilon, \, y^2x \rightarrow x, \, 
Y \rightarrow y, \, Yx \rightarrow yx, \, Xxy \rightarrow y, \, 
Yyx \rightarrow x, \, x^2 \rightarrow X, \, X^2 \rightarrow x, \, 
xyx \rightarrow y, \, Xy \rightarrow yx, \, Xyx \rightarrow xy, \,
x^2y \rightarrow yx, \, yX \rightarrow xy, \, 
yxy \rightarrow X, \, Yxy \rightarrow X, \, YX \rightarrow xy
\rangle$.
\end{center}
With the involutive complete rewrite system, we are now able to
uniquely reduce each word over the
alphabet $\{Y, X, y, x\}$ to one of the six
elements of $S_3$. To illustrate this, consider the word
$yXYx$. % $x^2yxy$.
Using the 10 element complete rewrite system $C'$
obtained by using the Knuth-Bendix algorithm,
there are several reduction paths for this word, as illustrated 
by the following diagram.
$$\xymatrix@R=5ex{
& yXYx \ar[ld]_{yX \rightarrow xy}
       \ar[rd]^{Y \rightarrow y} \\
xyYx \ar[rd]_{Y \rightarrow y}
&& yXyx \ar[ld]^{yX \rightarrow xy}
        \ar[rd]^{Xy \rightarrow yx} \\
& xy^2x \ar[rd]_{y^2 \rightarrow \varepsilon}
&& y^2x^2 \ar[ld]_{y^2 \rightarrow \varepsilon}
          \ar[rd]^{x^2 \rightarrow X} \\
& % x^2 \ar[d]_{x^2 \rightarrow X}
& x^2 \ar[rd]_{x^2 \rightarrow X}
&& y^2X \ar[ld]^{y^2 \rightarrow \varepsilon}
        \ar[rd]^{yX \rightarrow xy} \\
& && X && yxy \ar[d]^{yxy \rightarrow X} \\
&&&&& X
}$$
% $$\xymatrix@R=5ex{
% &&& x^2yxy \ar[lld]_{x^2 \rightarrow X}
%            \ar[d]^{xyx \rightarrow y}
%            \ar[rrd]^{yxy \rightarrow X}\\
% & Xyxy \ar[ld]_{yxy \rightarrow X}
%        \ar[rd]^{Xy \rightarrow yx}
% && xy^2 \ar[d]^{y^2 \rightarrow \varepsilon}
% && x^2X \ar[ld]_{x^2 \rightarrow X}
%         \ar[rd]^{xX \rightarrow \varepsilon} \\
% X^2 \ar[d]_{X^2 \rightarrow x}
% && yx^2y \ar[d]^{x^2 \rightarrow X}
% & x & X^2 \ar[d]_{X^2 \rightarrow x} && x \\
% x && yXy \ar[ld]_{Xy \rightarrow yx}
%          \ar[rd]^{yX \rightarrow xy} && x \\
% & y^2x \ar[d]_{y^2 \rightarrow \varepsilon}
% && xy^2 \ar[d]^{y^2 \rightarrow \varepsilon} \\
% & x    && x
% }$$
However, by involutively reducing the word
$yXYx$ % $x^2yxy$ 
with respect to the 19 element involutive
complete rewrite system $C''$, there is only
one reduction path, namely
$$\xymatrix{
yXYx \ar[d]^{Yx \rightarrow yx} \\
yXyx \ar[d]^{Xyx \rightarrow xy} \\
yxy  \ar[d]^{yxy \rightarrow X} \\
X
}$$
% $$\xymatrix{
% x^2yxy \ar[d]^{yxy \rightarrow X} \\
% x^2X   \ar[d]^{xX \rightarrow \varepsilon} \\
% x
% }$$
\end{example}

\subsection{Comparison of Divisions} \label{ExS4}

Following on from the $S_3$ example above, 
consider the basis $F :=
\{x^4 - 1, \, y^3 - 1, \, (xy)^2 - 1, \,
Xx - 1, \, xX - 1, \, Yy - 1, \, yY - 1\}$ over
the polynomial ring $\mathbb{Q}\langle Y, X, y, x\rangle$
corresponding to a monoid rewrite system for the
group $S_4$. With the monomial ordering being
DegLex, below we present some data collected
when, whilst using a prototype implementation of 
Algorithm \ref{noncom-inv} (as given in Appendix 
\ref{appB}), an Involutive Basis is computed for $F$ 
with respect to several different involutive divisions 
(the reduced DegLex Gr\"obner Basis for $F$ has 21 
elements).

\begin{remark}
The program was run using
FreeBSD 5.4 on an AMD Athlon XP 1800+ with
512MB of memory. 
\end{remark}
\begin{center}
\begin{tabular}{c|l||c|l}
Key & Involutive Division & Key & Involutive Division \\ \hline
1  & Left Division &
7  & Subword-Free Left Overlap Division \\
2  & Right Division &
8  & Right Overlap Division \\
3  & Left Overlap Division &
9  & Strong Right Overlap Division \\
4  & Strong Left Overlap Division &
10 & Two-Sided Right Overlap Division \\
5  & Two-Sided Left Overlap Division &
11 & Suffix-Only Right Overlap Division \\
6  & Prefix-Only Left Overlap Division &
12 & Subword-Free Right Overlap Division \\ \hline
\end{tabular}
\end{center}
\begin{center}
\begin{tabular}{c|c|c|c|c}
Division & Size of Basis & Number of
& Number of & Time \\
&& Prolongations & Involutive Reductions & \\ \hline
1  & 73 & 104 & 15947 & 0.77  \\
2  & 73 & 104 & 13874 & 0.74  \\
3  & 65 & 64  & 10980 & 8.62  \\
4  & 73 & 94  & 15226 & 23.14 \\
5  & 77 & 70  & 12827 & 16.04 \\
6  & 65 & 64  & 10980 & 8.97  \\
7  & 65 & 64  & 10980 & 7.13  \\
8  & 73 & 76  & 11046 & 13.27 \\
9  & 73 & 95  & 13240 & 26.16 \\
10 & 87 & 80  & 13005 & 24.53 \\
11 & 73 & 76  & 11046 & 13.40 \\
12 & 69 & 82  & 10458 & 9.52  \\ \hline
\end{tabular}
\end{center}

We note that the algorithm completes quickest with
respect to the global left or right divisions, as (i) we can
take advantage of the efficient involutive reduction with
respect to these divisions (see Section \ref{TGD}); and
(ii) the multiplicative variables for a particular
monomial with respect to these divisions is fixed
(each time the basis changes when using
one of the other local divisions, the multiplicative
variables have to be recomputed).
However, we also note that more prolongations are needed
when using the left or right divisions, so that, in the
long run, if we can devise an efficient way of finding
the multiplicative variables for a set of monomials with
respect to one of the local divisions, then the algorithm
could (perhaps) terminate more quickly than for the two global
divisions.

\section{Improvements to the Noncommutative Involutive Basis Algorithm}
\label{5point8}

\subsection{Efficient Reduction} \label{EffRed}

Conventionally, we divide a noncommutative polynomial $p$ with
respect to a set of polynomials $P$ using Algorithm
\ref{noncom-div}. In this algorithm, an important step is
to find out if a polynomial in $P$ divides one of the
monomials $u$ in the polynomial we are currently reducing,
stated as the condition `{\bf if} $(\LM(p_j) \mid u)$
{\bf then}' in Algorithm \ref{noncom-div}. One way of finding 
out if this condition is satisfied would be to execute the 
following piece of code, where $\alpha := \deg(u)$;
$\beta := \deg(\LM(p_j))$; and we note
that $\alpha-\beta+1$ operations are potentially
needed to find out if the condition is satisfied.
\begin{center}
\begin{algorithmic}
\STATE
$i = 1$;
\WHILE{($i \leqslant \alpha-\beta+1$)}
\IF{($\LM(p_j) == \SUB(u, i, i+\beta-1)$)}
\STATE
{\bf return} true;
\ELSE
\STATE
$i = i+1$;
\ENDIF
\ENDWHILE
\STATE
{\bf return} false;
\end{algorithmic}
\end{center}
When involutively dividing a polynomial $p$ with respect
to a set of polynomials $P$ and some involutive
division $I$, the corresponding problem
is to find out if some monomial $\LM(p_j)$ is
an {\it involutive} divisor of some monomial $u$.
At first glance, this problem seems more difficult
than the problem of finding out if $\LM(p_j)$ is a
conventional divisor of $u$, as it is not just 
sufficient to discover one way that $\LM(p_j)$ divides $u$ 
(as in the code above) --- we have to verify that
if we find a conventional divisor of $u$, then it is also an
involutive divisor of $u$. Naively, assuming that
thin divisors are being used, we could
solve the problem using the code shown below, code
that is clearly less efficient than the code for the
conventional case shown above.
\begin{center}
\begin{algorithmic}
\STATE
$i = 1$;
\WHILE{($i \leqslant \alpha-\beta+1$)}
\IF{($\LM(p_j) == \SUB(u, i, i+\beta-1)$)}
\IF{(($i == 1$) {\bf or} (($i > 1$) {\bf and}
($\SUB(u, i-1, i-1) \in \mathcal{M}^L_I(\LM(p_j), \LM(P))$))}
\IF{(($i == \alpha-\beta+1$) {\bf or}
(($i < \alpha-\beta+1$) {\bf and}
($\SUB(u, i+\beta, i+\beta) \in \mathcal{M}^R_I(\LM(p_j), \LM(P))$))}
\STATE
{\bf return} true;
\ENDIF
\ENDIF
\ELSE
\STATE
$i = i+1$;
\ENDIF
\ENDWHILE
\STATE
{\bf return} false;
\end{algorithmic}
\end{center}
However, for certain involutive divisions, it is possible
to take advantage of some of the properties of these
divisions in order to make it easier to discover
whether $\LM(p_j)$ is an involutive divisor of $u$. We have already
seen this in action in Section \ref{TGD}, where we saw that
$\LM(p_j)$ can only involutively divide $u$ with respect
to the left or right divisions if $\LM(p_j)$ is a suffix
or prefix of $u$ respectively.

Let us now consider an improvement to be used whenever (i)
an `overlap' division that assigns all variables to be
either left multiplicative or right multiplicative
is used (ruling out any `two-sided'
overlap divisions); and (ii) thick divisors are being used.
For the case of such an overlap division that assigns
all variables to be left multiplicative (for example
the left overlap division), the following
piece of code can be used to discover whether or not
$\LM(p_j)$ is an involutive divisor
of $u$ (note that a similar piece of code can be given
for the case of an overlap division assigning all variables
to be right multiplicative).
\begin{center}
\begin{algorithmic}
\STATE
$k = \alpha$; skip = 0;
\WHILE{($k \geqslant \beta+1$)}
\IF{($\SUB(u, k, k) \notin \mathcal{M}^R_I(\LM(p_j), \LM(P))$)}
\STATE
skip = $k$; $k = \beta$;
\ELSE
\STATE
$k = k-1$;
\ENDIF
\ENDWHILE \\[2mm]
\IF{(skip == 0)}
\STATE
$i = 1$;
\ELSE
\STATE
$i = \mathrm{skip}-\beta+1$;
\ENDIF \\[2mm]
\WHILE{($i \leqslant \alpha-\beta+1$)}
\IF{($\LM(p_j) == \SUB(u, i, i+\beta-1)$)}
\STATE
{\bf return} true;
\ELSE
\STATE
$i = i+1$;
\ENDIF
\ENDWHILE
\STATE
{\bf return} false;
\end{algorithmic}
\end{center}
We note that the final section of the code
(from `{\bf while} ($i \leqslant \alpha-\beta+1$)
{\bf do}' onwards) is identical to the code
for conventional reduction; the code before this
just chooses the initial value of $i$
(we rule out looking at certain subwords by analysing
which variables in $u$ are right nonmultiplicative
for $\LM(p_j)$). For example, if $u := xy^2xyxy$;
$\LM(p_j) := xyx$; and only the variable $x$ is right
nonmultiplicative for $p_j$, then in the
conventional case we need 4 subword comparisons
before we discover that $\LM(p_j)$ is a
conventional divisor of $u$; but in the involutive
case (using the code shown above) we only need 1
subword comparison before we discover that
$\LM(p_j)$ is an involutive divisor of $u$
(this is because the variable $\SUB(u, 6, 6) = x$
is right nonmultiplicative for $\LM(p_j)$,
leaving just two subwords of $u$ that are
potentially equal to $\LM(p_j)$ in such a way
that $\LM(p_j)$ is an involutive divisor of $u$).
\begin{center}
\begin{tabular}{ccc}
Conventional Reduction & \vspace*{5mm} &
Involutive Reduction \\
$\xymatrix @R=0.5pc @C=0.5pc{
      & x & y & y & x & y & x & y \\
i = 1 & x & y & x \\
i = 2 &   & x & y & x \\
i = 3 &   &   & x & y & x \\
i = 4 &   &   &   & x & y & x
}$
&&
$\xymatrix @R=0.5pc @C=0.5pc{
      & x & y & y & x & y & x & y \\
i = 4 &   &   &   & x & y & x
}$
\end{tabular}
\end{center}
Of course our new algorithm will not always `win'
in every case (for example if $u := xyx^2yxy$
and $\LM(p_j) := xyx$),
and we will always have the overhead from
having to determine the initial value of $i$,
but the impression should be that we have more
freedom in the involutive case to try these
sorts of tricks, tricks which may lead to
involutive reduction being more efficient
than conventional reduction.

\subsection{Improved Algorithms}

Just as Algorithm \ref{com-inv} was generalised
to give an algorithm for computing noncommutative
Involutive Bases in Algorithm \ref{noncom-inv}, it
is conceivable that other algorithms for
computing commutative Involutive Bases
(as seen for example in \cite{Gerdt05}) can be
generalised to the noncommutative case. Indeed,
in the source code given in Appendix \ref{appB},
a noncommutative version of an algorithm
found in \cite[Section 5]{Gerdt02} for computing commutative
Involutive Bases is given; we present below data
obtained by applying this new algorithm to our
$S_4$ example from Section \ref{ExS4} (the data
from Section \ref{ExS4} is given in brackets 
for comparison; we see that the new algorithm
generally analyses more prolongations
but performs less involutive reduction).
\begin{center}
\begin{tabular}{c|c|c|c|c}
Division & Size of Basis & Number of
& Number of & Time \\
&& Prolongations & Involutive Reductions & \\ \hline
1  & 73 (73) & 323 (104) & 875  (15947) & 0.72  (0.77)  \\
2  & 73 (73) & 327 (104) & 929  (13874) & 0.83  (0.74)  \\
3  & 70 (65) & 288 (64)  & 831  (10980) & 5.94  (8.62)  \\
4  & 73 (73) & 318 (94)  & 863  (15226) & 4.62  (23.14) \\
5  & 70 (77) & 288 (70)  & 831  (12827) & 5.79  (16.04) \\
6  & 70 (65) & 288 (64)  & 831  (10980) & 5.71  (8.97)  \\
7  & 69 (65) & 288 (64)  & 833  (10980) & 5.33  (7.13)  \\
8  & 68 (73) & 358 (76)  & 1092 (11046) & 28.51 (13.27) \\
9  & 73 (73) & 322 (95)  & 917  (13240) & 6.39  (26.16) \\
10 & 68 (87) & 358 (80)  & 1092 (13005) & 28.75 (24.53) \\
11 & 68 (73) & 358 (76)  & 1092 (11046) & 28.54 (13.40) \\
12 & 66 (69) & 364 (82)  & 1127 (10458) & 28.87 (9.52)  \\ \hline
\end{tabular}
\end{center}

% \subsection{Homogeneous Involutive Bases}
% 
% Start with $F$. Homogenise to get $F'$ (with respect to
% an extendible monomial ordering). Add in the
% homogenising polynomials' $H$. Problem: Need all the
% polynomials in $H$ to have every variable left and
% right multiplicative in order to be able to apply all
% polynomials in $H$ at all times. This is obviously not
% the case with the left, right and empty divisions;
% can find examples with overlap divisions that ensure
% that nonmultiplicative variables are assigned to
% members of $H$. Possible workaround: do not include
% $H$ in the homogenised basis but allow reductions
% `by $H$' anyway? However this may not be valid as
% we are reducing with respect to polynomials which
% are not in the basis and hence any Locally Involutive
% Basis may not be valid? Even if this problem is
% overcome we still have the problem of the procedure
% only working if the dehomogenised basis is autoreduced.
% But can still define the procedure and also explore
% whether divisions are extendible with respect to
% certain extendible monomial orderings (can construct a
% 4 x 12 table showing which involutive divisions
% are extendible with respect to which monomial orderings.

\subsection{Logged Involutive Bases}

A (noncommutative)
Logged Involutive Basis expresses each member of
an Involutive Basis in terms of members of the original
basis from which the Involutive Basis was computed.

\begin{defn}
Let $G = \{g_1, \hdots, g_p\}$ be
an Involutive Basis computed from
an initial basis $F = \{f_1, \hdots, f_m\}$. We say that $G$ is a
{\it Logged Involutive Basis}
\index{involutive basis!logged}
if, \index{logged involutive basis}
for each $g_i \in G$, we have an explicit expression of the form
$$g_i = \sum_{\alpha=1}^{\beta}
{\ell}_{\alpha}f_{k_{\alpha}}r_{\alpha},$$
where the $\ell_{\alpha}$ and the $r_{\alpha}$ are terms and
$f_{k_{\alpha}} \in F$ for all $1 \leqslant \alpha \leqslant \beta$.
\end{defn}

\begin{prop} \label{LogIBNC}
Let $F = \{f_1, \hdots, f_m\}$ be a finite basis
over a noncommutative polynomial ring.
If we can compute an Involutive Basis for $F$, then
it is always possible to compute a Logged Involutive Basis for $F$.
\end{prop}
\begin{pf}
Let $G = \{g_1, \hdots, g_p\}$ be an Involutive Basis
computed from the initial basis $F = \{f_1, \hdots, f_m\}$
using Algorithm \ref{noncom-inv}
(where $f_i \in R\langle x_1, \hdots, x_n\rangle$
for all $f_i \in F$).
If an arbitrary $g_i \in G$
is not a member of the original basis $F$, then
either $g_i$ is an involutively reduced prolongation,
or $g_i$ is obtained through the process of autoreduction.
In the former case, we can express $g_i$ in terms of
members of $F$ by substitution because either
$$g_i = x_jh - \sum_{\alpha=1}^{\beta}
\ell_{\alpha}h_{k_{\alpha}}r_{\alpha}$$
or $$g_i = hx_j - \sum_{\alpha=1}^{\beta}
\ell_{\alpha}h_{k_{\alpha}}r_{\alpha}$$
for a variable $x_j$; terms $\ell_{\alpha}$ and $r_{\alpha}$; and
polynomials $h$ and $h_{k_{\alpha}}$ which we already know how
to express in terms of members of $F$.
In the latter case,
$$g_i = h - \sum_{\alpha=1}^{\beta}
\ell_{\alpha}h_{k_{\alpha}}r_{\alpha}$$
for terms $\ell_{\alpha}, r_{\alpha}$ and polynomials
$h$ and $h_{k_{\alpha}}$ which we already know how to
express in terms of members of $F$, so it follows that we can again
express $g_i$ in terms of members of $F$.
\end{pf}

\begin{example}
Let $F := \{f_1, f_2\} = \{x^3 + 3xy - yx, \, y^2 + x\}$
generate an ideal over the polynomial
ring $\mathbb{Q}\langle x, y\rangle$; let the monomial ordering be
DegLex; and let the involutive division be the left division.
In obtaining an Involutive Basis for $F$ using Algorithm
\ref{noncom-inv}, a polynomial is added to $F$; $f_1$ is
involutively reduced during autoreduction; and then four
more polynomials are added to $F$,
giving an Involutive Basis $G :=
\{g_1, g_2, g_3, g_4, g_5, g_6, g_7\} =
\{x^3+2yx, \, y^2+x, \, xy-yx, \, y^2x+x^2, \,
xyx-yx^2, \, y^2x^2-2yx, \, xyx^2-2x^2\}$. 

The five new polynomials were obtained by
involutively reducing the prolongations
$f_2y$, $f_2x$, $g_3x$, $g_4x$ and $g_5x$ respectively.
\begin{eqnarray*}
f_2y & = & y^3+xy \\
& \xymatrix{\ar[r]_{\lhd}_(1){f_2} &} & y^3+xy - y(y^2+x) \\
& = & xy-yx; \\[2mm]
f_2x & = & y^2x+x^2; \\[2mm]
g_3x & = & xyx-yx^2; \\[2mm]
g_4x & = & y^2x^2+x^3 \\
& \xymatrix{\ar[r]_{\lhd}_(1){g_1} &} & y^2x^2+x^3-(x^3+2yx) \\
& = & y^2x^2-2yx; \\[2mm]
g_5x & = & xyx^2-yx^3 \\
& \xymatrix{\ar[r]_{\lhd}_(1){g_1} &} & xyx^2-yx^3+y(x^3+2yx) \\
& = & xyx^2+2y^2x \\
& \xymatrix{\ar[r]_{\lhd}_(1){g_4} &} & xyx^2+2y^2x-2(y^2x+x^2) \\
& = & xyx^2-2x^2.
\end{eqnarray*}
These reductions (plus the reduction
\begin{eqnarray*}
f_1 & \xymatrix{\ar[r]_{\lhd}_(1){g_3} &} & x^3+3xy-yx-3(xy-yx) \\
& = & x^3+2yx
\end{eqnarray*}
of $f_1$ performed during autoreduction after $g_3$
is added to $F$) enable us to give the following 
Logged Involutive Basis for $F$.
\begin{center}
\begin{tabular}{l|l} 
Member of $G$ & Logged Representation \\
\hline
$g_1 = x^3+2yx$    & $f_1 - 3f_2y + 3yf_2$ \\
$g_2 = y^2+x$      & $f_2$ \\
$g_3 = xy-yx$      & $f_2y-yf_2$ \\
$g_4 = y^2x+x^2$   & $f_2x$ \\
$g_5 = xyx-yx^2$   & $f_2yx-yf_2x$ \\
$g_6 = y^2x^2-2yx$ & $-f_1 + f_2x^2 + 3f_2y - 3yf_2$ \\
$g_7 = xyx^2-2x^2$ & $yf_1 + 3y^2f_2 + f_2yx^2 - 2f_2x
- yf_2x^2-3yf_2y$ \\ \hline
\end{tabular}
\end{center}
\end{example}

%
% Chapter 6
% Author: Gareth Evans
% Last Modified: 24th January 2006
%

\chapter{Gr\"obner Walks} \label{ChWalk}

When computing any Gr\"obner or Involutive Basis,
the monomial ordering that has been chosen is
a major factor in how long it will take for the algorithm to
complete. For example, consider the ideal $J$
generated by the basis
$F := \{-2x^3z+y^4+y^3z-x^3+x^2y, \,
2xy^2z + yz^3 + 2yz^2, \,
x^3y + 2yz^3 - 3yz^2 + 2\}$
over the polynomial
ring $\mathbb{Q}[x, y, z]$.
Using our test implementation of Algorithm
\ref{com-buch}, it takes less than a tenth of a second
to compute a Gr\"obner Basis for $F$ with respect
to the DegRevLex monomial ordering, but
90 seconds to compute a Gr\"obner Basis for
$F$ with respect to Lex.

The Gr\"obner Walk, \index{Gr\"obner walk} introduced 
\index{walk!Gr\"obner} by
Collart, Kalkbrener and Mall in \cite{CKMWalk},
forms part of a family of basis
conversion algorithms that can convert
Gr\"obner Bases with respect
to `fast' monomial orderings to
Gr\"obner Bases with respect to `slow' monomial
orderings (see Section \ref{BCAC} for a
brief discussion of other basis conversion algorithms).
This process is often quicker than computing
a Gr\"obner Basis for the `slow' monomial
ordering directly, as can be demonstrated by
stating that in our test implementation of the Gr\"obner
Walk, it only takes half a second to compute a Lex
Gr\"obner Basis for the basis $F$ defined above.

In this chapter, we will first recall the
theory of the (commutative) Gr\"obner Walk, based on
\cite{CKMWalk} and a paper \cite{AGK97} by
Amrhein, Gloor and K\"{u}chlin; the reader is
encouraged to read these papers in conjunction with
this Chapter. We then describe two generalisations of the
theory to give (i) a commutative Involutive Walk 
(due to Golubitsky \cite{Golub}); and (ii) noncommutative 
Walks between harmonious monomial orderings.

\section{Commutative Walks} \label{6point1}

To convert a Gr\"obner Basis with respect to one monomial ordering
to a Gr\"obner Basis with respect to another monomial ordering,
the Gr\"obner Walk works with the matrices associated
to the orderings. Fortunately, \cite{Robbiano85} and \cite{Weisp}
assert that any commutative monomial ordering has
an associated matrix, allowing the Gr\"obner Walk to
convert between any two monomial orderings.

\subsection{Matrix Orderings}

\begin{defn}
Let $m$ be a monomial over a polynomial ring
$R[x_1, \hdots, x_n]$ with associated multidegree
$(e^1, \hdots, e^n)$. If $\omega = (\omega^1, \hdots,
\omega^n)$ is an $n$-dimensional weight vector (where
$\omega^i \in \mathbb{Q}$ for all $1 \leqslant i \leqslant n$), we define
\index{$D$@$\deg_{\omega}$}
the {\it $\omega$-degree of $m$}, written
$\deg_{\omega}(m)$, to be the value
$$\deg_{\omega}(m) = (e^1 \times \omega^1) + (e^2 \times \omega^2)
+ \cdots + (e^n \times \omega^n).$$
\end{defn}

\begin{remark}
The $\omega$-degree of any term is equal
to the $\omega$-degree of the term's associated
monomial.
\end{remark}

\begin{defn}
Let $m_1$ and $m_2$ be two monomials over a polynomial
ring $R[x_1, \hdots, x_n]$ with associated
multidegrees $e_1 = (e_1^1, \hdots, e_1^n)$ and
$e_2 = (e_2^1, \hdots, e_2^n)$; and let $\Omega$ be
an $N\times n$ matrix. If $\omega_i$ denotes the
$n$-dimensional weight vector corresponding to the
$i$-th row of $\Omega$, then $\Omega$ determines a monomial
ordering as follows: $m_1 < m_2$ if
$\deg_{\omega_i}(m_1) < \deg_{\omega_i}(m_2)$ for some
$1 \leqslant i \leqslant N$ and
$\deg_{\omega_j}(m_1) = \deg_{\omega_j}(m_2)$ for all
$1 \leqslant j < i$.
\end{defn}

\begin{defn}
The corresponding matrices for the five monomial orderings 
defined in Section \ref{CMO} are
$$
\mathrm{Lex} = \left( \begin{array}{cccccc}
1 & 0 & 0 & \hdots & 0 & 0 \\
0 & 1 & 0 & \hdots & 0 & 0 \\
0 & 0 & 1 & \hdots & 0 & 0 \\
\vdots & \vdots & \vdots & \ddots & \vdots & \vdots \\
0 & 0 & 0 & \hdots & 1 & 0 \\
0 & 0 & 0 & \hdots & 0 & 1
\end{array}\right); \; \; \;
\mathrm{InvLex} = \left( \begin{array}{cccccc}
0 & 0 & 0 & \hdots & 0 & 1 \\
0 & 0 & 0 & \hdots & 1 & 0 \\
\vdots & \vdots & \vdots & \iddots & \vdots & \vdots \\
0 & 0 & 1 & \hdots & 0 & 0 \\
0 & 1 & 0 & \hdots & 0 & 0 \\
1 & 0 & 0 & \hdots & 0 & 0
\end{array}\right);
$$
$$
\mathrm{DegLex} = \left( \begin{array}{cccccc}
1 & 1 & 1 & \hdots & 1 & 1 \\
1 & 0 & 0 & \hdots & 0 & 0 \\
0 & 1 & 0 & \hdots & 0 & 0 \\
0 & 0 & 1 & \hdots & 0 & 0 \\
\vdots & \vdots & \vdots & \ddots & \vdots & \vdots \\
0 & 0 & 0 & \hdots & 1 & 0
\end{array}\right); \; \; \;
\mathrm{DegInvLex} = \left( \begin{array}{cccccc}
1 & 1 & 1 & \hdots & 1 & 1 \\
0 & 0 & 0 & \hdots & 0 & 1 \\
0 & 0 & 0 & \hdots & 1 & 0 \\
\vdots & \vdots & \vdots & \iddots & \vdots & \vdots \\
0 & 0 & 1 & \hdots & 0 & 0 \\
0 & 1 & 0 & \hdots & 0 & 0
\end{array}\right);
$$
$$
\mathrm{DegRevLex} = \left( \begin{array}{rrrrrr}
1 & 1  & 1  & \hdots & 1  & 1  \\
0 & 0  & 0  & \hdots & 0  & -1 \\
0 & 0  & 0  & \hdots & -1 & 0  \\
\vdots & \vdots & \vdots & \iddots & \vdots & \vdots \\
0 & 0  & -1 & \hdots & 0  & 0  \\
0 & -1 & 0  & \hdots & 0  & 0
\end{array}\right).
$$
\end{defn}

\begin{example}
Let $m_1 := x^2y^2z^2$ and $m_2 := x^2y^3z$ be two
monomials over the polynomial ring
$\mathcal{R} := \mathbb{Q}[x, y, z]$.
According to the matrix
$$
\left( \begin{array}{ccc}
1 & 1 & 1 \\
1 & 0 & 0 \\
0 & 1 & 0
\end{array}\right)
$$
representing the DegLex monomial ordering with respect to
$\mathcal{R}$, we can deduce that
$m_1 < m_2$ because
$\deg_{\omega_1}(m_1) = \deg_{\omega_1}(m_2) = 6$;
$\deg_{\omega_2}(m_1) = \deg_{\omega_2}(m_2) = 2$;
and $\deg_{\omega_3}(m_1) = 2 < \deg_{\omega_3}(m_2) = 3$.
\end{example}

\begin{defn} \label{initial-defn-ch6}
Given a polynomial $p$ and a weight vector $\omega$,
the \index{initial} {\it initial} of $p$ with respect to
$\omega$, written \index{$I$@$\init_{\omega}$}
$\init_{\omega}(p)$, is
the sum of those terms in $p$ that have
maximal $\omega$-degree. For example,
if $\omega = (0, 1, 1)$ and
$p = x^4 + xy^2z + y^3 + xz^2$, then
$\init_{\omega}(p) = xy^2z + y^3$.
\end{defn}

\begin{defn}
A weight vector $\omega$ is \index{compatible weight vector}
{\it compatible} with a monomial
ordering $O$ if,
given any polynomial $p = t_1
+ \cdots + t_m$ ordered in descending order with
respect to $O$, $\mathrm{deg}_{\omega}(t_1) \geqslant
\mathrm{deg}_{\omega}(t_i)$ holds for all
$1 < i \leqslant m$.
\end{defn}

\subsection{The Commutative Gr\"obner Walk Algorithm}

We present in Algorithm \ref{com-walk-grob} an
algorithm to perform the Gr\"obner Walk, modified
from an algorithm given in \cite{AGK97}.

\begin{algorithm}
\setlength{\baselineskip}{3.5ex}
\caption{The Commutative Gr\"obner Walk Algorithm}
\label{com-walk-grob}
\begin{algorithmic}
\vspace*{1.75mm}
\REQUIRE{A Gr\"obner Basis $G = \{g_1, g_2, \hdots, g_m\}$
         with respect to an admissible monomial ordering $O$ with
         an associated matrix $A$, where $G$ generates an
         ideal $J$ over a commutative polynomial ring
         $\mathcal{R} = R[x_1, \hdots, x_n]$.}
\ENSURE{A Gr\"obner Basis $H = \{h_1, h_2, \hdots, h_p\}$
        for $J$ with respect to an admissible monomial ordering $O'$
        with an associated matrix $B$.}
\vspace*{0.8mm}
\STATE
Let $\omega$ and $\tau$ be the weight vectors corresponding to
the first rows of $A$ and $B$; \\ % respectively; \\
Let $C$ be the matrix whose first row is equal to
$\omega$ and whose remainder is equal to the whole of
the matrix $B$; \\
$t = 0$; found = false;
\REPEAT
\STATE
Let $G' = \{\init_{\omega}(g_i)\}$ for all $g_i \in G$; \\
Compute a reduced Gr\"obner Basis $H'$ for $G'$ with respect
to the monomial ordering defined by the matrix $C$
(use Algorithms \ref{com-buch} and \ref{red-com});  \\
$H = \emptyset$; \\
\FOR{{\bf each} $h' \in H'$}
\STATE
Let $\sum_{i=1}^{j} p_ig'_i$ be the logged representation
of $h'$ with respect to $G'$
(where $g'_i \in G'$ and $p_i \in \mathcal{R}$),
obtained either through computing a Logged Gr\"obner Basis
above or by dividing $h'$ with respect to $G'$; \\
$H = H \cup \{\sum_{i=1}^{j} p_ig_i\}$, where
$\init_{\omega}(g_i) = g'_i$; \\
\ENDFOR
\STATE
Reduce $H$ with respect to $C$
(use Algorithm \ref{red-com}); \\
\IF{($t == 1$)}
\STATE
found = true;
\ELSE
\STATE
$t = \min (\{ s \mid \mathrm{deg}_{\omega(s)}(p_1) =
\mathrm{deg}_{\omega(s)}(p_i),
\mathrm{deg}_{\omega(0)}(p_1) \neq
\mathrm{deg}_{\omega(0)}(p_i),$ \\
$h = p_1 + \cdots + p_k \in H\} \cap (0, 1])$,
where $\omega(s) := \omega + s(\tau
- \omega)$ for $0 \leqslant s \leqslant 1$;
\ENDIF
\IF{($t$ is undefined)}
\STATE
found = true;
\ELSE
\STATE
$G = H$; $\omega = (1-t)\omega + t\tau$;
\ENDIF
\UNTIL{(found = true)}
\STATE
{\bf return} $H$;
\end{algorithmic}
\vspace*{0.75mm}
\end{algorithm}

\subsubsection{Some Remarks:}
\begin{itemize}
\item
In the first iteration of the {\bf repeat} \ldots
{\bf until} loop, $G'$ is a Gr\"obner Basis for
the ideal\footnote{The ideal $\init_{\omega}(J)$
is defined as follows: $p \in J$ if and only if
$\init_{\omega}(p) \in \init_{\omega}(J)$.}
$\init_{\omega}(J)$ with respect to the monomial
ordering defined by $C$, as $\omega$ is compatible
with $C$. During subsequent iterations of the same loop,
$G'$ is a Gr\"obner Basis for the ideal
$\init_{\omega}(J)$ with respect to the monomial ordering
used to compute $H$ during the previous iteration of the
{\bf repeat} \ldots {\bf until} loop,
as $\omega$ is compatible with this previous ordering.
\item
The fact that any set $H$ constructed by the {\bf for}
loop is a Gr\"obner Basis for $J$ with respect to the monomial
ordering defined by $C$ is proved in both \cite{AGK97}
and \cite{CKMWalk} (where you will also find proofs for the
assertions made in the previous paragraph).
\item
The section of code where we determine the value of $t$ is
where we determine the next step of the walk.
We choose $t$ to be the minimum value of $s$ in the
interval $(0, 1]$ such that, for some polynomial
$h \in H$, the $\omega$-degrees of $\LT(h)$ and
some other term in $h$ differ,
but the $\omega(s)$-degrees of the
same two terms are identical. We say that this
is the first point on the line segment between
the two weight vectors $\omega$
and $\tau$ where the initial of one of the
polynomials in $H$ {\it degenerates}.
\item
The success of the Gr\"obner Walk comes from the
fact that it breaks down a Gr\"obner Basis
computation into a series of smaller pieces,
each of which computes a Gr\"obner Basis for a
set of initials, a task that is usually quite
simple. There are still cases however where
this task is complicated and time-consuming,
and this has led to the development of so-called
{\it path perturbation} techniques that choose
`easier' paths on which to walk (see for
example \cite{AGK97} and \cite{Tran00}).
\end{itemize}

\subsection{A Worked Example}

\begin{example}
Let $F := \{xy-z, \; yz + 2x + z\}$ be a basis
generating an ideal $J$ over the polynomial ring
$\mathbb{Q}[x, y, z]$. Consider that we want to obtain the
Lex Gr\"obner Basis $H := \{2x + yz + z, \;
y^2z + yz + 2z\}$ for $J$ from the
DegLex Gr\"obner Basis $G := \{xy-z, \; yz + 2x + z, \;
2x^2 + xz + z^2\}$ for $J$ using the Gr\"obner Walk. Utilising
Algorithm \ref{com-walk-grob} to do this, we initialise
the variables as follows.
\begin{center}
$A = \left(
\begin{array}{ccc}
1 & 1 & 1 \\
1 & 0 & 0 \\
0 & 1 & 0
\end{array} \right)$;
$B = \left(
\begin{array}{ccc}
1 & 0 & 0 \\
0 & 1 & 0 \\
0 & 0 & 1
\end{array} \right)$;
$\omega = (1, 1, 1)$;
$\tau = (1, 0, 0)$; 
$C = \left(
\begin{array}{ccc}
1 & 1 & 1 \\
1 & 0 & 0 \\
0 & 1 & 0 \\
0 & 0 & 1
\end{array} \right)$;
$t = 0$; 
found = false.
\end{center}
Let us now describe what happens during each
pass of the {\bf repeat}\ldots {\bf until} loop
of Algorithm \ref{com-walk-grob}, noting that as
$A$ is equivalent to $C$ to begin with, nothing
substantial will happen during the first pass
through the loop.

{\bf Pass 1}
\begin{itemize}
\item
Construct the set of initials: $G' :=
\{g'_1, \, g'_2, \, g'_3\} = \{xy, \; yz, \;
2x^2 + xz + z^2\}$ (these are the terms
in $G$ that have maximal $(1,1,1)$-degree).
\item
Compute the Gr\"obner Basis $H'$ of $G'$ with respect to $C$.
\begin{eqnarray*}
\mathrm{S\mbox{-}pol}(g'_1, g'_2)
& = & \frac{xyz}{xy}(xy) - \frac{xyz}{yz}(yz) \\
        & = & 0; \\
\mathrm{S\mbox{-}pol}(g'_1, g'_3)
& = & \frac{x^2y}{xy}(xy) -
\frac{x^2y}{2x^2}(2x^2 + xz + z^2)  \\
& = & -\frac{1}{2}xyz - \frac{1}{2}yz^2 \\
& \rightarrow_{g'_1} & -\frac{1}{2}yz^2 \\
& \rightarrow_{g'_2} & 0; \\
\mathrm{S\mbox{-}pol}(g'_2, g'_3)
& = & 0 \; (\mbox{by Buchberger's First Criterion}).
\end{eqnarray*}
It follows that $H' = G'$.
\item
As $H' = G'$, $H$ will also be equal
to $G$, so that $H := \{h_1, \, h_2, \, h_3\} =
\{xy-z, \; yz + 2x + z, \;
2x^2 + xz + z^2\}$.
\item
Let
\begin{eqnarray*}
\omega(s) & := & \omega + s(\tau - \omega) \\
& = & (1,1,1) + s((1,0,0) - (1,1,1)) \\
& = & (1,1,1) + s(0, -1, -1) \\
& = & (1, 1-s, 1-s).
\end{eqnarray*}
To find the next value of $t$, we must find the
minimum value of $s$ such that the $\omega(s)$-degrees
of the leading term of a polynomial in $H$
and some other term in the same polynomial agree where
their $\omega$-degrees currently differ.

The $\omega$-degrees of the two
terms in $h_1$ differ, so we can seek a value of
$s$ such that
\begin{eqnarray*}
\deg_{\omega(s)}(xy) & = & \deg_{\omega(s)}(z) \\
1 + (1-s) & = & (1-s) \\
1 & = & 0 \; (\mathrm{inconsistent}).
\end{eqnarray*}
For $h_2$, we have two choices: either
\begin{eqnarray*}
\deg_{\omega(s)}(yz) & = & \deg_{\omega(s)}(x) \\
(1-s) + (1-s) & = & 1 \\
2-2s & = & 1 \\
s & = & \frac{1}{2};
\end{eqnarray*}
or
\begin{eqnarray*}
\deg_{\omega(s)}(yz) & = & \deg_{\omega(s)}(z) \\
(1-s) + (1-s) & = & (1-s) \\
(1-s) & = & 0 \\
s & = & 1.
\end{eqnarray*}
The $\omega$-degrees of all the terms in $h_3$
are the same, so we can ignore it.

It follows that the minimum value of $s$
(and hence the new value of $t$)
is $\frac{1}{2}$. As this value appears
in the interval $(0, 1]$, we set $G = H$;
set the new value of $\omega$ to be
$(1-\frac{1}{2})(1,1,1) + \frac{1}{2}(1,0,0) =
(1, \frac{1}{2}, \frac{1}{2})$
(and hence change $C$ to be the matrix
$\left(
\begin{array}{ccc}
1 & \frac{1}{2} & \frac{1}{2} \\
1 & 0 & 0 \\
0 & 1 & 0 \\
0 & 0 & 1
\end{array} \right)$); and embark upon
a second pass of the {\bf repeat}\ldots {\bf until} loop.
\end{itemize}

{\bf Pass 2}
\begin{itemize}
\item
Construct the set of initials: $G' :=
\{g'_1, \, g'_2, \, g'_3\} = \{xy, \; 2x+yz, \; 2x^2\}$
(these are the terms
in $G$ that have maximal $(1,\frac{1}{2},\frac{1}{2})$-degree).
\item
Compute the Gr\"obner Basis of $G'$ with respect to $C$.
\begin{eqnarray*}
\mathrm{S\mbox{-}pol}(g'_1, g'_2)
& = & \frac{xy}{xy}(xy) - \frac{xy}{2x}(2x + yz) \\
& = & -\frac{1}{2}y^2z =: g'_4; \\
\mathrm{S\mbox{-}pol}(g'_1, g'_3)
& = & \frac{x^2y}{xy}(xy)
- \frac{x^2y}{2x^2}(2x^2) \\
& = & 0; \\
\mathrm{S\mbox{-}pol}(g'_2, g'_3)
& = & \frac{x^2}{2x}(2x + yz) - \frac{x^2}{2x^2}(2x^2) \\
& = & \frac{1}{2}xyz \\
& \rightarrow_{g'_1} & 0; \\
\mathrm{S\mbox{-}pol}(g'_1, g'_4)
& = & \frac{xy^2z}{xy}(xy) -
\frac{xy^2z}{-\frac{1}{2}y^2z}\left(-\frac{1}{2}y^2z\right) \\
& = & 0; \\
\mathrm{S\mbox{-}pol}(g'_2, g'_4)
& = & 0 \; (\mbox{by Buchberger's First Criterion}); \\
\mathrm{S\mbox{-}pol}(g'_3, g'_4)
& = & 0 \; (\mbox{by Buchberger's First Criterion}).
\end{eqnarray*}
It follows that $G' = \{g'_1, \, g'_2, \, g'_3, \, g'_4\} =
\{xy, \; 2x+yz, \; 2x^2, \; -\frac{1}{2}y^2z\}$ is a
Gr\"obner Basis for $\init_{\omega}(J)$ with respect to $C$.

Applying Algorithm \ref{red-com} to $G'$, we can remove
$g'_1$ and $g'_3$ from $G'$ (because $\LM(g'_1) = y\times \LM(g'_2)$
and $\LM(g'_3) = x\times \LM(g'_2)$); we must also multiply
$g'_2$ and $g'_4$ by $\frac{1}{2}$ and $-2$ respectively
to obtain unit lead coefficients. This leaves us with the
unique reduced Gr\"obner Basis
$H' := \{h'_1, \, h'_2\} =
\{x+\frac{1}{2}yz, \; y^2z\}$ for $\init_{\omega}(J)$
with respect to $C$.
% This Gr\"obner Basis can be further reduced as
% follows:
% \begin{eqnarray*}
% h'_1 = xy & \rightarrow_{h'_2} & -\frac{1}{2}y^2z \\
% & \rightarrow_{h'_4} & 0; \\
% h'_3 = 2x^2 & \rightarrow_{h'_2} & -xyz \\
% & \rightarrow_{h'_2} & \frac{1}{2}y^2z^2 \\
% & \rightarrow_{h'_4} & 0.
% \end{eqnarray*}
% We therefore end up with
% $H' = \{h'_2, \, h'_4\} = \{2x+yz, \; -\frac{1}{2}y^2z\}$.
\item
We must now express the two elements of $H'$ in
terms of members of $G'$.
\begin{eqnarray*}
h'_1 = x + \frac{1}{2}yz & = & \frac{1}{2}g'_2; \\
h'_2 = y^2z & = & -2\left((xy) - \frac{1}{2}y(2x+yz)\right) \;
\mbox{(from the S-polynomial)} \\
& = & -2\left(g'_1 - \frac{1}{2}yg'_2\right).
\end{eqnarray*}
Lifting to the full polynomials, $h'_1$ lifts
to give the polynomial $h_1 := x + \frac{1}{2}yz + \frac{1}{2}z$; $h'_2$
lifts to give the polynomial $h_2 := -2((xy-z) - \frac{1}{2}y(2x + yz + z))
= -2xy+2z + 2xy +y^2z +yz =
y^2z + yz +2z$; and
we are left with the Gr\"obner Basis
$H := \{h_1, \, h_2\} =
\{x + \frac{1}{2}yz + \frac{1}{2}z,
\; y^2z +yz +2z\}$ for $J$ with respect to $C$.
\item
Let
\begin{eqnarray*}
\omega(s) & := & \omega + s(\tau - \omega) \\
& = & \left(1,\frac{1}{2},\frac{1}{2}\right) + s\left((1,0,0) -
\left(1,\frac{1}{2},\frac{1}{2}\right)\right) \\
& = & \left(1,\frac{1}{2},\frac{1}{2}\right) +
s\left(0, -\frac{1}{2}, -\frac{1}{2}\right) \\
& = & \left(1, \frac{1}{2}(1-s), \frac{1}{2}(1-s)\right).
\end{eqnarray*}
Finding the minimum value of $s$,
for $h_1$ we can have
\begin{eqnarray*}
\deg_{\omega(s)}(x) & = & \deg_{\omega(s)}(z) \\
1 & = & \frac{1}{2}(1-s) \\
s & = & -1 \; \mbox{(undefined: we must have $s \in (0, 1]$)}.
\end{eqnarray*}
Continuing with $h_2$, we can either have
\begin{eqnarray*}
\deg_{\omega(s)}(y^2z) & = & \deg_{\omega(s)}(yz) \\
3\left(\frac{1}{2}(1-s)\right) & = &
2\left(\frac{1}{2}(1-s)\right) \\
\frac{1}{2}(1-s) & = & 0 \\
s & = & 1;
\end{eqnarray*}
or
\begin{eqnarray*}
\deg_{\omega(s)}(y^2z) & = & \deg_{\omega(s)}(z) \\
3\left(\frac{1}{2}(1-s)\right) & = & \frac{1}{2}(1-s) \\
1-s & = & 0 \\
s & = & 1.
\end{eqnarray*}
It follows that the minimum value of $s$
(and hence the new value of $t$)
is $1$. As this value appears
in the interval $(0, 1]$, we set $G = H$;
set the new value of $\omega$ to be
$(1-1)(1, \frac{1}{2}, \frac{1}{2}) + 1(1,0,0) =
(1, 0, 0)$ (and hence change $C$ to be the matrix
$\left(
\begin{array}{ccc}
1 & 0 & 0 \\
1 & 0 & 0 \\
0 & 1 & 0 \\
0 & 0 & 1
\end{array} \right)
\equiv
\left(
\begin{array}{ccc}
1 & 0 & 0 \\
0 & 1 & 0 \\
0 & 0 & 1
\end{array} \right)
$); and embark upon a third (and final)
pass of the {\bf repeat}\ldots {\bf until} loop.
\end{itemize}

{\bf Pass 3}
\begin{itemize}
\item
Construct the set of initials: $G' :=
\{g'_1, \, g'_2\} = \{x, \; y^2z + yz +2z\}$
(these are the terms
in $G$ that have maximal $(1,0,0)$-degree).
\item
Compute the Gr\"obner Basis $H'$ of $G'$ with respect to $C$.
$$
\mathrm{S\mbox{-}pol}(g'_1, g'_2)
= 0 \; (\mbox{by Buchberger's First Criterion}).
$$
It follows that $H' = G'$.
\item
As $H' = G'$, $H$ will also be equal
to $G$, so that $H := \{h_1, \, h_2\} =
\{x + \frac{1}{2}yz + \frac{1}{2}z, \; y^2z + yz +2z\}$.
Further, as $t$ is now equal to $1$, we
have arrived at our target ordering (Lex)
and can return $H$ as the output Gr\"obner Basis, a basis
that is equivalent to the Lex Gr\"obner Basis
given for $J$ at the beginning of this example.
\end{itemize}
We can summarise the path taken during the walk
in the following diagram.
\begin{center}
\begin{picture}(0,0)%
\includegraphics{ch6d1.pstex}%
\end{picture}%
\setlength{\unitlength}{2605sp}%
\begingroup\makeatletter\ifx\SetFigFont\undefined%
\gdef\SetFigFont#1#2#3#4#5{%
  \reset@font\fontsize{#1}{#2pt}%
  \fontfamily{#3}\fontseries{#4}\fontshape{#5}%
  \selectfont}%
\fi\endgroup%
\begin{picture}(3960,5064)(901,-6583)
\put(4426,-3526){\makebox(0,0)[lb]{\smash{\SetFigFont{12}{14.4}{\rmdefault}{\mddefault}{\updefault}{\color[rgb]{0,0,0}$(1,1,1)$}%
}}}
\put(4426,-3211){\makebox(0,0)[lb]{\smash{\SetFigFont{12}{14.4}{\rmdefault}{\mddefault}{\updefault}{\color[rgb]{0,0,0}Pass 1}%
}}}
\put(3436,-5453){\makebox(0,0)[lb]{\smash{\SetFigFont{12}{14.4}{\rmdefault}{\mddefault}{\updefault}{\color[rgb]{0,0,0}Pass 3}%
}}}
\put(3436,-5768){\makebox(0,0)[lb]{\smash{\SetFigFont{12}{14.4}{\rmdefault}{\mddefault}{\updefault}{\color[rgb]{0,0,0}$(1,0,0)$}%
}}}
\put(3923,-4486){\makebox(0,0)[lb]{\smash{\SetFigFont{12}{14.4}{\rmdefault}{\mddefault}{\updefault}{\color[rgb]{0,0,0}Pass 2}%
}}}
\put(3923,-4801){\makebox(0,0)[lb]{\smash{\SetFigFont{12}{14.4}{\rmdefault}{\mddefault}{\updefault}{\color[rgb]{0,0,0}$(1,\frac{1}{2},\frac{1}{2})$}%
}}}
\put(3976,-2461){\makebox(0,0)[lb]{\smash{\SetFigFont{12}{14.4}{\familydefault}{\mddefault}{\updefault}{\color[rgb]{0,0,0}$y$}%
}}}
\put(4861,-6511){\makebox(0,0)[lb]{\smash{\SetFigFont{12}{14.4}{\familydefault}{\mddefault}{\updefault}{\color[rgb]{0,0,0}$x$}%
}}}
\put(901,-4636){\makebox(0,0)[lb]{\smash{\SetFigFont{12}{14.4}{\rmdefault}{\mddefault}{\updefault}{\color[rgb]{0,0,0}$\frac{1}{2}$}%
}}}
\put(1126,-1711){\makebox(0,0)[lb]{\smash{\SetFigFont{12}{14.4}{\familydefault}{\mddefault}{\updefault}{\color[rgb]{0,0,0}$z$}%
}}}
\put(2243,-3834){\makebox(0,0)[lb]{\smash{\SetFigFont{12}{14.4}{\rmdefault}{\mddefault}{\updefault}{\color[rgb]{0,0,0}1}%
}}}
\put(909,-4036){\makebox(0,0)[lb]{\smash{\SetFigFont{12}{14.4}{\rmdefault}{\mddefault}{\updefault}{\color[rgb]{0,0,0}1}%
}}}
\put(2919,-6151){\makebox(0,0)[lb]{\smash{\SetFigFont{12}{14.4}{\rmdefault}{\mddefault}{\updefault}{\color[rgb]{0,0,0}1}%
}}}
\put(1959,-5836){\makebox(0,0)[lb]{\smash{\SetFigFont{12}{14.4}{\rmdefault}{\mddefault}{\updefault}{\color[rgb]{0,0,0}$\frac{1}{2}$}%
}}}
\put(1671,-4351){\makebox(0,0)[lb]{\smash{\SetFigFont{12}{14.4}{\rmdefault}{\mddefault}{\updefault}{\color[rgb]{0,0,0}$\frac{1}{2}$}%
}}}
\end{picture}

\end{center}
\end{example}

\begin{algorithm}
\setlength{\baselineskip}{3.5ex}
\caption{The Commutative Involutive Walk Algorithm}
\label{com-walk-inv}
\begin{algorithmic}
\vspace*{2mm}
\REQUIRE{An Involutive Basis $G = \{g_1, g_2, \hdots, g_m\}$
         with respect to an involutive division $I$ and
         an admissible monomial ordering $O$ with
         an associated matrix $A$, where $G$ generates an
         ideal $J$ over a commutative polynomial ring
         $\mathcal{R} = R[x_1, \hdots, x_n]$.}
\ENSURE{An Involutive Basis $H = \{h_1, h_2, \hdots, h_p\}$
        for $J$ with respect to $I$ and an admissible
        monomial ordering $O'$
        with an associated matrix $B$.}
\vspace*{1mm}
\STATE
Let $\omega$ and $\tau$ be the weight vectors corresponding to
the first rows of $A$ and $B$; \\ % respectively; \\
Let $C$ be the matrix whose first row is equal to
$\omega$ and whose remainder is equal to the whole of
the matrix $B$; \\
$t = 0$; found = false;
\REPEAT
\STATE
Let $G' = \{\init_{\omega}(g_i)\}$ for all $g_i \in G$; \\
Compute an Involutive Basis $H'$ for $G'$ with respect
to the monomial ordering defined by the matrix $C$
(use Algorithm \ref{com-inv});  \\
$H = \emptyset$; \\
\FOR{{\bf each} $h' \in H'$}
\STATE
Let $\sum_{i=1}^{j} p_ig'_i$ be the logged representation
of $h'$ with respect to $G'$
(where $g'_i \in G'$ and $p_i \in \mathcal{R}$),
obtained either through computing a Logged Involutive Basis
above or by involutively dividing $h'$ with respect to $G'$; \\
$H = H \cup \{\sum_{i=1}^{j} p_ig_i\}$, where
$\init_{\omega}(g_i) = g'_i$; \\
\ENDFOR
\IF{($t == 1$)}
\STATE
found = true;
\ELSE
\STATE
$t = \min (\{ s \mid \mathrm{deg}_{\omega(s)}(p_1) =
\mathrm{deg}_{\omega(s)}(p_i),
\mathrm{deg}_{\omega(0)}(p_1) \neq
\mathrm{deg}_{\omega(0)}(p_i),$ \\
$h = p_1 + \cdots + p_k \in H\} \cap (0, 1])$,
where $\omega(s) := \omega + s(\tau
- \omega)$ for $0 \leqslant s \leqslant 1$;
\ENDIF
\IF{($t$ is undefined)}
\STATE
found = true;
\ELSE
\STATE
$G = H$; $\omega = (1-t)\omega + t\tau$;
\ENDIF
\UNTIL{(found = true)}
\STATE
{\bf return} $H$;
\end{algorithmic}
\vspace*{1mm}
\end{algorithm}

\subsection{The Commutative Involutive Walk Algorithm}

In \cite{Golub}, Golubitsky generalised the
Gr\"obner Walk technique to give a method for converting an
Involutive Basis with respect to one monomial
ordering to an Involutive Basis with respect to
another monomial ordering. Algorithmically, the way
in which we perform this {\it Involutive Walk}
\index{involutive walk} is 
\index{walk!involutive} virtually identical to the
way we perform the Gr\"obner Walk, as can be seen by
comparing Algorithms \ref{com-walk-grob} and
\ref{com-walk-inv}. The change however comes when
proving the correctness of the algorithm, as we
have to show that each $G'$ is an Involutive Basis
for $\init_{\omega}(J)$ and that each $H$ is an
Involutive Basis for $J$ (see \cite{Golub}
for these proofs).

\section{Noncommutative Walks} \label{6point2}

In the commutative case, any monomial ordering can be
represented by a matrix that provides a decomposition
of the ordering in terms of the rows of the matrix.
This decomposition is then utilised in the Gr\"obner
Walk algorithm when (for example) we use the first row
of the matrix to provide a set of initials for
a particular basis $G$
(cf. Definition \ref{initial-defn-ch6}).

In the noncommutative case, because monomials cannot be
represented by multidegrees, monomial orderings cannot
be represented by matrices. This seems to shut the door
on any generalisation of the Gr\"obner Walk to the
noncommutative case, as not only is there no first row
of a matrix to provide a set of initials, but no notion
of a walk between two matrices can be formulated
either.

Despite this, we note that in the commutative case, if
the first rows of the source and target matrices are
the same, then the Gr\"obner Walk will complete in one
pass of the algorithm, and all that is needed is
the first row of the source matrix to provide a set
of initials to work with.

Generalising to the noncommutative case,
it is possible that if we can find a way to decompose
a noncommutative monomial ordering to provide a set
of initials to work with, then a noncommutative
Gr\"obner Walk algorithm could complete
in one pass if the source and target monomial orderings
used the same method to compute sets of initials.

% Let us now consider one way of decomposing a
% noncommutative monomial ordering, based on breaking down
% an             ordering into a series of functions.

\subsection{Functional Decompositions}

Considering the monomial orderings defined in
Section \ref{NCMO}, we note that all the orderings
are defined step-by-step. For example,
the DegLex monomial ordering compares two monomials
by degree first, then by the first letter of each
monomial, then by the second letter, and so on.
This provides us with an opportunity to decompose
each monomial ordering into a series of steps or
functions, a decomposition we shall term a
{\it functional decomposition}.

\begin{defn}
An {\it ordering function} \index{ordering function}
is a function
$$\theta: \: M \longrightarrow \mathbb{Z}$$
that assigns an integer to any monomial $m \in M$,
where $M$ denotes the set of all monomials over a polynomial
ring $R\langle x_1, \hdots, x_n\rangle$.
We call the integer assigned by $\theta$ to $m$ the
$\theta$-value of $m$.
\end{defn}

\begin{remark}
% The integer an ordering function
% $\theta$ assigns to a monomial $m$ will be referred to
% as the $\theta$-value of $m$. Further,
The $\theta$-value of any term will be equal
to the $\theta$-value of the term's associated
monomial.
\end{remark}

\begin{defn}
A {\it functional decomposition} \index{functional
decomposition} $\Theta$ is a (possibly infinite) sequence
of ordering functions, written $\Theta = \{\theta_1,
\theta_2, \hdots\}.$
\end{defn}

\begin{defn}
Let $O$ be a monomial ordering; let $m_1$ and
$m_2$ be two arbitrary monomials such that
$m_1 < m_2$ with respect to $O$; and let
$\Theta = \{\theta_1, \theta_2, \hdots\}$ be a
functional decomposition. We say that
$\Theta$ defines $O$ if and
only if $\theta_i(m_1) < \theta_i(m_2)$ for some
$i \geqslant 1$ and $\theta_j(m_1) = \theta_j(m_2)$
for all $1 \leqslant j < i$.
\end{defn}

To describe the functional decompositions corresponding
to the monomial orderings defined in Section \ref{NCMO},
we first need the following definition.

\begin{defn}
Let $m$ be an arbitrary monomial over a polynomial
ring $R\langle x_1, \hdots, x_n\rangle$. The $i$-th
valuing function of $m$, written $\val_i(m)$,
\index{$V$@$\val_i$} is an ordering function that assigns
an integer to $m$ as follows.
$$\val_i(m) =
\begin{cases}
j & \text{if $\SUB(m, i, i) = x_j$ (where
$1 \leqslant j \leqslant n$).} \\
n+1 & \text{if $\SUB(m, i, i)$ is undefined.}
\end{cases}
$$
\end{defn}

Let us now describe the functional decompositions
corresponding to those monomial orderings defined
in Section \ref{NCMO} that are admissible.

\begin{defn}
The functional decomposition
$\Theta = \{\theta_1, \theta_2, \hdots\}$
corresponding to the DegLex monomial ordering
is defined (for an arbitrary monomial $m$) as follows.
$$\theta_i(m) =
\begin{cases}
\deg(m) & \text{if $i = 1$.} \\
n+1-\val_{i-1}(m) & \text{if $i > 1$.}
\end{cases}
$$
Similarly, we can define DegInvLex by
$$\theta_i(m) =
\begin{cases}
\deg(m) & \text{if $i = 1$.} \\
\val_{i-1}(m) & \text{if $i > 1$.}
\end{cases}
$$
and DegRevLex by
$$\theta_i(m) =
\begin{cases}
\deg(m) & \text{if $i = 1$.} \\
\val_{\deg(m)+2-i}(m) & \text{if $i > 1$.}
\end{cases}
$$
\end{defn}

\begin{example}
Let $m_1 := xyxz^2$ and $m_2 := xzyz^2$ be two
monomials over the polynomial ring $\mathbb{Q}\langle
x, y, z\rangle$. With respect to DegLex, we can work
out that $xyxz^2 > xzyz^2$, because $\theta_1(m_1) =
\theta_1(m_2)$ (or $\deg(m_1) = \deg(m_2)$);
$\theta_2(m_1) = \theta_2(m_2)$ (or $n+1-\val_1(m_1)
= n+1-\val_1(m_2)$, $3+1-1 = 3+1-1$); and
$\theta_3(m_1) > \theta_3(m_2)$ (or $n+1-\val_2(m_1)
> n+1-\val_2(m_2)$, $3+1-2 > 3+1-3)$.
Similarly, with respect to DegInvLex, we can work
out that $xyxz^2 < xzyz^2$ (because $\theta_3(m_1)
< \theta_3(m_2)$, or $2 < 3$); and with respect
to DegRevLex, we can work out that
$xyxz^2 < xzyz^2$ (because $\theta_4(m_1) <
\theta_4(m_2)$, or $1 < 2$).
\end{example}

\begin{defn}
Given a polynomial $p$ and an ordering function $\theta$,
the \index{initial} {\it initial} of $p$ with respect to
$\theta$, written \index{$I$@$\init_{\theta}$}
$\init_{\theta}(p)$, is
made up of those terms in $p$ that have
maximal $\theta$-value. For example,
if $\theta$ is the degree function and if
$p = x^4 + zxy^2 + y^3 + z^2x$, then
$\init_{\theta}(p) = x^4 + zxy^2$.
\end{defn}

\begin{defn}
Given an ordering function $\theta$,
a polynomial $p$ is said to be
\index{$T$@$\theta$-homogeneous}
{\it $\theta$-homogeneous} if
$p = \init_{\theta}(p)$.
\end{defn}

\begin{defn}
An ordering function $\theta$ is
\index{ordering function!compatible}
{\it compatible} \index{compatible ordering function} with a monomial
ordering $O$ if, given any polynomial $p = t_1
+ \cdots + t_m$ ordered in descending order with
respect to $O$, $\theta(t_1) \geqslant
\theta(t_i)$ holds for all $1 < i \leqslant m$.
\end{defn}

\begin{defn}
An ordering function $\theta$ is
\index{ordering function!extendible}
{\it extendible} \index{extendible ordering function} if, given any
$\theta$-homogeneous polynomial $p$,
% made up of terms that all have the same
% $\theta$-value, then
any multiple $upv$ of
$p$ by terms $u$ and $v$ is also
$\theta$-homogeneous.
% also contains only
% terms that have the same $\theta$-value.
\end{defn}

\begin{remark}
Of the ordering functions encountered so far, only
the degree function, $\val_1$ and\footnote{Think of 
$\val_{\deg(m)}$ as finding the value of the
final variable in $m$ (as opposed to $\val_1$ finding
the value of the first variable in $m$).} $\val_{\deg(m)}$
(for any given monomial $m$) are extendible.
\end{remark}

\begin{defn}
Two noncommutative monomial orderings $O_1$ and
$O_2$ are said to be {\it harmonious}
\index{monomial ordering!harmonious} if 
\index{harmonious monomial ordering} (i) there exist
functional decompositions
$\Theta_1 = \{\theta_{1_1}, \theta_{1_2}, \hdots\}$ and
$\Theta_2 = \{\theta_{2_1}, \theta_{2_2}, \hdots\}$
defining $O_1$ and $O_2$ respectively; and (ii)
the ordering functions $\theta_{1_1}$ and
$\theta_{2_1}$ are identical and extendible.
\end{defn}

\begin{remark}
The noncommutative monomial orderings DegLex,
DegInvLex and DegRevLex are all (pairwise)
harmonious.
\end{remark}

\subsection{The Noncommutative Gr\"obner Walk Algorithm
for Harmonious Monomial Orderings}

We present in Algorithm \ref{noncom-walk-grob} an
algorithm to perform a Gr\"obner Walk 
\index{walk!Gr\"obner} between two
harmonious noncommutative monomial orderings.
\index{Gr\"obner walk}

\begin{algorithm}
\setlength{\baselineskip}{3.5ex}
\caption{The Noncommutative Gr\"obner Walk Algorithm
for Harmonious Monomial Orderings}
\label{noncom-walk-grob}
\begin{algorithmic}
\vspace*{2mm}
\REQUIRE{A Gr\"obner Basis $G = \{g_1, g_2, \hdots, g_m\}$
         with respect to an admissible monomial ordering $O$ with
         an associated functional decomposition $A$,
         where $G$ generates an
         ideal $J$ over a noncommutative polynomial ring
         $\mathcal{R} = R\langle x_1, \hdots, x_n\rangle$.}
\ENSURE{A Gr\"obner Basis $H = \{h_1, h_2, \hdots, h_p\}$
        for $J$ with respect to an admissible monomial ordering $O'$
        with an associated functional decomposition $B$, where $O$ and $O'$
        are harmonious.}
\vspace*{1mm}
\STATE
Let $\theta$ be the ordering function corresponding to
the first ordering function of $A$; \\[1mm]
% $t = 0$; found = false; \\
% \REPEAT
% \STATE
Let $G' = \{\init_{\theta}(g_i)\}$ for all $g_i \in G$; \\
Compute a reduced Gr\"obner Basis $H'$ for $G'$ with respect
to the monomial ordering $O'$
(use Algorithms \ref{noncom-buch} and \ref{red-noncom});  \\
$H = \emptyset$; \\
\FOR{{\bf each} $h' \in H'$}
\STATE
Let $\sum_{i=1}^{j} \ell_ig'_ir_i$ be the logged representation
of $h'$ with respect to $G'$
(where $g'_i \in G'$ and the $\ell_i$ and the $r_i$ are terms),
obtained either through computing a Logged Gr\"obner Basis
above or by dividing $h'$ with respect to $G'$; \\
$H = H \cup \{\sum_{i=1}^{j} \ell_ig_ir_i\}$, where
$\init_{\theta}(g_i) = g'_i$; \\
\ENDFOR
\STATE
Reduce $H$ with respect to $O'$
(use Algorithm \ref{red-noncom}); \\[1mm]
% \IF{($t == 1$)}
% \STATE
% found = true;
% \ELSE
% \STATE
% FIND NEXT INITIAL FUNCTION
% \ENDIF
% \IF{($t$ is undefined)}
% \STATE
% found = true;
% \ELSE
% \STATE
% $G = H$; $\omega = (1-t)\omega + t\tau$;
% \ENDIF
% \UNTIL{(found = true)}
% \STATE
{\bf return} $H$;
\end{algorithmic}
\vspace*{1mm}
\end{algorithm}

Termination of Algorithm \ref{noncom-walk-grob} depends on the
termination of Algorithm \ref{noncom-buch} as used
(in Algorithm \ref{noncom-walk-grob})
to compute a noncommutative Gr\"obner Basis
for the set $G'$. The correctness of Algorithm \ref{noncom-walk-grob}
is provided by the following two propositions.

\begin{prop} \label{inGrob}
$G'$ is always a Gr\"obner Basis for
the ideal\footnote{The ideal $\init_{\theta}(J)$
is defined as follows: $p \in J$ if and only if
$\init_{\theta}(p) \in \init_{\theta}(J)$.}
$\init_{\theta}(J)$ with respect to
the monomial ordering $O$.
\end{prop}
\begin{pf}
Because $\theta$ is compatible with $O$ (by
definition), the S-polynomials involving members
of $G$ will be in one-to-one correspondence with the
S-polynomials involving members of $G'$, with the
same monomial being `cancelled' in each pair
of corresponding S-polynomials.

Let $p$ be an arbitrary S-polynomial involving
members of $G$ (with corresponding S-polynomial
$q$ involving members of $G'$).
Because $G$ is a Gr\"obner Basis for $J$ with respect to
$O$, $p$ will reduce to zero using $G$ by
the series of reductions
$$p \rightarrow_{g_{i_1}} p_1 \rightarrow_{g_{i_2}} p_2
\rightarrow_{g_{i_3}} \cdots \rightarrow_{g_{i_{\alpha}}} 0,$$
where $g_{i_j} \in G$ for all $1 \leqslant j \leqslant \alpha$.

{\bf Claim:}
$q$ will reduce to zero using $G'$ (and hence
$G'$ is a Gr\"obner Basis for $\init_{\theta}(J)$
with respect to $O$ by Definition 
\ref{grob-defn-noncom}) by the series of reductions
$$q \rightarrow_{\init_{\theta}(g_{i_1})} q_1
\rightarrow_{\init_{\theta}(g_{i_2})} q_2
\rightarrow_{\init_{\theta}(g_{i_3})} \cdots
\rightarrow_{\init_{\theta}(g_{i_{\beta}})} 0,$$
where $0 \leqslant \beta \leqslant \alpha$.

{\bf Proof of Claim:}
Let $w$ be the overlap word associated to the
S-polynomial $p$. If $\theta(\LM(p)) < \theta(w)$,
then because $\theta$ is extendible it is clear that $q =  0$,
and so the proof is complete.
Otherwise, we must have $q = \init_{\theta}(p)$,
and so by the compatibility of $\theta$ with $O$,
we can use the polynomial $\init_{\theta}(g_{i_1}) \in G'$
to reduce $q$ to give the polynomial $q_1$.
We now proceed by induction (if $\theta(\LM(p_1)) <
\theta(\LM(p))$ then $q_1 = 0$, \ldots), noting that
the process will terminate because 
$\init_{\theta}(p_{\alpha} = 0) = 0$.
\end{pf}

\begin{prop}
The set $H$ constructed by the {\bf for}
loop of Algorithm \ref{noncom-walk-grob}
is a Gr\"obner Basis for $J$ with respect
to the monomial ordering $O'$.
\end{prop}
\begin{pf}
By Definition \ref{grob-defn-noncom}, we can show that $H$
is a Gr\"obner Basis for $J$ by showing that all
S-polynomials involving members of $H$ reduce to zero
using $H$. Assume for a contradiction that an
S-polynomial $p$ involving members of $H$
does not reduce to zero using $H$, and
instead only reduces to a polynomial $q \neq 0$.

As all members of $H$ are members of the ideal $J$
(by the way $H$ was constructed as combinations of
elements of $G$), it is clear that $q$ is also
a member of the ideal $J$, as all we have done in
constructing $q$ is to reduce a combination of two
members of $H$ with respect to $H$. It follows
that the polynomial $\init_{\theta}(q)$ is a
member of the ideal $\init_{\theta}(J)$.

Because $H'$ is a Gr\"obner Basis for the ideal
$\init_{\theta}(J)$ with respect to $O'$, there must be a
polynomial $h' \in H'$ such that $h' \mid
\init_{\theta}(q)$. Let
$\sum_{i=1}^{j} \ell_ig'_ir_i$ be the logged representation
of $h'$ with respect to $G'$. Then it is clear that
$$\sum_{i=1}^{j} \ell_ig'_ir_i \mid \init_{\theta}(q).$$
However $\theta$ is compatible with $O'$, so that
$$\sum_{i=1}^{j} \ell_ig_ir_i \mid q.$$
It follows that there exists a polynomial $h \in H$
dividing our polynomial $q$, contradicting our initial
assumption.
\end{pf}

\subsection{A Worked Example}

\begin{example}
Let $F := \{x^2+y^2+8, \; 2xy+y^2+5\}$ be a basis
generating an ideal $J$ over the polynomial ring
$\mathbb{Q}\langle x, y\rangle$. Consider that we want to obtain the
DegLex Gr\"obner Basis $H := \{2xy+y^2+5, \; x^2+y^2+8, \;
5y^3-10x+37y, \; 2yx+y^2+5\}$ for $J$ from the
DegRevLex Gr\"obner Basis $G := \{2xy-x^2-3, \, y^2+x^2+8, \,
5x^3+6y+35x, \, 2yx-x^2-3\}$ for $J$ using the Gr\"obner Walk. Utilising
Algorithm \ref{noncom-walk-grob} to do this, we initialise
$\theta$ to be the degree function and we proceed as follows.
\begin{itemize}
\item
Construct the set of initials: $G' :=
\{g'_1, \, g'_2, \, g'_3, \, g'_4\} =
\{-x^2+2xy, \; x^2+y^2, \; 5x^3, \; -x^2+2yx\}$
(these are the terms in $G$ that have maximal degree).
% $\theta$-value).
\item
Compute the Gr\"obner Basis of $G'$ with respect to the
DegLex monomial ordering (for simplicity, we will not provide
details of those S-polynomials that reduce to zero or can be
ignored due to Buchberger's Second Criterion).
\begin{eqnarray*}
\mathrm{S\mbox{-}pol}(1, g'_1, 1, g'_2)
& = & (-x^2+2xy) - (-1)(x^2+y^2) \\
& = & 2xy+y^2 =: g'_5; \\
\mathrm{S\mbox{-}pol}(1, g'_1, 1, g'_4)
& = & (-1)(-x^2+2xy) - (-1)(-x^2+2yx) \\
& = & -2xy+2yx \\
& \rightarrow_{g'_5} & -2xy+2yx + (2xy+y^2) \\
& = & 2yx+y^2 =: g'_6; \\
\mathrm{S\mbox{-}pol}(y, g'_1, 1, g'_6)
& = & 2y(-x^2+2xy) - (-1)(2yx+y^2)x \\
& = & 4yxy+y^2x \\
& \rightarrow_{g'_5} & 4yxy+y^2x - 2y(2xy+y^2) \\
& = & y^2x-2y^3 \\
& \rightarrow_{g'_6} & y^2x-2y^3 -\frac{1}{2}y(2yx+y^2) \\
& = & -\frac{5}{2}y^3 =: g'_7.
\end{eqnarray*}
After $g'_7$ is added to the current basis, all
S-polynomials now reduce to zero, leaving the
Gr\"obner Basis $G' = \{g'_1, \, g'_2, \, g'_3, \, g'_4, \,
g'_5, \, g'_6, \, g'_7\} =
\{-x^2+2xy, \; x^2+y^2, \; 5x^3, \; -x^2+2yx, \;
2xy+y^2, \; 2yx+y^2, \; -\frac{5}{2}y^3\}$ for
$\init_{\theta}(J)$ with respect to $O'$.

Applying Algorithm \ref{red-noncom} to $G'$, we can remove
$g'_1$, $g'_2$ and $g'_3$ from $G'$ (because their lead
monomials are all multiplies of $\LM(g'_4)$); we must
multiply $g'_4$, $g'_5$, $g'_6$ and $g'_7$ by
$-1$, $\frac{1}{2}$, $\frac{1}{2}$ and $-\frac{2}{5}$
respectively (to obtain unit lead coefficients); and the
polynomial $g'_4$ can (then) be further reduced as follows.
\begin{eqnarray*}
g'_4 & = & x^2-2yx \\
& \rightarrow_{g'_6} & x^2-2yx +2\left(yx+\frac{1}{2}y^2\right) \\
& = & x^2+y^2.
\end{eqnarray*}
This leaves us with the unique reduced Gr\"obner Basis
$H' := \{h'_1, \, h'_2, \, h'_3, \, h'_4\}
= \{x^2+y^2, \; xy+\frac{1}{2}y^2, \; yx+\frac{1}{2}y^2, \; y^3\}$
for $\init_{\theta}(J)$ with respect to $O'$.
\item
We must now express the four elements of $H'$ in
terms of members of $G'$.
\begin{eqnarray*}
h'_1 = x^2+y^2 & = & g'_2; \\
h'_2 = xy+\frac{1}{2}y^2 & = & \frac{1}{2}(g'_1+g'_2) \;
\mbox{(from the S-polynomial)}; \\
h'_3 = yx+\frac{1}{2}y^2 & = & \frac{1}{2}(-g'_1+g'_4+(g'_1+g'_2)) \\
& = & \frac{1}{2}(g'_2+g'_4); \\
h'_4 = y^3 & = & -\frac{2}{5}\left(2y(g'_1) + (g'_2+g'_4)x
- 2y(g'_1+g'_2) - \frac{1}{2}y(g'_2+g'_4)\right) \\
& = & -\frac{2}{5}\left(g'_2x - \frac{5}{2}yg'_2 
+ g'_4x - \frac{1}{2}yg'_4\right).
\end{eqnarray*}
Lifting to the full polynomials,
we obtain the Gr\"obner Basis
$H := \{h_1, \, h_2, \, h_3, \, h_4\}$
as follows.
\begin{eqnarray*}
h_1 & = & g_2 \\
    & = & x^2+y^2+8; \\
h_2 & = & \frac{1}{2}(g_1+g_2) \\
    & = & \frac{1}{2}(-x^2+2xy-3+x^2+y^2+8) \\
    & = & xy+\frac{1}{2}y^2+\frac{5}{2}; \\
h_3 & = & \frac{1}{2}(g_2+g_4) \\
    & = & \frac{1}{2}(x^2+y^2+8-x^2+2yx-3) \\
    & = & yx+\frac{1}{2}y^2+\frac{5}{2}; \\
h_4 & = & -\frac{2}{5}\left(g_2x - \frac{5}{2}yg_2
          + g_4x - \frac{1}{2}yg_4\right) \\
    & = & -\frac{2}{5}\left(x^3+y^2x+8x-\frac{5}{2}yx^2
          -\frac{5}{2}y^3-20y\right. \\
    &   & \left.-x^3+2yx^2-3x+\frac{1}{2}yx^2
          -y^2x+\frac{3}{2}y\right) \\
    & = & y^3-2x+\frac{37}{5}y.
\end{eqnarray*}
The set $H$ does not reduce any further, so we return the
output DegLex Gr\"obner Basis $\{x^2+y^2+8, \,
xy+\frac{1}{2}y^2+\frac{5}{2}, \; yx+\frac{1}{2}y^2+\frac{5}{2}, \;
y^3-2x+\frac{37}{5}y\}$ for $J$, a basis that is
equivalent to the DegLex Gr\"obner Basis given for $J$
at the beginning of this example.
\end{itemize}
\end{example}

\subsection{The Noncommutative Involutive Walk Algorithm
for Harmonious Monomial Orderings}

We present in Algorithm \ref{noncom-walk-inv} an
algorithm to perform an Involutive Walk between two
\index{involutive walk}
harmonious \index{walk!involutive}
noncommutative monomial orderings.

\begin{algorithm}
\setlength{\baselineskip}{3.5ex}
\caption{The Noncommutative Involutive Walk Algorithm
for Harmonious Monomial Orderings}
\label{noncom-walk-inv}
\begin{algorithmic}
\vspace*{2mm}
\REQUIRE{An Involutive Basis $G = \{g_1, g_2, \hdots, g_m\}$
         with respect to an involutive division $I$ and an
         admissible monomial ordering $O$ with
         an associated functional decomposition $A$,
         where $G$ generates an
         ideal $J$ over a noncommutative polynomial ring
         $\mathcal{R} = R\langle x_1, \hdots, x_n\rangle$.}
\ENSURE{An Involutive Basis $H = \{h_1, h_2, \hdots, h_p\}$
        for $J$ with respect to $I$ and
        an admissible monomial ordering $O'$
        with an associated functional decomposition $B$, where
        $O$ and $O'$ are harmonious.}
\vspace*{1mm}
\STATE
Let $\theta$ be the ordering function corresponding to
the first ordering function of $A$; \\[1mm]
% $t = 0$; found = false; \\
% \REPEAT
% \STATE
Let $G' = \{\init_{\theta}(g_i)\}$ for all $g_i \in G$; \\
Compute an Involutive Basis $H'$ for $G'$ with respect
to $I$ and the monomial ordering $O'$
(use Algorithm \ref{noncom-inv});  \\
$H = \emptyset$; \\
\FOR{{\bf each} $h' \in H'$}
\STATE
Let $\sum_{i=1}^{j} \ell_ig'_ir_i$ be the logged representation
of $h'$ with respect to $G'$
(where $g'_i \in G'$ and the $\ell_i$ and the $r_i$ are terms),
obtained either through computing a Logged Involutive Basis
above or by involutively dividing $h'$ with respect to $G'$; \\
$H = H \cup \{\sum_{i=1}^{j} \ell_ig_ir_i\}$, where
$\init_{\theta}(g_i) = g'_i$; \\[1mm]
\ENDFOR
\STATE
% Reduce $H$ with respect to $C$
% (use Algorithm \ref{red-noncom}); \\[1mm]
% \IF{($t == 1$)}
% \STATE
% found = true;
% \ELSE
% \STATE
% FIND NEXT INITIAL FUNCTION
% \ENDIF
% \IF{($t$ is undefined)}
% \STATE
% found = true;
% \ELSE
% \STATE
% $G = H$; $\omega = (1-t)\omega + t\tau$;
% \ENDIF
% \UNTIL{(found = true)}
% \STATE
{\bf return} $H$;
\end{algorithmic}
\vspace*{1mm}
\end{algorithm}

Termination of Algorithm \ref{noncom-walk-inv} depends on the
termination of Algorithm \ref{noncom-inv} as used
(in Algorithm \ref{noncom-walk-inv})
to compute a noncommutative Involutive Basis
for the set $G'$. The correctness of Algorithm \ref{noncom-walk-inv}
is provided by the following two propositions.

\begin{prop} \label{inInv}
$G'$ is always an Involutive Basis for
the ideal $\init_{\theta}(J)$ with respect to $I$ and
the monomial ordering $O$.
\end{prop}
\begin{pf}
Let $p := ugv$ be an arbitrary multiple of a polynomial
$g \in G$ by terms $u$ and $v$.
Because $G$ is an Involutive Basis for $J$ with respect to
$I$ and $O$, $p$ will involutively reduce to zero using
$G$ by the series of involutive reductions
$$p \xymatrix{\ar[r]_I_(1){g_{i_1}} &} p_1
\xymatrix{\ar[r]_I_(1){g_{i_2}} &} p_2
\xymatrix{\ar[r]_I_(1){g_{i_3}} &} \cdots
\xymatrix{\ar[r]_I_(1){g_{i_{\alpha}}} &} 0,$$
where $g_{i_j} \in G$ for all $1 \leqslant j \leqslant \alpha$.

{\bf Claim:}
The polynomial $q := u\init_{\theta}(g)v$ will
involutively reduce to zero using $G'$ (and hence
$G'$ is an Involutive Basis for $\init_{\theta}(J)$
with respect to $I$ and $O$ by Definition \ref{IBNC})
by the series of involutive reductions
$$q
\xymatrix@C=3pc{\ar[r]_(0.4)I_(1){\init_{\theta}(g_{i_1})} &} q_1
\xymatrix@C=3pc{\ar[r]_(0.4)I_(1){\init_{\theta}(g_{i_2})} &} q_2
\xymatrix@C=3pc{\ar[r]_(0.4)I_(1){\init_{\theta}(g_{i_3})} &} \cdots
\xymatrix@C=3pc{\ar[r]_(0.4)I_(1){\init_{\theta}(g_{i_{\beta}})} &} 0,$$
where $1 \leqslant \beta \leqslant \alpha$.

{\bf Proof of Claim:}
Because $\theta$ is extendible, it is clear that
$q = \init_{\theta}(p)$. Further, because $\theta$
is compatible with $O$ (by
definition), the multiplicative variables of
$G$ and $G'$ with respect to $I$ will be identical,
and so it follows that because the polynomial $g_{i_1} \in G$
was used to involutively reduce $p$ to give the polynomial $p_1$,
then the polynomial $\init_{\theta}(g_{i_1}) \in G'$
can be used to involutively reduce $q$ to give the
polynomial $q_1$.

If $\theta(\LM(p_1)) < \theta(\LM(p))$,
then because $\theta$ is extendible it is clear that $q_1 = 0$,
and so the proof is complete.
Otherwise, we must have $q_1 = \init_{\theta}(p_1)$,
and so (again) by the compatibility of $\theta$ with $O$,
we can use the polynomial $\init_{\theta}(g_{i_2}) \in G'$
to involutively reduce $q_1$ to give the polynomial $q_2$.
We now proceed by induction (if $\theta(\LM(p_2)) <
\theta(\LM(p_1))$ then $q_2 = 0$, \ldots), noting that
the process will terminate because
$\init_{\theta}(p_{\alpha} = 0) = 0$.
\end{pf}

\begin{prop}
The set $H$ constructed by the {\bf for}
loop of Algorithm \ref{noncom-walk-inv}
is an Involutive Basis for $J$ with respect
to $I$ and the monomial ordering $O'$.
\end{prop}
\begin{pf}
By Definition \ref{IBNC}, we can show that $H$
is an Involutive Basis for $J$ by showing that any
multiple $upv$ of any polynomial $p \in H$ by any
terms $u$ and $v$ involutively reduces to zero 
using $H$. Assume for a contradiction that such a multiple
does not involutively reduce to zero using $H$, and
instead only involutively reduces to a
polynomial $q \neq 0$.

As all members of $H$ are members of the ideal $J$
(by the way $H$ was constructed as combinations of
elements of $G$), it is clear that $q$ is also
a member of the ideal $J$, as all we have done in
constructing $q$ is to reduce a multiple of a
polynomial in $H$ with respect to $H$. It follows
that the polynomial $\init_{\theta}(q)$ is a
member of the ideal $\init_{\theta}(J)$.

Because $H'$ is an Involutive Basis for the ideal
$\init_{\theta}(J)$ with respect to $I$ and $O'$, there must be a
polynomial $h' \in H'$ such that $h' \mid_I
\init_{\theta}(q)$. Let
$\sum_{i=1}^{j} \ell_ig'_ir_i$ be the logged representation
of $h'$ with respect to $G'$. Then it is clear that
$$\sum_{i=1}^{j} \ell_ig'_ir_i \mid_I \init_{\theta}(q).$$
However $\theta$ is compatible with $O'$ (in particular
the multiplicative variables for $H'$ and $H$ with respect
to $I$ and $O'$ will be identical), so that
$$\sum_{i=1}^{j} \ell_ig_ir_i \mid_I q.$$
It follows that there exists a polynomial $h \in H$
involutively dividing our polynomial $q$,
contradicting our initial assumption.
\end{pf}

\subsection{A Worked Example}

\begin{example} \label{appCch6}
Let $F := \{x^2+y^2+8, \; 2xy+y^2+5\}$ be a basis
generating an ideal $J$ over the polynomial ring
$\mathbb{Q}\langle x, y\rangle$. Consider that
we want to obtain the
DegRevLex Involutive Basis $H := \{2xy-x^2-3, \; y^2+x^2+8, \;
5x^3+6y+35x, \; 5yx^2+3y+10x, \; 2yx-x^2-3\}$
for $J$ from the DegLex Involutive Basis
$G := \{2xy+y^2+5, \; x^2+y^2+8, \;
5y^3-10x+37y, \; 5xy^2+5x-6y, \; 2yx+y^2+5\}$
for $J$ using the Involutive Walk,
where $H$ and $G$ are both Involutive Bases with respect
to the left division $\lhd$. Utilising
Algorithm \ref{noncom-walk-inv} to do this, we initialise
$\theta$ to be the degree function and we proceed as follows.
\begin{itemize}
\item
Construct the set of initials: $$G' :=
\{g'_1, \, g'_2, \, g'_3, \, g'_4, \, g'_5\} =
\{y^2+2xy, \; y^2+x^2, \; 5y^3, \; 5xy^2, \; y^2+2yx\}$$
(these are the terms in $G$ that have maximal degree).
% $\theta$-value).
\item
Compute the Involutive Basis of $G'$ with respect to $\lhd$ and the
DegRevLex monomial ordering. 
Step 1: autoreduce the set $G'$.
\begin{eqnarray*}
g'_1 & = & y^2+2xy \\
& \xymatrix{\ar[r]_{\lhd}_(1){g'_2} &} &
y^2+2xy-(y^2+x^2) \\
& = & 2xy-x^2 =: g'_6; \\
G' & = & (G'\setminus \{g'_1\})\cup \{g'_6\}; \\
% \end{eqnarray*}
% \begin{eqnarray*}
g'_2 & = & y^2+x^2 \\
& \xymatrix{\ar[r]_{\lhd}_(1){g'_5} &} &
y^2+x^2-(y^2+2yx) \\
& = & -2yx+x^2 =: g'_7; \\
G' & = & (G'\setminus \{g'_2\})\cup \{g'_7\}; \\
% \end{eqnarray*}
% \begin{eqnarray*}
g'_3 & = & 5y^3 \\
& \xymatrix{\ar[r]_{\lhd}_(1){g'_5} &} &
5y^3 - 5y(y^2+2yx) \\
& = & -10y^2x \\
& \xymatrix{\ar[r]_{\lhd}_(1){g'_7} &} &
-10y^2x - 5y(-2yx+x^2) \\
& = & -5yx^2 =: g'_8; \\
G' & = & (G'\setminus \{g'_3\})\cup \{g'_8\}; \\
% \end{eqnarray*}
% \begin{eqnarray*}
g'_4 & = & 5xy^2 \\
& \xymatrix{\ar[r]_{\lhd}_(1){g'_5} &} &
5xy^2 -5x(y^2+2yx) \\
& = & -10xyx \\
& \xymatrix{\ar[r]_{\lhd}_(1){g'_7} &} &
-10xyx - 5x(-2yx+x^2) \\
& = & -5x^3 =: g'_9; \\
G' & = & (G'\setminus \{g'_4\})\cup \{g'_9\}; \\
% \end{eqnarray*}
% \begin{eqnarray*}
g'_5 & = & y^2+2yx \\
& \xymatrix{\ar[r]_{\lhd}_(1){g'_7} &} &
y^2+2yx+(-2yx+x^2) \\
& = & y^2+x^2 =: g'_{10}; \\
G' & = & (G'\setminus \{g'_5\})\cup \{g'_{10}\}.
\end{eqnarray*}
Step 2: process the prolongations of the set
$G' = \{g'_6, \, g'_7, \, g'_8, \, g'_9, \, g'_{10}\}$.
Because all ten of these prolongations involutively
reduce to zero using $G'$, we are left with the Involutive
Basis $H' := \{h'_1, \, h'_2, \, h'_3, \, h'_4, \, h'_5\}
= \{2xy-x^2, \; -2yx+x^2, \; -5yx^2, \; -5x^3, \; y^2+x^2\}$ for
$\init_{\theta}(J)$ with respect to $\lhd$ and $O'$.
\item
We must now express the five elements of $H'$ in
terms of members of $G'$.
\begin{eqnarray*}
h'_1 = 2xy-x^2 & = & g'_1-g'_2 \;
\mbox{(from autoreduction)}; \\
h'_2 = -2yx+x^2 & = & g'_2-g'_5; \\
h'_3 = -5yx^2 & = & g'_3 -5yg'_5 - 5y(g'_2-g'_5) \\
& = & -5yg'_2+g'_3; \\
h'_4 = -5x^3 & = & g'_4-5xg'_5-5x(g'_2-g'_5) \\
& = & -5xg'_2+g'_4; \\
h'_5 = y^2+x^2 & = & g'_5+(g'_2-g'_5) \\
& = & g'_2.
\end{eqnarray*}
Lifting to the full polynomials,
we obtain the Involutive Basis
$H := \{h_1, \, h_2, \, h_3, \, h_4, \, h_5\}$
as follows.
\begin{eqnarray*}
h_1 & = & g_1-g_2 \\
    & = & (y^2+2xy+5)-(y^2+x^2+8) \\
    & = & 2xy-x^2-3; \\
h_2 & = & g_2-g_5 \\
    & = & (y^2+x^2+8)-(y^2+2yx+5) \\
    & = & -2yx+x^2+3; \\
h_3 & = & -5yg_2+g_3 \\
    & = & -5y(y^2+x^2+8)+(5y^3+37y-10x) \\
    & = & -5yx^2-3y-10x; \\
h_4 & = & -5xg_2+g_4 \\
    & = & -5x(y^2+x^2+8)+(5xy^2-6y+5x) \\
    & = & -5x^3-6y-35x; \\
h_5 & = & g_2 \\
    & = & y^2+x^2+8.
\end{eqnarray*}
We can now return the output
DegRevLex Involutive Basis
$H = \{2xy-x^2-3, \;
-2yx+x^2+3, \; -5yx^2-3y-10x, \; -5x^3-6y-35x, \; y^2+x^2+8\}$
for $J$ with respect to $\lhd$, a basis that is equivalent to
the DegRevLex Involutive Basis given for $J$
at the beginning of this example.
\end{itemize}
\end{example}

\subsection{Noncommutative Walks Between Any Two Monomial Orderings?}

Thus far, we have only been able to define a noncommutative
walk between two harmonious monomial orderings, where
we recall that the first ordering functions of the 
functional decompositions of the two
monomial orderings must be identical and extendible. For walks
between two arbitrary monomial orderings, the first ordering
functions need not be identical any more, but it is clear
that they must still be
extendible, so that (in an algorithm to perform such a walk)
each basis $G'$ is a Gr\"obner Basis
for each ideal $\init_{\theta}(J)$ (compare with the proofs
of Propositions \ref{inGrob} and \ref{inInv}). This
condition will also apply to any `intermediate' monomial ordering
we will encounter during the walk, but the challenge will be in
how to define these intermediate orderings, so that we
generalise the commutative concept of choosing a weight
vector $\omega_{i+1}$ on the line segment between two weight
vectors $\omega_i$ and $\tau$. %in order for $\omega_{i+1}$

\begin{openquestion} \label{oq4}
Is it possible to perform a noncommutative walk between two
admissible and extendible monomial orderings that are
not harmonious?
\end{openquestion}

\chapter{Conclusions} \label{ChConc}

\section{Current State of Play}

The goal of this thesis was to combine the
theories of noncommutative Gr\"obner Bases and
commutative Involutive Bases to give a theory of
noncommutative Involutive Bases. To
accomplish this, we started by surveying the
background theory in Chapters \ref{ChPr} to
\ref{ChCIB}, focusing our account on the various
algorithms associated with the theory. In particular,
we mentioned several improvements
to the standard algorithms, including how to
compute commutative Involutive Bases by
homogeneous methods, which required the
introduction of a new property (extendibility) of
commutative involutive divisions.

The theory of noncommutative Involutive Bases
was introduced in Chapter \ref{ChNCIB}, where we
described how to perform
noncommutative involutive reduction
(Definition \ref{NCID} and Algorithm \ref{noncom-inv-div});
introduced the notion of a noncommutative
involutive division (Definition \ref{noncom-div-defn});
described what is meant by a noncommutative
Involutive Basis (Definition \ref{IBNC}); and gave an
algorithm to compute noncommutative Involutive
Bases (Algorithm \ref{noncom-inv}). % , an implementation
% of which is given in Appendix \ref{appB}.
Several noncommutative involutive divisions were also
defined, each of which was shown to satisfy certain
properties (such as continuity) allowing the
deductions that all Locally Involutive Bases are
Involutive Bases; and that all Involutive Bases are
Gr\"obner Bases. % It follows that Algorithm \ref{noncom-inv}
% can be used as an alternative algorithm for
% computing noncommutative Gr\"obner Bases.

To finish, we partially generalised the theory of the
Gr\"obner Walk to the noncommutative case
in Chapter \ref{ChWalk},
yielding both Gr\"obner and Involutive Walks between
harmonious noncommutative monomial orderings.

\section{Future Directions}

As well as answering a few questions, the work in this thesis
gives rise to a number of new questions
we would like the answers to. Some of these questions
have already been posed as `Open Questions' in
previous chapters; we summarise below the content
of these questions.
\begin{itemize}
\item
% (Open Question \ref{oq1}) \\[1mm]
Regarding the procedure outlined in
Definition \ref{homprocIC} for computing an Involutive
Basis for a non-homogeneous basis by homogeneous
methods, if the set $G$ returned by the procedure
is not autoreduced, under what circumstances does
autoreducing $G$ result in obtaining a set
that is an Involutive Basis for the ideal generated by
the input basis $F$?
\item
% (Open Question \ref{oq2}) \\[1mm]
Apart from the empty, left and right divisions,
are there any other global noncommutative involutive
divisions of the following types:
\begin{enumerate}[(a)]
\item
strong and continuous;
\item
weak, continuous and Gr\"obner?
\end{enumerate}
\item
% (Open Question \ref{oq3}) \\[1mm]
Are there any conclusive noncommutative involutive divisions
that are also continuous and either strong or Gr\"obner?
\item
% (Open Question \ref{oq4}) \\[1mm]
Is it possible to perform a noncommutative walk between two
admissible and extendible monomial orderings that are
not harmonious?
\end{itemize}

In addition to seeking answers to the above questions,
there are a number of other directions we could take.
One area to explore would be the development of the
algorithms introduced in this thesis. For example, can
the improvements made to the involutive algorithms in
the commutative case, such as the {\it a priori} detection
of prolongations that involutively reduce to zero
(see \cite{Gerdt02}), be applied to the noncommutative
case? Also, can we develop multiple-object versions of
our algorithms, so that (for example) noncommutative Involutive
Bases for path algebras can be computed?

Implementations of any new or improved
algorithms would clearly build
upon the code presented in Appendix \ref{appB}. We
could also expand this code by implementing logged
versions of our algorithms; implementing efficient
methods for performing involutive reduction (as seen for example
in Section \ref{EffRed}); and implementing the algorithms
from Chapter \ref{ChWalk} for performing noncommutative
walks. These improved algorithms and implementations
could then be used (perhaps) to help judge the
relative efficiency and complexity of the involutive
methods versus the Gr\"obner methods.

\subsubsection{Applications}

As every noncommutative Involutive Basis is a 
noncommutative Gr\"obner Basis (at least for the
involutive divisions defined in this thesis),
applications for noncommutative Involutive Bases
will mirror those for noncommutative Gr\"obner
Bases. Some areas in which noncommutative Gr\"obner Bases have
already been used include operator theory; systems engineering
and linear control theory \cite{NCGB02}. Other areas in
noncommutative algebra which could also benefit from the theory
introduced in this thesis include term rewriting; Petri nets;
linear logic; quantum groups and coherence problems.

Further applications may come if we can
extend our algorithms to the multiple-object case.
It would be interesting (for example) to compare a
multiple-object algorithm to a (standard) one-object 
algorithm in cases where an Involutive Basis
for a multiple-object example can be computed
using the one-object algorithm by adding some extra
relations. This would tie in nicely with the
existing comparison between the commutative and
noncommutative versions of the Gr\"obner Basis
algorithm, where it has been noticed that although
commutative examples can be computed using
the noncommutative algorithm, taking this route may
in fact be less efficient than when using the
commutative algorithm to do the same computation.

% New applications in the field of noncommutative algebra.

% APPENDICES
\appendix
%
% Appendix A
% Author: Gareth Evans
% Last Modified: 30th January 2006
%

\chapter{Proof of Propositions 5.5.31 and 5.5.32}
\label{appA}
% \typeout{CHECK NUMBERS FOR PROPOSITIONS HAVE NOT CHANGED}
% \thispagestyle{empty}
% Test Numbers: \ref{alt-cts} and \ref{equiv}.

% \section{Proposition \ref{alt-cts}}
\section{Proposition 5.5.31}

{\bf (Proposition \ref{alt-cts})}
The two-sided left overlap division $\mathcal{W}$ is continuous.

\begin{pf}
Let $w$ be an arbitrary fixed monomial;
let $U$ be any set of monomials; and consider any sequence
$(u_1, \; u_2, \; \hdots, \; u_k)$ of monomials from $U$
($u_i \in U$ for all $1 \leqslant i \leqslant k$),
each of which is a conventional divisor of $w$
(so that $w = \ell_iu_ir_i$ for all $1 \leqslant i \leqslant k$,
where the $\ell_i$ and the $r_i$ are monomials).
For all $1 \leqslant i < k$, suppose that the monomial $u_{i+1}$
satisfies exactly one of the conditions (a) and (b) from
Definition \ref{ncc} (where
multiplicative variables are taken with respect to $\mathcal{W}$
over the set $U$). To show that $\mathcal{W}$ is continuous, we
must show that no two pairs $(\ell_i, r_i)$ and $(\ell_j, r_j)$ are
the same, where $i \neq j$.

Assume to the contrary that there are at least two identical pairs
in the sequence $$((\ell_1, r_1), \, (\ell_2, r_2), \, \hdots,
\, (\ell_k, r_k)),$$ so that we can choose two separate pairs
$(\ell_a, r_a)$ and $(\ell_b, r_b)$ from this sequence such that
$(\ell_a, r_a) = (\ell_b, r_b)$ and all the pairs
$(\ell_c, r_c)$ (for $a \leqslant c < b$) are different. We will now show
that such a sequence $((\ell_a, r_a), \, \hdots, \, (\ell_b, r_b))$
cannot exist.

To begin with, notice that for each monomial $u_{i+1}$ in the
sequence $(u_1, \hdots, u_k)$ of monomials
($1 \leqslant i < k$), if $u_{i+1}$ involutively divides a left
prolongation of the monomial $u_{i}$ (so that
$u_{i+1} \mid_{\mathcal{W}} (\SUFF(\ell_{i}, 1))u_{i}$), 
then $u_{i+1}$ must be
a prefix of this prolongation; if $u_{i+1}$
involutively divides a right prolongation of the monomial
$u_{i}$ (so that $u_{i+1} \mid_{\mathcal{W}} u_{i}(\PRE(r_{i}, 1))$),
then $u_{i+1}$ must be a suffix of this prolongation.
This is because in all other cases,
$u_{i+1}$ is either equal to $u_{i}$, in which case $u_{i+1}$
cannot involutively divide the (left or right) prolongation of
$u_{i}$ trivially; or $u_{i+1}$ is a subword of $u_{i}$, in which case
$u_{i+1}$ cannot involutively divide the (left or right)
prolongation of $u_{i}$ by definition of $\mathcal{W}$.
% We can therefore conclude that $b-a \geqslant 2$.

Following on from the above, we can deduce that $u_b$ is either
a suffix or a prefix of a prolongation of $u_{b-1}$, leaving the
following four cases, where $x^{\ell}_{b-1} = \SUFF(\ell_{b-1}, 1)$ and
$x^r_{b-1} = \PRE(r_{b-1}, 1)$.
\begin{center}
\begin{tabular}{ccc}
Case A ($\deg(u_{b-1}) \geqslant \deg(u_b)$)
&&
Case B ($\deg(u_{b-1}) + 1 = \deg(u_b)$) \\
$\xymatrix @R=0.5pc{
\ar@{-}[r]^*+{x^{\ell}_{b-1}} &
\ar@{<->}[rrrr]^*+{u_{b-1}} &&&& \\
\ar@{<->}[rrr]_*+{u_b} &&&
}$
& \hspace*{5mm} &
$\xymatrix @R=0.5pc{
\ar@{-}[r]^*+{x^{\ell}_{b-1}} & \ar@{<->}[rrr]^*+{u_{b-1}} &&& \\
\ar@{<->}[rrrr]_*+{u_b} &&&&
}$
\\[5mm]
Case C ($\deg(u_{b-1}) \geqslant \deg(u_b)$)
&&
Case D ($\deg(u_{b-1}) + 1 = \deg(u_b)$) \\
$\xymatrix @R=0.5pc{
\ar@{<->}[rrrr]^*+{u_{b-1}} &&&&
\ar@{-}[r]^*+{x^r_{b-1}} & \\
&& \ar@{<->}[rrr]_*+{u_b} &&&
}$
&&
$\xymatrix @R=0.5pc{
\ar@{<->}[rrr]^*+{u_{b-1}} &&& \ar@{-}[r]^*+{x^r_{b-1}} & \\
\ar@{<->}[rrrr]_*+{u_b} &&&&
}$
\end{tabular}
\end{center}
These four cases can all originate from one of the
following two cases (starting with a left prolongation or a
right prolongation), where $x^{\ell}_{a} = \SUFF(\ell_{a}, 1)$ and
$x^r_{a} = \PRE(r_{a}, 1)$.
\begin{center}
\begin{tabular}{ccc}
Case 1 && Case 2 \\
$\xymatrix @R=0.5pc{
\ar@{-}[r]^*+{x^{\ell}_{a}} & \ar@{<->}[rrr]^*+{u_{a}} &&&
}$
& \hspace*{5mm} &
$\xymatrix @R=0.5pc{
\ar@{<->}[rrr]^*+{u_{a}} &&& \ar@{-}[r]^*+{x^{r}_{a}} &
}$
\end{tabular}
\end{center}
So there are eight cases to deal with in total,
namely cases 1-A, 1-B, 1-C, 1-D, 2-A, 2-B, 2-C and 2-D.

We can immediately rule out cases 1-C and 2-A
because we can show that a particular variable is both
multiplicative and nonmultiplicative for monomial $u_a = u_b$
with respect to $U$, a contradiction.
In case 1-C, the variable is $x^{\ell}_a$: it has to be left
nonmultiplicative to provide a left prolongation for $u_a$,
and left multiplicative so that $u_b$ is an involutive divisor
of the right prolongation of $u_{b-1}$; in case 2-A,
the variable is $x^r_a$: it has to be right nonmultiplicative
to provide a right prolongation for $u_a$,
and right multiplicative so that $u_b$ is an involutive divisor
of the left prolongation of $u_{b-1}$.
We illustrate this in the
following diagrams by using a tick to denote a multiplicative
variable and a cross to denote a nonmultiplicative variable.
\begin{center}
\begin{tabular}{ccc}
Case 1-C && Case 2-A \\
$\xymatrix @R=1.5pc{
& \ar@{-}[r]^*+{x^{\ell}_{a}}_*+{\times} & \ar@{<->}[rrr]^*+{u_{a}} &&& \\
&&& \vdots \\
\ar@{<->}[rrrr]^*+{u_{b-1}} &&&& \ar@{-}[r]^*+{x^r_{b-1}} & \\
& \ar@{-}[r]^*+{x^{\ell}_{a}}_*+{\checkmark} 
& \ar@{<->}[rrr]^*+{u_b = u_a} &&&
}$
& \hspace*{5mm} &
$\xymatrix @R=1.5pc{
\ar@{<->}[rrr]^*+{u_{a}} &&& \ar@{-}[r]_*+{\times}^*+{x^r_a} & \\
&&& \vdots \\
\ar@{-}[r]^*+{x^{\ell}_{b-1}} &
\ar@{<->}[rrrr]^*+{u_{b-1}} &&&& \\
\ar@{<->}[rrr]^*+{u_b = u_a} &&& \ar@{-}[r]^*+{x^r_a}_*+{\checkmark} &
}$
\end{tabular}
\end{center}

For all the remaining cases, let us now consider how we may
construct a sequence $((\ell_a, r_a), \, \hdots, \,
(\ell_b, r_b) = (\ell_a, r_a))$. Because we know that
each $u_{c+1}$ is a prefix (or suffix) of a left 
(or right) prolongation of $u_c$ (where $a \leqslant c < b$),
it is clear that at some stage during the sequence,
some $u_{c+1}$ must be a proper suffix (or prefix) of a
prolongation, or else the degrees of the monomials in the
sequence $(u_a, \, \hdots)$ will strictly increase,
meaning that we can never encounter the same $(\ell, r)$
pair twice.
Further, the direction in which prolongations are taken
must change some time during the sequence, or else the
degrees of the monomials in one of the sequences
$(\ell_a, \, \hdots)$ and $(r_a, \, \hdots)$ will strictly
decrease, again meaning that we can never encounter the same
$(\ell, r)$ pair twice.
% the $(\ell, r)$ pairs will all be different.

A change in direction can only occur
if $u_{c+1}$ is equal to a prolongation of $u_c$, as
illustrated below.
\index{left prolongation turn}
\index{right prolongation turn}
\index{prolongation turn!left}
\index{prolongation turn!right}
\begin{center}
\begin{tabular}{ccc}
Left Prolongation Turn && Right Prolongation Turn \\
$\xymatrix @R=2pc{
\ar@{-}[r]^*+{x^{\ell}_{c}}_*+{\times} & \ar@{<->}[rrr]^*+{u_{c}} &&& \\
\ar@{<->}[rrrr]^*+{u_{c+1}} &&&& \ar@{-}[r]^*+{x^{r}_{c+1}}_*+{\times} &
}$
& \hspace*{5mm} &
$\xymatrix @R=2pc{
& \ar@{<->}[rrr]^*+{u_{c}} &&& \ar@{-}[r]_*+{\times}^*+{x^r_c} & \\
\ar@{-}[r]^*+{x^{\ell}_{c+1}}_*+{\times} & \ar@{<->}[rrrr]^*+{u_{c+1}} &&&&
}$
\end{tabular}
\end{center}
However, if no proper prefixes or suffixes are taken
during the sequence, it is clear that making left or
right prolongation turns will not affect the fact that the
degrees of the monomials in the
sequence $(u_a, \, \hdots)$ will strictly increase, once again
meaning that we can never encounter the same $(\ell, r)$
pair twice. It follows that our only course of action is
to make a (left or right) prolongation turn after
a proper prefix or a suffix of a prolongation has been taken. 
We shall call such prolongation turns {\it prefix}
or {\it suffix turns}.
\index{suffix turn}
\index{prefix turn}
\begin{center}
\begin{tabular}{ccc}
Prefix Turn && Suffix Turn \\
$\xymatrix @R=3pc{
& \ar@{-}[r]^*+{x^{\ell}_{c}}_*+{\times} & \ar@{<->}[rrrr]^*+{u_{c}} &&&& \\
\ar@{-}[r]^*+{x^{\ell}_{c+1}}_*+{\times} &
\ar@{<->}[rrr]^*+{u_{c+1}} &&& \ar@{-}[r]^*+{x_{c+2}^r}_*+{\checkmark} & \\
\ar@{<->}[rrrr]^*+{u_{c+2}} &&&& \ar@{-}[r]^*+{x_{c+2}^r}_*+{\times} &
}$
& \hspace*{5mm} &
$\xymatrix @R=3pc{
\ar@{<->}[rrrr]^*+{u_{c}} &&&& \ar@{-}[r]_*+{\times}^*+{x^r_c} & \\
& \ar@{-}[r]^*+{x^{\ell}_{c+2}}_*+{\checkmark} & \ar@{<->}[rrr]^*+{u_{c+1}}
&&& \ar@{-}[r]_*+{\times}^*+{x^r_{c+1}} & \\
& \ar@{-}[r]_*+{\times}^*+{x^{\ell}_{c+2}} & \ar@{<->}[rrrr]^*+{u_{c+2}} &&&&
}$
\end{tabular}
\end{center}

{\bf Claim:} It is impossible to perform a prefix turn when
$\mathcal{W}$ has been used to assign multiplicative variables.

{\bf Proof of Claim:} It is sufficient to
show that $\mathcal{W}$ cannot assign multiplicative
variables to $U$ as follows:
\begin{equation} \label{PreTurn}
x_c^{\ell} \notin \mathcal{M}^L_{\mathcal{W}}(u_c, U); \;
x^r_{c+2} \in \mathcal{M}^R_{\mathcal{W}}(u_{c+1}, U); \;
x^r_{c+2} \notin \mathcal{M}^R_{\mathcal{W}}(u_{c+2}, U).
\end{equation}
Consider how Algorithm \ref{inv-div-alt} can assign the variable
$x^r_{c+2}$ to be right
nonmultiplicative for monomial $u_{c+2}$. As things are set
up in the digram for the prefix turn, the only possibility
is that it is assigned due to the shown overlap between $u_c$ and $u_{c+2}$.
But this assumes that these two monomials actually overlap
(which won't be the case if $\deg(u_{c+1}) = 1$); that
$u_c$ is greater than or equal to $u_{c+2}$ with respect to the
DegRevLex monomial ordering (so any overlap assigns
a nonmultiplicative variable to $u_{c+2}$, not to $u_c$);
and that, by the time we come
to consider the prefix overlap between $u_c$ and
$u_{c+2}$ in Algorithm \ref{inv-div-alt}, the variable $x_c^{\ell}$ must
be left multiplicative for monomial $u_c$. But this final condition
ensures that Algorithm \ref{inv-div-alt} will terminate with
$x_c^{\ell}$ being left multiplicative for $u_c$, contradicting
Equation (\ref{PreTurn}). We therefore conclude that
the variable $x^r_{c+2}$ must be assigned right
nonmultiplicative for monomial $u_{c+2}$ via some other overlap. 

There are three possibilities for this overlap:
(i) there exists a monomial $u_d \in U$ such that
$u_{c+2}$ is a prefix of $u_d$; (ii) there exists a monomial
$u_d \in U$ such that $u_{c+2}$ is a subword of $u_d$; and (iii)
there exists a monomial $u_d \in U$ such that some prefix
of $u_d$ is equal to some suffix of $u_{c+2}$.
\begin{center}
\begin{tabular}{ccc}
Overlap (i) && Overlap (ii) \\
$\xymatrix @R=3pc @C=1.8pc{
\ar@{-}[r]^*+{x^{\ell}_{c+1}}_*+{\times} &
\ar@{<->}[rrr]^*+{u_{c+1}} &&& \ar@{-}[r]^*+{x_{c+2}^r}_*+{\checkmark} & \\
\ar@{<->}[rrrr]^*+{u_{c+2}} &&&& \ar@{-}[r]^*+{x_{c+2}^r}_*+{\times} & \\
\ar@{<->}[rrrrrrr]^*+{u_d} &&&&&&&
}$
& \vspace*{5mm} &
$\xymatrix @R=3pc @C=1.8pc{
& \ar@{-}[r]^*+{x^{\ell}_{c+1}}_*+{\times} &
\ar@{<->}[rrr]^*+{u_{c+1}} &&& \ar@{-}[r]^*+{x_{c+2}^r}_*+{\checkmark} & \\
& \ar@{<->}[rrrr]^*+{u_{c+2}} &&&& \ar@{-}[r]^*+{x_{c+2}^r}_*+{\times} & \\
\ar@{<->}[rrrrrrr]^*+{u_d} &&&&&&&
}$ \\
\multicolumn{3}{c}{Overlap (iii)} \\
\multicolumn{3}{c}{
$\xymatrix @R=3pc @C=1.8pc{
\ar@{-}[r]^*+{x^{\ell}_{c+1}}_*+{\times} &
\ar@{<->}[rrr]^*+{u_{c+1}} &&& \ar@{-}[r]^*+{x_{c+2}^r}_*+{\checkmark} & \\
\ar@{<->}[rrrr]^*+{u_{c+2}} &&&& \ar@{-}[r]^*+{x_{c+2}^r}_*+{\times} & \\
&& \ar@{<->}[rrrrr]^*+{u_d} &&&&&
}$
}
\end{tabular}
\end{center}
In cases (i) and (ii), the overlap shown between $u_{c+1}$ and $u_d$
ensures that Algorithm \ref{inv-div-alt} 
will always assign $x^r_{c+2}$ to be
right nonmultiplicative for monomial $u_{c+1}$,
contradicting Equation (\ref{PreTurn}). This leaves case (iii),
which we break down into two further subcases, dependent upon
whether $u_{c+1}$ is a prefix of $u_d$ or not. If $u_{c+1}$
is a prefix of $u_d$, then
Algorithm \ref{inv-div-alt} will again assign $x^r_{c+2}$ to be
right nonmultiplicative for $u_{c+1}$,
contradicting Equation (\ref{PreTurn}). Otherwise, assuming
that the shown overlap between $u_{c+2}$ and $u_d$
assigns $x^r_{c+2}$ to be right nonmultiplicative for
$u_{c+2}$ (so that the variable immediately to the left
of monomial $u_d$ must be left multiplicative), we must
again come to the conclusion that variable $x^r_{c+2}$ is
right nonmultiplicative for $u_{c+1}$
(due to the overlap between $u_{c+1}$ and $u_d$),
once again contradicting Equation (\ref{PreTurn}).

{\it Technical Point:} It is possible that several left
prolongations may occur between the monomials
$u_{c+1}$ and $u_{c+2}$ shown in the diagram for the
prefix turn, but, as long as no proper prefixes
are taken during this sequence (in which case
we potentially start another prefix turn), we
can apply the same proof as above
(replacing $c+2$ by $c+c'$) to show that we cannot
perform an {\it extended} prefix turn (as shown below)
with respect to $\mathcal{W}$.
\index{prefix turn!extended}
\index{extended prefix turn}
\begin{center}
\begin{tabular}{c}
Extended Prefix Turn\\
$\xymatrix @R=1.1pc{
&&&& & \ar@{-}[r]^*+{x^{\ell}_{c}}_*+{\times} 
& \ar@{<->}[rrrr]^*+{u_{c}} &&&& \\ \\
&&&& \ar@{-}[r]^*+{x^{\ell}_{c+1}}_*+{\times} &
\ar@{<->}[rrr]^*+{u_{c+1}} &&& \ar@{-}[r]^*+{x_{c+c'}^r}_*+{\checkmark} & \\
&&& \ar@{-}[r]_*+{\times} & \ar@{<->}[rrrr] &&&& \\
&& \ar@{-}[r]_*+{\times} & \ar@{<->}[rrrrr] &&&&& \\
&&&& \iddots \\
\ar@{-}[r]_*+{\times} & \ar@{<->}[rrrrrrr] &&&&&&& \\ \\
\ar@{<->}[rrrrrrrr]^*+{u_{c+c'}} &&&&&&&& 
\ar@{-}[r]^*+{x_{c+c'}^r}_*+{\times} &
}$
\end{tabular}
\end{center}
\hfill ${}_{\Box}$

Having ruled out prefix turns, we can now eliminate
cases 1-D, 2-C and 2-D because they require (i) a proper
prefix to be taken during the sequence (allowing
$\deg(r_{b-1}) = \deg(r_b)+1$); and (ii) the final
prolongation to be a right prolongation, ensuring
that a turn has to follow the proper prefix, and so
an (extended) prefix turn is required.

For Cases 1-A and 1-B, we start by taking a left
prolongation, which means that somewhere during
the sequence a proper suffix must be taken.
To do this, it follows that we must change the
direction that prolongations are taken. Knowing
that prefix turns are ruled out, we must therefore
turn by using a left prolongation turn, which
will happen after a finite number $a' \geqslant 1$
of left prolongations.
\begin{center}
$\xymatrix @R=1.1pc{
&&&& & \ar@{-}[r]^*+{x^{\ell}_{a}}_*+{\times} & 
\ar@{<->}[rrr]^*+{u_{a}} &&& \\ \\
&&&& \ar@{-}[r]^*+{x^{\ell}_{a+1}}_*+{\times} &
\ar@{<->}[rrrr]^*+{u_{a+1}} &&&& \\
&&& \ar@{-}[r]_*+{\times} & \ar@{<->}[rrrrr] &&&&& \\
&&&& \iddots \\
\ar@{-}[r]_*+{\times} & \ar@{<->}[rrrrrrrr] &&&&&&&& \\ \\
\ar@{<->}[rrrrrrrrr]^*+{u_{a+a'}} 
&&&&&&&&& \ar@{-}[r]^*+{x_{a+a'}^r}_*+{\times} &
}$
\end{center}
Considering how Algorithm \ref{inv-div-alt} assigns
the variable $x^r_{a+a'}$ to be right nonmultiplicative
for monomial $u_{a+a'}$, there are three possibilities:
(i) there exists a monomial $u_d \in U$ such that
$u_{a+a'}$ is a prefix of $u_d$; (ii) there exists a monomial
$u_d \in U$ such that $u_{a+a'}$ is a subword of $u_d$; and (iii)
there exists a monomial $u_d \in U$ such that some prefix
of $u_d$ is equal to some suffix of $u_{a+a'}$.
In each of these cases, there will be an overlap between
$u_a$ and $u_d$ that will ensure that Algorithm \ref{inv-div-alt}
also assigns the variable $x^r_{a+a'}$ to be right
nonmultiplicative for monomial $u_a$. This rules
out Case 1-A, as variable $x^r_{a+a'}$ must be
right multiplicative for monomial $u_b = u_a$ in order to
perform the final step of Case 1-A.

For Case 1-B, we must now make an (extended) suffix turn as
we need to finish the sequence prolongating to the left.
But, once we have done this, we must subsequently take
a proper prefix in order to ensure that 
$u_{b-1}$ is a suffix of $u_a = u_b$.
Pictorially, here is one way of accomplishing this,
where we note that any number of prolongations may occur between 
any of the shown steps.
\begin{center}
$\xymatrix @R=1.1pc{
\ar@{<->}[rrrrrr]^*+{u_{a+a'}} 
&&&&&& \ar@{-}[r]^*+{x_{a+a'}^r}_*+{\times} & \\
&& \ar@{-}[r]_*+{\checkmark} 
& \ar@{<->}[rrrr] &&&& \ar@{-}[r]_*+{\times} & \\
&& \ar@{-}[r]_*+{\times} & \ar@{<->}[rrrrr] &&&&& \\ \\
&& \ar@{<->}[rrrr]^*+{u_{a+a''}} &&&& \ar@{-}[r]_*+{\checkmark} &
}$
\end{center}
Once we have reached the stage where we are working with
a suffix of $u_a$, we may continue prolongating to the left
until we form the monomial $u_b = u_a$, seemingly providing
a counterexample to the proposition 
(we have managed to construct the same
($\ell, r$) pair twice). However, starting with the monomial
labelled $u_{a+a''}$ in the above diagram, if we follow
the sequence from $u_{a+a''}$ via left prolongations to
$u_b = u_a$, and then continue with the same
sequence as we started off with, we notice that by the
time we encounter the monomial $u_{a+a'}$ again,
an extended prefix turn has been made,
in effect meaning that the first prolongation of
$u_a$ we took right at the start of the sequence
was invalid.
\begin{center}
$\xymatrix @R=1.1pc{
&&&&& &&&& \ar@{-}[r]_*+{\times} &
\ar@{<->}[rr]^*+{u_{a+a''}} &&
\ar@{-}[r]_*+{\checkmark} & \\
&&&& &&&& \ar@{-}[r]_*+{\times} & \ar@{<->}[rrr] &&&& \\
&&&&&&&&& \iddots \\
&&&&& & \ar@{-}[r]^*(0.75)+{x^{\ell}_{b-1}}_*+{\times}
& \ar@{<->}[rrrrr]^*+{u_{b-1}} &&&&& \\ \\
&&&& & \ar@{-}[r]^*+{x^{\ell}_{a}}_*+{\times}
& \ar@{<->}[rrrrrr]^*+{u_b = u_{a}} &&&&&& \\ \\
&&&& \ar@{-}[r]^*+{x^{\ell}_{a+1}}_*+{\times} &
\ar@{<->}[rrrrrrr]^*+{u_{a+1}} &&&&&&& \\
&&& \ar@{-}[r]_*+{\times} & \ar@{<->}[rrrrrrrr] &&&&&&&& \\
&&&&&& \iddots \\
\ar@{-}[r]_*+{\times} & \ar@{<->}[rrrrrrrrrrr] &&&&&&&&&&& \\ \\
\ar@{<->}[rrrrrrrrrrrr]^*+{u_{a+a'}} &&&&&&&&&&&&
\ar@{-}[r]^*+{x_{a+a'}^r}_*+{\times} &
}$
\end{center}

This leaves Case 2-B. Here we start by taking a right
prolongation, meaning that somewhere during
the sequence a proper prefix must be taken.
To do this, it follows that we must change the
direction that prolongations are taken. There are
two ways of doing this: (i) by using
an (extended) suffix turn; (ii) by using a right prolongation turn.

In case (i), after performing the (extended) suffix turn, we need
to take a proper prefix so that the next monomial (say $u_c$) in
the sequence is a suffix of $u_a$; we then continue
by taking left prolongations until we form the
monomial $u_b = u_a$. This provides an
apparent counterexample to the proposition, but as for Case 1-B above,
by taking the right prolongation of $u_a$ the second
time around, we perform an extended prefix turn,
rendering the {\it first} right prolongation of $u_a$ invalid.
\newpage
\begin{center}
Case (i) \\
$\xymatrix @R=1.1pc{
&&& \ar@{-}[r]_*+{\times} &
\ar@{<->}[rr]^*+{u_{c}} &&
\ar@{-}[r]_*+{\checkmark} & \\
&& \ar@{-}[r]_*+{\times} & \ar@{<->}[rrr] &&&& \\
&&& \iddots \\
\ar@{-}[r]^*(0.75)+{x^{\ell}_{b-1}}_*+{\times}
& \ar@{<->}[rrrrr]^*+{u_{b-1}} &&&&& \\ \\
\ar@{<->}[rrrrrr]^*+{u_b = u_a}
&&&&&& \ar@{-}[r]^*+{x^r_a}_*+{\times} & \\ \\
}$
\end{center}
In case (ii), after we make a right prolongation
turn (which may itself occur
after a finite number of right prolongations),
we may now take the required proper prefix. But as
we are then required to take a proper suffix (in order to
ensure that we finish the sequence taking a left
prolongation), we need to make a turn. But as
this would entail making an (extended) prefix turn, we
conclude that case (ii) is also invalid.
\begin{center}
An Example of Case (ii) \\
$\xymatrix @R=1.1pc{
&& \ar@{<->}[rrrr]^*+{u_a} &&&& \ar@{-}[r]^*+{x^r_a}_*+{\times} & \\ \\
& \ar@{-}[r]^*+{x^{\ell}_{a+1}}_*+{\times} &
\ar@{<->}[rrrrr]^*+{u_{a+1}} &&&&& \\ \\
\ar@{-}[r]^*+{x^{\ell}_{a+2}}_*+{\times} &
\ar@{<->}[rrr]^*+{u_{a+2}} &&& \ar@{-}[r]_*+{\checkmark} & \\ \\
\ar@{<->}[rrrr]^*+{u_{a+3}} &&&& \ar@{-}[r]_*+{\times} &
}$
\end{center}
As we have now accounted for all eight possible
sequences, we can conclude that $\mathcal{W}$ is continuous.
\end{pf}

\newpage
% \section{Proposition \ref{equiv}}
\section{Proposition 5.5.32}

{\bf (Proposition \ref{equiv})}
The two-sided left overlap division $\mathcal{W}$
is a Gr\"obner involutive division.

\begin{pf}
We are required to show that if Algorithm
\ref{noncom-inv} terminates with $\mathcal{W}$ and some
arbitrary admissible monomial ordering $O$ as input,
then the Locally Involutive Basis $G$ it returns is a
noncommutative Gr\"obner Basis. By Definition \ref{grob-defn-noncom},
we can do this by showing that all S-polynomials involving
elements of $G$ conventionally reduce to zero using $G$.

Assume that $G = \{g_1, \hdots, g_p\}$ is sorted (by lead
monomial) with respect to
the DegRevLex monomial ordering (greatest first),
and let $U = \{u_1, \hdots, u_p\}
:= \{\LM(g_1), \hdots, \LM(g_p)\}$ be the set of leading monomials. 
Let $T$ be the
table obtained by applying Algorithm \ref{inv-div-alt} to $U$.
Because $G$ is a Locally Involutive Basis, every zero
entry $T(u_i, x_j^{\Gamma})$ ($\Gamma \in \{L, R\}$) in the table
corresponds to a prolongation $g_ix_j$ or $x_jg_i$ that involutively
reduces to zero.

Let $S$ be the set of S-polynomials involving elements of $G$, where
the $t$-th entry of $S$ ($1 \leqslant t \leqslant |S|$)
is the S-polynomial
$$s_t = c_{t}\ell_{t}g_ir_{t} - c'_{t}\ell'_{t}g_jr'_{t},$$
with $\ell_{t}u_ir_{t} = \ell'_{t}u_jr'_{t}$ being the
overlap word of the S-polynomial.
We will prove that every S-polynomial in $S$
conventionally reduces to zero using $G$.

Recall (from Definition \ref{ov-def})
that each S-polynomial in $S$ corresponds to a particular
type of overlap --- `prefix', `subword' or `suffix'. For the
purposes of this proof, let us now split the subword
overlaps into three further types  --- `left', `middle' and `right',
corresponding to the cases where a monomial $m_2$ is
a prefix, proper subword and suffix of a monomial $m_1$.
\begin{center}
\begin{tabular}{ccc}
Left & Middle & Right \\[1mm]
$\xymatrix @R=0.5pc{
\ar@{<->}[rrr]^*+{m_1} &&& \\
\ar@{<->}[rr]_*+{m_2} &&
}$
&
$\xymatrix @R=0.5pc{
\ar@{<->}[rrr]^*+{m_1} &&& \\
& \ar@{<->}[r]_*+{m_2} &
}$
&
$\xymatrix @R=0.5pc{
\ar@{<->}[rrr]^*+{m_1} &&& \\
& \ar@{<->}[rr]_*+{m_2} &&
}$
\end{tabular}
\end{center}
This classification provides
us with five cases to deal with in total,
which we shall process in the following order:
right, middle, left, prefix, suffix.

{\bf (1)} Consider an arbitrary entry $s_{t} \in S$
($1 \leqslant t \leqslant |S|$)
corresponding to a right overlap where the monomial
$u_j$ is a suffix of the monomial $u_i$. This means that
$s_{t} = c_tg_i - c'_t\ell'_tg_j$ for some $g_i, g_j \in G$,
with overlap word $u_i = \ell'_tu_j$. 
Let $u_i = x_{i_1}\hdots x_{i_{\alpha}}$;
let $u_j = x_{j_1}\hdots x_{j_{\beta}}$; and let $D = \alpha - \beta$.
$$\xymatrix @R=0.5pc {
u_i = & \ar@{-}[r]_{x_{i_1}} & \ar@{-}[r]_{x_{i_2}} & \ar@{--}[r]
& \ar@{-}[r]_{x_{i_D}} & \ar@{-}[r]_{x_{i_{D+1}}}
& \ar@{-}[r]_{x_{i_{D+2}}} & \ar@{--}[r]
& \ar@{-}[r]_{x_{i_{\alpha - 1}}} & \ar@{-}[r]_{x_{i_{\alpha}}} & \\
u_j = & &&&& \ar@{-}[r]_{x_{j_1}} & \ar@{-}[r]_{x_{j_2}} & \ar@{--}[r]
& \ar@{-}[r]_{x_{j_{\beta - 1}}} & \ar@{-}[r]_{x_{j_{\beta}}} &
}$$
Because $u_j$ is a suffix of $u_i$, 
it follows that $T(u_j, x^L_{i_{D}}) = 0$.
This gives rise to the prolongation
$x_{i_{D}}g_j$ of $g_j$. But we know that all prolongations
involutively reduce to zero ($G$ is a Locally Involutive Basis), 
so Algorithm \ref{noncom-inv-div} must find a monomial
$u_{k} = x_{k_1}\hdots x_{k_{\gamma}} \in U$ such that
$u_{k}$ involutively divides $x_{i_{D}}u_j$.
Assuming that $x_{k_{\gamma}} = x_{i_{\kappa}}$,
we can deduce that
any candidate for $u_k$ must be a suffix of $x_{i_{D}}u_j$
(otherwise $T(u_k, x^R_{i_{\kappa + 1}}) = 0$ because of the
overlap between $u_i$ and $u_k$).
But if $u_k$ is a suffix of $x_{i_{D}}u_j$, then we must have
$u_k = x_{i_{D}}u_j$ (otherwise $T(u_k, x^L_{i_{\alpha - \gamma}}) = 0$
again because of the overlap between $u_i$ and $u_k$).
We have therefore shown that there exists a monomial
$u_k = x_{k_1}\hdots x_{k_{\gamma}} \in U$ such that
$u_k$ is a suffix of $u_i$ and $\gamma = \beta + 1$.
$$\xymatrix @R=0.5pc {
u_i = & \ar@{-}[r]_{x_{i_1}} & \ar@{-}[r]_{x_{i_2}} & \ar@{--}[r]
& \ar@{-}[r]_{x_{i_D}} & \ar@{-}[r]_{x_{i_{D+1}}}
& \ar@{-}[r]_{x_{i_{D+2}}} & \ar@{--}[r]
& \ar@{-}[r]_{x_{i_{\alpha - 1}}} & \ar@{-}[r]_{x_{i_{\alpha}}} & \\
u_j = & &&&& \ar@{-}[r]_{x_{j_1}} & \ar@{-}[r]_{x_{j_2}} & \ar@{--}[r]
& \ar@{-}[r]_{x_{j_{\beta - 1}}} & \ar@{-}[r]_{x_{j_{\beta}}} & \\
u_k = & &&& \ar@{-}[r]_{x_{k_1}} & \ar@{-}[r]_{x_{k_2}}
& \ar@{-}[r]_{x_{k_3}} & \ar@{--}[r]
& \ar@{-}[r]_{x_{k_{\gamma - 1}}} & \ar@{-}[r]_{x_{k_{\gamma}}} &
}$$
In the case $D = 1$, it is clear that $u_k = u_i$, and so the first step
in the involutive reduction of the prolongation $x_{i_1}g_j$ of $g_j$
is to take away the multiple $(\frac{c_t}{c'_t})g_i$ of $g_i$
from $x_{i_1}g_j$ to leave the polynomial $x_{i_1}g_j
- (\frac{c_t}{c'_t})g_i = -(\frac{1}{c'_t})s_t$.
But as we know that all prolongations involutively reduce
to zero, we can conclude that the S-polynomial
$s_t$ conventionally reduces to zero.

For the case $D > 1$,
we can use the monomial $u_k$ together with Buchberger's Second Criterion
to simplify our goal of showing that the S-polynomial $s_t$ reduces
to zero. Notice that the monomial $u_k$
is a subword of the overlap word $u_i$ associated to $s_t$,
and so in order to show that $s_t$ reduces to zero,
all we have to do is to show that the two S-polynomials
$$s_u = c_ug_i - c'_{u}(x_{i_1}x_{i_2}\hdots x_{i_{D-1}})g_k$$
and
$$s_v = c_vg_k - c'_{v}x_{i_D}g_j$$
reduce to zero ($1 \leqslant u,v \leqslant |S|$).
But $s_v$ is an S-polynomial corresponding
to a right overlap of type $D = 1$ (because $\gamma-\beta = 1$),
and so $s_v$ reduces to zero. It remains
to show that the S-polynomial $s_u$
reduces to zero. But we can do this by using exactly the same argument
as above --- we can show that there exists a monomial
$u_{\pi} = x_{{\pi}_1}\hdots x_{{\pi}_{\delta}} \in U$
such that $u_{\pi}$ is a suffix of $u_i$ and $\delta = \gamma + 1$,
and we can deduce that the S-polynomial $s_u$ reduces to zero
(and hence $s_t$ reduces to $0$)
if the S-polynomial $$s_w = c_wg_i -
c'_{w}(x_{i_1}x_{i_2}\hdots x_{i_{D-2}})g_{\pi}$$ reduces to zero
($1 \leqslant w \leqslant |S|$).
By induction, there is a sequence
$\{u_{q_D}, u_{q_{D-1}}, \hdots, u_{q_2}\}$ of monomials
increasing uniformly in degree, so that
$s_t$ reduces to zero if the S-polynomial
$$s_{\eta} = c_{\eta}g_i - c'_{\eta}x_{i_1}g_{q_2}$$
reduces to zero ($1 \leqslant \eta \leqslant |S|$).
$$\xymatrix @R=0.5pc {
u_i = & \ar@{-}[r]_{x_{i_1}} & \ar@{-}[r]_{x_{i_2}} & \ar@{--}[r]
& \ar@{-}[r]_{x_{i_{D-1}}} & \ar@{-}[r]_{x_{i_D}}
& \ar@{-}[r]_{x_{i_{D+1}}} & \ar@{-}[r]_{x_{i_{D+2}}} & \ar@{--}[r]
& \ar@{-}[r]_{x_{i_{\alpha - 1}}} & \ar@{-}[r]_{x_{i_{\alpha}}} & \\
u_j = & &&&&& \ar@{-}[r]_{x_{j_1}} & \ar@{-}[r]_{x_{j_2}} & \ar@{--}[r]
& \ar@{-}[r]_{x_{j_{\beta - 1}}} & \ar@{-}[r]_{x_{j_{\beta}}} & \\
u_{q_{D}} = u_k = & &&&& \ar@{-}[r] & \ar@{-}[r] & \ar@{-}[r] & \ar@{--}[r]
& \ar@{-}[r] & \ar@{-}[r] & \\
u_{q_{D-1}} = u_{\pi} = & &&& \ar@{-}[r] 
& \ar@{-}[r] & \ar@{-}[r] & \ar@{-}[r]
& \ar@{--}[r] & \ar@{-}[r] & \ar@{-}[r] & \\
\vdots &&&&&&& \iddots \\
u_{q_{2}} = & & \ar@{-}[r] & \ar@{--}[r] & \ar@{-}[r] & \ar@{-}[r]
& \ar@{-}[r] & \ar@{-}[r] & \ar@{--}[r] & \ar@{-}[r] & \ar@{-}[r] & \\
}$$

But $s_{\eta}$ is always an S-polynomial corresponding
to a right overlap of type $D = 1$, and so $s_{\eta}$
reduces to zero --- meaning we can conclude that $s_t$
reduces to zero as well.

% {\bf (2) -- (5)} We refer to Appendix \ref{appA}.
% In summary, for the other types of overlap,
% we use the same iterative idea as above to show that
% a particular S-polynomial in a {\it chain} of
% S-polynomials reduces to zero. For a middle
% overlap, we will either encounter a right overlap
% or an overlap of `type $D = 1$' in our chain, cases
% that we have dealt with above. For the
% remaining types of overlap, we add the possibility of
% encountering a middle overlap in the chain.
  
{\bf (2)} Consider an arbitrary entry $s_{t} \in S$
($1 \leqslant t \leqslant |S|$)
corresponding to a middle overlap where the monomial
$u_j$ is a proper subword of the monomial $u_i$. This means that
$s_{t} = c_tg_i - c'_t\ell'_tg_jr'_t$ for some $g_i, g_j \in G$,
with overlap word $u_i = \ell'_tu_jr'_t$. 
Let $u_i = x_{i_1}\hdots x_{i_{\alpha}}$;
let $u_j = x_{j_1}\hdots x_{j_{\beta}}$; and choose $D$ such that
$x_{i_D} = x_{j_{\beta}}$.
$$\xymatrix @R=0.5pc @C=2.4pc {
u_i = & \ar@{-}[r]_{x_{i_1}} & \ar@{--}[r]
& \ar@{-}[r]_{x_{i_{D-\beta}}} & \ar@{-}[r]_{x_{i_{D-\beta+1}}}
& \ar@{-}[r]_{x_{i_{D-\beta+2}}} & \ar@{--}[r]
& \ar@{-}[r]_{x_{i_{D-1}}} 
& \ar@{-}[r]_{x_{i_{D}}} & \ar@{-}[r]_{x_{i_{D+1}}} &
\ar@{--}[r] & \ar@{-}[r]_{x_{i_{\alpha}}} & \\
u_j = & &&& \ar@{-}[r]_{x_{j_1}} & \ar@{-}[r]_{x_{j_2}} & \ar@{--}[r]
& \ar@{-}[r]_{x_{j_{\beta - 1}}} & \ar@{-}[r]_{x_{j_{\beta}}} &
}$$
Because $u_j$ is a proper subword of $u_i$, it follows that
$T(u_j, x^R_{i_{D+1}}) = 0$.
This gives rise to the prolongation
$g_jx_{i_{D+1}}$ of $g_j$. But we know that all prolongations
involutively reduce to zero, so there must exist a monomial
$u_{k} = x_{k_1}\hdots x_{k_{\gamma}} \in U$ such that
$u_{k}$ involutively divides $u_jx_{i_{D+1}}$.
Assuming that $x_{k_{\gamma}} = x_{i_{\kappa}}$,
any candidate for $u_k$ must be a suffix of $u_jx_{i_{D+1}}$
(otherwise $T(u_k, x^R_{i_{\kappa + 1}}) = 0$ because of the
overlap between $u_i$ and $u_k$). Unlike part (1) however,
we cannot determine the degree of $u_k$ (so that
$1 \leqslant \gamma \leqslant \beta+1$); we shall illustrate
this in the following diagram by using a squiggly line to
indicate that the monomial $u_k$ can begin anywhere 
(or nowhere if $u_k = x_{i_{D+1}}$) on the squiggly line.
$$\xymatrix @R=0.5pc @C=2.4pc {
u_i = & \ar@{-}[r]_{x_{i_1}} & \ar@{--}[r]
& \ar@{-}[r]_{x_{i_{D-\beta}}} & \ar@{-}[r]_{x_{i_{D-\beta+1}}}
& \ar@{-}[r]_{x_{i_{D-\beta+2}}} & \ar@{--}[r]
& \ar@{-}[r]_{x_{i_{D-1}}} 
& \ar@{-}[r]_{x_{i_{D}}} & \ar@{-}[r]_{x_{i_{D+1}}} &
\ar@{--}[r] & \ar@{-}[r]_{x_{i_{\alpha}}} & \\
u_j = & &&& \ar@{-}[r]_{x_{j_1}} & \ar@{-}[r]_{x_{j_2}} & \ar@{--}[r]
& \ar@{-}[r]_{x_{j_{\beta - 1}}} & \ar@{-}[r]_{x_{j_{\beta}}} & \\
u_k = & &&& \ar@{~}[rrrrr] &&&&& \ar@{-}[r]_{x_{k_{\gamma}}} &
}$$
We can now use the monomial $u_k$ 
together with Buchberger's Second Criterion
to simplify our goal of showing that the S-polynomial $s_t$ reduces
to zero. Notice that the monomial $u_k$
is a subword of the overlap word $u_i$ associated to $s_t$,
and so in order to show that $s_t$ reduces to zero,
all we have to do is to show that the two S-polynomials
$$s_u = c_ug_i - c'_{u}(x_{i_1}x_{i_2}\hdots
x_{i_{D+1-\gamma}})g_k(x_{i_{D+2}}\hdots x_{i_{\alpha}})$$
and\footnote{Technical point: if $\gamma \neq \beta+1$,
the S-polynomial $s_v$ could in fact appear as
$s_v = c_{v}g_jx_{i_{D+1}} - c'_v(x_{j_1}\hdots x_{i_{D+1-\gamma}})g_k$
and not as
$s_v = c_v(x_{j_1}\hdots x_{i_{D+1-\gamma}})g_k - c'_{v}g_jx_{i_{D+1}}$;
for simplicity we will treat both cases the same in the
proof as all that changes is the notation and the signs.}
$$s_v = c_v(x_{j_1}\hdots x_{i_{D+1-\gamma}})g_k - c'_{v}g_jx_{i_{D+1}}$$
reduce to zero ($1 \leqslant u,v \leqslant |S|$).

For the S-polynomial $s_v$, there are two cases to consider:
$\gamma = 1$, and $\gamma > 1$. In the former case,
because (as placed in $u_i$) the monomials
$u_j$ and $u_k$ do not overlap, we can use Buchberger's
First Criterion to say that the `S-polynomial' $s_v$ reduces to
zero (for further explanation, see the paragraph at the beginning of
Section \ref{BCNC}). In the latter case, note that $u_k$ is
the only involutive divisor of the prolongation
$u_jx_{i_{D+1}}$, as the existence of any suffix of $u_jx_{i_{D+1}}$
of higher degree than $u_k$ in $U$ will contradict the fact
that $u_k$ is an involutive divisor of $u_jx_{i_{D+1}}$; and
the existence of $u_k$ in $U$ ensures that any suffix of
$u_jx_{i_{D+1}}$ that exists in $U$ with a lower degree than
$u_k$ will not be an involutive divisor of $u_jx_{i_{D+1}}$. This means that
the first step of the involutive reduction of $g_jx_{i_{D+1}}$
is to take away the multiple
$(\frac{c_v}{c'_v})(x_{j_1}\hdots x_{i_{D+1-\gamma}})g_k$ of $g_k$
from $g_jx_{i_{D+1}}$ to leave the polynomial $g_jx_{i_{D+1}}
- (\frac{c_v}{c'_v})(x_{j_1}\hdots x_{i_{D+1-\gamma}})g_k
= -(\frac{1}{c'_v})s_v$.
But as we know that all prolongations involutively reduce
to zero, we can conclude that the S-polynomial
$s_v$ conventionally reduces to zero.

For the S-polynomial $s_u$, we note that if $D = \alpha-1$, then
$s_u$ corresponds to a right overlap, and so we know from part (1) that
$s_u$ conventionally reduces to zero. Otherwise, we proceed by
induction on the S-polynomial $s_u$ to produce a sequence
$\{u_{q_{D+1}}, u_{q_{D+2}}, \hdots, u_{q_{\alpha}}\}$ of monomials, so that
$s_u$ (and hence $s_t$) reduces to zero if the S-polynomial
$$s_{\eta} = c_{\eta}g_i - c'_{\eta}(x_{i_1}\hdots 
x_{i_{\alpha-\mu}})g_{q_{\alpha}}$$
reduces to zero ($1 \leqslant \eta \leqslant |S|$),
where $\mu = \deg(u_{q_{\alpha}})$.
$$\xymatrix @R=0.5pc @C=2.0pc {
u_i = & \ar@{-}[r]_{x_{i_1}} & \ar@{--}[r]
& \ar@{-}[r]_{x_{i_{D-\beta}}} & \ar@{-}[r]_{x_{i_{D-\beta+1}}}
& \ar@{--}[r] & \ar@{-}[r]_{x_{i_{D}}} & \ar@{-}[r]_{x_{i_{D+1}}} &
\ar@{-}[r]_{x_{i_{D+2}}} & \ar@{--}[r] & \ar@{-}[r]_{x_{i_{\alpha-1}}} &
\ar@{-}[r]_{x_{i_{\alpha}}} & \\
u_j = & &&& \ar@{-}[r]_{x_{j_1}} & \ar@{--}[r]
& \ar@{-}[r]_{x_{j_{\beta}}} & \\
u_{q_{D+1}} = u_k = & &&& \ar@{~}[rrr] &&& \ar@{-}[r]_{x_{k_{\gamma}}} & \\
u_{q_{D+2}} = & &&& \ar@{~}[rrrr] &&&& \ar@{-}[r] & \\
& &&&&&&& \ddots \\
u_{q_{\alpha}} = & &&& \ar@{~}[rrrrrrr] &&&&&&& \ar@{-}[r] &
}$$
But $s_{\eta}$ always corresponds to a right overlap, 
and so $s_{\eta}$ reduces to zero --- meaning we can conclude that $s_t$
reduces to zero as well.

{\bf (3)} Consider an arbitrary entry $s_{t} \in S$
($1 \leqslant t \leqslant |S|$)
corresponding to a left overlap where the monomial
$u_j$ is a prefix of the monomial $u_i$. This means that
$s_{t} = c_tg_i - c'_tg_jr'_t$ for some $g_i, g_j \in G$,
with overlap word $u_i = u_jr'_t$. Let $u_i = x_{i_1}\hdots x_{i_{\alpha}}$
and let $u_j = x_{j_1}\hdots x_{j_{\beta}}$.
$$\xymatrix @R=0.5pc {
u_i = & \ar@{-}[r]_{x_{i_1}} & \ar@{-}[r]_{x_{i_2}} & \ar@{--}[r] &
\ar@{-}[r]_{x_{i_{\beta-1}}} & \ar@{-}[r]_{x_{i_{\beta}}} &
\ar@{-}[r]_{x_{i_{\beta+1}}} & \ar@{--}[r] &
\ar@{-}[r]_{x_{i_{\alpha-1}}} & \ar@{-}[r]_{x_{i_{\alpha}}} & \\
u_j = & \ar@{-}[r]_{x_{j_1}} & \ar@{-}[r]_{x_{j_2}} & \ar@{--}[r] &
\ar@{-}[r]_{x_{j_{\beta-1}}} & \ar@{-}[r]_{x_{j_{\beta}}} &
}$$
Because $u_j$ is a prefix of $u_i$, it follows that
$T(u_j, x^R_{i_{\beta+1}}) = 0$.
This gives rise to the prolongation
$g_jx_{i_{\beta+1}}$ of $g_j$. But we know that all prolongations
involutively reduce to zero, so there must exist a monomial
$u_{k} = x_{k_1}\hdots x_{k_{\gamma}} \in U$ such that
$u_{k}$ involutively divides $u_jx_{i_{\beta+1}}$.
Assuming that $x_{k_{\gamma}} = x_{i_{\kappa}}$,
any candidate for $u_k$ must be a suffix of $u_jx_{i_{\beta+1}}$
(otherwise $T(u_k, x^R_{i_{\kappa + 1}}) = 0$ because of the
overlap between $u_i$ and $u_k$).
$$\xymatrix @R=0.5pc {
u_i = & \ar@{-}[r]_{x_{i_1}} & \ar@{-}[r]_{x_{i_2}} & \ar@{--}[r] &
\ar@{-}[r]_{x_{i_{\beta-1}}} & \ar@{-}[r]_{x_{i_{\beta}}} &
\ar@{-}[r]_{x_{i_{\beta+1}}} & \ar@{--}[r] &
\ar@{-}[r]_{x_{i_{\alpha-1}}} & \ar@{-}[r]_{x_{i_{\alpha}}} & \\
u_j = & \ar@{-}[r]_{x_{j_1}} & \ar@{-}[r]_{x_{j_2}} & \ar@{--}[r] &
\ar@{-}[r]_{x_{j_{\beta-1}}} & \ar@{-}[r]_{x_{j_{\beta}}} & \\
u_k = & \ar@{~}[rrrrr] &&&&& \ar@{-}[r]_{x_{k_{\gamma}}} &
}$$
If $\alpha = \gamma$,
then it is clear that $u_k = u_i$,
and so the first step in the involutive reduction of the
prolongation $g_jx_{i_{\alpha}}$
is to take away the multiple
$(\frac{c_t}{c'_t})g_i$ of $g_i$
from $g_jx_{i_{\alpha}}$ to leave the polynomial $g_jx_{i_{\alpha}}
- (\frac{c_t}{c'_t})g_i = -(\frac{1}{c'_t})s_t$.
But as we know that all prolongations involutively reduce
to zero, we can conclude that the S-polynomial
$s_t$ conventionally reduces to zero.

Otherwise, if $\alpha > \gamma$, we can now use the monomial $u_k$
together with Buchberger's Second Criterion
to simplify our goal of showing that the S-polynomial $s_t$ reduces
to zero. Notice that the monomial $u_k$
is a subword of the overlap word $u_i$ associated to $s_t$,
and so in order to show that $s_t$ reduces to zero,
all we have to do is to show that the two S-polynomials
$$s_u = c_ug_i - c'_{u}(x_{i_1}\hdots
x_{i_{\beta+1-\gamma}})g_k(x_{i_{\beta+2}}\hdots x_{i_{\alpha}})$$
and
$$s_v = c_v(x_{i_1}\hdots x_{i_{\beta+1-\gamma}})g_k 
- c'_{v}g_jx_{i_{\beta+1}}$$
reduce to zero ($1 \leqslant u,v \leqslant |S|$).

The S-polynomial $s_v$ reduces to zero by comparison with
part (2). For the S-polynomial $s_u$,
first note that if $\alpha = \beta+1$, then
$s_u$ corresponds to a right overlap, and so we know from part (1) that
$s_u$ conventionally reduces to zero. Otherwise, if $\gamma \neq
\beta+1$, then $s_u$ corresponds to a middle overlap, and so
we know from part (2) that $s_u$ conventionally reduces to zero.
This leaves the case where $s_u$ corresponds to another left
overlap, in which case we proceed by induction on $s_u$, eventually
coming across either a middle overlap or a right overlap because
we move one letter at a time to the right after each
inductive step.
$$\xymatrix @R=0.5pc {
u_i = & \ar@{-}[r]_{x_{i_1}} & \ar@{-}[r]_{x_{i_2}} & \ar@{--}[r] &
\ar@{-}[r]_{x_{i_{\beta-1}}} & \ar@{-}[r]_{x_{i_{\beta}}} &
\ar@{-}[r]_{x_{i_{\beta+1}}} & \ar@{-}[r]_{x_{i_{\beta+2}}} & \ar@{--}[r] &
\ar@{-}[r]_{x_{i_{\alpha-1}}} & \ar@{-}[r]_{x_{i_{\alpha}}} & \\
u_j = & \ar@{-}[r]_{x_{j_1}} & \ar@{-}[r]_{x_{j_2}} & \ar@{--}[r] &
\ar@{-}[r]_{x_{j_{\beta-1}}} & \ar@{-}[r]_{x_{j_{\beta}}} & \\
u_k = & \ar@{~}[rrrrr] &&&&& \ar@{-}[r]_{x_{k_{\gamma}}} & \\
      & \ar@{~}[rrrrrr] &&&&&& \ar@{-}[r] & \\
      &&&&&& \ddots \\
      & \ar@{~}[rrrrrrrrr] &&&&&&&&& \ar@{-}[r] & \\
}$$

{\bf (4 and 5)} In Definition \ref{ov-def}, we defined a prefix
overlap to be an overlap where, given two monomials $m_1$ and
$m_2$ such that $\deg(m_1) \geqslant \deg(m_2)$, a prefix of
$m_1$ is equal to a suffix of $m_2$; suffix overlaps were defined
similarly. If we drop the condition on the degrees
of the monomials, it is clear that every suffix
overlap can be treated as a prefix overlap (by swapping the
roles of $m_1$ and $m_2$); this allows us to deal with the
case of a prefix overlap only.

Consider an arbitrary entry $s_{t} \in S$
($1 \leqslant t \leqslant |S|$)
corresponding to a prefix overlap where a prefix of the monomial
$u_i$ is equal to a suffix of the monomial $u_j$. This means that
$s_{t} = c_t\ell_tg_i - c'_tg_jr'_t$ for some $g_i, g_j \in G$,
with overlap word $\ell_tu_i = u_jr'_t$.
Let $u_i = x_{i_1}\hdots x_{i_{\alpha}}$;
let $u_j = x_{j_1}\hdots x_{j_{\beta}}$; and choose $D$ such that
$x_{i_D} = x_{j_{\beta}}$.
$$\xymatrix @R=0.5pc @C=2.4pc{
u_i = &&&& \ar@{-}[r]_{x_{i_1}} & \ar@{--}[r]
& \ar@{-}[r]_{x_{i_{D}}}
& \ar@{-}[r]_{x_{i_{D+1}}} & \ar@{--}[r]
& \ar@{-}[r]_{x_{i_{\alpha - 1}}} & \ar@{-}[r]_{x_{i_{\alpha}}} & \\
u_j = & \ar@{-}[r]_{x_{j_1}} & \ar@{--}[r] & \ar@{-}[r]_{x_{j_{\beta-D}}}
& \ar@{-}[r]_{x_{j_{\beta-D+1}}} & \ar@{--}[r]
& \ar@{-}[r]_{x_{j_{\beta}}} &
}$$
By definition of $\mathcal{W}$, at least one of
$T(u_i, x^L_{j_{\beta-D}})$ and $T(u_j, x^R_{i_{D+1}})$ is equal to zero.
\begin{itemize}
\item Case $T(u_j, x^R_{i_{D+1}}) = 0$. \\[2mm]
Because we know that the prolongation
$g_jx_{i_{D+1}}$ involutively reduces to zero, there must exist a monomial
$u_{k} = x_{k_1}\hdots x_{k_{\gamma}} \in U$ such that
$u_{k}$ involutively divides $u_jx_{i_{D+1}}$.
This $u_k$ must be a suffix of $u_jx_{i_{D+1}}$
(otherwise, assuming that $x_{k_{\gamma}} = x_{j_{\kappa}}$,
we have $T(u_k, x^R_{i_{D+1}}) = 0$ if $\gamma = \beta$
(because of the overlap between $u_i$ and $u_k$);
$T(u_k, x^L_{j_{\beta-\gamma}}) = 0$ if
$\gamma < \beta$ and $\kappa = \beta$ (because of the overlap between
$u_j$ and $u_k$); and $T(u_k, x^R_{j_{\kappa + 1}}) = 0$
if $\gamma < \beta$ and $\kappa < \beta$ 
(again because of the overlap between $u_j$ and $u_k$)).
$$\xymatrix @R=0.5pc @C=2.4pc{
u_i = &&&& \ar@{-}[r]_{x_{i_1}} & \ar@{--}[r]
& \ar@{-}[r]_{x_{i_{D}}}
& \ar@{-}[r]_{x_{i_{D+1}}} & \ar@{--}[r]
& \ar@{-}[r]_{x_{i_{\alpha - 1}}} & \ar@{-}[r]_{x_{i_{\alpha}}} & \\
u_j = & \ar@{-}[r]_{x_{j_1}} & \ar@{--}[r] & \ar@{-}[r]_{x_{j_{\beta-D}}}
& \ar@{-}[r]_{x_{j_{\beta-D+1}}} & \ar@{--}[r]
& \ar@{-}[r]_{x_{j_{\beta}}} & \\
u_k = & \ar@{~}[rrrrrr] &&&&&& \ar@{-}[r]_{x_{k_{\gamma}}} &
}$$
Let us now
use the monomial $u_k$ together with Buchberger's Second Criterion
to simplify our goal of showing that the S-polynomial $s_t$ reduces
to zero. Because $u_k$ is a subword of the overlap word $\ell_tu_i$
associated to $s_t$, in order to show that $s_t$ reduces to zero,
all we have to do is to show that the two S-polynomials
$$s_u =
\begin{cases}
c_{u}(x_{k_{1}}\hdots x_{j_{\beta-D}})g_i
- c'_{u}g_k(x_{i_{D+2}}\hdots x_{i_{\alpha}})
& \text{if $\gamma > D+1$} \\
c_{u}g_i - c'_{u}\ell'_ug_k(x_{i_{D+2}}\hdots x_{i_{\alpha}})
& \text{if $\gamma \leqslant D+1$}
\end{cases}
$$
% $$s_u = c_{u}(x_{k_{1}}\hdots x_{j_{\beta-D}})g_i
% - c'_{u}\ell'_ug_k(x_{i_{D+2}}\hdots x_{i_{\alpha}})$$
and
$$s_v = c_{v}g_jx_{i_{D+1}} 
- c'_{v}(x_{j_1}\hdots x_{j_{\beta+1-\gamma}})g_k$$
reduce to zero ($1 \leqslant u,v \leqslant |S|$).

The S-polynomial $s_v$ reduces to zero by comparison with
part (2). For the S-polynomial $s_u$, first note that if $\alpha = D+1$,
then either $u_k$ is a suffix of $u_i$,
$u_i$ is a suffix of $u_k$, or $u_k = u_i$; it follows that
$s_u$ reduces to zero trivially if $u_k = u_i$, and $s_u$ reduces to
zero by part (1) in the other two cases.

If however $\alpha \neq D+1$,
then either $s_u$ is a middle overlap (if $\gamma < D+1$), a left overlap
(if $\gamma = D+1$), or another prefix overlap.
% (if $\gamma = \beta+1$, $\alpha = \beta$).
The first two cases can be handled by parts (2) and (3)
respectively; the final case is handled by induction,
where we note that after each step of the induction, the value
$\alpha+\beta-2D$ strictly decreases (regardless of
which case $T(u_j, x^R_{i_{D+1}}) = 0$
or $T(u_i, x^L_{j_{\beta-D}}) = 0$ applies), so we are
guaranteed at some stage to find an overlap that is not
a prefix overlap, enabling us to verify that the
S-polynomial $s_t$ conventionally reduces to zero.

\item Case $T(u_i, x^L_{j_{\beta-D}}) = 0$. \\[2mm]
Because we know that the prolongation
$x_{j_{\beta - D}}g_i$ involutively reduces to zero, 
there must exist a monomial
$u_{k} = x_{k_1}\hdots x_{k_{\gamma}} \in U$ such that
$u_{k}$ involutively divides $x_{j_{\beta - D}}u_i$.
This $u_k$ must be a prefix of $x_{j_{\beta - D}}u_i$
(otherwise, assuming that $x_{k_{\gamma}} = x_{i_{\kappa}}$,
we have $T(u_k, x^L_{j_{\beta-D}}) = 0$ if $\gamma = \alpha$
(because of the overlap between $u_j$ and $u_k$);
$T(u_k, x^L_{i_{\kappa-\gamma}}) = 0$ if $\gamma < \alpha$ and
$\kappa = \alpha$ (because of the overlap between
$u_i$ and $u_k$); and $T(u_k, x^R_{i_{\kappa + 1}}) = 0$
if $\gamma < \alpha$ and $\kappa < \alpha$ 
(again because of the overlap between $u_i$ and $u_k$)).
$$\xymatrix @R=0.5pc @C=2.4pc{
u_i = &&&& \ar@{-}[r]_{x_{i_1}} & \ar@{--}[r]
& \ar@{-}[r]_{x_{i_{D}}}
& \ar@{-}[r]_{x_{i_{D+1}}} & \ar@{--}[r]
& \ar@{-}[r]_{x_{i_{\alpha - 1}}} & \ar@{-}[r]_{x_{i_{\alpha}}} & \\
u_j = & \ar@{-}[r]_{x_{j_1}} & \ar@{--}[r] & \ar@{-}[r]_{x_{j_{\beta-D}}}
& \ar@{-}[r]_{x_{j_{\beta-D+1}}} & \ar@{--}[r]
& \ar@{-}[r]_{x_{j_{\beta}}} & \\
u_k = &&& \ar@{-}[r]_{x_{k_1}} & \ar@{~}[rrrrrrr] &&&&&&&
}$$

Let us now
use the monomial $u_k$ together with Buchberger's Second Criterion
to simplify our goal of showing that the S-polynomial $s_t$ reduces
to zero. Because $u_k$ is a subword of the overlap word $\ell_tu_i$
associated to $s_t$, in order to show that $s_t$ reduces to zero,
all we have to do is to show that the two S-polynomials
$$s_u = c_{u}x_{k_{1}}g_i - c'_ug_k(x_{i_{\gamma}}\hdots x_{i_{\alpha}})$$
and
$$s_v =
\begin{cases}
c_{v}g_j(x_{i_{D+1}}\hdots x_{k_{\gamma}})
- c'_{v}(x_{j_1}\hdots x_{j_{\beta-D-1}})g_k
& \text{if $\gamma > D+1$} \\
c_{v}g_j - c'_{v}(x_{j_1}\hdots x_{j_{\beta-D-1}})g_kr'_v
& \text{if $\gamma \leqslant D+1$}
\end{cases}
$$
% $$s_v = c_{v}g_jr_v - c'_{v}(x_{j_1}\hdots x_{j_{\beta-D-1}})g_kr'_v$$
reduce to zero ($1 \leqslant u,v \leqslant |S|$).

The S-polynomial $s_u$ reduces to zero by comparison with
part (2). For the S-polynomial $s_v$, first note that if $\beta-D = 1$,
then either $u_k$ is a prefix of $u_j$,
$u_j$ is a prefix of $u_k$, or $u_k = u_j$; it follows that
$s_v$ reduces to zero trivially if $u_k = u_j$, and $s_v$ reduces to
zero by part (3) in the other two cases.

If however $\beta-D \neq 1$,
then either $s_v$ is a middle overlap (if $\gamma < D+1$), a right overlap
(if $\gamma = D+1$), or another prefix overlap.
% (if $\gamma = \beta+1$, $\alpha = \beta$).
The first two cases can be handled by parts (2) and (1)
respectively; the final case is handled by induction,
where we note that after each step of the induction, the value
$\alpha+\beta-2D$ strictly decreases (regardless of
which case $T(u_j, x^R_{i_{D+1}}) = 0$
or $T(u_i, x^L_{j_{\beta-D}}) = 0$ applies), so we are
guaranteed at some stage to find an overlap that is not
a prefix overlap, enabling us to verify that the
S-polynomial $s_t$ conventionally reduces to zero.
\end{itemize}
\end{pf}

% \begin{remark}
% Part (1) of the above proof also appears in Chapter \ref{ChNCIB}.
% \end{remark}

%
% Appendix B
% Author: Gareth Evans
% Last Modified: 10th September 2005
%

\chapter{Source Code} \label{appB}

In this Appendix, we will present ANSI C
source code for an initial implementation of the
noncommutative Involutive Basis algorithm
(Algorithm \ref{noncom-inv}), together with
an introduction to {\sf AlgLib}, a set of
ANSI C libraries providing data
types and functions that serve as building
blocks for the source code.

\section{Methodology}

A problem facing anyone wanting to implement mathematical ideas is
the choice of language or system in which to do the implementation.
The decision depends on the task at hand.  If all that is required is
a convenient environment for prototyping ideas, a symbolic computation
system such as Maple \cite{Maple05}, 
Mathematica \cite{Mathematica04} or MuPAD \cite{MuPAD} may suffice.  Such
systems have a large collection of mathematical data types, functions
and algorithms already present; tools that will not be available in a
standard programming language.  There is however always a price to pay
for convenience.  These common systems are all interpreted and use
a proprietary programming syntax, making it it difficult to use other
programs or libraries within a session.  It also makes such systems less
efficient than the execution of compiled programs.

The {\sf AlgLib} libraries can be said to 
provide the best of both worlds, as
they provide data types, functions and algorithms to allow programmers
to more easily implement certain mathematical algorithms (including
the algorithms described in this thesis) in the ANSI C programming
language. For example, {\sf AlgLib} contains the {\sf FMon} \cite{FMon01}
and {\sf FAlg} \cite{FAlg01} libraries, respectively
containing data types and functions to perform computations
in the free monoid on a set of symbols
and the free associative algebra on a set of symbols.
Besides the benefit of the efficiency of compiled programs, 
the strict adherence to ANSI C makes
programs written using the libraries highly portable.

% What is the best way of implementing a computationally
% intense mathematical algorithm? Do we use a symbolic
% computation system such as Maple \cite{Maple05}
% or Mathematica \cite{Mathematica04},
% or do we start from scratch using a programming language
% such as C or Java?
%
% Certainly it is easier to start with a symbolic
% computation system, as a large collection of mathematical
% data types, functions and algorithms will already
% be present, tools that will not be
% available in a standard programming language. However,
% algorithms implemented in a symbolic computation
% system may not be as efficient as those implemented
% using a programming language, as (for example) the code may be
% interpreted or many other superfluous
% tools may be loaded into memory.
%
% The {\sf AlgLib} libraries can be said to provide
% the best of both worlds, as they provide % the
% data types, functions and algorithms to allow % that allow a
% programmers to more easily implement certain mathematical algorithms
% (including the algorithms described in this thesis)
% in the ANSI C programming language. For
% example, {\sf AlgLib} contains the {\sf FMon} \cite{FMon01}
% and {\sf FAlg} \cite{FAlg01} libraries, respectively
% containing data types and functions to perform computations
% in the free monoid on a set of symbols
% and the free associative algebra on a set of symbols.

\subsection{MSSRC}

{\sf AlgLib} is supplied by MSSRC \cite{MSSRC}, a company whose
Chief Scientist is Prof. Larry Lambe, an honorary
professor at the University of Wales, Bangor. For an
introduction to MSSRC, we quote the following passage
from \cite{LarLet}.

\begin{quote}
Multidisciplinary Software Systems Research Corporation (MSSRC)
was conceived as a company devoted to furthering the long-term
effective use of mathematics and mathematical computation.
MSSRC researches, develops, and markets advanced mathematical
tools for engineers, scientists, researchers, educators,
students and other serious users of mathematics. These tools
are based on providing levels of power, productivity
and convenience far greater than existing tools while
maintaining mathematical rigor at all times. The company also
provides computer education and training.

MSSRC has several lines of ANSI C libraries for providing
mathematical support for research and implementation of
mathematical algorithms at various levels of complexity.
No attempt is made to provide the user of these libraries
with any form of Graphical User Interface (GUI). All
components are compiled ANSI C functions which represent
various mathematical operations from basic (adding,
subtracting, multiplying polynomials, etc.) to advanced
(operations in the free monoid on an arbitrary number of
symbols and beyond). In order to use the libraries
effectively, the user must be expert at ANSI C programming,
e.g., in the style of Kernighan and Richie \cite{KandR}
and as such, they are not suited for the casual user.
This does not imply in any way that excellent user
interfaces for applications of the libraries cannot be
supplied or are difficult to implement by well
experienced programmers.

The use of MSSRC's libraries has been reported in a
number of places such as \cite{M1}, \cite{M2}, \cite{M3},
\cite{M4} and elsewhere.
\end{quote}

\subsection{AlgLib}

To give a taste of how 
{\sf AlgLib} has been used to implement
the algorithms considered in this thesis, consider one
of the basic operations of these algorithms, the task of
subtracting two polynomials to yield a third polynomial
(an operation essential for computing an S-polynomial).
In ordinary ANSI C, there is no data type for a polynomial,
and certainly no function for subtracting two
polynomials; {\sf AlgLib} however does supply
these data types and functions, both in the commutative and
noncommutative cases. For example, the {\sf AlgLib} data type for a 
noncommutative polynomial is an {\it FAlg},
and the {\sf AlgLib} function for subtracting two such polynomials
is the function {\it fAlgMinus}. It follows that we
can write ANSI C code for subtracting two 
noncommutative polynomials,
as illustrated below where we subtract the polynomial
$2b^2+ab+4b$ from the polynomial $2\times(b^2+ba+3a)$.
\lstset{
language=C,
frame=single,
framesep=5pt,
columns=fullflexible,
basicstyle=\scriptsize,
commentstyle={\it \textcolor[rgb]{0.4,0.4,0.4} },
stringstyle=\ttfamily,
morekeywords={FAlg,  QInteger}
}

{\bf Source Code}
\begin{lstlisting}
# include <fralg.h>

int 
main( argc, argv )
int argc;
char *argv[];
{
  // Define Variables
  FAlg p, q, r;
  QInteger two;

  // Set Monomial Ordering (DegLex)
  theOrdFun = fMonTLex;

  // Initialise Variables
  p = parseStrToFAlg("b^2 + b*a + 3*a");
  q = parseStrToFAlg("2*b^2 + a*b + 4*b");
  two = parseStrToQ("2");

  // Perform the calculation and display the result on screen
  r = fAlgMinus( fAlgScaTimes( two, p ), q );
  printf("2*(%s) - (%s) = %s\n", fAlgToStr( p ), fAlgToStr( q ), fAlgToStr( r ) );

  return EXIT_SUCCESS;
}
\end{lstlisting}

{\bf Program Output}
\lstset{
language=,
columns=fullflexible,
basicstyle=\scriptsize,
commentstyle={\it \textcolor[rgb]{0.4,0.4,0.4} },
stringstyle=\ttfamily,
morekeywords={FAlg,  QInteger}
}
\begin{lstlisting}
ma6:mssrc-aux/thesis> fAlgMinusTest
2*(b^2 + b a + 3 a) - (2 b^2 + a b + 4 b) = 2 b a - a b -4 b + 6 a
ma6:mssrc-aux/thesis>
\end{lstlisting}

\section{Listings} \label{Bpoint2}

Our implementation of the noncommutative Involutive Basis
algorithm is arranged as follows: {\it involutive.c} is the main
program, dealing with all the input and output and calling
the appropriate routines; the `{\it \_functions}' files contain
all the procedures and functions used by the program;
and {\it README} describes how to use the program, including
what format the input files should take and what the
different options of the program are used for.

In more detail, {\it arithmetic\_functions.c} contains
functions for dividing a polynomial by its
(coefficient) greatest common divisor and for converting
user specified generators to ASCII generators (and vice-versa);
{\it file\_functions.c} contains all the functions
needed to read and write polynomials and variables to
and from disk; {\it fralg\_functions.c} contains
functions for monomial orderings, polynomial division
and reduced Gr\"obner Bases computation;
{\it list\_functions.c} contains some extra functions
needed to deal with displaying, sorting and
manipulating lists;
and {\it ncinv\_functions.c} contains all the
involutive routines, for example the Involutive Basis algorithm
itself and associated functions for determining
multiplicative variables and for performing autoreduction.
\vfill
%\begin{doublespace}
\begin{center}
\begin{tabular}{lp{12.5cm}}
\multicolumn{2}{l}{{\large {\bf Contents}}} \\[3mm]
\ref{appBf1}  & README                  \dotfill\ \pageref{appBf1} \\[2mm]
\ref{appBf2}  & arithmetic\_functions.h \dotfill\ \pageref{appBf2} \\[2mm]
\ref{appBf3}  & arithmetic\_functions.c \dotfill\ \pageref{appBf3} \\[2mm]
\ref{appBf4}  & file\_functions.h       \dotfill\ \pageref{appBf4} \\[2mm]
\ref{appBf5}  & file\_functions.c       \dotfill\ \pageref{appBf5} \\[2mm]
\ref{appBf6}  & fralg\_functions.h      \dotfill\ \pageref{appBf6} \\[2mm]
\ref{appBf7}  & fralg\_functions.c      \dotfill\ \pageref{appBf7} \\[2mm]
\ref{appBf8}  & list\_functions.h       \dotfill\ \pageref{appBf8} \\[2mm]
\ref{appBf9}  & list\_functions.c       \dotfill\ \pageref{appBf9} \\[2mm]
\ref{appBf10} & ncinv\_functions.h      \dotfill\ \pageref{appBf10} \\[2mm]
\ref{appBf11} & ncinv\_functions.c      \dotfill\ \pageref{appBf11} \\[2mm]
\ref{appBf12} & involutive.c            \dotfill\ \pageref{appBf12}
\end{tabular}
\end{center}
%\end{doublespace}

\lstset{
language=,
frame=,
columns=flexible,
basicstyle=\scriptsize,
commentstyle={\it \textcolor[rgb]{0.4,0.4,0.4} },
}

\vspace*{-24pt}
\subsection{README} \label{appBf1}
\begin{lstlisting}
****************************************************************
* HOW TO USE THE INVOLUTIVE PROGRAM - QUICK GUIDE *
****************************************************************

NAME
       involutive - Computes Noncommutative Involutive Bases for ideals.

SYNOPSIS
       involutive [OPTION]... [FILE]...

DESCRIPTION
       Here are the options for the program.

       -a
            e.g. > involutive -d -a file.in
            Optimises the lexicographical ordering according to
            the frequency of the variables in the input basis
            (most frequent = lexicographically smallest).

       -c(n)
            e.g. > involutive -c2 file.in
            Chooses which involutive algorithm to use.
            n is a required number between 1 and 2.

            1: *DEFAULT * Gerdt's Algorithm
            2: Seiler's Algorithm

       -d
            e.g. > involutive -d file.in
            Allows the user to calculate a DegLex
            Involutive Basis for the basis in file.in.

       -e(n)
            e.g. > involutive -e2 -s2 file.in
            Allows the user to select the type of Overlap
            Division to use. n is a required number between
            1 and 5. Note: Must be used with either the
            -s1 or -s2 options.

            Left Overlap Division:

                   A          B              C              D
                -----  -------  ---------  ------
            -----x        ----x           ----x            x---

            1: * DEFAULT * A, B, C (weak, Gr\"obner)
            2: A, B, C, Strong (strong if used with -m2)
            3: A, B, C, D (weak, Gr\"obner)
            4: A, B (weak, Gr\"obner)
            5: A (weak, Gr\"obner)

            Right Overlap Division:

                A              B            C             D
            -----        -------  ---------  ------
              x-----           x---       x----       ---x

            1: * DEFAULT * A, B, C (weak, Gr\"obner)
            2: A, B, C, Strong (strong if used with -m2)
            3: A, B, C, D (weak, Gr\"obner)
            4: A, B (weak, Gr\"obner)
            5: A (weak, Gr\"obner)

       -f
            e.g. > involutive -f file.in
            Removes any fractions from the input basis.

       -l
            e.g. > involutive -l file.in
            Allows the user to calculate a Lex
            Involutive Basis for the basis in file.in.
            Warning: program may go into an infinite loop
            (Lex is not an admissible monomial ordering).

       -m(n)
            e.g. > involutive -m2 file.in
            Selects which method of deciding whether a monomial
            involutively divides another monomial is used.
            n is a required number between 1 and 2.

            1: * DEFAULT * 1st letters on left and right (thin divisor)
            2: All letters on left and right (thick divisor)

       -o(n)
            e.g. > involutive -o2 file.in
            Allows the user to select how the basis is sorted
            during the algorithm. n is a required number between
            1 and 3.

            1: * DEFAULT * DegRevLex Sorted
            2: No Sorting
            3: Sorting by Main Ordering

       -p
            e.g. > involutive -l -p file.in
            An interactive Ideal Membership Problem Solver.
            There are two ways the solver can be used:
            either a file containing a list of polynomials
            (e.g. x*y-z;
                  x^2-z^2+y^2; ) can be given, or the
            polynomials can be input
            manually (e.g. x*y-z). The solver tests to see
            whether the Involutive Basis computed in the
            algorithm reduces the polynomials given to zero.

       -r   * DEFAULT *
            e.g. > involutive -r file.in
            Allows the user to calculate a DegRevLex
            Involutive Basis for the basis in file.in.

       -s(n)
            e.g. > involutive -s2 file.in
            Allows the user to select the type of Involutive
            Basis to calculate. n is a required number between
            1 and 5. Note: If an `Overlap' Division is selected,
            the type of Overlap Division can be chosen with
            the -e(n) option.

            1: Left Overlap Division (local, cts, see -e option)
            2: Right Overlap Division (local, cts, see -e option)
            3: * DEFAULT * Left Division (global, cts, strong)
            4: Right Division (global, cts, strong)
            5: Empty Division (global, cts, strong)

       -v(n)
            e.g. > involutive -v3 file.in
            Changes the amount of information given out by the
            program (i.e. the `verbosity' of the program).
            n is a number between 0 and 9. Rough Guide:

            0: Silent (no output given).
            1: * DEFAULT *
            2: Returns Number of Reductions Carried Out,
               Prints Out Every Polynomial Found
            3: More Autoreduction Information,
               Prolongation Information
            4: More Details of Steps Taken in Algorithm
            5: More Global Division Information
            6: Step-by-Step Reduction, Overlap Information
            7: Shows Multiplicative Grids
            8: More Overlap Division Information
            9: All Other Information

       -w
            e.g. > involutive -w file.in
            Allows the user to calculate an Involutive Basis
            for the basis in file.in using the Wreath
            Product Monomial Ordering.

       -x
            e.g. > involutive -x file.in
            Ignores any prolongations of degree greater than or 
            equal to 2d, where d is a value determined by the degree 
            of the largest degree lead monomial in the current minimal basis.
            Warning: May not return a valid Involutive Basis
            (only a valid Gr\"obner Basis).

FILE FORMATS

       There is one file format for the input basis:

       IDEALS:
       --------------
       x; y; z;
       x*y - z;
       2*x + y*z + z;
       --------------

       First line = List of variables in order. In the
                    above, x; y; z; represents x > y > z.
       Remaining lines = Polynomial generators (which must be
                         terminated by semicolons).

OUTPUT

       As output, the program provides a reduced Gr\"obner Basis and
       an Involutive Basis for the input ideal (if it can calculate it).

       For the following, assume that our input basis was given as file.in.
       * If a DegRevLex Gr\"obner Basis is calculated, it is stored as file.drl.
       * If a DegLex Gr\"obner Basis is calculated, it is stored as file.deg.
       * If a Lex Gr\"obner Basis is calculated, it is stored as file.lex.
       * If a Wreath Product Gr\"obner Basis is calculated, it is stored as file.wp.

       The Involutive Basis is given as <Gr\"obner Basis>.inv.
       For example, if a DegLex Involutive Basis is calculated,
       it is stored as file.deg.inv.

       Note that the program has the ability to recognise the .in suffix and
       replace it with .drl, .deg, .lex or .wp as necessary.
       If your input file does not have a .in suffix then the program will
       simply append the appropriate suffix onto the end of the file name.
       For example, using the command
       > involutive FILE
       we obtain file.drl if FILE = file.in
       and obtain e.g. file.other.drl if FILE = file.other.
\end{lstlisting}

\lstset{
language=C,
columns=fullflexible,
basicstyle=\scriptsize,
commentstyle={\it \textcolor[rgb]{0.4,0.4,0.4} },
stringstyle=\ttfamily,
numbers=left,
numberstyle=\textcolor[rgb]{0.53,0.53,0.53},
numbersep=5pt,
morekeywords={
String, Pointer, ULong, UInt, Length, FAlg, FILE, Short,
Short, Long, Bool, StrList, IntList, FAlgPair, FAlgList,
FAlgPairList,  QInteger, FMon, FMonList, FMonPairList,
FMonPair, FMonTriple, IntFMonList, IntFMonPairList,
QIntegerList, QPolyList, QPoly, CritPairList, ULPolList,
ULMonom, Integer
}
}

\subsection{arithmetic\_functions.h} \label{appBf2}
\begin{lstlisting}
/*
 * File: arithmetic_functions.h
 * Author: Gareth Evans
 * Last Modified: 29th September 2004
 */

// Initialise file definition
# ifndef ARITHMETIC_FUNCTIONS_HDR
# define ARITHMETIC_FUNCTIONS_HDR

// Include MSSRC Libraries
# include <fralg.h>

//
// Numerical Functions
//

// Returns the numerical value of a 3 letter word
ULong ASCIIVal( String );
// Returns the 3 letter word of a numerical value
String ASCIIStr( ULong );
// Returns the monomial corresponding to the 3 letter word of a numerical value
FMon ASCIIMon( ULong );

//
// QInteger Functions
//

// Calculate Alternative LCM of 2 QIntegers
QInteger AltLCMQInteger( QInteger, QInteger );

//
// FAlg Functions
//

// Divides the input FAlg by its common GCD
FAlg findGCD( FAlg );
// Returns maximal degree of lead term for the given FAlgList
ULong maxDegree( FAlgList );
// Returns the position of the smallest LM(g) in the given FAlgList
ULong fAlgListLowest( FAlgList );

# endif // ARITHMETIC_FUNCTIONS_HDR
\end{lstlisting}

\subsection{arithmetic\_functions.c} \label{appBf3}
\begin{lstlisting}
/*
 * File: arithmetic_functions.c
 * Author: Gareth Evans
 * Last Modified: 11th February 2005
 */

/*
 * ===================
 * Numerical Functions
 * ===================
 */

/*
 * Function Name: ASCIIVal
 *
 * Overview: Returns the numerical value of a 3 letter word
 *
 * Detail: Given a String containing 3 letters from the set
 * {A, B, ..., Z}, this function returns the numerical
 * value of the String according to the following rule:
 * AAA = 1, AAB = 2, ..., AAZ = 26, ABA = 27, ABB = 28,
 * ..., ABZ = 52, ACA = 53, ...
 *
 */
ULong
ASCIIVal( word )
String word;
{
  ULong back = 0;

  // Add on 17576*value of 1st letter (A = 0, B = 1, ...)
  back = back + 17576*( (ULong)( (int)word[0] - (int)'A' ) );
  // Add on 26*value of 2nd letter (A = 0, B = 1, ...)
  back = back + 26*( (ULong)( (int)word[1] - (int)'A' ) );
  // Add on the value of the 3rd letter (A = 1, B = 2, ...)
  back = back + (ULong)( (int)word[2] - (int)'A' + 1 );

  return back;
}

/*
 * Function Name: ASCIIStr
 *
 * Overview: Returns the 3 letter word of a numerical value
 *
 * Detail: Given a ULong, this function returns the
 * 3 letter String corresponding to the following rule:
 * 1 = AAA, 2 = AAB, ..., 26 = AAZ, 27 = ABA, 28 = ABB,
 * ..., 52 = ABZ, 53 = ACA, ...
 *
 */
String
ASCIIStr( number )
ULong number;
{
  String back = strNew();
  int i = 0, j = 0, k;

  // Take away multiples of 26^2 to get the first letter
  while( number > 17576 )
  {
    i++;
    number = number - 17576;
  }

  // Take away multiples of 26 to get the second letter
  while( number > 26 )
  {
    j++;
    number = number - 26;
  }

  // We are now left with the third letter
  k = (int) number - 1;

  // Convert the numbers to a String
  sprintf( back, "%c%c%c", (char)( (int)'A' + i ),
                           (char)( (int)'A' + j ),
                           (char)( (int)'A' + k ) );

  // Return the three letters
  return back;
}

/*
 * Function Name: ASCIIStr
 *
 * Overview: Returns the monomial corresponding to the 
 *           3 letter word of a numerical value
 *
 * Detail: Given a ULong, this function returns the
 * monomial corresponding to the following rule:
 * 1 = AAA, 2 = AAB, ..., 26 = AAZ, 27 = ABA, 28 = ABB,
 * ..., 52 = ABZ, 53 = ACA, ...
 *
 */
FMon
ASCIIMon( number )
ULong number;
{
  // Obtain the String corresponding to the input
  // number and change it to an FMon
  return parseStrToFMon( ASCIIStr( number ) );
}

/*
 * ==================
 * QInteger Functions
 * ==================
 */

/*
 * Function Name: AltLCMQInteger
 *
 * Overview: Calculates an 'alternative' LCM of 2 QIntegers
 *
 * Detail: Given two QIntegers a = an/ad and b = bn/bd,
 * this function calculates the LCM given
 * by alt_lcm(a, b) = (a*b)/(alt_gcd(a, b))
 * = (an*bn*ad*bd)/(ad*bd*gcd(an, bn)*gcd(ad, bd))
 * = (an*bn)/(gcd(an, bn)*gcd(ad, bd)).
 *
 */
QInteger
AltLCMQInteger( a, b )
QInteger a, b;
{
  Integer an = a -> num,
          ad = a -> den,
          bn = b -> num,
          bd = b -> den;

  return qDivide( zToQ( zTimes( an, bn ) ),
                  zToQ( zTimes( zGcd( an, bn ), zGcd( ad, bd ) ) ) );
}

/*
 * ==============
 * FAlg Functions
 * ==============
 */

/*
 * Function Name: findGCD
 *
 * Overview: Divides the input FAlg by its common GCD
 *
 * Detail: Given an FAlg, this function divides the
 * polynomial by its common GCD so that the output
 * polynomial g cannot be written as g = cg', where
 * g' is a polynomial and c is an integer, c > 1.
 *
 */
FAlg
findGCD( input )
FAlg input;
{
  FAlg output = input, process = input;
  QInteger coef;
  Integer GCD = zOne, numerator, denominator;
  Bool first = 0, allNeg = qLess( fAlgLeadCoef( input ), qZero() );

  if( (ULong) fAlgNumTerms( input ) == 1 ) // If poly has just 1 term
  {
    // Return that term with a unit coefficient
    return fAlgMonom( qOne(), fAlgLeadMonom( input ) );
  }
  else // Poly has more than 1 term
  {
    while( process ) // Go through each term
    {
      coef = fAlgLeadCoef( process );    // Read the lead coefficient
      numerator = coef -> num;           // Break the coefficient down
      denominator = coef -> den;         // into a numerator and a denominator
      process = fAlgReductum( process ); // Get ready to look at the next term

      if( zIsOne( denominator ) != (Bool) 1 ) // If we encounter a fraction
      {
        return input; // We cannot divide through by a GCD so just return the input
      }
      else // The coefficient was an integer
      {
        if( first == 0 ) // If this is the first term
        {
          first = (Bool) 1;
          GCD = numerator; // Set the GCD to be the current numerator
        }
        else // Recursively calculate the GCD
          GCD = zGcd( GCD, numerator );
      }
    }

    if( zLess( GCD, zZero ) == (Bool) 1 ) // If the GCD is negative
      GCD = zNegate( GCD ); // Negate the GCD
    if( zLess( zOne, GCD ) == (Bool) 1 ) // If the GCD is > 1
      output = fAlgZScaDiv( output, GCD ); // Divide the poly by the GCD
  }

  if( allNeg == (Bool) 1 ) // If the original coefficient was negative
    return fAlgZScaTimes( zMinusOne, output ); // Return the negated polynomial
  else
    return output;
}

/*
 * Function Name: maxDegree
 *
 * Overview: Returns maximal degree of lead term for the given FAlgList
 *
 * Detail: Given an FAlgList, this function calculates the degree
 * of the lead term for each element of the list and returns
 * the largest value found.
 *
 */
ULong
maxDegree( input )
FAlgList input;
{
  ULong test, output = 0;

  while( input ) // For each polynomial in the list
  {
    // Calculate the degree of the lead monomial
    test = fMonLength( fAlgLeadMonom( input -> first ) );
    if( test > output ) output = test;
    input = input -> rest; // Advance the list
  }

  // Return the maximal value
  return output;
}

/*
 * Function Name: fAlgListLowest
 *
 * Overview: Returns the position of the smallest LM(g) in the given FAlgList
 *
 * Detail: Given an FAlgList, this function looks at all the leading
 * monomials of the elements in the list and returns the position of
 * the smallest lead monomial with respect to the monomial ordering
 * currently being used.
 *
 */
ULong
fAlgListLowest( input )
FAlgList input;
{
  ULong output = 0, i, len = fAlgListLength( input );
  FMon next, lowest;

  if( input ) // Assume the 1st lead monomial is the smallest to begin with
  {
    lowest = fAlgLeadMonom( input -> first );
    output = 1;
  }
  for( i = 1; i < len; i++ ) // For the remaining polynomials
  {
    input = input -> rest;
    // Extract the next lead monomial
    next = fAlgLeadMonom( input -> first );

    // If this lead monomial is smaller than the current smallest
    if( theOrdFun( next, lowest ) == (Bool) 1 )
    {
      // Make this lead monomial the smallest
      output = i+1;
      lowest = fAlgScaTimes( qOne(), next );
    }
  }

  // Return position of smallest lead monomial
  return output;
}

/*
 * ===========
 * End of File
 * ===========
 */
\end{lstlisting}

\subsection{file\_functions.h} \label{appBf4}
\begin{lstlisting}
/*
 * File: file_functions.h
 * Author: Gareth Evans
 * Last Modified: 14th July 2004
 */

// Initialise file definition
# ifndef FILE_FUNCTIONS_HDR
# define FILE_FUNCTIONS_HDR

// Include MSSRC Libraries
# include <fralg.h>

// MAXLINE denotes the length of the longest allowable line in a file
# define MAXLINE 5000

//
// Low Level File Handling Functions
//

// Read a line from a file; return length
int getLine( FILE *, char[], int );
// Pick an integer from a list such as "2, 5, 6,"
int intFromStr( char[], int, int * );
// Pick a variable from a list such as "a; b; c;"
String variableFromStr( char[], int, int * );
// Pick an FMon from a list such as "a; b; c;"
FMon fMonFromStr( char[], int, int * );
// Pick an FAlg from a string such as "x*y - z;"
FAlg fAlgFromStr( char[], int, int * );

//
// High Level File Reading Functions
//

// Routine to read an FMonList from the first line of a file
FMonList fMonListFromFile( FILE * );
// Routine to read an FAlgList from a file
FAlgList fAlgListFromFile( FILE * );

//
// High Level File Writing Functions
//

// Writes an FMon (in parse format) followed by a semicolon to a file
void fMonToFile( FILE *, FMon );
// Writes an FMonList to a file on a single line
void fMonListToFile( FILE *, FMonList );

//
// File Name Modification Functions
//

// Appends ".drl" onto a string (except in special case "*.in")
String appendDotDegRevLex( char[] );
// Appends ".deg" onto a string (except in special case "*.in")
String appendDotDegLex( char[] );
// Appends ".lex" onto a string (except in special case "*.in")
String appendDotLex( char[] );
// Appends ".wp" onto a string (except in special case "*.in")
String appendDotWP( char[] );
// Calculates the length of an input string
int filenameLength( char[] );

# endif // FILE_FUNCTIONS_HDR
\end{lstlisting}

\subsection{file\_functions.c} \label{appBf5}
\begin{lstlisting}
/*
 * File: file_functions.c
 * Author: Gareth Evans
 * Last Modified: 16th August 2004
 */

/*
 * ==================================
 * Low Level File Handling Functions
 * (Used in the high level functions)
 * ==================================
 */

/*
 * Function Name: getLine
 *
 * Overview: Read a line from a file; return length
 *
 * Detail: Given a file _infil_, we read the first line
 * of the file, placing the contents into the string _s_.
 * The third parameter _lim_ determines the maximum length
 * of any line to be returned (when we call the function
 * this is usually MAXLINE); the returned integer tells
 * us the length of the line we have just read.
 *
 * Known Issues: The length of a line is sometimes returned
 * incorrectly when a file saved in Windows is used
 * on a UNIX machine. Resave your file in UNIX.
 */
int
getLine( infil, s, lim )
FILE *infil;
char s[];
int lim;
{
  int c, i;

  /*
   * Place characters in _s_ as long as (1) we do not exceed _lim_ number of
   * characters; (2) the end of the file is not encountered; (3) the end of the
   * line is not encountered.
   */
  for( i = 0; ( i < lim-1 ) && ( ( c = fgetc(infil) ) != -1 ) && ( c != (int)'\n' ); i++ )
  {
    s[i] = (char)c;
  }
  if( c == (int)'\n' ) // if the for loop was terminated due to reaching end of line
  {
    s[i] = (char)c; // add the newline character to our string
    i++;
  }
  s[i] = '\0'; // '\0' is the null character

  return i-1; // The -1 is used to compensate for the null character
}

/*
 * Function Name: intFromStr
 *
 * Overview: Pick an integer from a list such as "2, 5, 6,"
 *
 * Detail: Starting from position _j_ in a string _s_,
 * read in an integer and return it. Note that the integer
 * in the string must be terminated with a comma and that
 * the sign of the integer is taken into account.
 * Once the integer has been read, place the position we
 * have reached in the string in the variable _pk_.
 */
int
intFromStr( s, j, pk )
char s[];
int j, *pk;
{
  char c;
  int n = 0, sign = 1, k = j;
  c = s[k];

  // Traverse through any empty space
  while( c == ' ' )
  {
    k++;
    c = s[k];
  }

  // If a sign is present, process it
  if( c == '+' )
  {
    k++;
    c = s[k];
  }
  else if( c == '-' )
  {
    sign = -1;
    k++;
    c = s[k];
  }

  // Until a comma is encountered (signalling the
  // end of the integer)
  while( c != ',' )
  {
    if( ( c >= '0' ) && ( c <= '9' ) )
    {
      n = 10*n + (int)(c - '0'); // the "- '0'" is needed to get the correct integer
    }
    else
    {
      printf("Error: Incorrect Input in File (%c is not a number).\n", c);
      exit( EXIT_FAILURE );
    }
    k++;
    c = s[k];
  }
  *pk = k+1; // return the finishing position

  /*
   * Note: In this function we return *pk = k+1 and not *pk = k as
   * in subsequent functions because this function has a slightly
   * different structure due to having to deal with the + and -
   * characters at the beginning of the string.
   */
    
  return sign*n; // return the integer
}

/*
 * Function Name: variableFromStr
 *
 * Overview: Pick a variable from a list such as "a; b; c;"
 *
 * Detail: Starting from position _j_ in a string _s_,
 * read in a String and return it. Note that the String
 * in the string must be terminated with a semicolon.
 * Once the String has been read, place the position we
 * have reached in the string in the variable _pk_.
 */
String
variableFromStr( s, j, pk )
char s[];
int j, *pk;
{
  char c = ' ';
  int i = 0, k = j;
  String back = strNew(), concat;

  sprintf( back, "" ); // Initialise back

  // Until a semicolon is encountered
  while( c != ';' )
  {
    c = s[k]; // Pick a character from the string

    // If a semicolon was encountered
    if( c == ';' )
    {
      concat = strNew();
      sprintf( concat, "%c", '\0' );
      back = strConcat( back, concat ); // Finish with the null character
    }
    else if( c != ' ' )
    {
      concat = strNew();
      sprintf( concat, "%c", c );
      // Transfer character to output String
      if( i == 0 ) back = strCopy( concat );
      else back = strConcat( back, concat );
      i++;
    }
    k++;
  }
  *pk = k; // Place finish position in the variable _pk_

  return back; // Return the String
}

/*
 * Function Name: fMonFromStr
 *
 * Overview: Pick an FMon from a list such as "a; b; c;"
 *
 * Detail: Starting from position _j_ in a string _s_,
 * read in an FMon and return it. Note that the FMon
 * in the string must be terminated with a semicolon.
 * Once the FMon has been read, place the position we
 * have reached in the string in the variable _pk_.
 */
FMon
fMonFromStr( s, j, pk )
char s[];
int j, *pk;
{
  char c = ' ', a[MAXLINE];
  int i = 0, k = j;
  FMon back;

  // Until a semicolon is encountered
  while( c != ';' )
  {
    c = s[k]; // Pick a character from the string

    // If we have found a semicolon
    if( c == ';' )
    {
      a[i] = '\0'; // Finish the string with the null character
    }
    else
    {
      a[i] = c; // Continue to process...
      i++;
    }
    k++;
  }
  *pk = k; // Place the finish position in the variable _pk_

  back = parseStrToFMon( a ); // Convert the string to an FMon
  return back; // Return the FMon
}

/*
 * Function Name: fAlgFromStr
 *
 * Overview: Pick an FAlg from a string such as "x*y - z;"
 *
 * Detail: Starting from position _j_ in a string _s_,
 * read in an FAlg and return it. Note that the FAlg
 * in the string must be terminated with a semicolon.
 * Once the FAlg has been read, place the position we
 * have reached in the string in the variable _pk_.
 */
FAlg
fAlgFromStr( s, j, pk )
char s[];
int j, *pk;
{
  char c = ' ', a[MAXLINE];
  int i = 0, k = j;
  FAlg back;

  // Until a semicolon is encountered
  while( c != ';' )
  {
    c = s[k]; // Read a character from the string

    // If a semicolon is encountered
    if( c == ';' )
    {
      a[i] = '\0'; // Finish with the null character
    }
    else
    {
      a[i] = c; // Continue to process...
      i++;
    }
    k++;
  }
  *pk = k; // Place the finish position in the variable _pk_

  back = parseStrToFAlg( a ); // Convert the string to an FAlg
  return back; // Return the FAlg
}

/*
 * =================================
 * High Level File Reading Functions
 * =================================
 */

/*
 * Function Name: fMonListFromFile
 *
 * Overview: Routine to read an FMonList from the first line of a file
 *
 * Detail: Given an input file, this function
 * reads the first line of the file and returns
 * the semicolon separated FMonList found on that line.
 * For example, if the input is a list such as a; b; A; B;
 * then the output is the FMonList (a, b, A, B).
 */
FMonList
fMonListFromFile( infil )
FILE *infil;
{
  FMon w;
  FMonList words = fMonListNul;
  char s[MAXLINE];
  int j = 0, k = 0, len = 0;

  // Get the first line of the file and its length
  len = getLine( infil, s, MAXLINE );

  // While there are more FMons to be found
  while( j < len )
  {
    w = fMonFromStr( s, j, &k ); // Obtain an FMon
    j = k; // Set the next starting position
    words = fMonListPush( w, words ); // Construct the list
  }

  // Return the list - note that we must reverse the list
  // because it has been read in reverse order.
  return fMonListFXRev( words );
}

/*
 * Function Name: fAlgListFromFile
 *
 * Overview: Routine to read an FAlgList from a file
 *
 * Detail: Given an input file, this function
 * takes each line of the file in turn, pushing one FAlg from
 * each line onto an FAlgList. This process is
 * continued until there are no more lines in the file
 * to process. For example, if the input is a list such as
 *
 * 2*x - 4*y;
 * 5*x*y;
 * 4 + 5*x + 60*y;
 *
 * then the output is the FAlgList
 * (2x-4y, 5xy, 4+5x+60y).
 */
FAlgList
fAlgListFromFile( infil )
FILE *infil;
{
  FAlg entry;
  FAlgList back = fAlgListNul;
  char s[MAXLINE];
  int j = 0, k = 0, len;

  // Get the first line of the file
  len = getLine( infil, s, MAXLINE );

  // While there are still lines to process
  while( len > 0 )
  {
    entry = fAlgFromStr( s, j, &k ); // Obtain an FAlg from a line
    back = fAlgListPush( entry, back ); // Push the FAlg onto the list
    len = getLine( infil, s, MAXLINE ); // Get a new line
  }

  // Return the list - note that we must reverse the list
  // because it has been read in reverse order.
  return fAlgListFXRev( back );
}

/*
 * =================================
 * High Level File Writing Functions
 * =================================
 */

/*
 * Function Name: fMonToFile
 *
 * Overview: Writes an FMon (in parse format) followed by a semicolon to a file
 *
 * Detail: Given an input file and an FMon, this function
 * writes the FMon to file in parse format followed by a semicolon.
 */
void
fMonToFile( infil, w )
FILE *infil;
FMon w;
{
  FMon wM;
  ULong length;

  // If the FMon is non-empty
  if ( fMonEqual( w, fMonOne() ) != (Bool) 1 )
  {
    // While there are letters left in the FMon
    while ( w )
    {
      wM = fMonLeadPowFac( w ); // Obtain a factor
      fprintf( infil, "%s", fMonToStr( wM ) ); // Write the factor to file
      length = fMonLength( wM );
      w = fMonSuffix( w, fMonLength( w ) - length );
      if ( fMonEqual( w, fMonOne() ) != (Bool) 1 )
      {
        // In parse format, to separate variables we use an asterisk
        fprintf( infil, "*" );
      }
    }
    fprintf( infil, ";" ); // At the end write a semicolon to file
  }
  else // Just write a semicolon to file
  {
    fprintf( infil, ";" );
  }
}

/*
 * Function Name: fMonListToFile
 *
 * Overview: Writes an FMonList to a file on a single line
 *
 * Detail: Given an input file and an FMonList, this function
 * writes the list to file as l1; l2; l3; ...
 */
void
fMonListToFile( infil, L )
FILE *infil;
FMonList L;
{
  ULong i, length = fMonListLength( L );

  // For each element of the list
  for( i = 1; i <= length; i++ )
  {
    // Write an FMon to file
    fMonToFile( infil, L -> first );

    // If there are more FMons left to look at
    if( i < length )
    {
      fprintf( infil, " " ); // Provide a space between elements
    }
    else // else terminate the line
    {
     fprintf( infil, "\n" );
    }
    L = L -> rest;
  }
}

/*
 * ================================
 * File Name Modification Functions
 * ================================
 */

/*
 * Function Name: appendDotDegRevLex
 *
 * Overview: Appends ".drl" onto a string (except in special case "*.in")
 *
 * Detail: Given an input character array, this function
 * appends the String ".drl" onto the end of the character array.
 * In the special case that the input ends with ".in", the function
 * replaces the ".in" with ".drl".
 */
String
appendDotDegRevLex( input )
char input[];
{
  int length = (int) strlen( input );
  String back = strNew();

  // First check for .in at the end of the file name
  if ( input[length-1] == 'n' & input[length-2] == 'i' & input[length-3] == '.' )
  {
    input[length-2] = 'd';
    input[length-1] = 'r';
    sprintf( back, "%s%s", input, "l" );
  }
  else // Just append with ".drl"
  {
    sprintf( back, "%s%s", input, ".drl" );
  }

  return back;
}

/*
 * Function Name: appendDotDegLex
 *
 * Overview: Appends ".deg" onto a string (except in special case "*.in")
 *
 * Detail: Given an input character array, this function
 * appends the String ".deg" onto the end of the character array.
 * In the special case that the input ends with ".in", the function
 * replaces the ".in" with ".deg".
 */
String
appendDotDegLex( input )
char input[];
{
  int length = (int) strlen( input );
  String back = strNew();

  // First check for .in at the end of the file name
  if ( input[length-1] == 'n' & input[length-2] == 'i' & input[length-3] == '.' )
  {
    input[length-2] = 'd';
    input[length-1] = 'e';
    sprintf( back, "%s%s", input, "g" );
  }
  else // Just append with ".deg"
  {
    sprintf( back, "%s%s", input, ".deg" );
  }

  return back;
}

/*
 * Function Name: appendDotLex
 *
 * Overview: Appends ".lex" onto a string (except in special case "*.in")
 *
 * Detail: Given an input character array, this function
 * appends the String ".lex" onto the end of the character array.
 * In the special case that the input ends with ".in", the function
 * replaces the ".in" with ".lex".
 */
String
appendDotLex( input )
char input[];
{
  int length = (int) strlen( input );
  String back = strNew();

  // First check for .in at the end of the file name
  if ( input[length-1] == 'n' & input[length-2] == 'i' & input[length-3] == '.' )
  {
    input[length-2] = 'l';
    input[length-1] = 'e';
    sprintf( back, "%s%s", input, "x" );
  }
  else // Just append with ".lex"
  {
    sprintf( back, "%s%s", input, ".lex" );
  }

  return back;
}

/*
 * Function Name: appendDotWP
 *
 * Overview: Appends ".wp" onto a string (except in special case "*.in")
 *
 * Detail: Given an input character array, this function
 * appends the String ".wp" onto the end of the character array.
 * In the special case that the input ends with ".in", the function
 * replaces the ".in" with ".wp".
 */
String
appendDotWP( input )
char input[];
{
  int length = (int) strlen( input );
  String back = strNew();

  // First check for .in at the end of the file name
  if ( input[length-1] == 'n' & input[length-2] == 'i' & input[length-3] == '.' )
  {
    input[length-2] = 'w';
    input[length-1] = 'p';
    sprintf( back, "%s", input );
  }
  else // Just append with ".wp"
  {
    sprintf( back, "%s%s", input, ".wp" );
  }

  return back;
}

/*
 * Function Name: filenameLength
 *
 * Overview: Calculates the length of an input string
 *
 * Detail: Given an input character array, this function
 * finds the length of that character array
 */
int
filenameLength( s )
char s[];
{
  int i = 0;

  while( s[i] != '\0' ) i++;

  return i;
}

/*
 * ===========
 * End of File
 * ===========
 */
\end{lstlisting}

\subsection{fralg\_functions.h} \label{appBf6}
\begin{lstlisting}
/*
 * File: fralg_functions.h
 * Author: Gareth Evans
 * Last Modified: 10th August 2005
 */

// Initialise file definition
# ifndef FRALG_FUNCTIONS_HDR
# define FRALG_FUNCTIONS_HDR

// Include MSSRC Libraries
# include <fralg.h>

// Include System Libraries
# include <limits.h>

// Include *_functions Libraries
# include "list_functions.h"
# include "arithmetic_functions.h"

//
// External Variables Required
//

extern ULong nRed;          // Stores how many reductions have been performed
extern int nOfGenerators,   // Holds the number of generators
           pl;              // Holds the "Print Level"

//
// Functions Defined in fralg_functions.c
//

//
// Ordering Functions
//

// Returns 1 if 1st arg <_{Lex} 2nd arg
Bool fMonLex( FMon, FMon );
// Returns 1 if 1st arg <_{InvLex} 2nd arg
Bool fMonInvLex( FMon, FMon );
// Returns 1 if 1st arg <_{DegRevLex} 2nd arg
Bool fMonDegRevLex( FMon, FMon );
// Returns 1 if 1st arg <_{WreathProduct} 2nd arg
Bool fMonWreathProd( FMon, FMon );

//
// Alphabet Manipulation Functions
//

// Substitutes ASCII generators for original generators in a list of polynomials
FAlgList preProcess( FAlgList, FMonList );
// Substitutes original generators for ASCII generators in a given polynomial
String postProcess( FAlg, FMonList );
// As above but gives back its output in parse format
String postProcessParse( FAlg, FMonList );
// Adjusts the original generator order (1st arg) according to frequency of generators in 2nd arg
FMonList alphabetOptimise( FMonList, FAlgList );

//
// Polynomial Manipulation Functions
//

// Returns all possible ways that 2nd arg divides 1st arg; 3rd arg = is division possible?
FMonPairList fMonDiv( FMon, FMon, Short * );
// Returns the first way that 2nd arg divides 1st arg; 3rd arg = is division possible?
FMonPairList fMonDivFirst( FMon, FMon, Short * );
// Finds all possible overlaps of 2 FMons
FMonPairList fMonOverlaps( FMon, FMon );
// Returns the degree-based initial of a polynomial
FAlg degInitial( FAlg );
// Reverses a monomial
FMon fMonReverse( FMon );

//
// Groebner Basis Functions
//

// Returns the normal form of a polynomial w.r.t. a list of polynomials
FAlg polyReduce( FAlg, FAlgList );
// Minimises a given Groebner Basis
FAlgList minimalGB( FAlgList );
// Reduces each member of a Groebner Basis w.r.t. all other members
FAlgList reducedGB( FAlgList );
// Tests whether a given FAlg reduces to 0 using the given FAlgList
Bool idealMembershipProblem( FAlg, FAlgList );

# endif // FRALG_FUNCTIONS_HDR
\end{lstlisting}

\subsection{fralg\_functions.c} \label{appBf7}
\lstinputlisting{fralg_functions.c}

\subsection{list\_functions.h} \label{appBf8}
\begin{lstlisting}
/*
 * File: list_functions.h
 * Author: Gareth Evans
 * Last Modified: 9th August 2005
 */

// Initialise file definition
# ifndef LIST_FUNCTIONS_HDR
# define LIST_FUNCTIONS_HDR

// Include MSSRC Libraries
# include <fralg.h>

//
// External Variables Required
//

extern int pl; // Holds the "Print Level"

//
// Display Functions
//

// Displays an FMonList in the format l1\n l2\n l3\n...
void fMonListDisplay( FMonList );
// Displays an FMonList in the format l1 > l2 > l3...
void fMonListDisplayOrder( FMonList );
// Displays an FMonPairList in the format (l1, l2)\n (l3, l4)\n...
void fMonPairListMultDisplay( FMonPairList );
// Displays an FAlgList in the format p1\n p2\n p3\n...
void fAlgListDisplay( FAlgList );

//
// List Extraction Functions
//

// Returns the ith member of an FMonList (i = 1st arg)
FMon fMonListNumber( ULong, FMonList );
// Returns the ith member of an FMonPairList (i = 1st arg)
FMonPair fMonPairListNumber( ULong, FMonPairList );
// Returns the ith member of an FAlgList (i = 1st arg)
FAlg fAlgListNumber( ULong, FAlgList );

//
// List Membership Functions
//

// Does the FAlg appear in the FAlgList? (1 = yes)
Bool fAlgListIsMember( FAlg, FAlgList );

//
// List Position Functions
//

// Gives position of 1st appearance of FMon in FMonList
ULong fMonListPosition( FMon, FMonList );
// Gives position of 1st appearance of FAlg in FAlgList
ULong fAlgListPosition( FAlg, FAlgList );
   
//
// Sorting Functions
//

// Swaps 2 elements in arrays of ULongs and FMons
void alphabetArraySwap( ULong[], FMon[], ULong, ULong );
// Sorts an array of ULongs (largest first) and applies the same changes to the FMon array
void alphabetArrayQuickSort( ULong[], FMon[], ULong, ULong );
// Swaps 2 elements in arrays of FAlgs and FMons
void fAlgArraySwap( FAlg[], FMon[], ULong, ULong );
// Sorts an array of FAlgs using DegRevLex (largest first)
void fAlgArrayQuickSortDRL( FAlg[], FMon[], ULong, ULong );

// Sorts an array of FAlgs using theOrdFun (largest first)
void fAlgArrayQuickSortOrd( FAlg[], FMon[], ULong, ULong );
// Sorts an FAlgList (largest first)
FAlgList fAlgListSort( FAlgList, int );
// Swaps 2 elements in arrays of FMons, ULongs and ULongs
void multiplicativeArraySwap( FMon[], ULong[], ULong[], ULong, ULong );
// Sorts input data to OverlapDiv w.r.t. DegRevLex (largest first)
void multiplicativeQuickSort( FMon[], ULong[], ULong[], ULong, ULong );

//
// Insertion Sort Functions
//

// Insert into list according to DegRevLex
FAlgList fAlgListDegRevLexPush( FAlg, FAlgList );
// As above, but also returns the insertion position
FAlgList fAlgListDegRevLexPushPosition( FAlg, FAlgList, ULong * );
// Insert into list according to the current monomial ordering
FAlgList fAlgListNormalPush( FAlg, FAlgList );
// As above, but also returns the insertion position
FAlgList fAlgListNormalPushPosition( FAlg, FAlgList, ULong * );

//
// Deletion Functions
//

// Removes the (1st arg)-th element from the list
FMonList fMonListRemoveNumber( ULong, FMonList );
// Removes the (1st arg)-th element from the list
FMonPairList fMonPairListRemoveNumber( ULong, FMonPairList );
// Removes the (1st arg)-th element from the list
FAlgList fAlgListRemoveNumber( ULong, FAlgList );

//
// Normalising Functions
//

// Removes any fractions found in the FAlgList by scalar multiplication
FAlgList fAlgListRemoveFractions( FAlgList );

# endif // LIST_FUNCTIONS_HDR
\end{lstlisting}

\subsection{list\_functions.c} \label{appBf9}
\lstinputlisting{list_functions.c}

\subsection{ncinv\_functions.h} \label{appBf10}
\begin{lstlisting}
/*
 * File: ncinv_functions.h
 * Author: Gareth Evans
 * Last Modified: 6th July 2005
 */

// Initialise file definition
# ifndef NCINV_FUNCTIONS_HDR
# define NCINV_FUNCTIONS_HDR

// Include MSSRC Libraries
# include <fralg.h>

//
// External Variables Required
//

extern ULong nOfProlongations, // Stores the number of prolongations calculated
             nRed;             // Stores how many reductions have been performed
extern int degRestrict,        // Determines whether of not prolongations are restricted by degree
           EType,              // Stores the type of Overlap Division
           IType,              // Stores the involutive division used
           nOfGenerators,      // Holds the number of generators
           pl,                 // Holds the "Print Level"
           SType,              // Determines how the basis is sorted
           MType;              // Determines involutive division method

//
// Functions Defined in ncinv_functions.c
//

//
// Overlap Functions
//

// Returns the union of (non-)multiplicative variables (1st arg) and a generator (2nd arg)
FMon multiplicativeUnion( FMon, FMon );
// Does the generator (1st arg) appear in the list of multiplicative variables (2nd arg)?
int fMonIsMultiplicative( FMon, FMon );
// Does the 1st arg appear as a subword in the 2nd arg (yes (1)/no (0))
int fMonIsSubword( FMon, FMon );
// Is the 1st arg a subword of the 2nd arg; if so, return start pos in 2nd arg
ULong fMonSubwordOf( FMon, FMon, ULong );

// Returns size of smallest overlap of type (suffix of 1st arg = prefix of 2nd arg)
ULong fMonPrefixOf( FMon, FMon, ULong, ULong );
// Returns size of smallest overlap of type (prefix of 1st arg = suffix of 2nd arg)
ULong fMonSuffixOf( FMon, FMon, ULong, ULong );

//
// Multiplicative Variables Functions
//

// Returns no ('empty') multiplicative variables
void EMultVars( FMon, ULong *, ULong * );
// All variables left mult., no variables right mult.
void LMultVars( FMon, ULong *, ULong * );
// All variables right mult., no variables left mult.
void RMultVars( FMon, ULong *, ULong * );
// Returns local overlap-based multiplicative variables
FMonPairList OverlapDiv( FAlgList );

//
// Polynomial Reduction and Basis Completion Functions
//

// Reduces 1st arg w.r.t. 2nd arg (list) and 3rd arg (vars)
FAlg IPolyReduce( FAlg, FAlgList, FMonPairList );
// Autoreduces an FAlgList recursively until no more reductions are possible
FAlgList IAutoreduceFull( FAlgList );
// Implements Seiler's original algorithm for computing locally involutive bases
FAlgList Seiler( FAlgList );
// Implements Gerdt's advanced algorithm for computing locally involutive bases
FAlgList Gerdt( FAlgList );

# endif // NCINV_FUNCTIONS_HDR
\end{lstlisting}

\subsection{ncinv\_functions.c} \label{appBf11}
\lstinputlisting{ncinv_functions.c}

\subsection{involutive.c} \label{appBf12}
\lstinputlisting{involutive.c}

%
% Appendix C
% Author: Gareth Evans
% Last Modified: 20th September 2005
%

\chapter{Program Output} \label{appC}

\lstset{
language=,
frame=single,
framesep=5pt,
columns=fullflexible,
basicstyle=\scriptsize,
numbers=none,
}

In this Appendix, we provide sample sessions
showing how the program given in Appendix \ref{appB}
can be used to compute noncommutative Involutive Bases with
respect to different involutive divisions and
monomial orderings.

\section{Sample Sessions}

\subsection{Session 1: Locally Involutive Bases}

{\bf Task:} If $F  := \{x^2y^2 - 2xy^2 + x^2, \; x^2y - 2xy\}$
generates an ideal $J$ over the polynomial ring
$\mathbb{Q}\langle x, y\rangle$, compute a Locally
Involutive Basis for $F$ with respect to the strong left
overlap division $\mathcal{S}$; thick divisors;
and the DegLex monomial ordering.

{\bf Origin of Example:} Example \ref{appCex1}.

{\bf Input File:}
\begin{lstlisting}
x; y;
x^2*y^2 - 2*x*y^2 + x^2;
x^2*y - 2*x*y;

\end{lstlisting}

{\bf Plan:} Apply the program given in Appendix \ref{appB}
to the above file, using the `-c2' option to select
Algorithm \ref{noncom-inv}; the `-d' option to select
the DegLex monomial ordering; the \linebreak
`-m2' option to select
thick divisors; and the `-e2' and `-s1' options
to select the strong left overlap division.

{\bf Program Output:}
\begin{lstlisting}
ma6:mssrc-aux/thesis> time involutive -c2 -d -e2 -m2 -s1 thesis1.in

*** NONCOMMUTATIVE INVOLUTIVE BASIS PROGRAM (LOCAL DIVISION) ***

Using the DegLex Ordering with x > y

Polynomials in the input basis:
x^2 y^2 - 2 x y^2 + x^2,
x^2 y - 2 x y,
[2 Polynomials]

Computing an Involutive Basis...
Added Polynomial #3 to Basis...
Added Polynomial #4 to Basis...
Autoreduction reduced the basis to size 3...
Added Polynomial #4 to Basis...
Autoreduction reduced the basis to size 3...
Added Polynomial #4 to Basis...
Added Polynomial #5 to Basis...
...Involutive Basis Computed.

Here is the Involutive Basis
((Left, Right) Multiplicative Variables in Brackets):
x y^2 x, (x y, 1),
x y^2, (x y, y),
x y x, (x y, 1),
x y, (x y, 1),
x^2, (x y, 1),
[5 Polynomials]

Computing the Reduced Groebner Basis...
...Reduced Groebner Basis Computed.

Here is the Reduced Groebner Basis:
x y,
x^2,
[2 Polynomials]

Writing Reduced Groebner Basis to Disk... Done.
Writing Involutive Basis to Disk... Done.

0.000u 0.007s 0:00.15 0.0%      0+0k 0+2io 16pf+0w
ma6:mssrc-aux/thesis>
\end{lstlisting}

{\bf Output File:}
\begin{lstlisting}
x; y;
x*y^2*x; (x y, 1);
x*y^2; (x y, y);
x*y*x; (x y, 1);
x*y; (x y, 1);
x^2; (x y, 1);
\end{lstlisting}

\subsection{Session 2: Involutive Complete Rewrite Systems}

{\bf Task:} If $F := \{x^3-1, \; y^2-1, \; (xy)^2-1, \;
Xx-1, \; xX-1, \; Yy-1, \; yY-1\}$ generates an ideal
$J$ over the polynomial ring
$\mathbb{Q}\langle Y, X, y, x\rangle$, compute an
Involutive Basis for $F$ with respect to the left division
$\lhd$ and the DegLex monomial ordering.

{\bf Origin of Example:} Example \ref{S3} ($F$ corresponds to a
monoid rewrite system for the group $S_3$; we want to
compute an involutive complete rewrite system for
$S_3$).

{\bf Input File:}
\begin{lstlisting}
Y; X; y; x;
x^3 - 1;
y^2 - 1;
(x*y)^2 - 1;
X*x - 1;
x*X - 1;
Y*y - 1;
y*Y - 1;

\end{lstlisting}

{\bf Plan:} Apply the program given in Appendix \ref{appB}
to the above file, using the `-c2' option to select
Algorithm \ref{noncom-inv} and the `-d' option to select
the DegLex monomial ordering (the left division
is selected by default).

{\bf Program Output:}
\begin{lstlisting}
ma6:mssrc-aux/thesis> time involutive -c2 -d thesis2.in

*** NONCOMMUTATIVE INVOLUTIVE BASIS PROGRAM (GLOBAL DIVISION) ***

Using the DegLex Ordering with Y > X > y > x

Polynomials in the input basis:
x^3 - 1,
y^2 - 1,
x y x y - 1,
X x - 1,
x X - 1,
Y y - 1,
y Y - 1,
[7 Polynomials]

Computing an Involutive Basis...
Added Polynomial #8 to Basis...
Added Polynomial #9 to Basis...
Added Polynomial #10 to Basis...
Added Polynomial #11 to Basis...
Added Polynomial #12 to Basis...
Added Polynomial #13 to Basis...
Autoreduction reduced the basis to size 11...
Added Polynomial #12 to Basis...
Added Polynomial #13 to Basis...
Added Polynomial #14 to Basis...
Added Polynomial #15 to Basis...
Added Polynomial #16 to Basis...
Added Polynomial #17 to Basis...
Added Polynomial #18 to Basis...
Added Polynomial #19 to Basis...
Added Polynomial #20 to Basis...
Added Polynomial #21 to Basis...
Added Polynomial #22 to Basis...
Added Polynomial #23 to Basis...
Autoreduction reduced the basis to size 19...
Added Polynomial #20 to Basis...
Autoreduction reduced the basis to size 19...
Added Polynomial #20 to Basis...
Added Polynomial #21 to Basis...
Added Polynomial #22 to Basis...
Added Polynomial #23 to Basis...
Added Polynomial #24 to Basis...
Added Polynomial #25 to Basis...
Added Polynomial #26 to Basis...
Added Polynomial #27 to Basis...
Added Polynomial #28 to Basis...
Added Polynomial #29 to Basis...
Added Polynomial #30 to Basis...
Added Polynomial #31 to Basis...
Added Polynomial #32 to Basis...
Added Polynomial #33 to Basis...
Added Polynomial #34 to Basis...
Added Polynomial #35 to Basis...
Added Polynomial #36 to Basis...
Added Polynomial #37 to Basis...
Added Polynomial #38 to Basis...
Added Polynomial #39 to Basis...
Added Polynomial #40 to Basis...
Autoreduction reduced the basis to size 29...
Added Polynomial #30 to Basis...
Autoreduction reduced the basis to size 19...
...Involutive Basis Computed.

Here is the Involutive Basis
((Left, Right) Multiplicative Variables in Brackets):
y^2 - 1, (Y X y x, 1),
X x - 1, (Y X y x, 1),
x X - 1, (Y X y x, 1),
Y y - 1, (Y X y x, 1),
y^2 x - x, (Y X y x, 1),
Y - y, (Y X y x, 1),
Y x - y x, (Y X y x, 1),
X x y - y, (Y X y x, 1),
Y y x - x, (Y X y x, 1),
x^2 - X, (Y X y x, 1),
X^2 - x, (Y X y x, 1),
x y x - y, (Y X y x, 1),
X y - y x, (Y X y x, 1),
X y x - x y, (Y X y x, 1),
x^2 y - y x, (Y X y x, 1),
y X - x y, (Y X y x, 1),
y x y - X, (Y X y x, 1),
Y x y - X, (Y X y x, 1),
Y X - x y, (Y X y x, 1),
[19 Polynomials]

Computing the Reduced Groebner Basis...
...Reduced Groebner Basis Computed.

Here is the Reduced Groebner Basis:
y^2 - 1,
X x - 1,
x X - 1,
Y - y,
x^2 - X,
X^2 - x,
x y x - y,
X y - y x,
y X - x y,
y x y - X,
[10 Polynomials]

Writing Reduced Groebner Basis to Disk... Done.
Writing Involutive Basis to Disk... Done.

0.105u 0.000s 0:00.16 62.5%     197+727k 0+2io 0pf+0w
ma6:mssrc-aux/thesis> 
\end{lstlisting}

{\bf Output File:}
\begin{lstlisting}
Y; X; y; x;
y^2 - 1; (Y X y x, 1);
X*x - 1; (Y X y x, 1);
x*X - 1; (Y X y x, 1);
Y*y - 1; (Y X y x, 1);
y^2*x - x; (Y X y x, 1);
Y - y; (Y X y x, 1);
Y*x - y*x; (Y X y x, 1);
X*x*y - y; (Y X y x, 1);
Y*y*x - x; (Y X y x, 1);
x^2 - X; (Y X y x, 1);
X^2 - x; (Y X y x, 1);
x*y*x - y; (Y X y x, 1);
X*y - y*x; (Y X y x, 1);
X*y*x - x*y; (Y X y x, 1);
x^2*y - y*x; (Y X y x, 1);
y*X - x*y; (Y X y x, 1);
y*x*y - X; (Y X y x, 1);
Y*x*y - X; (Y X y x, 1);
Y*X - x*y; (Y X y x, 1);
\end{lstlisting}

\subsection{Session 3: Noncommutative Involutive Walks}

{\bf Task:} If $G' := \{y^2+2xy, \; y^2+x^2, \;
5y^3, \; 5xy^2, \; y^2+2yx\}$
generates an ideal $J$ over the polynomial ring
$\mathbb{Q}\langle x, y\rangle$, compute an
Involutive Basis for $G'$ with respect to the left
division $\lhd$ and the DegRevLex monomial ordering.

{\bf Origin of Example:} Example \ref{appCch6} ($G'$ corresponds
to a set of initials in the noncommutative Involutive Walk
algorithm; we want to compute an Involutive Basis $H'$
for $G'$).

{\bf Input File:}
\begin{lstlisting}
x; y;
y^2 + 2*x*y;
y^2 + x^2;
5*y^3;
5*x*y^2;
y^2 + 2*y*x;

\end{lstlisting}

{\bf Plan:} Apply the program given in Appendix \ref{appB}
to the above file, using the `-c2' option to select
Algorithm \ref{noncom-inv} (the DegRevLex monomial ordering
and the left division are selected by default).

{\bf Program Output:}
\begin{lstlisting}
ma6:mssrc-aux/thesis> time involutive -c2 thesis3.in

*** NONCOMMUTATIVE INVOLUTIVE BASIS PROGRAM (GLOBAL DIVISION) ***

Using the DegRevLex Ordering with x > y

Polynomials in the input basis:
y^2 + 2 x y,
y^2 + x^2,
5 y^3,
5 x y^2,
y^2 + 2 y x,
[5 Polynomials]

Computing an Involutive Basis...
...Involutive Basis Computed.

Here is the Involutive Basis
((Left, Right) Multiplicative Variables in Brackets):
2 y x - x^2, (x y, 1),
y x^2, (x y, 1),
x^3, (x y, 1),
2 x y - x^2, (x y, 1),
y^2 + x^2, (x y, 1),
[5 Polynomials]

Computing the Reduced Groebner Basis...
...Reduced Groebner Basis Computed.

Here is the Reduced Groebner Basis:
2 y x - x^2,
x^3,
2 x y - x^2,
y^2 + x^2,
[4 Polynomials]

Writing Reduced Groebner Basis to Disk... Done.
Writing Involutive Basis to Disk... Done.

0.005u 0.000s 0:00.07 0.0%      0+0k 0+2io 0pf+0w
ma6:mssrc-aux/thesis>
\end{lstlisting}

{\bf More Verbose Program Output:} (we select the
`-v3' option to obtain more information about the
autoreduction that occurs at the start of the algorithm).
\begin{lstlisting}
ma6:mssrc-aux/thesis> time involutive -c2 -v3 thesis3.in

*** NONCOMMUTATIVE INVOLUTIVE BASIS PROGRAM (GLOBAL DIVISION) ***

Using the DegRevLex Ordering with x (AAB) > y (AAA)

Polynomials in the input basis:
AAA^2 + 2 AAB AAA,
AAA^2 + AAB^2,
5 AAA^3,
5 AAB AAA^2,
AAA^2 + 2 AAA AAB,
[5 Polynomials]

Computing an Involutive Basis...
Autoreducing...
Looking at element p = AAA^2 + 2 AAA AAB of basis
Reduced p to AAB AAA - AAA AAB
Looking at element p = 5 AAB AAA^2 of basis
Reduced p to AAB AAA AAB
Looking at element p = AAB AAA - AAA AAB of basis
Reduced p to AAB AAA - AAA AAB
Looking at element p = 5 AAA^3 of basis
Reduced p to AAA^2 AAB
Looking at element p = AAB AAA AAB of basis
Reduced p to AAB AAA AAB
Looking at element p = AAB AAA - AAA AAB of basis
Reduced p to AAB AAA - AAA AAB
Looking at element p = AAA^2 + AAB^2 of basis
Reduced p to 2 AAA AAB - AAB^2
Looking at element p = AAA^2 AAB of basis
Reduced p to AAA AAB^2
Looking at element p = 2 AAA AAB - AAB^2 of basis
Reduced p to 2 AAA AAB - AAB^2
Looking at element p = AAB AAA AAB of basis
Reduced p to AAB^3
Looking at element p = AAA AAB^2 of basis
Reduced p to AAA AAB^2
Looking at element p = 2 AAA AAB - AAB^2 of basis
Reduced p to 2 AAA AAB - AAB^2
Looking at element p = AAB AAA - AAA AAB of basis
Reduced p to 2 AAB AAA - AAB^2
Looking at element p = AAB^3 of basis
Reduced p to AAB^3
Looking at element p = AAA AAB^2 of basis
Reduced p to AAA AAB^2
Looking at element p = 2 AAA AAB - AAB^2 of basis
Reduced p to 2 AAA AAB - AAB^2
Looking at element p = AAA^2 + 2 AAB AAA of basis
Reduced p to AAA^2 + AAB^2
Looking at element p = 2 AAB AAA - AAB^2 of basis
Reduced p to 2 AAB AAA - AAB^2
Looking at element p = AAB^3 of basis
Reduced p to AAB^3
Looking at element p = AAA AAB^2 of basis
Reduced p to AAA AAB^2
Looking at element p = 2 AAA AAB - AAB^2 of basis
Reduced p to 2 AAA AAB - AAB^2
Analysing AAA AAB...
Adding Right Prolongation by variable #0 to S...
Adding Right Prolongation by variable #1 to S...
Analysing AAA AAB^2...
Adding Right Prolongation by variable #0 to S...
Adding Right Prolongation by variable #1 to S...
Analysing AAB^3...
Adding Right Prolongation by variable #0 to S...
Adding Right Prolongation by variable #1 to S...
Analysing AAB AAA...
Adding Right Prolongation by variable #0 to S...
Adding Right Prolongation by variable #1 to S...
Analysing AAA^2...
Adding Right Prolongation by variable #0 to S...
Adding Right Prolongation by variable #1 to S...
...Involutive Basis Computed.
Number of Prolongations Considered = 0

Here is the Involutive Basis
((Left, Right) Multiplicative Variables in Brackets):
2 AAA AAB - AAB^2, (all, none),
AAA AAB^2, (all, none),
AAB^3, (all, none),
2 AAB AAA - AAB^2, (all, none),
AAA^2 + AAB^2, (all, none),
[5 Polynomials]

Computing the Reduced Groebner Basis...

Looking at element p = 2 AAA AAB - AAB^2 of basis
Reduced p to 2 AAA AAB - AAB^2

Looking at element p = AAB^3 of basis
Reduced p to AAB^3

Looking at element p = 2 AAB AAA - AAB^2 of basis
Reduced p to 2 AAB AAA - AAB^2

Looking at element p = AAA^2 + AAB^2 of basis
Reduced p to AAA^2 + AAB^2
Number of Reductions Carried out = 34
...Reduced Groebner Basis Computed.

Here is the Reduced Groebner Basis:
2 AAA AAB - AAB^2,
AAB^3,
2 AAB AAA - AAB^2,
AAA^2 + AAB^2,
[4 Polynomials]

Writing Reduced Groebner Basis to Disk... Done.
Writing Involutive Basis to Disk... Done.

0.000u 0.005s 0:00.04 0.0%      0+0k 0+2io 0pf+0w
ma6:mssrc-aux/thesis>
\end{lstlisting}

{\bf Output File:}
\begin{lstlisting}
x; y;
2*y*x - x^2; (x y, 1);
y*x^2; (x y, 1);
x^3; (x y, 1);
2*x*y - x^2; (x y, 1);
y^2 + x^2; (x y, 1);
\end{lstlisting}

\subsection{Session 4: Ideal Membership}

{\bf Task:} If $F := \{x+y+z-3, \; x^2+y^2+z^2-9, \;
x^3+y^3+z^3-24\}$ generates an ideal $J$ over the polynomial ring
$\mathbb{Q}\langle x, y, z\rangle$, are the
polynomials $x+y+z-3$; $x+y+z-2$; $xz^2+yz^2-1$;
$zyx+1$ and $x^{10}$ members of $J$?

{\bf Input File:}
\begin{lstlisting}
x; y; z;
x + y + z - 3;
x^2 + y^2 + z^2 - 9;
x^3 + y^3 + z^3 - 24;

\end{lstlisting}

{\bf Plan:} To solve the ideal membership problem for the
five given polynomials, we first need to obtain a Gr\"obner or
Involutive Basis for $F$. We shall do this by applying
the program given in Appendix \ref{appB} to compute
an Involutive Basis for $F$ with respect to the
DegLex monomial ordering and the right division $\rhd$ (this
requires the `-d' and `-s4' options respectively).
Once the Involutive Basis has been computed (which then
allows the program to compute the unique reduced
Gr\"obner Basis $G$ for $F$), we can start an ideal
membership problem solver (courtesy of the `-p' option)
which allows us to type in a polynomial $p$ and
find out whether or not $p$ is a member of $J$
(the program reduces $p$ with respect to $G$,
testing to see whether or not a zero remainder is obtained).

{\bf Program Output:}
\lstinputlisting{thesis4.out}

{\bf Output File:}
\lstinputlisting{thesis4.deg.inv}

% BIBLIOGRAPHY
\clearpage
\phantomsection
\addcontentsline{toc}{chapter}{Bibliography}
\bibliographystyle{ncib}
\bibliography{ncib}

% INDEX
\clearpage
\phantomsection
\addcontentsline{toc}{chapter}{Index}
\printindex

\end{document}